\documentclass[reqno, 11pt]{amsbook}
\usepackage[utf8]{inputenc}
\usepackage[T1]{fontenc}
\usepackage{lmodern}
\usepackage{epsfig,amsmath, amsthm ,amsfonts,latexsym,stmaryrd, amssymb}
\usepackage[abs]{overpic}
\usepackage{simpler-wick}
\usepackage{caption}
\usepackage{subcaption}
\usepackage{tikz-cd}
\usepackage{dsfont}
\usepackage{mathtools}
\usepackage{tikz-cd}
\usetikzlibrary{shapes,arrows}
\usepackage{comment}
\usepackage{enumitem}
\usepackage{marginnote}
\usepackage[all,cmtip]{xy}
\usepackage[backend=biber, maxnames=4, style=alphabetic]{biblatex}
\usepackage{amssymb}
\usepackage{movie15}
\usepackage{animate}
\usepackage{hyperref}
\usepackage{tikz}

\newcommand{\gsim}{\raisebox{-0.13cm}{~\shortstack{$>$ \\[-0.07cm]
      $\sim$}}~}

\tikzcdset{scale cd/.style={every label/.append style={scale=#1},
    cells={nodes={scale=#1}}}}

\theoremstyle{plain}

\newtheorem{theorem}{Theorem}[section]

\newtheorem{dummy}{anything}[section]
\newtheorem{thm}{Theorem}[section]

\newtheorem{lemma}[dummy]{Lemma}

\newtheorem{proposition}{Proposition}[section]
\newtheorem*{prop}{Proposition}
\newtheorem{corollary}[dummy]{Corollary}
\newtheorem{conjecture}[dummy]{Conjecture}

\theoremstyle{definition}
\newtheorem{definition}[dummy]{Definition}

\newtheorem{example}[dummy]{Example}

\newtheorem{remark}[dummy]{Remark}

\newtheorem*{assumptions}{Assumptions}
\theoremstyle{remark}
\newcommand{\N}{\mathbb{N}}
\newcommand{\Z}{\mathbb{Z}}
\newcommand{\R}{\mathbb{R}}
\newcommand{\C}{\mathbb{C}}

\newcommand{\M}{\mathcal{M}}

\newcommand{\J}{\mathcal{J}}
\newcommand{\X}{X_{\Ha}}
\newcommand{\wind}{\operatorname{wind}}
\newcommand{\A}{\mathbf{A}}

\newcommand{\Ha}{\mathcal{H}}

\def\br#1{\left\{#1\right\}}
\def\bp#1{\left(#1\right)}

\def\fl#1{\lfloor#1\rfloor}
\def\ceil#1{\lceil#1\rceil}
\def\abs#1{\left|#1\right|}

\def\R{\mathbb{R}}
\def\Z{\mathbb{Z}}
\def\C{\mathbb{C}}
\def\N{\mathbb{N}}

\DeclareFontFamily{U}{mathx}{}
\DeclareFontShape{U}{mathx}{m}{n}{<-> mathx10}{}
\DeclareSymbolFont{mathx}{U}{mathx}{m}{n}
\DeclareMathAccent{\widecheck}{0}{mathx}{"71}

\usepackage{graphicx}
\usepackage{epstopdf}

\AppendGraphicsExtensions{.gif}

\begin{document}

\title[The symplectic geometry of the three-body problem]{The symplectic geometry of the three-body problem \\ \small On Floer theory, contact topology, symplectic dynamics and astrodynamics}

\author{Agustin Moreno}






\maketitle

\tableofcontents

\section{Introduction}

The current book grew out of a set of lecture notes \cite{M22}, that the author wrote in 2020 for a graduate mini-course aimed at graduate students in mathematics in UDELAR (Uruguay), remotely from the Mittag--Leffler Institute in Sweden. It is the result of 7 years of research on a classical conundrum which has grown very close to the author's heart, and which has been around since the times of Newton, Kepler, Poincar\'e, and so many other great scientists. Namely, this will be a book centered around the circular, restricted three-body problem, or CR3BP for short. This is the dynamical system obtained when three point-masses are left to interact with each other according to Newton's second law of gravitation, and moreover one of them is assumed negligible when compared to the other two (the \emph{primaries)}. The latter are further assumed to move in circles around their center of mass. Despite the centuries, and the simplifications made, this is still a poorly understood problem today, and unravelling its secrets is an astounding open challenge. The aim of this book is then, starting from basic material, move at quick strides towards some of the current research in this subject.

The CR3BP is not only interesting from a theoretical point of view (and indeed very large portions of the modern scientific discourse can be traced to this problem alone), but also from a practical perspective, due to its deep connections to astronomy and space exploration. Namely, the CR3BP is the most basic model approximating the motion of a spacecraft under the influence of a Planet--Moon system. This is a modern interpretation: unlike the times of Newton, when space travel was but a wild opium dream, in the current day and age, when mission proposals to remote regions of our expanding Universe are common currency, the CR3BP is one of the preeminent models used for spacecraft trajectory design. In the context of astrodynamics, the CR3BP is then the theoretical starting point supporting many high-fidelity (or \emph{ephemeris}) numerical studies which go into actual mission proposals. While finding trajectories that meet the requirements of an actual mission is a very complicated art, the families of periodic orbits found in the CR3BP, as well as the stable/unstable manifolds of some of them, can be used as building blocks for designing the desired trajectories, and to transfer between them.

The difference between this book and others is the approach, perspective, and scope. Firstly, the emphasis is on the \emph{spatial} case of the CR3BP (where the small mass moves in three-dimensional space), as opposed to the \emph{planar} case (where the small mass moves in the plane). While the planar problem has been extensively studied since the times of Poincaré, as it is a lower-dimensional problem and hence more tractable, the spatial problem is more physically meaningful and amenable to applications, e.g.\ to space exploration. The price to pay is the high dimension of the system (a six dimensional phase-space), which e.g.\ renders visualization harder, and imposes the need of global and higher-dimensional topological methods. Secondly, we have chosen to present the material from the vantage point and perspective of modern symplectic geometry. This is a currently very active field of research, which has been developed in earnest only in the last 40+ years, since the introduction of the notion of pseudo-holomorphic curves due to Gromov \cite{Gro85}, the development of Floer theory shortly after, the accompanying work of Giroux in contact topology \cite{Gir02}, and the invention of the framework of symplectic dynamics by Hofer \cite{BH}. 

Throughout the book, we will restrict our attention exclusively to the \emph{low-energy} and \emph{near-primary} dynamics, as this is the setup in which modern methods from symplectic and contact geometry can be made to bear on the problem (see Theorem \ref{thm:contact_type}). We have chosen to focus on intuition as opposed to formality, and have attempted to keep technicalities to a minimum, in the hope of getting to the research material as quickly as possible. The treatment will therefore be rather concise, adding references where the details here omitted can be consulted. While the intended audience is mostly pure mathematicians (graduate students and researchers), the second part of this book may be of interest to applied scientists with interest in Hamiltonian systems and bifurcations of periodic orbits, e.g.\ engineers. 

\subsection{Organization of the book} 

The book is split into two inter-related but complementary parts. Part I deals with the purely theoretical aspects of the problem, whereas Part II deals with the aspects that point towards applications (which we call the ``practical'' aspects, although some of the engineers that the author works with would also call them ``theoretical''). Part I is heavily based on the author's collaboration with Otto van Koert \cite{MvK20a,MvK20b}, on the author's paper \cite{M20}, on the author's paper with Arthur Limoge \cite{LM25}, on the author's conversations with Rohil Prasad \cite{Pr24}, on the author's work with Bahar Acu \cite{AM18}, and with Francesco Ruscelli \cite{MR23}. Part II draws mostly from the author's collaboration with Urs Frauenfelder, Otto van Koert, Dayung Koh, Bhanu Kumar and Cengiz Aydin \cite{FM,FMb,FKM,AFvKKM}, and is complemented with numerical work carried out by Otto van Koert, Dayung Koh, Cengiz Aydin, and Bhanu Kumar.

\medskip

\textbf{Part I: theoretical aspects.} Chapter \ref{ch:basic_notions} introduces the basic notions from symplectic and Hamiltonian dynamics, and its odd-dimensional counterpart contact geometry and Reeb dynamics. These are the geometries underlying classical mechanics, and the rest of the book will be expressed in this language. 

Chapter \ref{ch:celestial_mechanics} discusses the main problems from celestial mechanics that we will be interested in, from the very general $n$-body problem to the more tractable CR3BP, which is the main focus of the book, as well as limit cases like Hill's lunar problem and the rotating Kepler problem (RKP). We also discuss collision regularization, a classical mathematical artifact by which binary collisions between bodies may be continued.

In Chapter \ref{ch:open_books_and_dynamics}, we discuss the notion of open book decompositions from a topological and dynamical point of view, in particular introducing the notion of a \emph{global hypersurface of section} and touching upon Giroux correspondence, and discussing the main examples which will later appear in the CR3BP. We also include three digressions aimed at illustrating the role of open books in contact and symplectic topology. 

In Chapter \ref{ch:contact_geometry_in_the_CR3BP} we start in earnest with the modern approach to the CR3BP. After giving a historical account (as the author's bias understands them), we arrive at the advent of the modern methods of contact geometry in the CR3BP, whose starting point is Theorem \ref{thm:contact_type} from \cite{AFvKP,CJK18}. This chapter is based on the collaboration of the author with Otto van Koert \cite{MvK20a}, and with Bahar Acu \cite{AM}. The main points are:

\medskip

\begin{itemize}
    \item \textbf{(Open books in the CR3BP)} Existence of adapted open book decompositions for the spatial CR3BP in the low-energy range (Theorem \ref{thm:openbooks});
    \medskip
    \item \textbf{(Hamiltonian return maps)} Existence of Hamiltonian return maps reducing the continuous spatial dynamics to a discrete dynamics in dimension $4$ (Theorem \ref{thm:returnmap});
    \medskip
    \item \textbf{(Iterated picture)} Introducing the structure of an \emph{iterated planar} contact manifold on the low-energy energy levels sets (Theorem \ref{thm:IP});
    \medskip
    \item \textbf{(Rotating Kepler problem)} An explicit study of the return map in the integrable limit case of the RKP (Theorem \ref{thm:integrablecase});    
\end{itemize}
We also give a digression addressing the technicality that the symplectic form degenerates at the boundary of a global hypersurface of section, by defining the notion of a \emph{degenerate} Liouville domain (phenomenon which also arises in the setting of \emph{billiards}).

In Chapter \ref{ch:Floer_homology}, we give an overview of different flavors of Floer homology which we will need (Hamiltonian, Lagrangian, wrapped, local). It is meant as a reference chapter, although we will not attempt to provide proofs, as this is by now a standard subject in symplectic geometry. We will focus on basic definitions and uses rather than rigor, and give references where appropriate. 

Chapter \ref{ch:Ham_twist_maps} deals with the fixed point theory of what we call \emph{Hamiltonian twist maps}, based on the collaboration of the author with Otto van Koert \cite{MvK20b}. Following Poincaré's approach to the problem of finding periodic orbits in the planar problem, once we found the global section, we wish to prove an abstract fixed-point theorem for the return map. This chapter addresses this problem, the main points being:

\medskip

\begin{itemize}
    \item \textbf{(A generalized Poincar\'e--Birkhoff theorem)} A fixed-point theorem for Hamiltonian twist maps generalizing the classical Poincar\'e--Birkhoff theorem (Theorem \ref{thm:GBP}), aimed at the existence problem for periodic orbits;

    \medskip
    
    \item \textbf{(A relative Poincar\'e--Birkhoff theorem)} A fixed-point theorem (Theorem \ref{thm:chords}) aimed at the existence problem for chords between Lagrangians.
\end{itemize}

The relative version is inspired by the observation that the \emph{collision locus} $L$ is a Lagrangian in the global hypersurface of section given by Theorem \ref{thm:openbooks}, and chords in $L$ (i.e.\ points in $L\cap f(L)$ with $f$ the return map) correspond to consecutive collision orbits, i.e.\ the small mass collides with a primary once, and another time in the future. While this only make sense through the artifact of regularization, these orbits may be perturbed to actual orbits which pass close to the primaries, and therefore may be used as gravity assists (or \emph{flybys}) used to reach another target. We should emphasize that the twist condition as we introduced it suffers from several shortcomings (see Remark \ref{rk:twist_condition}), and adaptations of the above fixed-point theorems will likely be needed before applying them to the CR3BP. With this in mind, we also include three digressions: the first one explains how a given Hamiltonian twist map indeed arises as the return map for some adapted flow; the second one gives alternative definitions of the twist condition, which might be more adapted to the setup of the CR3BP but for which no fixed-point theorem is apparent; and the last one discusses an example of Morrison \cite{Morr82} of a Hamiltonian map on the ball with no interior fixed points, as well as the outlook concerning the study of Hamiltonian maps on Liouville domains.

Chapter \ref{ch:symp_dynamics} gives an basic introduction to \emph{symplectic dynamics}, a framework introduced by Hofer in order to address old problems but to also ask new questions, at the interface of symplectic geometry and dynamical systems. The exposition will include the basics of the theory of pseudo-holomorphic curves in symplectizations, Hofer--Wysocki--Zehnder's groundbreaking paper \cite{HWZ98}, and the (still work in progress) Siefring intersection theory in higher-dimensions. This discussion is aimed at making the author's paper \cite{M20} accessible, which fits into the scope of symplectic dynamics, and associates to the (low-energy, near-primary) spatial dynamics of the CR3BP a dynamics on $S^3$. 

But most importantly, the last section of this chapter (which also falls under the umbrella of symplectic dynamics) is a very brief introduction to the theory of \emph{feral curves} as introduced by Fish--Hofer \cite{FH23} and enhanced by Prasad \cite{Pr24}. In this section, we give a completely novel application to the \emph{planar} CR3BP, i.e.\ a rigorous proof of the following statement:

\begin{itemize}
    \item (\textbf{Existence of invariant subsets with dense union for the planar problem}) We prove existence (Theorem \ref{thm:goal}), for any fixed mass ratio and any fixed energy below $H(L_2)$ (i.e.\ below the second Lagrange point), for the unregularized and near-primary dynamics, of infinitely many proper, distinct, either closed or properly embedded, invariant subsets whose union is dense in the corresponding level set. 
\end{itemize}

In other words, this is a density result (and a direct application of the theory of feral curves) stating that the planar CR3BP dynamics is enormously rich. The energy constraint is due to the fact that the level sets can be regularized to be compact (near the primaries) only up to $L_2$. As far as the author understands, given the fact that this result is completely non-perturbative and relies on very novel machinery, \emph{there is no other rigorous proof of an existence result of these characteristics available in the literature.}  

\medskip

\textbf{Part II: practical aspects.} The second part of the book deals with material which is closer to applications. Chapter \ref{ch:symp_data_analysis} introduces a ``symplectic toolkit'' designed to study periodic orbits, their bifurcations in families, and their stability, with emphasis on symmetric orbits. The basic notions are the \emph{B-signs} \cite{FM}, the \emph{GIT-sequence} \cite{FM}, the CZ-indices, and the Floer numerical invariants (defined as the Euler characteristics of various local Floer homology groups). 

Chapter \ref{ch:num_work} contains numerical work, namely:

\medskip

\begin{itemize}
    \item \textbf{(Bifurcation graphs)} bifurcation graphs for various systems of interest (Hill's lunar problem, Saturn-Enceladus, Jupiter-Europa, Earth-Moon), produced by Cengiz Aydin; and

\medskip
    
    \item \textbf{(GIT plots)} numerical plots in the GIT sequence, produced by Dayung Koh; and

\medskip
    

    
    \item \textbf{(Earth--Moon system :Halo orbits)} We summarize a numerical investigation using the symplectic toolkit for periodic orbit families in the full Earth–Moon CR3BP. Near the Moon, prograde, retrograde, and Halo orbits are considered, discovering previously-unknown orbit families linking them together through bifurcations and singularities, confirming a 1968 conjecture of Broucke. This is based on work by Bhanu Kumar \cite{KM25}. New families connecting to the Halo orbits are found and described (starting from planar orbits studied by Broucke \cite{Br68}), with plots produced by Cengiz Aydin and Bhanu Kumar. 
\end{itemize}

\medskip


\subsection{Acknowledgments} This book draws heavily from my ongoing collaboration with Otto van Koert, Urs Frauenfelder, Dayung Koh, Cengiz Aydin, and Bhanu Kumar. Much of what appears in these pages is due to their insights, and so I am very grateful to them for the work that they have diligently put into what have quickly become very fruitful years of interactions. Let us hope for more to come. 

I am grateful to several people from whom I learned so much over the years. To name a few, in no specific order: Helmut Hofer, Dennis Sullivan, Dan Scheeres, Chris Wendl, Kai Cieliebak, Peter Sarnak, Richard Montgomery, Umberto Hryniewicz, Peter Albers, Sergei Tabachnikov, Lei Zhao, Connor Jackman, Jo Nelson, Julian Chaidez, Sobhan Seyfaddini, Ed Belbruno, Michael Hutchings, Vini Ramos, Janko Latschev, Gabriel Paternain, Georgios Dimitroglou Rizell, Rohil Prasad, Richard Siefring, Ezequiel Maderna, Alejandro Passeggi, Rafael Potrie, Bhanu Kumar, and my students Arthur Limoge, Favio Pir\'an, Francesco Ruscelli and Aidan Latona.

I would like to thank the warm hospitality of the Institute of Advanced Study in Princeton, where several of the ideas in this book where brewed while I was a member, as well as the Mittag--Leffler Institute in Sweden, where the lecture notes in which this book is based on where first conceived, while I was a fellow.

The author is supported by the Sonderforschungsbereich TRR 191 Symplectic Structures in Geometry, Algebra and Dynamics, funded by the DFG (Projektnummer 281071066 – TRR 191), by the DFG under Germany's Excellence Strategy EXC 2181/1 - 390900948 (the Heidelberg STRUCTURES Excellence Cluster), and by the Air Force Office of Scientific Research under award number FA8655-24-1-7012.

\part{Theoretical aspects}

\chapter{Basic notions}\label{ch:basic_notions}

This chapter is devoted to the basic concepts underlying the general principles of classical mechanics. In particular, we will focus on the modern language of symplectic and contact geometry, in which we will express the rest of the book.

\smallskip

\section{Symplectic geometry and Hamiltonian dynamics} 

\subsection{Symplectic geometry} Roughly speaking, symplectic geometry is the geometry of phase-space (where one keeps track of position and velocities of classical particles, and so it is a theory in even dimensions). Formally, a \emph{symplectic manifold} is a pair $(M,\omega)$, where $M$ is a smooth manifold with $\dim(M)=2n$ even, and $\omega \in \Omega^2(M)$ is a two-form (the \emph{symplectic form}) satisfying:
\begin{itemize}
    \item\textbf{(closedness)} $d\omega=0$;
    \item\textbf{(non-degeneracy)} $\omega^n=\omega\wedge\dots\wedge \omega \in \Omega^{2n}(M)$ is nowhere-vanishing, and hence a volume form. Equivalently, the map
    $$
    \mathfrak{X}(M)\rightarrow \Omega^1(M)
    $$
    $$
    X\mapsto i_X\omega=\omega(X,\cdot)
    $$
    is a linear isomorphism, where $\mathfrak{X}(M)$ denotes the space of smooth vector fields on $M$.
\end{itemize}

Note that symplectic manifolds are always orientable. We assume that $M$ is always oriented by the orientation induced by the symplectic form.

\begin{example}(From classical mechanics).$\;$
\begin{itemize}
    \item\textbf{(Phase-space)} $(\mathbb{R}^{2n},\omega_{std})$, where, writing $(q,p)\in \mathbb{R}^{2n}=\mathbb{R}^n\oplus \mathbb{R}^n$ ($q=$position, $p=$momenta), we have
    $$
    \omega_{std}=-d\lambda_{std}=dq\wedge dp,
    $$
    where $\lambda_{std}=pdq$ is the standard \emph{Liouville form}. Here we use the shorthand notation $dq\wedge dp=\sum_{i=1}^ndq_i\wedge dp_i$ and similarly $pdq=\sum_{i=1}^n p_idq_i$.
    \medskip
    
    \item\textbf{(cotangent bundles)} $(T^*Q,\omega_{std})$, where $Q$ is a closed $n$-manifold, and $\omega_{std}$ is defined invariantly as
    $$
    \omega_{std}=-d\lambda_{std},
    $$
    with
    $$
    (\lambda_{std})_{(q,p)}(\eta)=p(d_{(q,p)}\pi(\eta)),
    $$
    also called the standard Liouville form. Here, $q$ is a point in the base, and $p$ a covector in $T_qQ^*$, and
    $$
    \pi:T^*Q\rightarrow Q
    $$
    is the natural projection to the base. Note that phase-space corresponds to the case $Q=\mathbb{R}^{n}$.

    If $Q$ is equipped with a Riemannian metric, we denote the \emph{co-disk} bundle $\mathbb D^*Q=\{(q,p)\in T^*Q: \vert p \vert \leq 1\}$, endowed with the restriction of $\omega_{std}$, and which has boundary the \emph{unit cotangent bundle} $S^*Q=\{(q,p)\in T^*Q: \vert p \vert = 1\}$.
    
\end{itemize}    
    
\end{example}

\begin{example}(From complex algebraic/Kähler geometry).$\;$
\begin{itemize}
    \item(Projective varieties) Complex projective space $\mathbb{C}P^n$ admits a natural symplectic form, called the \emph{Fubini-Study} form $\omega_{FS}$, defined as follows. Let
    $$K:\mathbb{C}^{n}\rightarrow \mathbb{R}$$
    $$K(z)=\mbox{log}\left(1+\sum_{i=1}^{n}\vert z_i\vert^2\right).$$
    In homogeonous coordinates $(\zeta_0:\dots:\zeta_n)$ for $\mathbb{C}P^n$, let $U_\alpha=\{(\zeta_0:\dots:\zeta_n):\zeta_\alpha \neq 0\}$ and $$\varphi_\alpha:  U_\alpha \rightarrow \mathbb{C}^n,$$ $$\varphi_\alpha(\zeta_0:\dots:\zeta_n)=\left(\frac{\zeta_0}{\zeta_i},\dots,\frac{\zeta_{i-1}}{\zeta_i},\frac{\zeta_{i+1}}{\zeta_i},\dots,\frac{\zeta_{n}}{\zeta_i}\right)=(z_1^\alpha,\dots,z_n^\alpha)$$ be the standard affine chart around $(0:\dots:1:\dots:0)$. Let $K_\alpha=K \circ \varphi_\alpha$, and define
    $$
    \omega_\alpha=\sqrt{-1}\partial \overline{\partial} K_\alpha=\sum_{i,j=1}^n h_{ij}(z^\alpha)dz^\alpha_i\wedge d\overline{z}_j^\alpha.
    $$
    Here, one computes
    $$
    h_{ij}(z^\alpha)=\frac{\delta_{ij}\left(1+\sum_{i=1}^{n}\vert z^\alpha_i\vert^2\right)-z_i^\alpha \overline{z}_j^\alpha}{\left(1+\sum_{i=1}^{n}\vert z^\alpha_i\vert^2\right)^2}
    $$
    One checks that on overlaps $U_\alpha \cap U_\beta$, we have $\omega_\alpha=\omega_\beta$, and so we get a well-defined global $\omega_{FS}$ so that $\omega_{FS}\vert_{U_\alpha}=\omega_\alpha$. 
    The $K_\alpha$ are what is called a local K\"ahler potential (or plurisubharmonic function) for the Fubini-Study form. Every algebraic/analytic projective variety inherits a symplectic form via restriction of the ambient Fubini-study form.
    
    \medskip
    
    \item\textbf{(Affine varieties: Stein manifolds)} The standard complex affine space $\mathbb{C}^n$ carries the standard symplectic form via the identification $\mathbb{C}^n=\mathbb{R}^{2n}$, which in complex notation is
    $$
    \omega_{std}=\frac{\sqrt{-1}}{2}\sum_{i=1}^n dz_i\wedge d\overline{z}_j=:\frac{\sqrt{-1}}{2}dz\wedge d\overline{z}=-d\lambda_{std}
    $$
    with $\lambda_{std}=\frac{\sqrt{-1}}{4}(zd\overline{z}-\overline{z}dz)$. This admits the standard plurisubharmonic function
    $$
    f_{std}(z)=\vert z \vert^2,
    $$
    i.e.\ $\omega_{std}=\sqrt{-1}\partial \overline{\partial} f_{std}$. This function is exhausting (i.e.\ $\{z:f(z)\leq c\}$ is compact for every $c \in \mathbb{R}$), and is a Morse function (with a unique critical point at the origin). 
    
    By analogy to the projective case, a Stein manifold $X$ is a properly embedded complex submanifold of $\mathbb{C}^n$, endowed with the restriction of the standard symplectic form, the standard complex structure $i$, and the standard plurisubharmonic function. One may further assume (after a small perturbation) that $f_{std}$ defines a Morse function on $X$. 
    \end{itemize}
    \end{example}
    
    The above examples (projective and affine) are all instances of K\"ahler manifolds, i.e.\ the symplectic form is suitably compatible with an integrable complex structure, and with a Riemannian metric. One way to obtain Stein manifolds from projective varieties is to remove a collection of generic hyperplane sections, i.e.\ the intersection of the variety with the zero sets of generic homogeneous polynomials of degree $1$. The Liouville form (i.e.\ the primitive of the resulting symplectic form), depends on the number of sections.
\medskip

A general important feature of symplectic manifolds (or rather the reason for their existence) is that they are locally modelled on phase-space:

\begin{thm}[\textbf{Darboux's theorem for symplectic manifolds}] If $p\in (M,\omega)$ is an arbitrary point in a symplectic manifold, we can find local charts centered at $p$, so that $(M,\omega)$ is isomorphic to standard phase-space $(\mathbb{R}^{2n},\omega_{std})$ in this local chart.
\end{thm}
The notion of isomorphism we use above is the obvious one: two symplectic manifolds $(M_1,\omega_1)$ and $(M_2,\omega_2)$ are \emph{symplectomorphic} if there exists a diffeomorphism $f: M_1\rightarrow M_2$ satisfying $f^*\omega_2=\omega_1$. In particular, a symplectomorphism preserves volume, i.e.\ $f^*\omega_2^n=\omega_1^n$. Darboux's theorem is usually interpreted as saying that, unlike in Riemannian geometry where the curvature is a local isometry invariant, there are no local invariants for symplectic manifolds (they locally all look the same). A source of symplectomorphisms on cotangent bundles are the \emph{physical transformations}, i.e.\ those induced by a diffeomorphism on the base $f: Q_1\rightarrow Q_2$, given by
$$
f_*:T^*Q_1\rightarrow T^*Q_2,
$$
$$
f_*(q,p)=(f(q),(d^*_qf)^{-1}(p)).
$$
An important class of submanifolds of a given symplectic manifold consists of the \emph{Lagrangian} submanifolds, i.e.\ half-dimensional manifolds $L^n\subset M^{2n}$ satisfying $\omega\vert_L\equiv 0$. Standard examples of such are the zero section $Q\subset T^*Q$, the cotangent fiber $T_q^*Q\subset T^*Q$, the graph of a closed $1$-form $\alpha:Q\rightarrow T^*Q$, and $\mathbb RP^n\subset \mathbb CP^n$. More generally, a submanifold $N\subset M$ is \emph{isotropic} if $\omega\vert_N\equiv 0$, i.e.\ $TN\subset TN^\omega=\{v \in TM: \omega(v,w)=0 \mbox{ for } w\in TN\}$ (the symplectic complement). It is \emph{co-isotropic} if $TN^\omega\subset TN$. Lagrangians correspond to those which are co-isotropic and isotropic, i.e.\ $TL=TL^\omega$. A simple lemma from linear algebra implies that the dimension of an isotropic submanifold is at most $n=\dim(M)/2$, whereas the dimension of a co-isotropic is at least $n=\dim(M)/2$.

\subsection{Hamiltonian dynamics.} From a dynamical perspective, symplectic manifolds are the natural geometric space where one can study Hamiltonian dynamics, via the \emph{Hamiltonian formalism}. On a cotangent bundle $T^*Q$, the idea is to model the motion of a particle moving along the manifold $Q$, subject to the principle of least action associated to a given physical problem. 

In general, we start with a symplectic manifold $(M,\omega)$, and a \emph{Hamiltonian} $H:M\rightarrow \mathbb{R}$, which is simply a function (which we assume $C^1$, say), thought of as the \emph{energy} function of the mechanical system. The symplectic form implicitly defines a vector field $X_H \in \mathfrak{X}(M)$ (the \emph{Hamiltonian vector field} or \emph{Hamiltonian gradient} of $H$) via the equation
$$
i_{X_H}\omega=dH.
$$
Note that this uniquely defines $X_H$ due to non-degeneracy of $\omega$. The above equation is the global, invariant version for the following. 

\begin{example}\textbf{(Fundamental example: Hamilton equations)} Whenever $(M,\omega)=(\mathbb{R}^{2n},\omega_{std})$, we have
$$
X_H=\left(\frac{\partial H}{\partial p},-\frac{\partial H}{\partial q}\right)=\frac{\partial H}{\partial p} \partial_q -\frac{\partial H}{\partial q}\partial_p.
$$
In other words, a solution $x(t)=(q(t),p(t))$ to the ODE $\dot x(t)=X_H(x(t))$ is precisely a solution to the Hamilton equations
$$
\left\{\begin{array}{cc}
    \dot q=& \frac{\partial H}{\partial p}  \\
    \dot p= & -\frac{\partial H}{\partial q} 
\end{array}\right.
$$
\end{example}

By Darboux's theorem, we see that, locally, solutions to the Hamiltonian flow are solutions to the above.

\medskip

More invariantly, we consider the Hamiltonian flow $\phi_t^H: M \rightarrow M$ generated by $H$, i.e.\ the unique solution to the equations
$$
\phi_0^H=\mbox{id},\; \frac{d}{dt}\phi_t^H=X_H \circ \phi_t^H.
$$
This flow can be thought of as a symmetry of the symplectic manifold, since it preserves the symplectic form:
$$
\frac{d}{dt}(\phi_t^H)^*\omega=(\phi_t^H)^*\mathcal{L}_{X_H}\omega=(\phi_t^H)^*(i_{X_H}d\omega+di_{X_H}\omega)=(\phi_t^H)^*d^2H=0,
$$
and so $(\phi_t^H)^*\omega=(\phi_0^H)^*\omega=\omega$ for every $t$. A symplectomorphism $f:(M,\omega)\rightarrow (M,\omega)$ is called \emph{Hamiltonian} whenever $f=\phi_H^1$ is the time-$1$ map of a Hamiltonian flow. Hamiltonian maps then preserve volume (which is a way of stating Liouville's theorem from classical mechanics).

\begin{example}\textbf{(Simple harmonic oscillator)} The simple harmonic oscillator is given by the Hamilton flow of $H:\mathbb{R}^2\rightarrow \mathbb{R}$, $H(q,p)=\frac{p^2}{2m}+\frac{m\omega^2x^2}{2},$ where $\omega=\sqrt{\frac{k}{m}}$ is the angular frequency, $k$ is the spring constant, $m$ is the mass of a classical particle with position $x$ and momenta $p$. 
\end{example}

\begin{remark}
The Hamiltonian usually also depends on time. We have assumed for simplicity that it does not, i.e.\ it is autonomous. We will see that this will hold for the simplified versions of the three body problem we will consider, i.e.\ the restricted case.
\end{remark}

In the above symplectic formalism, it is a fairly straightforward matter to write down the fundamental conservation of energy principle (in the autonomous case):

\begin{thm}[\textbf{Conservation of energy}] Assume $H$ is autonomous. Then
$$
dH(X_H)=0.
$$
In other words, the level sets $H^{-1}(c)$ are invariant under the Hamiltonian flow.
\end{thm}

This is also usually written down using the \emph{Poisson bracket} $\{F,G\}=dF(X_G)=-dG(X_H)$ as
$$
\{H,H\}=0,
$$
which is another way of saying that $H$ is preserved under the Hamiltonian flow of itself, or that $H$ is a conserved quantity (or integral) of the motion. The proof fits in one line:
$$
dH(X_H)=i_{X_H}\omega(X_H)=\omega(X_H,X_H)=0,
$$
since $\omega$ is skew-symmetric.

\subsection{Periodic orbits and monodromy.} 

Given a $2n$-dimensional symplectic manifold $(M,\omega)$ and a Hamiltonian $H:M\rightarrow \mathbb R$, with Hamiltonian flow $\phi^H_t: M\rightarrow M$, a \emph{periodic orbit} $x \in C^\infty(S^1,M)$ is a solution of the ODE
$$\partial_t x(t)=T \cdot X_H(x(t)), \quad t \in S^1,$$
where $T$ is a positive real number, the \emph{period} of the periodic orbit. Equivalently, $$x(t)=\phi^{T\cdot\; t}_H(x(0)).$$
Denoting $x_0=x(0)$, the differential 
$$M_{x}:=d_{x_0}\phi^H_T: T_{x_0}M \to T_{x_0}M$$
is a linear symplectic map of the symplectic vector space
$(T_{x_0}M,\omega_{x_0})$, i.e.\ $$M_{x}^* \omega_{x_0}=\omega_{x_0}.$$ The map $M_x$ is called the \emph{monodromy}. After choosing a basis of $(T_{x_0}M,\omega_{x_0})$ in which the symplectic form is standard (this is called a \emph{symplectic basis}), $M_{x}$ is becomes a symplectic $2n\times 2n$--matrix, i.e.\ it satisfies the equation
$$
M_{x}^t J M_{x}=J.
$$
Here, $J=\left(\begin{array}{cc}
    0 & -\mathds 1 \\
   \mathds 1  & 0
\end{array}\right)$ is the standard almost complex structure on $\mathbb R^{2n}$, satisfying $J^2=-\mathds 1$. We denote by $$Sp(2n)=\{ M \in M_{2n\times 2n}(\mathbb R): M^t J M=J\}$$ the space of symplectic matrices (the \emph{symplectic group}). Note that the choice of a different point $x_0$ along the orbit $x$ changes the monodromy up to \emph{symplectic} conjugation, i.e.\ up to conjugating with a symplectic matrix. Therefore $M_x$ is strictly speaking an element of $Sp(2n)/Sp(2n)$, where $Sp(2n)$ acts on itself by conjugation.

\medskip

It is an easy exercise to show that if $M$ is symplectic, and $\mu \in \mathbb C$ is an eigenvalue of $M$, then so are $\overline{\mu},1/\mu, 1/\overline{\mu}$. Then we have the following possibilities for $\mu$:
\medskip
\begin{itemize}
    \item\textbf{($\mathcal{P}$, parabolic)} $\mu=\pm 1$, in which case it has even multiplicity;
    \medskip
    \item \textbf{($\mathcal{E}$, elliptic)} $\vert \mu \vert=1$, in which case it comes as an elliptic pair $\mu,\overline{\mu}=1/\mu$;
    \medskip
    \item \textbf{($\mathcal{H}^+$, positive hyperbolic)} $\mu \in \mathbb R$, $\mu>0$, $\mu \neq 1$, in which case both $\mu,1/\mu$ are positive;
    \item \textbf{($\mathcal{H}^-$, negative hyperbolic)} $\mu \in \mathbb R$, $\mu<0$, $\mu \neq -1$, in which case both $\mu,1/\mu$ are negative;
    \item\textbf{($\mathcal{N}$, complex quadruple)} $\mu\notin S^1\cup \;\mathbb R$, in which case it comes in a quadruple $\mu,\overline{\mu},1/\mu, 1/\overline{\mu}$.
\end{itemize}

\medskip

Note that if $H$ is time-independent then $1$ appears twice as a \emph{trivial} eigenvalue of $M_x$, as $X_H$ is a corresponding eigenvector of $M_x$, and the spectrum of $M_x$ satisfies the above symmetries. We can ignore these if we consider the \emph{reduced} monodromy matrix $M_x^{red}\in Sp(2n-2)$, obtained by fixing the energy and dropping the direction of the flow, i.e.\ 
$$
M_x^{red}: T_{x_0}H^{-1}(c)/\langle X_H(x_0)\rangle\rightarrow T_{x_0}H^{-1}(c)/\langle X_H(x_0)\rangle.
$$
This map preserves a symplectic form $\overline{\omega}$ on $T_{x_0}H^{-1}(c)/\langle X_H(x_0)\rangle$ defined by \emph{symplectic reduction} (i.e.\ satisfying $\pi^*\overline{\omega}=i^*\omega$, where $i: H^{-1}(c)\hookrightarrow M$ is the inclusion, and $\pi: T_{x_0}H^{-1}(c) \rightarrow T_{x_0}H^{-1}(c)/\langle X_H(x_0)\rangle$ is the quotient map).

\begin{definition} $\;$   
\begin{itemize}
    \item A \emph{Floquet multiplier} of $x$ is an eigenvalue of $M_x$, which is not one of the trivial eigenvalues (i.e.\ an eigenvalue of $M_x^{red}$).
    \item An orbit is \emph{non-degenerate} if $1$ does not appear among its Floquet multipliers. 
    \item An orbit is \emph{stable} if all its Floquet multipliers are semi-simple and lie in the unit circle.
\end{itemize}
\end{definition} 

\subsection{Symmetries} The role of symmetry in physics has been prominent since the work of Emmy Noether. We will be interested, in what follows, in $\mathbb Z_2$ symmetries, i.e.\ involutions.  

\medskip

An \emph{involution} is a map $\rho:(M,\omega)\rightarrow (M,\omega)$ satisfying $\rho^2=id$, and it is symplectic or anti-symplectic if $\rho^*\omega=\pm\omega$ respectively. Its \emph{fixed-point locus} is $\mathrm{Fix}(\rho)=\{x\in M:\rho(x)=x\}$, which is a symplectic submanifold of $M$ in the symplectic case, and a Lagrangian submanifold of $M$ in the anti-symplectic case. An anti-symplectic or symplectic involution $\rho$ is a \emph{symmetry} of the system if $H\circ \rho=H.$ A periodic orbit $x$ is \emph{symmetric} with respect to an anti-symplectic involution $\rho$ if $\rho(x(-t))=x(t)$ for all $t$. The \emph{symmetric points} of the symmetric orbit $x$ are the two intersection points of $x$ with $\mathrm{Fix}(\rho)$, i.e.\
$$x\big(0\big),\,\,x\big(\tfrac{T}{2}\big) \in \mathrm{Fix}(\rho).$$
In particular, half of the symmetric periodic orbit is a Hamiltonian chord (i.e.\ trajectory) from
$\mathrm{Fix}(\rho)$ to itself. Hence we can think of a symmetric periodic orbit in two ways,
either as a closed string, or as an open string from the Lagrangian $\mathrm{Fix}(\rho)$ to itself.

The monodromy matrix of a symmetric orbit at a symmetric point is a \emph{Wonenburger} matrix, i.e.\ it satisfies
\begin{equation*}\label{symsymp}
M=M_{A,B,C}=\left(\begin{array}{cc}
A & B\\
C & A^t
\end{array}\right)\in Sp(2n),
\end{equation*}
where
\begin{equation}\label{eq:Wonenburger}
B=B^t,\quad C=C^t,\quad AB=BA^t,\quad
A^tC=CA,\quad A^2-BC=id,
\end{equation}
equations which ensure that $M$ is symplectic. The eigenvalues of $M$ are determined by those of the first block $A$ (see \cite{FM}):
\begin{itemize}
    \item If $\lambda$ is an eigenvalue of $M$ then its \emph{stability index} $a(\lambda)=\frac{1}{2}(\lambda + 1/\lambda)$ is an eigenvalue of $A$. 
    \item If $a$ is an eigenvalue of $A$ then $\lambda(a)=a+\sqrt{a^2-1}$ is an eigenvalue of $M$, for any choice of complex square root.
\end{itemize}

Note that in order to write the monodromy matrix in Wonenburger form, we implicitly chose a basis for $\mathrm{Fix}(\rho)$ at a symmetric point of the orbit (and extended it to a symplectic basis). A different choice of basis amounts to acting with an invertible matrix $R \in GL_n(\mathbb{R})$, via
\begin{equation*}\label{act}
R_*\big(A,B,C\big)=\Big(RAR^{-1},RBR^t,(R^t)^{-1}CR^{-1}\Big),
\end{equation*}
i.e.,\ $M_{A,B,C}$ is replaced by $M_{R_*(A,B,C)}$. We denote the space of Wonenburger matrices by
$$
Sp^\mathcal{I}(2n)=\{M_{A,B,C}: A,B,C \mbox{ satisfy } (\ref{eq:Wonenburger})\},
$$
which comes with the above action of $GL_n(\mathbb R)$. 

By a beautiful result of Wonenburger, every symplectic matrix $M\in Sp(2n)$ can be written as a product of two linear anti-symplectic involutions, i.e.\ $M=I_1 I_2$. From this, it is straightforward to derive the following fact (see \cite{FM}):

\begin{thm}\label{thm:Wonenburger}
    Every symplectic matrix $M\in Sp(2n)$ is symplectically conjugated to a Wonenburger matrix.
\end{thm}

In other words, the natural map
$$
Sp^\mathcal{I}(2n)/GL_n(\mathbb R) \rightarrow Sp(2n)/Sp(2n),
$$
$$
[M_{A,B,C}] \mapsto [M_{A,B,C}],
$$
is surjective.

In the presence of a symmetric periodic orbit, the above algebraic fact has a geometric interpretation: the monodromy matrix at each point of the orbit (a symplectic matrix) is symplectically conjugated via the linearized flow to the monodromy matrix at any of the symmetric points of the orbit (a Wonenburger matrix). The above discussion is the starting point for the \emph{GIT sequence} \cite{FM}, which will be discussed in Chapter \ref{ch:symp_data_analysis}.

\subsection{Monodromy splittings} In the presence of a symplectic symmetry, periodic orbits lying in the symplectic fixed-point locus have monodromy matrices which split into components. Namely, if $\sigma$ is a sympletic symmetry of the Hamiltonian $H$, and $x$ is a periodic orbits with $x_0=x(t_0)\in \mathrm{Fix}(\sigma)$ and period $T$, consider the splitting
$$
T_{x_0}M=E_1\oplus E_{-1}=T_{x_0}\mathrm{Fix}(\sigma)\oplus E_{-1}
$$
into $\pm 1$ eigenspaces of $d_{x_0}\sigma$, which are symplectically orthogonal. Since $\sigma$ commutes with the Hamiltonian flow, the monodromy $M_x$ leaves the splitting invariant, i.e.\ as a matrix it is of the form
$$
M_x=\left(\begin{array}{cc}
   M_p  & 0 \\
    0 & M_s
\end{array} \right)
$$
for symplectic matrices $M_p,M_s$. Moreover, reducing the matrix is also compatible with this splitting, i.e.\
$$
M^{red}_x=\left(\begin{array}{cc}
   M^{red}_p  & 0 \\
    0 & M_s
\end{array} \right),
$$
where $M_p^{red}$ is the reduction of $M_p$.

\subsection{Compatible almost complex structures} An \emph{almost complex structure} on an even dimensional manifold $M$ is $J\in \text{End}(TM)$ satisfying $J^2=-\mathds 1$. Given a symplectic form $\omega$, an almost complex structure $J$ is \emph{compatible} with $\omega$ if
\begin{itemize}
    \item $\omega$ is $J$-invariant, i.e.\ $\omega(J\cdot,J\cdot)=\omega(\cdot,\cdot)$; and
    \item $g=\omega(\cdot,J\cdot)$ is a Riemannian metric on $M$.
\end{itemize}

By a well-known result of Gromov, the space of almost complex structures compatible with a given symplectic form is non-empty and contractible (see e.g.\ \cite{MS17}).

\section{Contact geometry and Reeb dynamics} 

\subsection{Contact geometry.} Contact geometry is, roughly speaking, the odd-dimensional analogue of symplectic geometry, and arises on level sets of Hamiltonians satisfying a suitable convexity assumption (see Prop. \ref{prop:contacttype}). Formally, a \emph{(strict) contact manifold} is a pair $(X,\alpha)$, where $X$ is a smooth manifold with $\dim(X)=2n-1$ odd, and $\alpha \in \Omega^1(X)$ is a $1$-form (the \emph{contact form}) satisfying the \emph{contact} condition:
$$
\alpha\wedge d\alpha^{n-1}\neq 0 \mbox{ is nowhere-vanishing, and hence a volume form}.
$$
Contact manifolds are therefore orientable (see Remark \ref{rk:coorientable} below). The codimension-$1$ distribution $\xi=\ker \alpha\subset TM$ (a choice of hyperplane at each tangent space, varying smoothly with the point), is called the \emph{contact structure} or \emph{contact distribution}, and $(M,\xi)$ is a \emph{contact manifold}. 

\begin{example}$\;$
\begin{itemize}
    \item \textbf{(standard)} The standard contact form on $\mathbb{R}^{2n-1}=\mathbb{R}\oplus \mathbb{R}^{n-1}\oplus \mathbb{R}^{n-1}\ni (z,q,p)$ is
    $$
    \alpha_{std}=dz-pdq,
    $$
    where we again use the short-hand notation $pdq=\sum_{i=1}^np_idq_i$.
    \medskip
    \item \textbf{(First-jet bundles)} Given a manifold $Q$, its first-jet bundle $J^1(Q)\rightarrow Q$, by definition, has total space the collection of all possible first-derivatives of maps $f: Q\rightarrow \mathbb{R}$. The fiber over $q$ is as all possible tuples $(q,f(q),d_qf)$, and so $J^1(Q)\cong \mathbb{R}\times T^*Q$. It carries the natural contact form 
    $$
    \alpha=dz+\lambda_{std},
    $$
    where $z$ is the coordinate on the first factor, and $\lambda_{std}$ is the standard Liouville form on $T^*Q$; note that the standard contact form corresponds to the case $Q=\mathbb{R}^{n-1}$.
    \medskip
    \item \textbf{(contactization)} More generally: If $(M,\omega=d\lambda)$ is an exact symplectic manifold, then its \emph{contactization} is
    $$
    (\mathbb{R}\times M,dz+\lambda),
    $$
    where $z$ is the coordinate in the first factor.
\end{itemize}
\end{example}

The contact condition should be thought of as a \emph{maximally non-integrability} condition, as follows. Recall the following theorem from differential geometry:

\begin{thm}[\textbf{Frobenius' theorem}]
If $\alpha \wedge d\alpha\equiv 0$, then $\xi=\ker \alpha\subset TM$ is integrable. That is, there are codimension-$1$ submanifolds whose tangent space is $\xi$.
\end{thm}

The condition in Frobenius' theorem is equivalent to $d\alpha\vert_\xi\equiv 0$. The contact condition is the extreme opposite of the above: $d\alpha\vert_\xi>0$ is symplectic, i.e.\ non-degenerate. In fact:

\begin{prop}
If $Y\subset (X,\xi)$ is a submanifold of a $(2n-1)$-dimensional contact manifold so that $TY\subset \xi$ (i.e.\ $Y$ is \emph{isotropic}), then $\dim(Y)\leq n-1$. 
\end{prop}

The isotropic submanifolds of maximal dimension $n-1$ are called \emph{Legendrians}. The analogous theorem of Darboux in the contact category is the following.

\begin{thm}[\textbf{Darboux's theorem for contact manifolds}] If $p \in (X,\alpha)$ is an arbitrary point in a strict contact manifold, we can find a local chart $U\cong \mathbb{R}^{2n-1}$ centered at $p$, so that $\alpha\vert_U=\alpha_{std}$.
\end{thm}

\subsection{Reeb dynamics} Whereas a contact manifold is a geometric object, a \emph{strict} contact manifold is a dynamical one, as we shall see below. Note first that the choice of contact form for a contact structure $\xi$ on $X$ is not unique: if $\alpha$ is such a choice, then $\nu\alpha$ is also, for any smooth positive function $\nu: X\rightarrow \mathbb R$, $\nu>0$. This is in fact the only ambiguity.

Given a contact form $\alpha$, it defines an autonomous dynamical system on $X$, generated by the \emph{Reeb vector field} $R_\alpha\in \mathfrak{X}(X)$. This is defined implicitly via:
\begin{itemize}
    \item $i_{R_\alpha}d\alpha=0$;
    \item $\alpha(R_\alpha)=1$.
\end{itemize}
To understand the above, note that, since $d\alpha\vert_\xi$ is symplectic, the kernel of $d\alpha$ is the $1$-dimensional distribution $TX/\xi\subset TX$. This is trivialized (as a real line bundle) via a choice of contact form, which also gives it an orientation induced from the one on $M$. The Reeb vector field then lies in this $1$-dimensional distribution; the second condition normalizes it so that it points precisely in the positive direction with respect to the co-orientation.
We emphasize that the Reeb vector field depends significantly on the contact form, and not the contact structure; different choices give, in general, very different dynamical systems.

\begin{remark}\label{rk:coorientable}
There are also examples of contact manifolds which are not globally co-orientable (e.g.\ the space of contact elements); we will not be concerned with those.
\end{remark}

The Reeb flow $\varphi_t$ has the property that it preserves the geometry in a strict way, i.e. it is a \emph{strict contactomorphism}. This means that $\varphi_t^*\alpha=\alpha$, or in other words, the Reeb vector field generates a (strict) local symmetry of the (strict) contact manifold. This fact easily follows from the Cartan formula:
$$
\frac{d}{dt}\varphi_t^*\alpha=\varphi_t^*(di_{R_\alpha}\alpha+i_{R_\alpha}d\alpha)=\varphi_t^*(d(1)+0)=0,
$$
and so $\varphi_t^*\alpha=\varphi_0^*\alpha=\alpha.$ 

More generally, a (not necessarily strict) contactomorphism is a diffeomorphism $f$ such that $f^*(\xi)=\xi$, or $f^*\alpha=\nu \alpha$ for some strictly positive smooth function $\nu$.

\medskip

\subsection{The bridge} The fundamental relationship between symplectic and contact geometry lies in the following. If the symplectic form $\omega=d\lambda$ is exact (which can only happen if the symplectic manifold is open, by Stokes' theorem), then we have a \emph{Liouville} vector field $V$, defined implicitly via
$$
i_V\omega=\lambda,
$$
where we again use non-degeneracy of $\omega$. To understand this vector field, consider $\varphi_t$ the flow of $V$. The Cartan formula implies
$$
\frac{d}{dt}\varphi_t^*\omega=\varphi_t^*(di_V\omega+i_Vd\omega)=\varphi_t^*(d\lambda)=\varphi_t^*\omega,
$$
and so, integrating, we get
$$
\varphi_t^*\omega=e^t\omega.
$$
Taking the top wedge power of this equation: $\varphi_t^*\omega^n=e^{nt}\omega^n$, and we see that the symplectic volume grows exponentially along the flow of $V$, i.e.\ $\varphi_t$ is a \emph{symplectic dilation}.

Assume that $X\subset (M,\omega=d\lambda)$ is a co-oriented codimension-$1$ submanifold, and the Liouville vector field is positively transverse to $X$. Then we obtain a volume form on $X$ by contraction:
$$
0<i_V\omega^n\vert_X=ni_V\omega\wedge \omega^{n-1}\vert_X=n\lambda\wedge d\lambda^{n-1}\vert_X=n\alpha \wedge d\alpha^{n-1},
$$
where $\alpha=\lambda\vert_X$.
We have proved:
\begin{proposition}\label{prop:contacttype}
If $\omega=d\lambda$, and the associated Liouville vector field $V$ is positively transverse to $X$, then
$(X,\alpha=\lambda\vert_X=i_V\omega\vert_X)$ is a strict contact manifold.
\end{proposition}

A hypersurface $X$ as in the above proposition is then called \emph{contact-type}. The most relevant example to keep in mind, is when $X=H^{-1}(c)$ is the level set of a Hamiltonian (in fact, locally this is always the case). In this situation:

\begin{proposition} If $X=H^{-1}(c)$ is contact-type, then the Reeb dynamics on $X$ is a positive reparametrization of the Hamiltonian dynamics of $H$.
\end{proposition}

\begin{figure}
    \centering
    \includegraphics[width=0.7\linewidth]{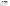}
    \caption{The fundamental relationship between contact and symplectic geometry is summarized here.}
    \label{fig:contactsymp}
\end{figure}

In other words, \emph{Reeb dynamics on contact-type Hamiltonian level sets is dynamically equivalent to Hamiltonian dynamics.} See Figure \ref{fig:contactsymp} for an abstract sketch.

\begin{example}$\;$
\begin{itemize}
    \item\textbf{(star-shaped domains)} Assume that $X\subset \mathbb{R}^{2n}$ is \emph{star-shaped}, i.e.\ it bounds a compact domain $D$ containing the origin, and the radial vector field $V=\frac{1}{2}(q\partial_q+p\partial_p)=r\partial_r$ is positively transverse to $X$ (with the boundary orientation). Since $V$ is precisely the Liouville vector field associated to $\lambda_0=\frac{1}{2}(qdp-pdq)$, every star-shaped domain is contact-type.
    \medskip
    \item \textbf{(standard contact form on $S^3$)} As a particular case, let $$S^3=\{z\in \mathbb{R}^4:\vert z\vert=1\}\subset \mathbb{R}^4$$ be the round $3$-sphere. Then $S^3=H^{-1}(1/2)$, where $H:\mathbb{R}^4\rightarrow \mathbb{R}$, $H(z)=\frac{1}{2}\vert z\vert^2$, and it is star-shaped. Writing $z=(z_1,z_2)=(x_1,y_1,x_2,y_2)$, the radial vector field 
    $$
    V=\frac{1}{2}r\partial_r=\frac{1}{2}(x_1\partial_{x_1}+y_1\partial_{y_1}+x_2\partial_{x_2}+y_2\partial_{y_2})
    $$
    is Liouville and induces the contact form 
    $$
    \alpha=i_V\omega_{std}\vert_{S^3}=\lambda_{std}\vert_{S^3}=\frac{1}{2}(x_1dy_1-y_1dx_1+x_2dy_2-y_2dx_2)\vert_{S^3}
    $$
    on $S^3$ whose Reeb vector field is
    $$
    R_\alpha=2(x_1\partial_{y_1}-y_1\partial_{x_1}+x_2\partial_{y_2}-y_2\partial_{x_2}).
    $$
    Its Reeb flow is, in complex coordinates, $\varphi_t(z_1,z_2)= e^{2\pi it}(z_1,z_2)$, whose orbits are precisely the fibers of the Hopf fibration $S^3\ni (z_1,z_2) \mapsto [z_1:z_2]\in \mathbb{C}P^1$. In particular, this flow is periodic, and all orbits have the same period. 
    
    The Hopf fibration $\pi: S^3\rightarrow S^2=\mathbb{C}P^1$ is an example of what is usually called a \emph{prequantization bundle}, i.e.\ the contact form $\alpha$ is a connection form whose curvature form on the base is symplectic. In other words, $d\alpha=i\pi^*\omega_{FS}$ for a symplectic form $\omega_{FS}$ on $S^2$, and its Reeb orbits are the $S^1$-fibers (here, $\omega_{FS}$ is the Fubini-Study metric on $\mathbb{C}P^1$, and the line bundle associated to the principal $S^1$-bundle $\pi$ is $\mathcal{O}(1)\rightarrow \mathbb{C}P^1$).
    \medskip
    
    \item \textbf{(ellipsoids)} Given $a,b>0$, define the \emph{ellipsoid}
    $$
    E(a,b)=\left\{(z_1,z_2)\in \mathbb{C}^2: \frac{\pi\vert z_1 \vert^2}{a}+  \frac{\pi\vert z_2 \vert^2}{b}\leq 1\right\},
    $$
    a star-shaped domain. The restriction of the symplectic form $\omega_{std}$ is a symplectic form on $E(a,b)$, and its boundary $\partial E(a,b)$ inherits a contact form $\lambda_{std}\vert_{\partial E(a,b)}$ whose Reeb flow is
    $$
    \varphi_t(z_1,z_2)=(e^{2\pi i a t}z_1,e^{2\pi i b t}z_2).
    $$
    In particular, if $a,b$ are rationally independent, then this Reeb flow has only two periodic orbits, passing through the points $z_1=0$, or $z_2=0$. If $a=b$, $E(a,a)$ is the unit ball, and we recover the Hopf flow along the standard $S^3=\partial E(a,a)$. 
    \medskip
    \item\textbf{(Unit cotangent bundle and geodesic flows)} Given a manifold $Q$, choose a Riemannian metric on $TQ$ (which induces a metric on $T^*Q$), and consider its unit cotangent bundle
    $$
    S^*Q=\{(q,p)\in T^*Q: \vert p \vert=1\}.
    $$
    We have $S^*Q=H^{-1}(1/2)$, where $H:T^*Q\rightarrow \mathbb{R}$, $H(q,p)=\frac{\vert p \vert^2}{2}$ is the kinetic energy Hamiltonian. The radial vector field $V=p\partial_p$ on each fiber is the Liouville vector field associated to $-\lambda_{std}$, and is positively transverse to $S^*Q$. It follows that $\alpha_{std}:=-\lambda_{std}\vert_{S^*Q}$ is a contact form, and $(S^*Q,\xi_{std}=\ker \alpha_{std})$ is called the standard contact structure on $S^*Q$. Its Reeb dynamics is the (co)geodesic flow. We see that \emph{a geodesic flow is a particular case of a Reeb flow.}
    \medskip
    \item{\textbf{(Fiberwise star-shaped domains)}} More generally, if $Q$ is a Riemannian manifold, a domain $D\subset T^*Q$ such that the radial vector field $V=p\partial_p$ is transverse to $\partial D$ is called \emph{fiberwise} star-shaped. It inherits a contact structure as in the previous example. 
    
\end{itemize}
\end{example}

\subsection{Symplectization.} Given a contact form $\alpha$ on $X$, its \emph{symplectization} is the symplectic manifold
$$
(\mathbb{R}\times X, \omega=d(e^t\alpha)).
$$
The Liouville vector field is $V=\partial_t$, which is positively transverse to all slices $\{t\}\times X$, where it induces the contact form $i_V\omega=e^t\alpha$. Note that the Reeb dynamics is the same in each slice (i.e.\ it is only rescaled by a constant positive multiple). In fact, the symplectization is the ``universal neighbourhood'' for every contact-type hypersurface:
\begin{proposition}
Let $X\subset (M,\omega)$ be a contact-type hypersurface, with $\omega=d\lambda$ exact near $X$. Then we can find sufficiently small $\epsilon>0$, and an embedding
$$
\Phi:(-\epsilon,\epsilon)\times X\hookrightarrow M, 
$$
so that $\Phi^*\omega=d(e^t\alpha)$ where $\alpha=\lambda\vert_X$.
\end{proposition}

In other words, a contact manifold is always contact-type in \emph{some} symplectic manifold, and vice-versa. We can summarize this discussion in the following motto: \emph{contact geometry is $\mathbb{R}$-invariant symplectic geometry.}

\begin{remark}
One also calls the symplectic manifold $(\mathbb{R}\times X, \omega=d(r\alpha))$ the symplectization of $\alpha$; this is related to the above by the obvious change of coordinates $r=e^t$. We shall use the two interchangeably. Note that $X=\{t=0\}=\{r=1\}$.
\end{remark}

\subsection{Weinstein and Liouville manifolds}
We now discuss an important class of symplectic manifolds, introduced by Eliashberg and Gromov \cite{EG}, where both contact and symplectic geometry, as well as Morse theory, are intertwined. We follow Cieliebak--Eliashberg's definition \cite{CE12}.

\begin{definition} A \emph{Weinstein manifold} is a tuple $(W, \omega, X,\varphi)$, where 
    \begin{itemize}
        \item $(W,\omega)$ is a symplectic manifold,
        \item $\varphi: W \rightarrow \mathbb R$ is an exhausting generalized Morse function,
        \item  $X$ is a complete vector field which is Liouville for $\omega$ and gradient-like for $\varphi$.
    \end{itemize}
    
\end{definition}

Here, a function $\varphi: W \rightarrow \mathbb R$ is \emph{exhausting} if it is proper (i.e.\
preimages of compact sets are compact) and bounded from below. It \emph{Morse} if all its critical points are nondegenerate, and \emph{generalized Morse} if its critical points are either nondegenerate or embryonic, where the latter means that there exist local coordinates $x_1, \dots, x_m$ near the critical point $p$ where the function $\varphi$ coincides with the time $t=0$ function $\varphi_0$ in the birth–death family

$$
\varphi_t(x)=\varphi_t(p) \pm tx_1+x_1^3-\sum_{i=2}^k x_i^2 + \sum_{i=k+1}^m x_j^2.
$$

A vector field $X$ is \emph{complete} if its flow exists for all times. It is \emph{gradient-like} for a function $\varphi$ if $d\varphi(X) \geq \delta(\vert X \vert^2 + \vert d\varphi \vert^2)$, for some positive function $\delta: W \rightarrow \mathbb R^+$ (norms are taken with respect to any Riemannian metric on $W$). Away from critical points this just means $d\varphi(X) > 0$, whereas critical points $p$ of $\varphi$ agree with zeroes of $X$, and $p$ is nondegenerate (embryonic) as a critical point of $\varphi$ iff it is nondegenerate (embryonic) as
a zero of $X$. Here a zero $p$ of a vector field $X$ is embryonic if $X$ agrees near $p$,
up to higher order terms, with the gradient of a function having $p$ as an embryonic critical point.

The compatibility of the Liouville structure with the Morse function implies that stable manifolds of critical points are isotropic, and unstable manifolds, co-isotropic. This imposes a strong topological condition: the index of all critical points is at most $n=\dim(W)/2$, and therefore $W$ is homotopy equivalent to a $CW$ complex of half its dimension. In particular, if $\dim(W)\geq 4$, its boundary is connected.

A Weinstein \emph{cobordism} is a tuple $(W,\omega, X,\varphi)$ where $(W,\omega)$ is a compact symplectic manifold with contact-type boundary $\partial W=\partial_+W\sqcup -\partial _-W$, i.e.\ $X$ is Liouville and is inwards-pointing along $\partial_-W$ and outwards-pointing along $\partial_+W$, and $\varphi$ is generalized Morse but where the condition on exhausting is replaced by asking that $\partial_\pm W$ be a regular level sets of $\varphi$. A \emph{Weinstein domain} is then a Weinstein cobordism with $\partial_-W=\emptyset$. 

A Weinstein manifold is \emph{finite-type} if the Morse function $\varphi$ has finitely many critical points. Therefore one can find a large value $t\in \mathbb R$ such that $W_{cpt}=\{x\in W:\varphi(x)\leq t\}$ is compact and contains all the critical points, and so $W$ is the \emph{completion} of the Weinstein domain $W_{cpt}$, i.e.\ 
$$
W=W_{cpt}\cup_{\partial W_{cpt}} [0,+\infty)\times \partial W_{cpt} 
$$
obtained by attaching the symplectization of the contact manifold $\partial W_{cpt}$ to the boundary of the domain $W_{cpt}$.

Stein manifolds are all Weinstein, with the plurisubharmonic function playing the role of the Morse function. The fact that up to deformation the converse is also true is a deep result of Eliashberg (see \cite{CE12} for all details on this story).    

\emph{A Liouville manifold} is a more relaxed notion than that of a Weinstein manifold, i.e.\ it is a tuple $(W,\omega,X)$ with $X$ a complete Liouville vector field for the symplectic form $\omega$ on $W$. A \emph{Liouville} cobordism and \emph{Liouville domain} are defined analogously, without the conditions on the existence of a Morse function as above. Therefore Weinstein manifolds/domains are Liouville manifolds/domains. The converse is not true, as e.g.\ there exist examples of Liouville domains with disconnected contact-type boundary (see \cite{M91, Mi95, G95, MNW}). See also Section \ref{sec:Liouville} for more background on Liouville domains and manifolds.
 
Weinstein domains can be thought of as being obtained by performing a sequence of handle attachments on the ball, by a construction originally introduced by Weinstein \cite{W91}. In other words, Weinstein domains are handlebodies, where the index of the handles is always at most half the dimension of the manifold.

\chapter{Celestial mechanics}\label{ch:celestial_mechanics}

In this chapter, we introduce the basic problems from celestial mechanics that we will be interested in. The treatment will be brief, as our main interest lies in the chapters to come, and moreover this subject is so classical that the number of references is large. The main character of the story is the CR3BP, to which we will devote more time.

\section{The n-body problem} The setup of the classical $n$-body problem consists of $n$ bodies in $\mathbb{R}^3$, viewed as point-like masses, subject to the gravitational interactions between them, which are governed by Newton's laws of motion. Given initial positions and velocities, the problem consists in predicting the future positions and velocities of the bodies. The understanding of the resulting dynamical system an outstanding open problem. In what follows, we will briefly discuss the general case, and restrict our attention to the simplified circular, restricted case.

Let $q_i\in \mathbb R^3$ be the position vector of the $i$-th mass $m_i$. By Newton's law, we derive the equations of motion to be
$$
m_i \ddot q_i=\sum_{\substack{j=1\\j\neq i}}^n \frac{Gm_im_j(q_j-q_i)}{\Vert q_j - q_i\Vert^3}=-\frac{\partial U}{\partial q_i},
$$
where $G$ is the gravitational constant, and $U$ is the potential energy
$$
U=-\sum_{1\leq i < j \leq n} \frac{Gm_im_j}{\Vert q_i-q_j\Vert}.
$$
If the momentum is defined as $p_i=m_i\dot q_i$, then the Hamiltonian describing these equations is
$$
H=T+U,
$$
where $T$ is the kinetic energy
$$
T=\sum_{i=1}^n\frac{\Vert p_i \Vert^2}{2m_i}.
$$
The problem can be reduced via integrals of motion by appealing to its symmetries. Translational symmetry implies that the center of mass
$$
C=\frac{\sum_{i=1}^n m_iq_i}{\sum_{i=1}^nm_i} 
$$
moves in a straight line, i.e.\ $C(t)=Lt+C_0$, and $L,C_0$ are constants of motion which give six integrals. Rotational symmetry implies that the total angular momentum
$$
A=\sum_{i=1}^nq_i\times p_i
$$
is constant, which gives three more integrals. The last integral is the energy $H$. In total, there are always ten integrals of motion.

\section{Kepler problem} The two-body problem is the most basic model in celestial mechanics. The solutions to this problem can be described by conics, perhaps one of the most beautiful connections between geometry and the laws of nature. As this is the starting point for any study in mechanics, let us briefly revisit this age old problem.  

\subsection{Kepler's laws of planetary motion} Published between 1609 and 1619, the laws of planetary motion were empirically derived by Kepler, from the astronomical observations of his mentor Tycho Brahe. They serve as the most basic description of the motion of Planets around the Sun. They are classically expressed as follows.

\medskip

\begin{itemize}
    \item A planet's motion traces an ellipse, with the Sun at one of the two foci;

\medskip
    
    \item A line segment joining a planet and the Sun sweeps out equal areas during equal intervals of time;

\medskip

    \item The square of the orbital period $T$ is proportional to the cube of the length of the semi-major axis $a$, i.e.\ $T^2\propto a^3$.
\end{itemize}

In the language of Newtonian mechanics, if $r_1,r_2$ denote the position vectors of the masses $m_1,m_2$, we set $r=r_2-r_1$, for which the equations of motion are
$$
\ddot{r}=-\frac{Gm}{\Vert r\Vert^3}r=\nabla U,
$$
where $m=m_1+m_2$, and the potential is $U=\frac{Gm}{\Vert r \Vert}$.

In the language of Hamiltonian dynamics, the Kepler problem is described by the Hamiltonian system
$$
K:T^*(\mathbb R^3\backslash\{0\})\rightarrow \mathbb R,
$$
$$
K(q,p)=\frac{\Vert p \Vert^2}{2} - \frac{1}{\Vert q \Vert}.
$$

As this Hamiltonian is autonomous, it is preserved under its flow.
As the potential is a central force (i.e.\ depends only on $\Vert q \Vert$), angular momentum is also conserved. This implies that the motion always lies in a plane.

\begin{figure}
    \centering
    \includegraphics{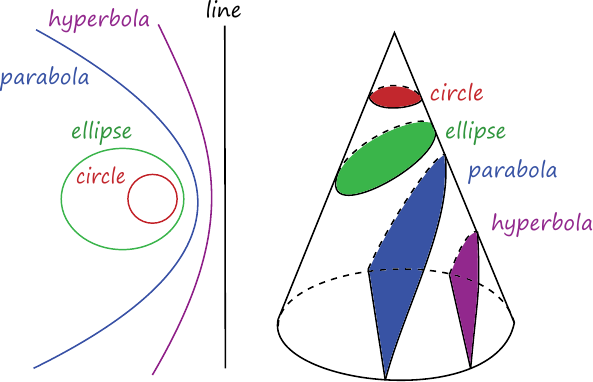}
    \caption{Solutions to the Kepler problem are conics.}
    \label{fig:conics}
\end{figure}

The geometry of the orbits depends on the energy. If $K<0$, the (periodic) orbits are ellipses. If $K=0$, we obtain parabolas, and of $K>0$, hyperbolas. There are also collision orbits, which degenerate into straight lines. In polar coordinates $(r,f)$ centered at one of the foci, the general solution has the form
$$
r=\frac{a(1-e^2)}{1+e\cos f},
$$
where $a$ is the length of the semi-major axis, and $e<1$ is the eccentricity (so that $e=0$ corresponds to circles). The case $e>1$ corresponds to hyperbolas.

The third Kepler law can be equivalently expressed by the fact that the period of a Kepler ellipse depends only on the energy $K<0$ and is given by the formula
$$
T=T(K)=\frac{\pi}{2(-K)^{3/2}}.
$$

\section{The circular restricted three-body problem (CR3BP)}

The CR3BP is a simplification of the general $n$-body problem, for $n=3$, and where the focus is only on the dynamics of one of the masses, which is assumed negligibile by comparison. Concretely, we consider three bodies: Earth (E), Moon (M) and Satellite (S), with masses $m_E, m_M, m_S$ (of course these names may be replaced by, say, Jupiter, Europa, asteroid, respectively). One has the following cases and assumptions.

\medskip

\begin{itemize}
    \item\textbf{(Restricted case)} $m_S=0$, i.e.\ the Satellite is negligible when compared with the \emph{primaries} E and M);

\medskip
    
    \item\textbf{(Circular assumption)} Each primary moves in a circle, centered around the common center of mass of the two (as opposed to general ellipses);

\medskip
    
    \item\textbf{(Planar case)} S moves in the ecliptic plane containing the primaries;

\medskip
    
    \item\textbf{(Spatial case)} The planar assumption is dropped, and S is allowed to move in three-dimensional space.
\end{itemize}

The restricted problem then consists in understanding the dynamics of the trajectories of the Satellite, whose motion is affected by the primaries, but not vice-versa. We denote the \emph{mass ratio} by  $\mu = \frac{m_M}{m_E+ m _M} \in [0,1],$ and we normalize so that $m_E+ m _M=1$, and so $\mu=m_M$ can be thought of as the mass of the Moon. 

In a suitable inertial plane spanned by the $E$ and $M$, the position of the Earth becomes
$$E(t) = (\mu \cos(t),\mu \sin(t),0),$$ and the position of the Moon is
$$M(t) = (-(1-\mu) \cos(t),-(1-\mu)\sin(t),0).$$ The time-dependent Hamiltonian whose Hamiltonian dynamics we wish to study is then
$$
H_t: T^*(\mathbb{R}^3\backslash\{E(t), M(t)\})\rightarrow \mathbb{R}
$$
$$
H_t(q,p)=\frac{1}{2}\Vert  p\Vert^2 - \frac{\mu}{\Vert q-M(t)\Vert } - \frac{1-\mu}{\Vert q-E(t)\Vert },
$$
i.e.\ the sum of the kinetic energy plus the two gravitational potentials associated to each primary. Note that this Hamiltonian is time-dependent. To remedy this, we choose rotating coordinates, in which both primaries are at rest; the price to pay is the appearance of angular momentum term in the Hamiltonian which represents the centrifugal and Coriolis forces in the rotating frame. Namely, we undo the rotation of the frame, and assume that the positions of Earth and Moon are $E=(\mu,0,0),$ $M=(-1+\mu,0,0)$. After this (time-dependent) change of coordinates, which is just the Hamiltonian flow of $L=p_1q_2-p_2q_1$, the Hamiltonian becomes
$$
H: \mathbb{R}^3\backslash\{E,M\}\times \mathbb{R}^3 \rightarrow \mathbb{R}
$$
$$
H(q,p)=\frac{1}{2}\Vert p\Vert^2 - \frac{\mu}{\Vert q-M\Vert } - \frac{1-\mu}{\Vert q-E\Vert } +p_1q_2-p_2q_1,
$$
and in particular is \emph{autonomous}. By preservation of energy, this means that it is a preserved quantity of the Hamiltonian motion. The planar problem is the subset $\{p_3=q_3=0\},$ which is clearly invariant under the Hamiltonian dynamics.

\begin{figure}
    \centering
    \includegraphics[width=0.6\linewidth]{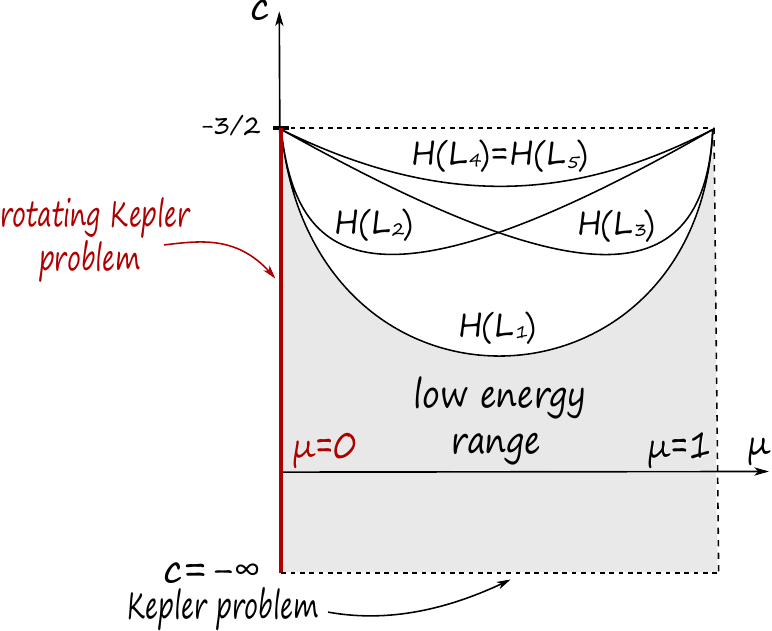}
    \caption{The critical values of $H$.}
    \label{fig:critvalues}
\end{figure}

The Hamiltonian $H$ is invariant under the anti-symplectic involutions 
$$
\rho_1: (q_1,q_2,q_3,p_1,p_2,p_3) \mapsto (q_1,-q_2,q_3,-p_1,p_2,-p_3),
$$
$$
\overline{\rho}_1: (q_1,q_2,q_3,p_1,p_2,p_3) \mapsto (q_1,-q_2,-q_3,-p_1,p_2,p_3),
$$
with corresponding Lagrangian fixed-point loci given by 

$$
L_1=\mbox{Fix}(\rho_1)=\{q_2=p_1=p_3=0\},
$$
$$
\overline{L}_1=\mbox{Fix}(\overline{\rho}_1)=\{q_2=q_3=p_1=0\}.
$$ 
Their composition $\sigma=\rho_1 \circ \overline{\rho}_1$ is symplectic, and corresponds to reflection along the \emph{ecliptic} $\{q_3=0\}$, having fixed point locus the planar problem.

\begin{figure}
    \centering
    \includegraphics[width=0.7\linewidth]{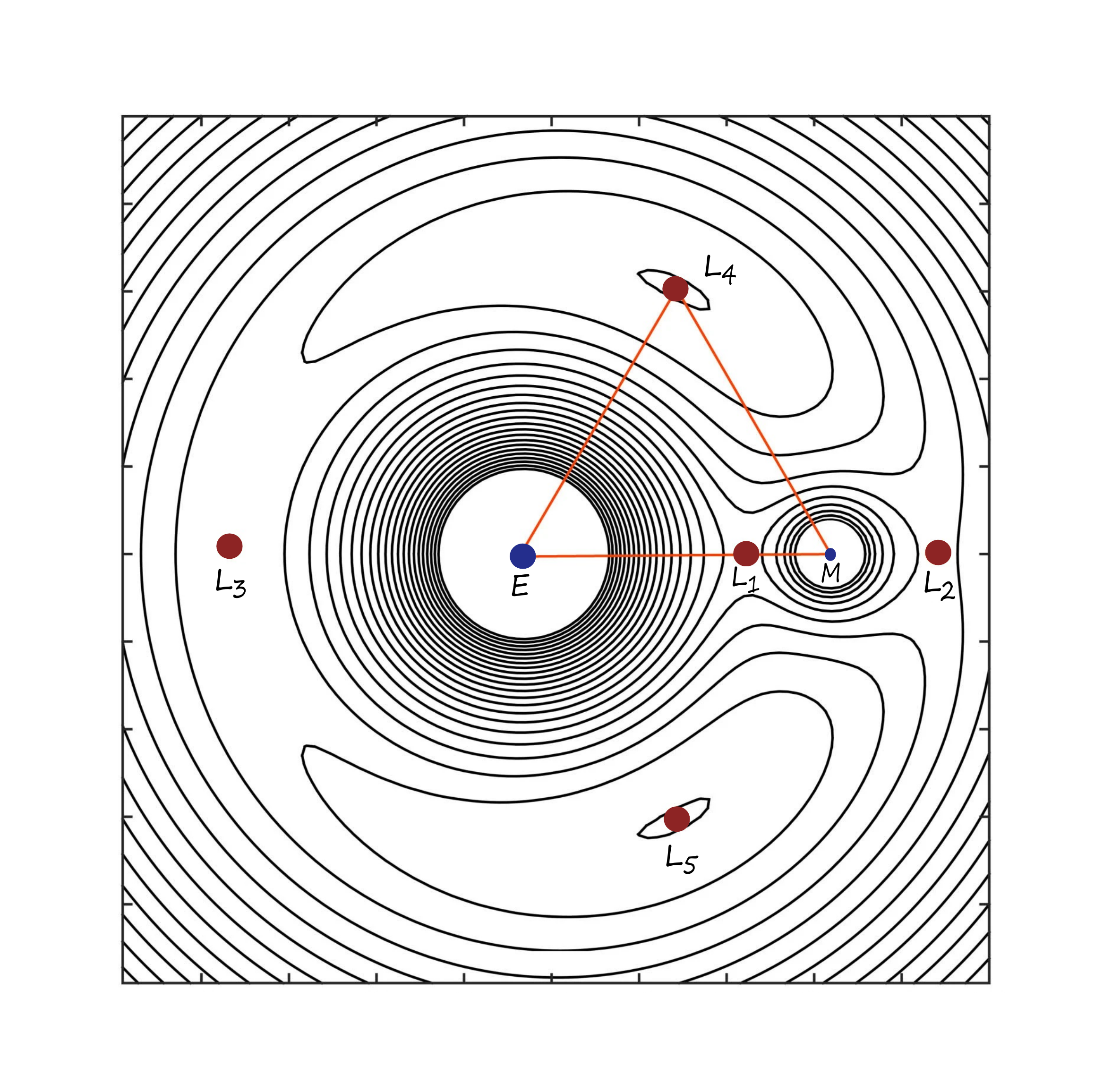}
    \caption{The Hill regions and the Lagrange points for the planar problem.}
    \label{fig:Hillregions}
\end{figure}

As computed by Euler and Lagrange, there are precisely five critical points of $H$, called the \emph{Lagrangian points} $L_i,$ $i=1,\dots,5$, ordered so that $H(L_1)<H(L_2)< H(L_3)<H(L_4)=H(L_5)$ (in the case $\mu <1/2$; if $\mu=1/2$ we further have $H(L_2)=H(L_3)$). See Figure \ref{fig:critvalues}. $L_1,L_2,L_3$, all saddle points, lie in the axis between Earth and Moon (they are the \emph{collinear} Euler points). $L_1$ lies between the latter, while $L_2$ on the opposite side of the Moon, and $L_3$ on the opposite side of the Earth. The others, $L_4$, $L_5$, are maxima, and are called the \emph{triangular} Lagrangian points, as they form equilateral triangles. For $c \in \mathbb{R}$, consider the energy hypersurface $\Sigma_c=H^{-1}(c)$. If 
$$
\pi: \mathbb{R}^3\backslash\{E,M\}\times \mathbb{R}^3 \rightarrow  \mathbb{R}^3\backslash\{E,M\}, \; \pi(q,p)=q, 
$$
is the projection onto the position coordinate, we define the \emph{Hill's region} of energy $c$ as 
$$
\mathcal{K}_c=\pi(\Sigma_c) \in \mathbb{R}^3\backslash\{E,M\}.
$$
This is the region in space where the Satellite of energy $c$ is allowed to move. If $c < H(L_1)$ lies below the first critical energy value, then $\mathcal{K}_c$ has three connected components: a bounded one around the Earth, another bounded one around the Moon, and an unbounded one. Namely, if the Satellite starts near one of the primaries, and has low energy, then it stays near the primary also in the future. The unbounded region corresponds to asteroids which stay away from the primaries. Denote the first two components by $\mathcal{K}_c^E$ and $\mathcal{K}_c^M$, as well as $\Sigma_c^E=\pi^{-1}(\mathcal{K}_c^E)\cap \Sigma_c$, $\Sigma_c^M=\pi^{-1}(\mathcal{K}_c^M)\cap \Sigma_c$, the components of the corresponding energy hypersurface over the bounded components of the Hill region. As $c$ crosses the first critical energy value, the two connected components $\mathcal{K}_c^E$ and $\mathcal{K}_c^M$ get glued to each other into a new connected component $\mathcal{K}_c^{E,M}$, which topologically is their connected sum. Then, the Satellite in principle has enough energy to transfer between Earth and Moon. In terms of Morse theory, crossing critical values corresponds precisely to attaching handles, so similar handle attachments occur as we sweep through the energy values until the Hill region becomes all of position space. See Figure \ref{fig:Hillregions}.

\section{Collision regularization} The $5$-dimensional energy hypersurfaces are non-compact, due to collisions of the massless body $S$ with one of the primaries, i.e.\ when $q=M$ or $q=E$. Note that the Hamiltonian becomes singular at collisions because of the gravitaional potentials, and conservation of energy implies that the momenta necessarily explodes whenever $S$ collides (i.e.\ $p=\infty$). Fortunately, there are ways to regularize the dynamics even after collision. Intuitively, the effect is: whenever $S$ collides with a primary, it bounces back to where it came from, and hence we continue the dynamics beyond the catastrophe. More formally, one is looking for a compactification of the energy hypersurface, which may be viewed as the level set of a new Hamiltonian on another symplectic manifold, in such a way that the Hamiltonian dynamics of the compact, regularized level set is a \emph{reparametrization} of the original one (time is forgotten under regularization).

Two body collisions can be regularized via Moser's recipe. This consists in interchanging position and momenta, and compactifying by adding a point at infinity corresponding to collisions (where the velocity explodes). The bounded components $\Sigma_c^E$ and $\Sigma_c^M$ (for $c < H(L_1)$), as well as $\Sigma_c^{E,M}$ (for $c \in (H(L_1),H(L_1)+\epsilon$)), are thus compactified to compact manifolds $\overline{\Sigma}_c^E$, $\overline{\Sigma}_c^M$, and $\overline{\Sigma}_c^{E,M}$. The first two are diffeomorphic to $S^*S^3=S^3 \times S^2$, and should be thought of as level sets in two different copies of $(T^*S^3,\omega_{std})$ of a suitable regularized Hamiltonian $Q: T^*S^3\rightarrow \mathbb{R}$. The fiber of the level sets $\overline{\Sigma}_c^E$, $\overline{\Sigma}_c^M$ over (a momenta) $p \in S^3$ is a $2$-sphere of allowable positions $q$ in order to have fixed energy. If $p=\infty$ is the North pole, the fiber, called the \emph{collision locus}, is the result of a real blow-up at a primary, i.e.\ we add all possible ``infinitesimal'' positions nearby (which one may think of as all unit directions in the tangent space of the primary). On the other hand, $\overline{\Sigma}_c^E$ is a copy of $S^*S^3 \# S^*S^3$, which can be understood in terms of handle attachments along a critical point of index $1$. In the planar problem, the situation is similar: we obtain copies of $S^*S^2=\mathbb{R}P^3$ and $\mathbb{R}P^3\#\mathbb{R}P^3$.

\begin{figure}
    \centering
    \includegraphics[width=0.8 \linewidth]{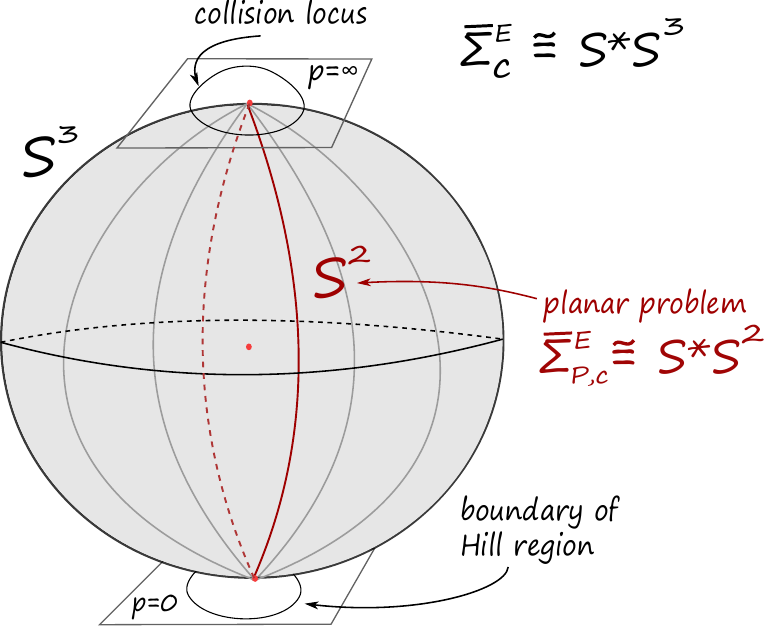}
    \caption{In Moser regularization near the Earth, we add a Legendrian sphere of collisions at the North pole (for fixed energy). The planar problem, which also contains collisions, is an invariant subset.}
    \label{fig:moserreg}
\end{figure}

Another classical way of regularizing collisions is due to Levi--Civita, which works only for the planar problem. This can be viewed as a dynamics on $S^3$, which doubly covers the Moser regularization on $\mathbb RP^3$ (and similarly for the connected sum on $S^3\# S^3$).

\begin{remark}[Global regularization] An important point to keep in mind (especially for Section \ref{sec:feral} below) is that \emph{every} level set for the CR3BP can be regularized to be the level set of a Hamiltonian in $T^*S^3$ (for the spatial problem) and in $T^*S^2$ (for the planar problem), although the formulas take a different form. We refer to Lemma 7.5.1 in \cite{FvK18} for the planar case. The level set is not necessarily contact-type (cf.\ Theorem \ref{thm:contact_type} below).
\end{remark}

In terms of formulas, regularization can be done as follows. 

\subsection{Stark--Zeeman systems.} We will only do the subcritical case $c<H(L_1)$. By restricting the Hamiltonian to the Earth or Moon component, we can view the three-body problem as a \emph{Stark--Zeeman} system, which is a more general class of mechanical systems. 

To define such systems in general, consider a twisted symplectic form
$$
\omega=d\vec p\wedge d\vec q +\pi^*\sigma_B,
$$
with $\sigma_B=\frac{1}{2}\sum B_{ij}dq_i\wedge dq_j$ a $2$-form on the position variables (a \emph{magnetic} term, which physically represents the presence of an electromagnetic field, as in Maxwell's equations), and $\pi(q,p)=q$ the projection to the base. A {\em Stark--Zeeman system} for such a symplectic form is a Hamiltonian of the form
$$
H(\vec q, \vec p)=\frac{1}{2} \Vert \vec p \Vert^2+V_0(\vec q)+V_1(\vec q),
$$
where $V_0(\vec q)=-\frac{g}{\Vert \vec q\Vert}$ for some positive coupling constant $g$, and $V_1$ is an extra potential.\footnote{In this section, we will use the symbol $\vec{\phantom{~}}\;$ for vectors in $\R^3$ to make our formulas for Moser regularization simpler. We will use the convention that $\xi\in \R^4$ has the form $(\xi_0,\vec \xi)$.}

We will make two further assumptions.
\begin{assumptions}

\begin{enumerate}$\;$
\item[(A1)] We assume that the magnetic field is exact with primitive $1$-form $\vec A$.
Then with respect to $d\vec p \wedge d\vec q$ we can write
$$
H(\vec q, \vec p)=\frac{1}{2}\Vert \vec p+\vec A(\vec q) \Vert^2 +V_0(\vec q) +V_1(\vec q).
$$
\item[(A2)] We assume that $\vec A(\vec q)=(A_1(q_1,q_2),A_2(q_1,q_2),0)$, and that the potential satisfies that symmetry $V_1(q_1,q_2,-q_3)=V_1(q_1,q_2,q_3)$.
\end{enumerate}
\end{assumptions}
Observe that these assumptions imply that the planar problem, defined as the subset $\{(\vec q,\vec p):q_3=p_3=0\}$, is an invariant set of the Hamiltonian flow. Indeed, we have
\begin{equation}\label{Hamvf}
\dot q_3 = \frac{\partial H}{\partial p_3}= p_3,
\text{ and }
\dot p_3 = -\frac{\partial H}{\partial q_3}
= -\frac{g q_3}{\Vert \vec q\Vert^3} -\frac{\partial V_1}{\partial q_3}.
\end{equation}
Both these terms vanish on the subset $q_3=p_3=0$ by noting that the symmetry implies that $\frac{\partial V_1}{\partial q_3}\big|_{q_3=0}=0$.

For non-vanishing $g$, Stark--Zeeman systems have a singularity corresponding to two-body collisions, which we will regularize by Moser regularization.
To do so, we will define a new Hamiltonian $Q$ on $T^*S^3$ whose dynamics correspond to a reparametrization of the dynamics of $H$.
We will describe the scheme for energy levels $H=c$, which we need to \emph{fix} a priori (i.e.\ the regularization is not in principle for all level sets at once). Define the intermediate Hamiltonian
$$
K(\vec q,\vec p):=(H(\vec q, \vec p)-c)\Vert\vec q\Vert.
$$
For $\vec q\neq 0$, this function is smooth, and its Hamiltonian vector field equals
$$
X_K= \Vert\vec q\Vert\cdot X_H +(H-c) X_{\Vert\vec q\Vert}. 
$$
We observe that $X_K$ is a multiple of $X_H$ on the level set $K=0$.
Writing out $K$ gives
$$
K=
\left( \frac{1}{2}(\Vert \vec p\Vert^2+1)-(c+1/2) +\langle\vec p,\vec A \rangle +\frac{1}{2}\Vert\vec A\Vert^2+V_1(\vec q)
\right) \Vert\vec q\Vert
-g
. 
$$
\textbf{Stereographic projection.} We now substitute with the stereographic coordinates. The basic idea is to switch the role of momentum and position in the $\vec q,\vec p$-coordinates, and use the $\vec p$-coordinates as position coordinates in $T^*\R^n$ (for any $n$), where we think of $\R^n$ as a chart for $S^n$.
We set 
$$
\vec x=- \vec p, \;\;\vec y=\vec q.
$$ We view $T^*S^n$ as a symplectic submanifold of $T^*\R^{n+1}$, via
$$
T^*S^n
=
\{
(\xi,\eta)\in T^*\mathbb{R}^{n+1} \vert\; \Vert \xi\Vert^2=1,\;
\langle \xi, \eta \rangle =0
\}.
$$
Let $N=(1,0,\dots,0)\in S^n$ be the north pole. To go from $T^*S^n\backslash T^*_N S^n$ to $T^*\R^n$ we use the stereographic projection, given by
\begin{equation}
\label{eq:regularized_to_unregularized}
\begin{split}
\vec x &= \frac{\vec \xi}{1-\xi_0} \\
\vec y &= \eta_0 \vec \xi +(1-\xi_0) \vec \eta .
\end{split}
\end{equation}

To go from $T^*\R^n$ to $T^*S^n\backslash T^*_NS^n$, we use the inverse given by
\begin{equation}
\label{eq:unregularized_to_regularized}
\begin{split}
\xi_0 &= \frac{\Vert \vec x\Vert ^2-1}{\Vert \vec x\Vert ^2+1} \\
\vec \xi &= \frac{2 \vec x}{\Vert \vec x\Vert ^2+1} \\
\eta_0 &= \langle \vec x,\vec y \rangle \\
\vec \eta &= \frac{\Vert \vec x\Vert ^2+1}{2}\vec y 
- \langle \vec x, \vec y \rangle \vec x.
\end{split}
\end{equation}

These formulas imply the following identities
$$
\frac{2}{\Vert \vec x\Vert ^2+1}= 1-\xi_0, \;\;
\Vert \vec y\Vert =\frac{2\Vert\eta \Vert}{\Vert \vec x\Vert^2+1}
=(1-\xi_0) \Vert \eta\Vert 
$$
which allows us to simplify the expression for $K$. Setting $n=3$, we obtain a Hamiltonian $\tilde K$ defined on $T^*S^3$, given by
\[
\begin{split}
\tilde K
&= 
\left( 
\frac{1}{1-\xi_0}
-(c+1/2) 
-\frac{1}{1-\xi_0} \langle \vec \xi , \vec A(\xi,\eta) \rangle
+
\frac{1}{2}\Vert \vec A(\xi,\eta)\Vert ^2
+
V_1(\xi,\eta)
\right)
(1-\xi_0) \Vert \eta \Vert -g 
\\
& =\Vert \eta \Vert
\left(
1-(1-\xi_0)(c+1/2)
-\langle \vec \xi, \vec A(\xi,\eta) \rangle
+
(1-\xi_0)\left(\frac{1}{2}\Vert \vec A(\xi,\eta)\Vert ^2
+
V_1(\xi,\eta)
\right)
\right)
-g
\end{split}
\]
Put
\begin{equation}
\label{eq:Stark_Zeeman_f}
\begin{split}
f(\xi,\eta)
&=
1+(1-\xi_0)\left(-(c+1/2)
+
\frac{1}{2}\Vert \vec A(\xi,\eta)\Vert ^2
+
V_1(\xi,\eta)
\right)
-\langle \vec \xi, \vec A(\xi,\eta) \rangle
\\
&=1+(1-\xi_0) b(\xi,\eta) +M(\xi,\eta)
\end{split}
\end{equation}
where 
$$
b(\xi,\eta)=-(c+1/2)
+
\frac{1}{2}\Vert \vec A(\xi,\eta)\Vert ^2
+
V_1(\xi,\eta)
$$ 
$$
M(\xi,\eta)=-\langle \vec \xi, \vec A(\xi,\eta) \rangle
$$
Note that the collision locus corresponds to $\xi_0=1$, i.e.\ the cotangent fiber over $N$. The notation is supposed to suggest that $(1-\xi_0)b(\xi,\eta)$ vanishes on the collision locus and $M$ is associated with the magnetic term; it is not the full magnetic term, though. We then have that 
$$
\tilde K=\Vert \eta \Vert f(\xi,\eta)-g.$$ To obtain a smooth Hamiltonian, we define the Hamiltonian
$$
Q(\xi,\eta):= \frac{1}{2} f(\xi,\eta)^2 \Vert \eta \Vert^2.
$$
The dynamics on the level set $Q=\frac{1}{2} g^2$ are a reparametrization of the dynamics of $\tilde K=0$, which in turn correspond to the dynamics of $H=c$.
\begin{remark}
We have chosen this form to stress that $Q$ is a deformation of the Hamiltonian describing the geodesic flow on the round sphere, which is given by level sets of the Hamiltonian
$$
Q_0 =\frac{1}{2} \Vert \eta \Vert^2.
$$
This is the dynamics that one obtains in the regularized Kepler problem (the two-body problem; see below), corresponding to the Reeb dynamics of the contact form given by the standard Liouville form. 
\end{remark}

\subsection{Formula for the CR3BP} Since the CR3BP is our main interest, we now give the explicit formula for this problem.
By completing the squares, we obtain
$$
H(\vec q, \vec p)=\frac{1}{2}\left((p_1+q_2)^2+(p_2-q_1)^2+p_3^2\right)-\frac{\mu}{\Vert \vec q - \vec m\Vert}-\frac{1-\mu}{\Vert \vec q -\vec e\Vert} -\frac{1}{2}(q_1^2+q_2^2).
$$
This is then a Stark--Zeeman system with primitive 
$$\vec A=(q_2,-q_1,0),$$ coupling constant $g=\mu$,
and potential 
\begin{equation}\label{potential}
V_1(\vec q)=-\frac{1-\mu}{\Vert \vec q -\vec e\Vert} -\frac{1}{2}(q_1^2+q_2^2),
\end{equation}
both of which satisfy Assumptions (A1) and (A2).

After a computation, we obtain
\begin{equation}\label{eq:f3bp}
f(\xi,\eta)=
1+ \left( 1-\xi_{{0}} \right)  \left( -(c+1/2)+ \xi_{{2}}\eta_{{1}}-\xi_{{1}}
\eta_{{2}} \right) -\xi_{{2}} \mu-\frac {\mu(1-\xi_{{0}}) }{\Vert \vec \eta (1-\xi_0)+\vec \xi \eta_0 + \vec m - \vec e \Vert},
\end{equation}
and we have
\begin{equation}\label{eq:b}
b(\xi,\eta)= -(c+1/2) -\frac {\mu}{\Vert \vec \eta (1-\xi_0)+\vec \xi \eta_0 + \vec m - \vec e \Vert}
\end{equation}
\begin{equation}\label{eq:M1}
M(\xi,\eta)=(1-\xi_0) (\xi_{{2}}\eta_{{1}}-\xi_{{1}}
\eta_{{2}}) - \xi_2\mu.
\end{equation}

\subsection{Levi-Civita regularization.}

We follow the exposition in \cite{FvK}. Consider the map
$$
\mathcal{L}:\mathbb{C}^2\backslash(\mathbb{C}\times \{0\})\rightarrow T^*\mathbb{C}\backslash \mathbb{C},
$$
$$
(u,v)\mapsto \left(\frac{u}{\overline{v}},2v^2\right),
$$
where we view $\mathbb{C}\subset T^*\mathbb{C}$ as the zero section. Using $\mathbb{C}$ as a chart for $S^2$ via the stereographic projection along the north pole, this map extends to a map
$$
\mathcal{L}: \mathbb{C}^2\backslash \{0\} \rightarrow T^*S^2\backslash S^2,
$$
which is a degree $2$ cover. Writing $(p,q)$ for coordinates on $T^*\mathbb{C}=\mathbb{C}\times \mathbb{C}$ (this is the \emph{opposite} to the standard convention, and comes from the Moser regularization), the Liouville form on $T^*\mathbb{C}$ is $\lambda=q_1dp_1+q_2dp_2$, with associated Liouville vector field $X=q_1\partial_{q_1}+q_2\partial_{q_2}$. One checks that
$$
\mathcal{L}^*\lambda=2(v_1du_1-u_1dv_1+v_2du_2-u_2dv_2),
$$
whose derivative is the symplectic form
$$
\omega=d\lambda=4(dv_1\wedge du_1+dv_2\wedge du_2).
$$
Note that $\lambda$ and $\omega$ are \emph{different} from the standard Liouville and symplectic forms (resp.) on $\mathbb{C}^2$. However, the associated Liouville vector field defined via $i_V\omega=\lambda$ coincides with the standard Liouville vector field
$$
V=\frac{1}{2}(u_1\partial_{u_1}+u_2\partial_{u_2}+v_1\partial_{v_1}+v_2\partial_{v_2}),
$$
and we have $\mathcal{L}^*X=V$. We conclude:

\begin{lemma}
A closed hypersurface $\Sigma\subset T^*S^2$ is fiber-wise star-shaped if and only if $\mathcal{L}^{-1}(\Sigma)\subset \mathbb{C}^2\backslash \{0\}$ is star-shaped.
\end{lemma}

Note that $\Sigma\cong S^*S^2\cong \mathbb{R}P^3$, and $\mathcal{L}^{-1}(\Sigma)\cong S^3$, and so $\mathcal{L}$ induces a two-fold cover between these two hypersurfaces.

\subsection{Kepler problem.} We now work out the Moser and Levi-Civita regularizations of the Kepler problem at energy $-\frac{1}{2}$. Recall that its Hamiltonian is given by
$$
E: T^*(\mathbb{R}^2\backslash\{0\})\rightarrow \mathbb{R},
$$
$$
E(q,p)=\frac{1}{2}\Vert p \Vert^2-\frac{1}{\Vert q \Vert}.
$$
The result of Moser regularization is the Hamiltonian
$$
K(p,q)=\frac{1}{2}\left(\Vert q \Vert \left(E(-q,p)+\frac{1}{2}\right)+1\right)^2=\frac{1}{2}\left(\frac{1}{2}\left(\Vert p \Vert^2+1\right)\Vert q \Vert\right)^2.
$$
This is the kinetic energy of the ``momentum'' $q$, with respect to the round metric, viewed in the stereographic projection chart. It follows that its Hamiltonian flow is the round geodesic flow. Moreover, we have $$X_K\vert_{E^{-1}(-1/2)}(p,q)=\Vert q \Vert X_E\vert_{E^{-1}(-1/2)}(-q,p),$$ so that the Kepler flow is a reparametrization of the round geodesic flow.

To understand the Levi-Civita regularization, we consider the shifted Hamiltonian $H=E+\frac{1}{2}$ (which has the same Hamiltonian dynamics). After substituing variables via the Levi-Civita map $\mathcal{L}$, we obtain
$$
H(u,v)=\frac{\Vert u \Vert^2}{2\Vert v \Vert^2}-\frac{1}{2\Vert v \Vert^2}+\frac{1}{2}.
$$
We then consider the Hamiltonian
$$
Q(u,v)=\Vert v \Vert^2H(u,v)=\frac{1}{2}(\Vert u \Vert^2+\Vert v\Vert^2-1).
$$
The level set $Q^{-1}(0)=H^{-1}(0)$ is the $3$-sphere, and the Hamiltonian flow of $Q$, a reparametrization of that of $H$, is the flow of two uncoupled harmonic oscillators. This is precisely the Hopf flow. We summarize this discussion in the following.

\begin{proposition}\label{prop:doublecover}
The Moser regularization of the Kepler problem is the geodesic flow on $S^2$. Its Levi-Civita regularization is the Hopf flow on $S^3$, i.e.\ the double cover of the geodesic flow on $S^2$ (cf.\ Rk.\ \ref{rk:lift}). 
\end{proposition}

\section{The rotating Kepler problem (RKP)} The rotating Kepler problem (RKP) is the boundary case of the CR3BP corresponding to $\mu=0$, i.e.\ there is no Moon anymore, and the Earth is now lying in the origin (but the frame is still rotating). This is a completely integrable system, whose Hamiltonian is explicitly given by
$$
H:T^*(\mathbb R^3\backslash \{0\})\rightarrow \mathbb R,
$$
$$
H=K+L,
$$
where $K(q,p)=\frac{\Vert p \Vert^2}{2}$ is the kinetic energy, and $L(q,p)=p_1q_2-p_2q_1$ is the angular momentum term. One easily checks that
$$
\{H,L\}=\{H,K\}=\{L,K\}=0,
$$
and therefore the Hamiltonian flows are related via
$$
\phi^H_t=\phi^K_t \circ \phi^L_t=\phi^L_t \circ \phi^K_t.
$$
In particular, $K<0$ and $T=T(K)=\frac{\pi}{2(-K)^{3/2}}$ is the period of a Kepler ellipse of energy $K$ (given by Kepler's third law), then, unless it is a circle, a trajectory through $(q,p)$ is periodic for the RKP if and only if a \emph{resonance} condition is satisfied, i.e.\
$$
T(K(q,p))=\frac{k}{l} 2\pi,
$$
for some coprime $l,k \in \mathbb Z$. See \cite{AFFvK13} for more details on periodic orbits in the RKP. The remaining integral for $H$, besides $K$ and $L$, is the last component of the Laplace--Runge--Lenz vector, see \cite{FvK}. 

The Moser regularization of the RKP has the expression
$$
Q(\xi,\eta)=\frac{1}{2}f(\xi,\eta)^2\Vert \eta \Vert^2,
$$
where 
$$
f(\xi,\eta)=1+(1-\xi_0)(-c-1/2+\xi_2\eta_1-\xi_1\eta_2),
$$
is obtained from Equation (\ref{eq:f3bp}) by putting $\mu=0$.

\section{Hill's lunar problem} 

Hill’s lunar problem \cite{H77} is a limit case of the CR3BP where the first primary is much larger than the second one, the second primary has is much larger than the satellite, and the satellite moves very close to the second primary. The Hamiltonian describing the spatial problem is
$$
H: T^*(\mathbb R^3\backslash\{0\})\rightarrow \mathbb R
$$
$$
H(q.p)=\frac{\Vert p \Vert^2}{2}-\frac{1}{\Vert q \Vert}+q_1p_2-q_2p_1-q_1^2+\frac{1}{2}q_2^2+\frac{1}{2}q_3^2.
$$
The planar problem is obtained by setting $q_3=p_3=0$. We refer to \cite{FvK} for a derivation of the Hamiltonian in the planar case, and to \cite{BFvK}, for the spatial case. Roughly speaking, the Hill lunar problem is ``very close'' to the RKP, as the quadratic terms in the Hamiltonian may be suitably understood as a perturbation.

\subsection{Symmetries of Hill's lunar problem}\label{sec:symmetries_Hill} We now present the results of \cite{Ay23}, which completely characterizes the linear symmetries of the lunar problem.

The planar problem is invariant under
the the two commuting linear anti-symplectic involutions
$$
\rho_1 : T^*\mathbb R^2 \rightarrow  T^*\mathbb R^2, (q, p) \mapsto (q_1, -q_2, -p_1, p_2),
$$
$$
\rho_2 : T^*\mathbb R^2 \rightarrow T^*\mathbb R^2, (q, p)\mapsto (-q_1, q_2, p_1, -p_2),
$$
satisfying $\rho_1\circ \rho_1=\rho_1\circ \rho_2=-id$ (symplectic), and so generating the Klein group $$\Sigma_2=\langle \rho_1,\rho_2\rangle =\mathbb Z_2 \oplus \mathbb Z_2.$$ Note that $\rho_1,\rho_2$ are respectively the physical transformations induced by reflection along the $q_2$-axis and $q_1$-axis. For the planar problem, we have the following.

\begin{thm}[Aydin \cite{Ay23}, \textbf{planar problem}] The group of linear involutions which are symplectic or anti-symplectic symmetries of the planar Hill lunar problem is precisely $\Sigma_2=\mathbb Z_2 \oplus \mathbb Z_2$.
\end{thm}

The spatial problem is invariant under the symplectic involution
$$
\sigma: T^*\mathbb R^3\rightarrow T^*\mathbb R^3,
$$
$$
(q_1, q_2, q_3, p_1, p_2, p_3) \mapsto (q_1, q_2, -q_3, p_1, p_2, -p_3),
$$
which is induced by reflection along the \emph{ecliptic} $\{q_3=0\}$. The planar problem is then the restriction of the spatial problem to the (symplectic) fixed-point locus
$$
\mathrm{Fix}(\sigma)=\{q_3=p_3=0\}.
$$
Further linear symplectic symmetries are $-\sigma$ and $\pm id$, where $-\sigma$ is induced by a rotation around the $q_3$-axis by $\pi$. We also have four anti-symplectic involutions

\medskip

\begin{tabular}{c|c}
   $\rho_1(q, p)=(q_1, -q_2, q_3, -p_1, p_2, -p_3)$  &  \text{ induced by reflection at the }$q_1q_3$\text{-plane }\\
    $\rho_2(q, p)=(-q_1, q_2, q_3, p_1, -p_2, -p_3)$ &  \mbox{ induced by reflection at the $q_2q_3$-plane} \\
   $\overline{\rho}_1(q, p)=(q_1, -q_2, -q_3, -p_1, p_2, p_3)$  &  \mbox{ induced by rotation around the $q_1$-axis by $\pi$}\\
   $\overline{\rho}_2(q, p)=(-q_1, q_2, -q_3, p_1, -p_2, p_3)$   & \mbox{ induced by rotation around the $q_2$-axis by $\pi$}\\
\end{tabular}

\medskip

Moreover, we have
$$
\rho_1\vert_{\mathrm{Fix}(\sigma)}=\overline{\rho}_1\vert_{\mathrm{Fix}(\sigma)} \mbox{ and } \rho_2\vert_{\mathrm{Fix}(\sigma)}=\overline{\rho}_2\vert_{\mathrm{Fix}(\sigma)}
$$
coincide with the $\rho_1,\rho_2$ for the planar problem, respectively. These eight symmetries form a group
$$
\{\pm id, \pm \sigma, \rho_1,\rho_2,\overline{\rho}_1,\overline{\rho}_2\}=\langle \rho_1,\rho_2,\sigma\rangle \cong \mathbb Z_2\oplus \mathbb Z_2 \oplus \mathbb Z_2=\Sigma_3.
$$
For the spatial problem, we have the following.

\begin{thm}[Aydin \cite{Ay23}, \textbf{spatial problem}] The group of linear involutions which are symplectic or anti-symplectic symmetries of the spatial Hill lunar problem is precisely $\Sigma_3=\mathbb Z_2 \oplus \mathbb Z_2 \oplus \mathbb Z_3$.
\end{thm}

\chapter{Open books and dynamics}\label{ch:open_books_and_dynamics}

The contents of this chapter lie at the intersection of topology and dynamics. The main character is the notion of open book decompositions, which is purely topological. The way in which it interacts with dynamics is encapsulated in the notion of a global hypersurface of section, which is a higher-dimensional version of the more classical notion of a global surface of section due to Poincar\'e. The emphasis is on examples, in particular those which arise in mechanics. We will include three digressions aimed at illustrating the use of open books in modern contact and symplectic topology.

\section{Open book decompositions} We have the following fundamental notion from smooth topology.

\begin{definition}[\textbf{Open book decomposition}] Let $M$ be a closed manifold. A (concrete) open book decomposition on $M$ is a fibration $\pi: M\backslash B \rightarrow S^1$, where $B \subset M$ is a closed, codimension-$2$ submanifold with trivial normal bundle. We further assume that $\pi(b,r,\theta)=\theta$ along some collar neighbourhood $B\times \mathbb{D}^2\subset M$, where $(r,\theta)$ are polar coordinates on the disk factor.
\end{definition}

\begin{figure}
    \centering
    \includegraphics{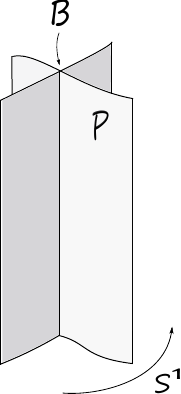}
    \caption{A neighbourhood of the binding look precisely like the pages of an open book, whose front cover has been glued to its back cover via some gluing map (the monodromy).}
    \label{fig:openbookabstract}
\end{figure}

Note that collar neighbourhoods of $B$ exist, since they are trivializations of its normal bundle. $B$ is called the \emph{binding}, and the closure of the fibers $P_\theta=\overline{\pi^{-1}}(\theta)$ are called the \emph{pages}, which satisfy $\partial P_\theta=B$ for every $\theta$. We usually denote a concrete open book by the pair $(\pi,B)$. See Figure \ref{fig:openbookabstract}.

The above concrete notion also admits an abstract version, as follows. Given the data of a typical page $P$ (a manifold with boundary $B$), and a diffeomorphism $\varphi: P \rightarrow P$ with $\varphi=id$ in a neighbourhood of $B$, we can abstractly construct a manifold
$$
M:=\mathbf{OB}(P,\varphi):=B\times \mathbb{D}^2\bigcup_\partial P_\varphi,
$$
where $P_\varphi=P\times [0,1]\backslash (x,0)\sim (\varphi(x),1)$ is the associated mapping torus. By gluing the obvious fibration $P_\varphi \rightarrow S^1$ with the angular map $(b,r,\theta)\mapsto \theta$ defined on $B\times \mathbb{D}^2$, we see that this abstract notion recovers the concrete one. Reciprocally, every concrete open book can also be recast in abstract terms, where the choices are unique up to isotopy. However, while the two notions are equivalent from a topological perspective, it is important to make distinctions between the abstract and the concrete versions for instance when studying dynamical systems adapted to the open books (as we shall do below), since dynamics is of course very sensitive to isotopies.

\begin{figure}
    \centering
    \includegraphics[width=0.7 \linewidth]{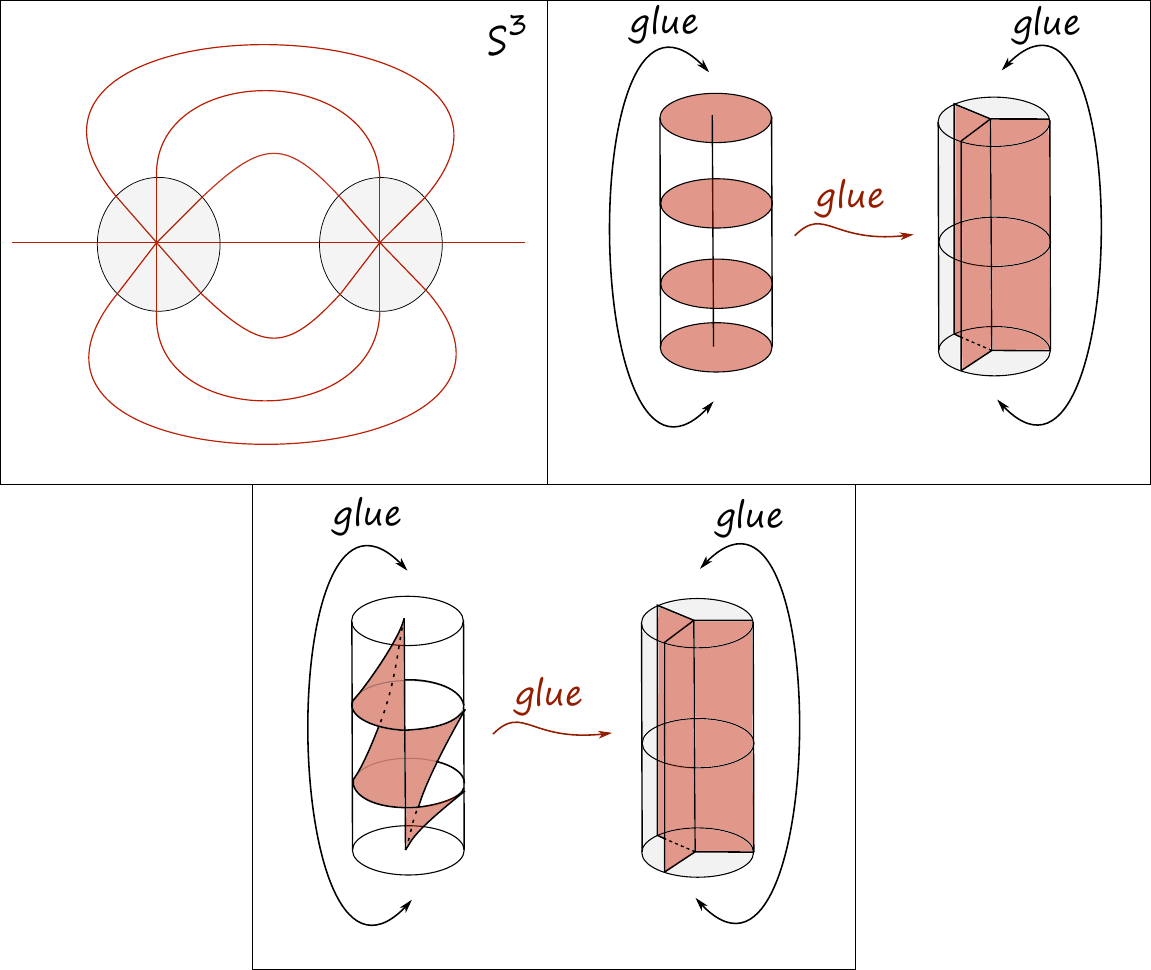}
    \caption{The disk-like pages of the trivial open book in $S^3$ (above) are obtained by gluing two foliations on two solid tori; similarly for its stabilized version (below), whose pages are annuli. Here we use the genus 1 Heegaard splitting for $S^3$.}
    \label{fig:openbooksons3}
\end{figure}

\begin{example}\label{ex}$\;$
\begin{itemize}
    \item\textbf{(Trivial open book)} Since the relative mapping class group of $\mathbb{D}^2$ is trivial, the only possible monodromy for an open book with disk-like pages is $S^3=\mathbf{OB}(\mathbb{D}^2,\mathds{1})$. Viewing $S^3=\{(z_1,z_2) \in \mathbb{C}^2: \vert z_1 \vert^2+\vert z_2 \vert^2=1\}$, let $B=\{z_1=0\}\subset S^3$ be the binding (the unknot). The concrete version is e.g.\ $\pi:S^3\backslash B\rightarrow S^1$, $\pi(z_1,z_2)=\frac{z_1}{\vert z_1\vert}$. See Figure \ref{fig:openbooksons3}.
    
    \medskip
    
    \item\textbf{(Stabilized version)} We also have $S^3=\mathbf{OB}(\mathbb D^*S^1,\tau)$, where $\tau$ is the positive Dehn twist along the zero section $S^1$ of the annulus $\mathbb D^*S^1$. A concrete version is $\pi:S^3\backslash L\rightarrow S^1$, $\pi(z_1,z_2)=\frac{z_1z_2}{\vert z_1z_2\vert}$, where $L=\{z_1z_2=0\}$ is the Hopf link. This is the positive stabilization of the trivial open book, an operation which does not change the manifold (see below). See Figure \ref{fig:openbooksons3}.
    
    \medskip
    
    \item\textbf{(Milnor fibrations)} More generally, let $f:\mathbb{C}^2\rightarrow \mathbb{C}$ be a polynomial which vanishes at the origin, and has no singularity in $S^3$ except perhaps the origin. Then $\pi_f:S^3\backslash B_f\rightarrow S^1$, $\pi_f(z_1,z_2)=\frac{f(z_1,z_2)}{\vert f(z_1,z_2)\vert}$, $B_f=\{f(z_1,z_2)=0\}\cap S^3$, is an open book for $S^3$, called the \emph{Milnor fibration} of the hypersurface singularity $(0,0)$. The link $B_f$ is the \emph{link} of the singularity, and the binding of the open book, whereas the page is called the Milnor \emph{fibre}. If $f$ has no critical point at $(0,0)$, then $B_f$ is necessarily the unknot.
    
    \medskip
    
    \item \textbf{(A trivial product)} We have $S^1\times S^2=\mathbf{OB}(\mathbb D^*S^1,\mathds{1})$. This can be easily seen by removing the north and south poles of $S^2$ (whose $S^1$-fibers become the binding), and projecting the resulting manifold $\mathbb D^*S^1\times S^1$ to the second factor.
    
    \medskip
    
    \item \textbf{(Some lens spaces)} We have $\mathbb{R}P^3=\mathbf{OB}(\mathbb D^*S^1,\tau^2)$, as follows from taking the quotient of the stabilized open book in $S^3$ via the double cover $S^3\rightarrow \mathbb{R}P^3$. More generally, for $p\geq 1$, we have $L(p,p-1)=\mathbf{OB}(\mathbb D^*S^1,\tau^p)$, and for $p\leq 0$, $L(p,1)=\mathbf{OB}(\mathbb D^*S^1,\tau^{p})$. Here, $L(p,q)=S^3/\mathbb{Z}_p,$ is the lens space, where the generator $\zeta=e^{\frac{2\pi i}{p}}\in \mathbb{Z}_p$ acts via $\zeta\cdot (z_1,z_2)=(\zeta. z_1, \zeta^q.z_2)$. For $p=0,1,2$, we recover the above examples.

\end{itemize}
\end{example}

In general, we have the following important result from smooth topology, which says that the open book construction achieves all closed, odd-dimensional manifolds:

\begin{thm}[Alexander ($\dim=3$), Winkelnkemper (simply-connected, $\dim\geq 7$), Lawson ($\dim \geq 7$), Quinn ($\dim\geq 5$)]
If $M$ is closed and odd-dimensional, then $M$ admits an open book decomposition.
\end{thm}

\subsection{Open books in contact topology} So far, we have discussed open books in terms of smooth topology. We now tie it with contact geometry, via the fundamental work of Emmanuel Giroux, which basically shows that contact manifolds can be studied from a purely topological perspective. One therefore usually speaks of the field contact \emph{topology}, when the object of study is the contact manifold itself (as opposed e.g.\ to a Reeb dynamical system on the contact manifold).

If $M$ is oriented and endowed with an open book decomposition, then the natural orientation on the circle induces an orientation on the pages, which in turn induce the boundary orientation on the binding. The fundamental notion is the following.

\begin{definition}[Giroux] Let $(M,\xi)$ be an oriented contact manifold, and $(\pi,B)$ an open book decomposition on $M$. Then $\xi$ is \emph{supported} by the open book if one can find a positive contact form $\alpha$ for $\xi$ (called a \emph{Giroux form}) such that:

\begin{enumerate}
    \item $\alpha_B:=\alpha\vert_B$ is a positive contact form for $B$;
    \item $d\alpha\vert_P$ is a positive symplectic form on the interior of every page $P$.
\end{enumerate}
Here, the a priori orientations on binding and pages are the ones described above. Also, by a \emph{positive} contact form, we mean a contact form $\alpha$ on $M^{2n-1}$ such that the orientation induced by the volume form $\alpha \wedge d\alpha^{n-1}$ coincides with the given orientation on $M$.
\end{definition}

The above conditions are equivalent to:
\begin{itemize}
    \item[(1)'] $R_\alpha\vert_B$ is tangent to $B$;
    \item[(2)'] $R_\alpha$ is positively transverse to the interior of every page.
\end{itemize}

In the above situation, $(B,\xi_B=\ker \alpha_B)$ is a codimension-$2$ contact submanifold, i.e.\ $\xi_B=\xi\vert_B$. 
\begin{thm}[Giroux \cite{Gir02}]
Every open book decomposition supports a unique isotopy class of contact structures. Any contact structure admits a supporting open book decomposition with Weinstein page.
\end{thm}
Here, two contact structures are isotopic if they can be joined by a smooth path $\xi_t$ of contact structures. An important result in contact geometry is \emph{Gray's stability}, which says that isotopic contact structures are \emph{contactomorphic}, i.e.\ there exists a diffeomorphism which carries one to the other. One may further assume that the pages in the above theorem are Stein manifolds, as discussed above (which are in particular \emph{Weinstein}, i.e.\ the Liouville vector field is pseudo-gradient for a Morse function). One may unequivocally use $\mathbf{OB}(P,\varphi)$ to denote the unique isotopy class of contact structures that this open book supports; we write $\mathbf{OB}(P,\varphi)=(M,\xi)$. 

Giroux's original result is actually much stronger in dimension $3$, since it moreover states that the supporting open book is unique up to a suitable notion of \emph{positive stabilization}, which can be thought of as two cancelling surgeries which therefore smoothly do not change the ambient manifold. In arbitrary dimension, this procedure consists of choosing a regular Lagrangian $n$-disk $D$ inside the $2n$-dimensional page $P$ with Legendrian boundary in $\partial P$, attaching an $n$-handle $H$ along the attaching sphere $S^{n-1}\cong \partial D \subset \partial P$, and considering the Lagrangian sphere $S\cong S^n$ obtained by gluing $D$ with the core of $H$. One then replaces the monodromy $\varphi$ with $\varphi \circ \tau_S$, where $\tau_S$ is the right-handed Dehn--Seidel twist along $S$ (an exact symplectomorphism defined by Arnold in dimension $4$ in \cite{A95} and extended by Seidel to higher-dimensions --see e.g.\ \cite{Sei00}--, and which is a generalization of the classical Dehn twist on the annulus). In abstract notation:
$$
\mathbf{OB}(P,\varphi)\rightsquigarrow \mathbf{OB}(P\cup H,\varphi \circ \tau_S).
$$
The handle attachent on the page can be seen as an index $n$ surgery on $M^{2n+1}$, whereas composing with the monodromy adds a cancelling index $n+1$ surgery, so that $\mathbf{OB}(P,\varphi)\cong \mathbf{OB}(P\cup H,\varphi \circ \tau_S)$. Note that if $P$ is a surface then $D$ is simply a properly embedded arc in $P$, and $\tau_S$ is the right-handed Dehn twist along the loop $S$.

\begin{thm}[Giroux's correspondence \cite{Gir02}]
If $\dim(M)=3$, there is a 1:1 correspondence
$$
\{\mathrm{contact}\;\mathrm{structures}\}\Big/\mathrm{isotopy}\longleftrightarrow \{\mathrm{open}\;\mathrm{books}\}\Big/\mathrm{positive}\;\mathrm{stabilization}
$$
\end{thm}

This bijection is why one talks about Giroux's \emph{correspondence}, which reduces the topological study of contact manifolds to the topological study of open books. Let us emphasize that in the above result only the contact \emph{structure} is fixed, and the contact form (and hence the dynamics) is auxiliary; Giroux's result is \emph{not} dynamical, but rather topological/geometrical.

The analogous general uniqueness statement in higher-dimensions has only very recently been established, based on the very recent developments in higher-dimensional convex hypersurface theory as initiated by Honda--Huang \cite{HH}:

\begin{thm}[Breen--Honda--Huang \cite{BHH}]
    Any two supporting Weinstein open book decompositions are stably equivalent.
\end{thm}

Here, two Weinstein open book decompositions are \emph{stably equivalent} if they are related by a sequence of positive stabilizations and destabilizations, conjugations of the monodromy, and Weinstein homotopies (where the homotopies allow for the appearance of Morse, birth-death, and swallowtail type critical points; see \cite{BHH}). The above result is the last missing piece which allows us to talk about Giroux's correspondence in arbitrary dimensions.

\section{Global hypersurfaces of section} From a dynamical point of view, one wishes to adapt the underlying topology to the given dynamics, rather than vice-versa. We therefore make the following.
\begin{definition}
Given a flow $\varphi_t:M\rightarrow M$ of an autonomous vector field on an odd-dimensional closed oriented manifold $M$ carrying a concrete open book decomposition $(\pi,B)$, we say that the open book is \emph{adapted to the dynamics} if:
\begin{itemize}
    \item $B$ is $\varphi_t$-invariant;
    \item $\varphi_t$ is positively transverse to the interior of each page;
    \item for each $x\in M\backslash B$ and $P$ a page, then the orbit of $x$ intersects the interior of $P$ in the future, and in the past, i.e.\ there exists $\tau^+(x)>0$ and $\tau^-(x)<0$ such that $\varphi_{\tau^\pm(x)}(x)\in \mbox{int}(P)$. 
\end{itemize}

\end{definition}

Note that the third condition actually follows from the second one, since we require it for every page and these foliate the complement of $B$. If $\varphi_t$ is a Reeb flow, then the above is equivalent to asking that the (given) contact form is a Giroux form for the (auxiliary) open book. It follows from the definition, that each page is a \emph{global hypersurface of section}, defined as follows:

\begin{definition}(Global hypersurface of section) A \emph{global hypersurface of section} for an autonomous flow $\varphi_t$ on a manifold $M$ is a codimension-$1$ submanifold $P\subset M$, whose boundary (if non-empty) is flow-invariant, whose interior is transverse to the flow, and such that the orbit of every point in $M\backslash \partial P$ intersects the interior of $P$ in the future and past.
\end{definition}

\subsection{Poincar\'e return map.} Given a global hypersurface of section $P$ for a flow $\varphi_t$, this induces a Poincar\'e return map, defined as
$$
f:\mbox{int}(P)\rightarrow \mbox{int}(P),\;f(x)=\varphi_{\tau(x)}(x),
$$
where $\tau(x)=\min\{t >0: \varphi_t(x)\in \mbox{int}(P)\}$. This is clearly a diffeomorphism. And, by construction, periodic points of $f$ (i.e.\ points $p$ for which $f^k(p)=p$ for some $k\geq 1$) are in 1:1 correspondence with closed \emph{spatial} orbits (those which are not fully contained in the binding).

Moreover, in the case of a Reeb dynamics we have:
\begin{proposition}\label{prop:symplecto}
If $\varphi_t$ is the Reeb flow of a contact form $\alpha$, and $P$ is a global hypersurface of section with induced return map $f$, then $\omega=d\alpha\vert_P=d\lambda$, with $\lambda=\alpha\vert_P$, is a symplectic form on $\mbox{int}(P)$, and 
$$
f: (\mbox{int}(P),\omega)\rightarrow(\mbox{int}(P),\omega)
$$
is a symplectomorphism, i.e.\ $f^*\omega=\omega$.
\end{proposition}
In fact, $f$ is an \emph{exact} symplectomorphism, which means that $f^*\lambda=\lambda+d\tau$ for some smooth function $\tau$ (i.e.\ the return time). Differentiating this equation, we obtain $f^*\omega=\omega$. In dimension $2$, a symplectic form is just an area form, and so the above proposition simply says that the return map is area-preserving.

The proof is quite simple: $\omega$ is symplectic precisely because the Reeb vector field, which spans the kernel of $d\alpha$, is transverse to the interior of $P$ (note, however, that it is degenerate at $\partial P$). For $x \in \mbox{int}(P)$, $v \in T_xP$, we have
$$
d_xf(v)=d_x\tau(v)R_{\alpha}(f(x))+d_x\varphi_{\tau(x)}(v).
$$
Using that $\varphi_t$ satisfies $\varphi_t^*\alpha=\alpha$, we obtain
\begin{equation}
    \begin{split}
        (f^*\lambda)_x(v)&=\alpha_{f(x)}(d_xf(v))\\
        &=d_x\tau(v)+(\varphi_{\tau(x)}^*\alpha)_x(v)\\
        &=d_x\tau(v)+\lambda_x(v).\\
    \end{split}
\end{equation}
Therefore 
\begin{equation}\label{Tlambda}
f^*\lambda=d\tau + \lambda,
\end{equation}
which proves the proposition. 

\begin{remark}
In general, the return map might not necessarily extend to the boundary, and indeed there are many examples on which this doesn't hold; this is a delicate issue which usually relies on analyzing the linearized flow equation along the normal direction to the boundary.
\end{remark}

\begin{remark}[Monodromy $\neq$ return map] We wish to emphasize the often puzzling fact that the monodromy of an open book should \emph{not} be confused with the return map of some adapted Reeb flow. First of all, the return map (a \emph{dynamical} object encoding the dynamics) is a map, while the monodromy (a \emph{topological} object encoding the underlying manifold) is strictly speaking an isotopy class of maps relative boundary. Moreover, the return map, as opposed to the monodromy, might \emph{not} necessarily be the identity near or at the boundary (and in most interesting cases it is not). Even more crucially, while the monodromy can be made to preserve a symplectic form on the page (with \emph{infinite} volume), this is different from that preserved by the return map, which has \emph{finite} volume and degenerates at the boundary. The two forms are related in that the former is a completion of a truncation of the latter, however; see App.\ B in \cite{MvK20a} for details. 
    
\end{remark}

Let us discuss two simple but important examples of open books supporting a Reeb dynamics.

\begin{example} Let us discuss two important but simple examples of open books supporting a Reeb dynamics.

\medskip

\begin{itemize}

\item \textbf{(Hopf flow)} The trivial open book on $S^3$, as well as its stabilized version, are both adapted to the Hopf flow. The return map is the identity in both cases.

\medskip

\item \textbf{(Ellipsoids)} 
More generally, the trivial and stabilized open books on $S^3$ are adapted to the Reeb dynamics of every ellipsoid $E(a,b)$. In the trivial case, the return map on each page is the rotation by angle $2\pi\frac{a}{b}$; and in the stabilized case, we get a map of the annulus which rotates the two boundary components in the same direction (i.e.\ it is \emph{not} a twist map, and therefore the classical Poincaré--Birkhoff theorem does not apply).
\end{itemize}

\end{example}

\begin{figure}
    \centering
    \includegraphics{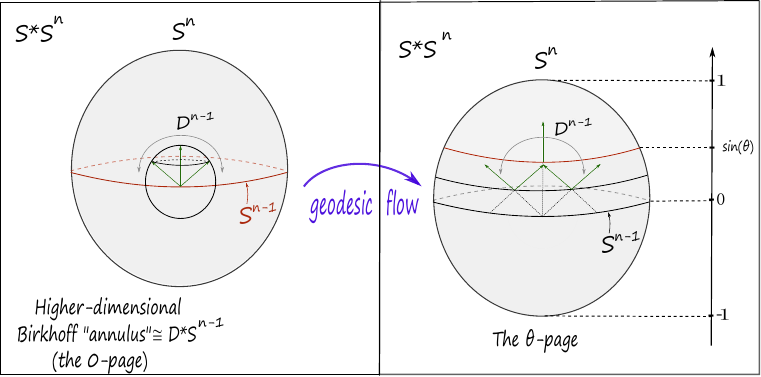}
    \caption{The geodesic open book for $S^*S^n$.}
    \label{fig:birkhoffannulus}
\end{figure}

\section{Open books in mechanics}

We now discuss open books that naturally arise in classical mechanical systems, including the CR3BP.

\subsection{Geodesic flow on $S^*S^n$, and the geodesic open book.} \label{sec:geodesic} We write $$T^*S^{n}=\left\{(\xi,\eta) \in T^*\mathbb{R}^{n+1} =\mathbb{R}^{n+1}\oplus \mathbb{R}^{n+1}:\Vert \xi \Vert=1,\; \langle \xi,\eta \rangle =0\right\}.$$ The Hamiltonian for the geodesic flow is $Q= \frac{1}{2} \Vert \eta \Vert^2\vert_{T^*S^{n}}$ with Hamiltonian vector field
$$
X_Q=\eta \cdot \partial_{\xi} -\xi \cdot \partial_{\eta}.
$$
This is the Reeb vector field of the standard Liouville form $\lambda_{std}$ on the energy hypersurface $\Sigma=Q^{-1}(\frac{1}{2})=S^*S^{n}$. We have the invariant set
$$
B:=\{ (\xi_0,\ldots,\xi_n;\eta_0,\ldots, \eta_n) \in \Sigma~|~\xi_n=\eta_n=0\}=S^*S^{n-1}.
$$
Define the circle-valued map
$$
\pi_g: \Sigma \setminus B \longrightarrow S^1,
\quad
(\xi_0,\ldots,\xi_n;\eta_0,\ldots, \eta_n)
\longmapsto \frac{\eta_n +i \xi_n}{\Vert \eta_n+i\xi_n\Vert}.
$$
This is a concrete open book on $S^*S^n$, which we shall refer to as the \emph{geodesic} open book. The page $\xi_n=0$ and $\eta_n>0$, i.e.\ the fiber over $1\in S^1$, corresponds to a higher-dimensional version of the famous \emph{Birkhoff annulus} (when $n=2$), and is a copy of $\mathbb{D}^*S^{n-1}$. Indeed, it consists of those (co)-vectors whose basepoint lies in the equator, and which point upwards to the upper-hemisphere. See Figure \ref{fig:birkhoffannulus}.

We then consider the angular form 
$$
\omega_g =d\pi_g=\frac{\eta_n d\xi_n -\xi_n d \eta_n}{\xi_n^2 +\eta_n^2}.
$$
We see that $\omega_g(X_Q)=1>0$, away from $B$. This means that $(B,\pi_g)$ is a supporting open book for $\Sigma$ and the pages of $\pi_g$ are global hypersurfaces of section for $X_Q$. In fact, all of its pages are obtained from the Birkhoff annulus by flowing with the geodesic flow. In terms of the contact structure $\xi_{std}=\ker \lambda_{std}$, this open book corresponds to the abstract open book $(S^*S^{n},\xi_{std})=\mathbf{OB}(\mathbb{D}^*S^{n-1},\tau^2)$ supporting $\xi_{std}$. Here, $\tau: \mathbb{D}^*S^{n-1}\rightarrow \mathbb{D}^*S^{n-1}$ is the Dehn--Seidel twist. For $n=2$, we re-obtain the open book $\mathbb{R}P^3=S^*S^2=\mathbf{OB}(\mathbb{D}^*S^1,\tau^2)$. 

This is the abstract open book which will be relevant for the CR3BP; see Theorem \ref{thm:openbooks} below.

\subsection{Double cover of $S^*S^2$.}
We focus on $n=2$, and consider $$S^*S^2=\{(\xi,\eta)\in T^*\mathbb{R}^3: \Vert \xi \Vert=\Vert \eta \Vert=1, \langle \xi, \eta \rangle=0\},$$ the unit cotangent bundle of $S^2$, with canonical projection $\pi_0: S^*S^2 \rightarrow S^2$, $\pi_0(\xi,\eta)=\xi$. It is easy to see that the map
$$
\Phi: S^*S^2\rightarrow SO(3),
$$
$$
\Phi(\xi,\eta)=(\xi,\eta,\xi \times \eta),
$$
is a diffeomorphism, where we view $\xi,\eta$ as column vectors, and so $S^*S^2\cong SO(3)\cong \mathbb{R}P^3$. The projection $\pi_0$ on $SO(3)$ becomes $\pi_0(A)=A(e_1)$, i.e.\ the first column of the matrix $A \in SO(3)$. We have $\pi_1(S^*S^2)=\mathbb{Z}_2$, generated by the $S^1$-fiber. By definition, the double cover of $SO(3)$ is the Spin group Spin$(3)$, which can be constructed as follows. Consider the quaternions
$$
\mathbb{H}=\{a+bi+cj+dk: a,b,c,d \in \mathbb{R}\},
$$
with $i^2=j^2=k^2=-1$, $ij=k, jk=i, ki=j$. We identify $S^3=Sp(1):=\{q \in \mathbb{H}: \Vert q \Vert=1\},$ and $\mathbb{R}^3=\mbox{Im}(\mathbb{H})=\langle i,j,k\rangle$ the set of purely imaginary quaternions. The conjugate of $q=a+bi+cj+dk$ is $\overline{q}=a-bi-cj-dk$. We then define
$$
p: S^3 \rightarrow SO(3),
$$
$$
p(q)(v)=\overline{q} v q,
$$
where $v \in \mbox{Im}(\mathbb{H})=\mathbb{R}^3.$ We have $\Vert \overline{q} v q \Vert=\Vert q \Vert^2 \Vert v\Vert =\Vert v \Vert,$ and $p(q)$ is seen to preserve orientation, so indeed $p(q)\in SO(3)$. Clearly $p(-q)=p(q)$, and the map $p$ is in fact a double cover, so that $S^3=\mbox{Spin}(3)$.

Identifying $i$ with $e_1$, we have $\pi_0(p(q))=p(q)(i)=\overline{q}iq$. A short computation gives
$$
\overline{q}iq=(a+bi+cj+dk)^*i(a+bi+cj+dk)=(a^2+b^2-c^2-d^2)i+2(bc-ad)j+2(ac+bd)k.
$$
On the other hand, the Hopf map may be defined as the map $$\pi: S^3 \rightarrow S^2,\;\pi(z_1,z_2)=(\vert z_1 \vert ^2-\vert z_2 \vert^2, 2\mbox{Re}z_1\overline{z_2},2\mbox{Im}z_1\overline{z_2}),$$ where we view $S^3=\{(z_1,z_2)\in \mathbb{C}^2; \vert z_1\vert^2+\vert z_2 \vert^2=1\}$ and $S^2\subset \mathbb{R}^3$. Writing $q=a+bi+cj+dk=z_1+z_2 j$, i.e.\ $z_1=a+ib$, $z_2=c+id$, one can easily check that 
$$
(\vert z_1 \vert ^2-\vert z_2 \vert^2,2\mbox{Re}z_1\overline{z_2},2\mbox{Im}z_1\overline{z_2})=(a^2+b^2-c^2-d^2,2(bc-ad), 2(ac+bd)).
$$
We have proved the following.
\begin{proposition}\label{prop:Hopf}
The Hopf fibration is the fiber-wise double cover of the canonical projection $\pi_0$, i.e.\ we have a commutative diagram
\begin{center}
    \begin{tikzcd}
        S^1 \arrow[hookrightarrow, d] \arrow[r, "z\mapsto z^2"] & S^1 \arrow[hookrightarrow, d] \\
        S^3=\mbox{Spin}(3)\arrow[r, "p" ] \arrow[d, "\pi"]&
        SO(3)=S^*S^2 \arrow[d, "\pi_0"]\\
        S^2\arrow[equal, r] & 
        S^2\\
    \end{tikzcd}
\end{center}
\end{proposition}

\subsection{Magnetic flows and quaternionic symmetry.} In this section, we expose the beautiful construction of \cite{AG18 } (to which we refer the reader for further details here omitted), relating the quaternions with Reeb flows on $S^3$, as double covers of magnetic flows on $S^*S^2$. 

On $S^2$, consider an area form $\sigma$ (the \emph{magnetic field}), and the \emph{twisted} symplectic form $\omega_\sigma$, defined on $T^*S^2$ via
$$
\omega_\sigma=\omega_{std}-\pi_0^*\sigma,
$$
where $\pi_0: T^*S^2\rightarrow S^2$ is the natural projection. Fixing a metric $g$ on $S^2$, the Hamiltonian flow of the kinetic Hamiltonian $H(q,p)=\frac{\Vert p \Vert^2}{2}$, computed with respect to $\omega_\sigma$, is called the \emph{magnetic flow} of $(g,\sigma)$. Note that $\sigma=0$ corresponds to the geodesic flow of $g$. Physically, the magnetic flow models the motion of a particle on $S^2$ subject to a magnetic field (the terminology comes from Maxwell's equations, which can be recast in this language). From now on, we fix $\sigma$ to be the standard area form on $S^2$, with total area $4\pi$, and $g$ the standard metric with constant Gaussian curvature $1$.

On $S^*S^2$, we can choose a connection $1$-form $\alpha$ satisfying $d\alpha=\pi^*\sigma$, which is a contact form (usually called a \emph{prequantization form}). We identify $T^*S^2\backslash S^2$ with $\mathbb{R}^+\times S^*S^2$, and denoting by $r\in \mathbb{R}^+$ the radial coordinate, we have the associated symplectization form $d(r\alpha)$. Consider the $S^1$-family of symplectic forms
$$
\omega_\theta=\cos \theta \;d(r\alpha)+\sin \theta \;d(r\alpha_{std}),\; \theta \in \mathbb{R}/2\pi \mathbb{Z},
$$
defined on $\mathbb{R}^+\times S^*S^2=T^*S^2\backslash S^2$, where $d(r\alpha_{std})=\omega_{std}$. The Hamiltonian flow of the kinetic Hamiltonian $H$, with respect to $\omega_\theta$, and along $r=1$, is easily seen to be the magnetic flow of $(g,-\cot \theta\cdot \sigma)$ up to constant reparametrization. In particular, for $\theta=\pi/2 \mbox{ mod } \pi$, we obtain the geodesic flow, whose orbits are great circles; for other values of $\theta$ the strength of the magnetic field increases, and the orbits become circles of smaller radius with an increasing left drift. For $\theta=0 \mbox{ mod } \pi,$ the circles become points and the flow rotates the fibers of $S^*S^2$, i.e.\ this is the magnetic flow with ``infinite'' magnetic field.

We now construct the double covers of these magnetic flows on $S^3$, using the hyperk\"ahler structure on $\mathbb{H}=\mathbb{R}^4=\mathbb{C}^2$. We view $S^3$ as the unit sphere in $\mathbb{H}$. Every unit vector 
$$
c=c_1 i +c_2 j +c_2 k \in S^2 \subset \mathbb{R}^3
$$
may be viewed as a complex structure on $\mathbb{H}$, i.e.\ $c^2=-\mathds{1}$. Denoting the radial coordinate on $\mathbb{R}^4$ by $r$, we obtain an $S^2$-family of contact forms on $S^3$ given by
$$
\alpha_c=-2dr \circ c\vert_{TS^2}, \; c \in S^2.
$$
The Reeb vector field of $\alpha_c$ is $R_c=\frac{1}{2}c\partial_r$. Note that $\alpha_i$ is the standard contact form on $S^3$, whose Reeb orbits are the Hopf fibers.

We then consider the quaternionic action of $S^3$ on itself, given by
$$
l_a: S^3 \rightarrow S^3
$$
$$
u \mapsto au,
$$
for $a \in S^3$. Recall that we also have the action of $S^3$ on $S^2$ via the $SO(3)$-action of the previous section, i.e.\ $a\cdot c=p(a)(c)=ac\overline{a} \in S^2,$ for $a\in  S^3,$ $c\in S^2$, and $p: S^3 \rightarrow SO(3)$ the Spin group double cover. One checks directly that $(l_a)_* \alpha_c=\alpha_{ac\overline{a}}=\alpha_{a\cdot c}$. In particular, $(l_a)_*\alpha_i=\alpha_{\pi(a)}$, where $\pi$ is the Hopf fibration.

On the other hand, the stabilizer of $i \in S^2$ under the $S^3$-action is the circle 
$$
\mbox{Stab}(i)=\{\cos(\varphi)+i\sin(\varphi):\varphi \in S^1\}\cong S^1\subset S^3.
$$
The action of an element in this subgroup on $S^3$ then fixes $\alpha_i$, but reparametrizes its Reeb orbits, i.e.\ rotates the Hopf fibers. We then consider an $S^1$-subgroup $\{a_\theta\}\subset S^3$ of unit quaternions which are transverse to this stabilizer, intersecting it only at the identity, given by
$$
a_\theta=\cos(\theta/2)+k\sin(\theta/2), \;\theta\in[0,\pi]
$$
for which
$$
\pi(a_\theta)=a_\theta i \overline{a}_\theta=i\cos\theta +j\sin\theta .
$$
Define $$\alpha_\theta:=\alpha_{\pi(a_\theta)}=\cos \theta \;\alpha_i + \sin \theta\; \alpha_j,$$ with Reeb vector field $R_\theta:=R_{\pi(a_\theta)}$. One further checks that
$$
\alpha_\theta=p^*(\cos \theta\; \alpha + \sin \theta \;\alpha_{std}),
$$
and so
$$
\widetilde \omega_\theta:=d\alpha_\theta=p^*\omega_\theta\vert_{S^*S^2}
$$
is the double cover of the twisted symplectic form $\omega_\theta$ along the unit cotangent bundle (alternatively, we can also think of $\widetilde \omega_\theta$ as being defined on $\mathbb{R}^4\backslash\{0\}=\mathbb{R}^+\times S^3$ as the symplectization of $\alpha_\theta$). We have obtained:

\begin{thm}[\cite{AGZ18}] There are contact forms $\alpha_i,\alpha_j$ and an $S^1$-action on $S^3$, sending $\alpha_i$ to contact forms $\alpha_\theta=\cos \theta\;\alpha_i+\sin \theta\;\alpha_j$, $\theta \in S^1$, such that the Reeb flow of $\alpha_\theta$ doubly covers the magnetic flow of $\omega_\theta$.
\end{thm}

\begin{remark}\label{rk:lift}
Note that for $\theta=0$, corresponding to the infinite magnetic flow, this reduces to the statement of Proposition \ref{prop:Hopf}. For $\theta=\pi/2$, this says that we can lift the geodesic flow on $S^2$ to (a rotated version of) the Hopf flow. Of course, this statement depends on choices; we could have arranged that the lift is precisely the Hopf flow by changing our choice of coordinates.
\end{remark}

\subsection{The magnetic open book decompositions.}

\begin{figure}
    \centering
    \includegraphics[width=0.5 \linewidth]{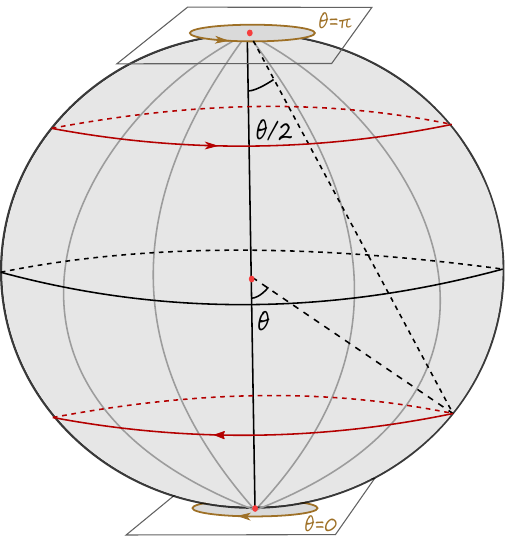}
    \caption{The binding of the magnetic open book $\overline{p}_\theta$ (in red), consisting of two circles of latitude $\theta$ and $\pi-\theta$, doubly covered by two Reeb orbits of $\alpha_\theta$. At $\theta=\pi$ the action of $a_\pi$ maps the Hopf fiber over a point to the Hopf fiber over its antipodal (cf.\ \cite[Fig.\ 1]{AGZ18}).}
    \label{fig:magneticbinding}
\end{figure}

We now tie the previous discussion with open book decompositions. We have seen that the geodesic open book on $S^*S^2$ is constructed in such a way that it is adapted to the geodesic flow of the round metric. On the other hand, by considering the action on $S^3$ of the subgroup $\{a_\theta\}\subset S^3$ of the previous section, we obtain an $S^1$-family $\{p_\theta: S^3\backslash a_\theta(L) \rightarrow S^1\}$ of open book decompositions on $S^3$ (here, $L$ is the Hopf link). These are respectively adapted to the Reeb dynamics of $\alpha_\theta$, and start from the stabilized open book $p_0$ on $S^3$ (adapted to $\alpha_i$ by the example discussed above); they are all just rotations of each other. 

Note that Proposition \ref{prop:Hopf}, the push-forward of $p_0$ under the Hopf map, i.e.\ $\overline{p}_0:=\pi_*(p_0)=p_0 \circ \pi^{-1}: S^*S^2\backslash B_0\rightarrow S^1$ where $B_0$ is the disjoint union of the unit cotangent fibers over the north and south poles $N,S$ in $S^2$ (i.e.\ the image of the Hopf link under $\pi$), is adapted to the infinite magnetic flow. The pages are cylinders obtained as follows: $S^*S^2\backslash B_0\cong ((-1,1)\times S^1)\times S^1$ is a trivial bundle over $S^2\backslash\{N,S\}\cong (-1,1)\times S^1$ (the Euler class of $S^*S^2$ is $-2$), and $\overline{p}_0$ is the trivial fibration.

The push-forward $\overline{p}_\theta=\pi_*(p_\theta): S^*S^2\backslash B_\theta \rightarrow S^1$ is then an open book decomposition on $S^*S^2$, which coincides with the geodesic open book at $\theta=\pi$. The binding $B_\theta$ consists of two magnetic geodesics for $\omega_\theta$; see Figure \ref{fig:magneticbinding}. We call any element of the family $\{\overline{p}_\theta\}$, a \emph{magnetic open book decomposition}. 

\section{Digression: open books and Heegaard splittings} A $3$-dimensional genus $g$ (orientable) handlebody $H_g$ is the $3$-manifold with boundary resulting by taking the boundary connected sum of $g$ copies of the solid $2$-torus $S^1\times \mathbb{D}^2$ (here, we set $H_0=B^3$ the $3$-ball). $H_g$ can also be obtained by attaching a sequence of $g$ $1$-handles to $B^3$. Its boundary is $\Sigma_g$, the orientable surface of genus $g$. A \emph{Heegaard splitting} of genus $g$ of a closed $3$-manifold $X$ is a decomposition 
$$
X=H_g\bigcup_f H_g^\prime,
$$
where $f: \Sigma_g=\partial H_g \rightarrow \Sigma_g=\partial H_g^\prime$ is a homeomorphism of the boundary of two copies of $H_g$. The surface $\Sigma_g$ is called the splitting surface. Different choices of $f$ in the mapping class group of $\Sigma_g$ give, in general, different $3$-manifolds. In fact, it is a fundamental theorem of $3$-dimensional topology that every closed $3$-manifold admits a Heegaard splitting. We have also touched upon another structural result for $3$-manifolds: namely, that every closed $3$-manifold admits an open book decomposition. Let us then discuss how to induce a Heegaard splitting from an open book.

Starting from a concrete open book decomposition $M\backslash B\rightarrow S^1=\mathbb{R}/\mathbb{Z}$ of abstract type $M=\mathbf{OB}(P,\varphi)$, we obtain a Heegaard splitting via 
$$
H_g=\pi^{-1}([0,1/2])\cup B,\;H_g^\prime=\pi^{-1}([1/2,1])\cup B,
$$
where the splitting surface $\Sigma_g=P_0\cup_B P_{1/2}$ is the double of the page $P_0=\pi^{-1}(0)$, obtained by gluing $P_0$ to its ``opposite'' $P_{1/2}=\pi^{-1}(1/2)$. The gluing map $f$ is simply given by $\varphi$ on $P_0$, and the identity on $P_{1/2}$. Stabilizing the open book translates into a stabilization of the Heegaard splitting. 

This shows that the Heegaard diagram thus induced is rather special, since the gluing map is trivial on ``half'' of the splitting surface. In fact, not every Heegaard splitting arises this way, as is easy to see (e.g.\ the lens spaces are precisely the $3$-manifolds with Heegaard splittings of genus $1$, but only the lens spaces discussed in Example \ref{ex} arise from an open book with annulus page, since its relative mapping class group is generated by the Dehn twist).

\section{Digression: open books and Lefschetz fibrations/pencils} We now explore some further interplay between symplectic and algebraic geometry.

\begin{definition}[\textbf{Lefschetz fibration}] Let $M$ be a compact, connected, oriented, smooth $4$-manifold with boundary. A \emph{Lefschetz fibration} on $M$ is a smooth map $\pi : M \rightarrow S$, where $S$ is a compact, connected, oriented surface with boundary, such that each critical point $p$ of $\pi$ lies in the interior of $M$ and has a local complex coordinate chart $(z_1,z_2)\in \mathbb{C}^2$ centered at $p$ (and compatible with the orientation of $M$), together with a local complex coordinate $z$ near $\pi(p)$, such that $\pi(z_1,z_2)=z_1^2+z_2^2$ in this chart.

\end{definition}

\begin{figure}
    \centering
    \includegraphics[width=0.6 \linewidth]{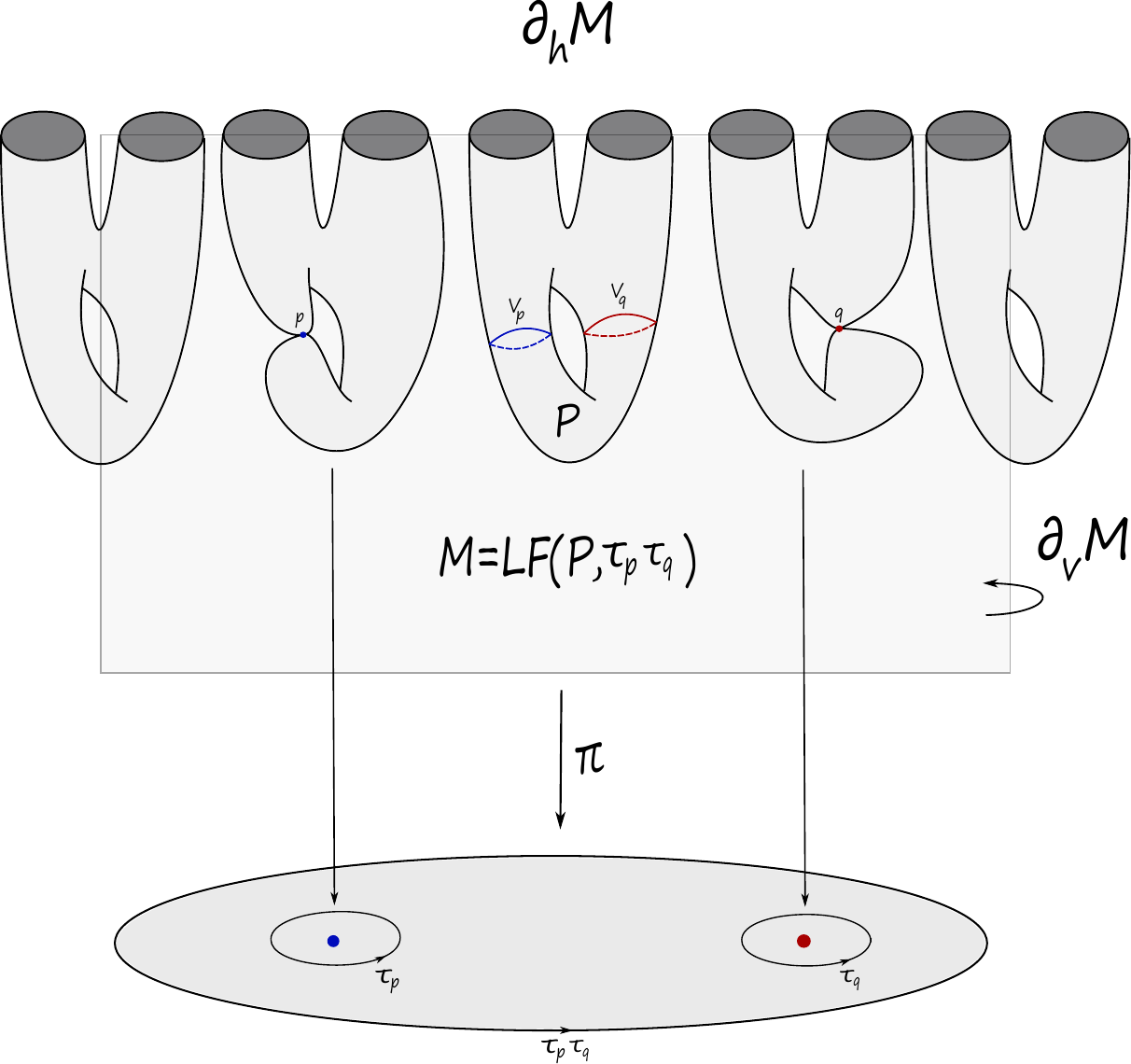}
    \caption{The Lefschetz fibration $\mathbf{LF}(P,\tau_p\tau_q)$ over $\mathbb{D}^2$.}
    \label{fig:LF}
\end{figure}

In other words, each critical point has a local (complex) Morse chart, and is therefore non-degenerate. We then have finitely many critical points due to compactness of $M$. One may also (up to perturbation of $\pi$) assume that there is a single critical point on each fiber of $\pi$. The regular fibers are connected oriented surfaces with boundary, whereas the singular fibers are immersed oriented surfaces with a transverse self-intersection (or node). This singularity is obtained from nearby fibers by pinching a closed curve (the \emph{vanishing cycle}) to a point. See Figure \ref{fig:LF}.

The boundary of a Lefschetz fibration splits into two pieces:
$$
\partial M=\partial_h M\cup \partial_v M, 
$$
where 
$$
\partial_h M=\bigcup_{b\in S}\partial \pi^{-1}(b),\;\partial_v M=\pi^{-1}(\partial S).
$$
By construction, $\partial_h M$ is a circle fibration over $S$, and $\partial_v M$ is a surface fibration over $\partial S$. If we focus on the case $S=\mathbb{D}^2$, the two-disk, denoting the regular fiber $P$ and $B=\partial P$, we necessarily have that $\partial_hM$ is trivial as a fibration, and $\partial_v M$ is the mapping torus $P_\phi$ of some monodromy $\phi: P\rightarrow P$. Therefore
$$
\partial M= \partial_h M \cup \partial_v M= B\times \mathbb{D}^2 \bigcup P_\phi=\mathbf{OB}(P,\phi).
$$
Now, the monodromy $\phi$ is not arbitrary, since orientations here play a crucial role. While every element in the symplectic mapping class group of a surface is a product of powers of Dehn twists along some simple closed loops, it turns out that $\phi$ is necessarily a product of \emph{positive} powers of Dehn twists (once orientations are all fixed). In fact, $\phi=\prod_{p \in \mbox{crit}(\pi)}\tau_p,$ where $\tau_p=\tau_{V_p}$ is the positive (or right-handed) Dehn twist along the corresponding vanishing cycle $V_p\cong S^1\subset P$. This can be algebraically encoded via the monodromy representation $$\rho: \pi_1(\mathbb{D}^2\backslash \mbox{critv}(\pi))\rightarrow \mbox{MCG}(P,\partial P),$$ where $\mbox{critv}(\pi)=\{x_1,\dots,x_n\}$, $x_i=\pi(p_i)$, is the finite set of critical values of $\pi$. We have $$\pi_1(\mathbb{D}^2\backslash \{x_1,\dots,x_n\})=\langle g_\partial,g_1,\dots,g_n: g_\partial=\prod_{i=1}^n g_i\rangle,$$ where $g_i$ is a small loop around $x_i$ and $g_\partial=\partial \mathbb{D}^2$, and $\rho$ is defined via $\rho(g_i)=\tau_{V_{p_i}}$. 

\begin{figure}
    \centering
    \includegraphics[width=0.6 \linewidth]{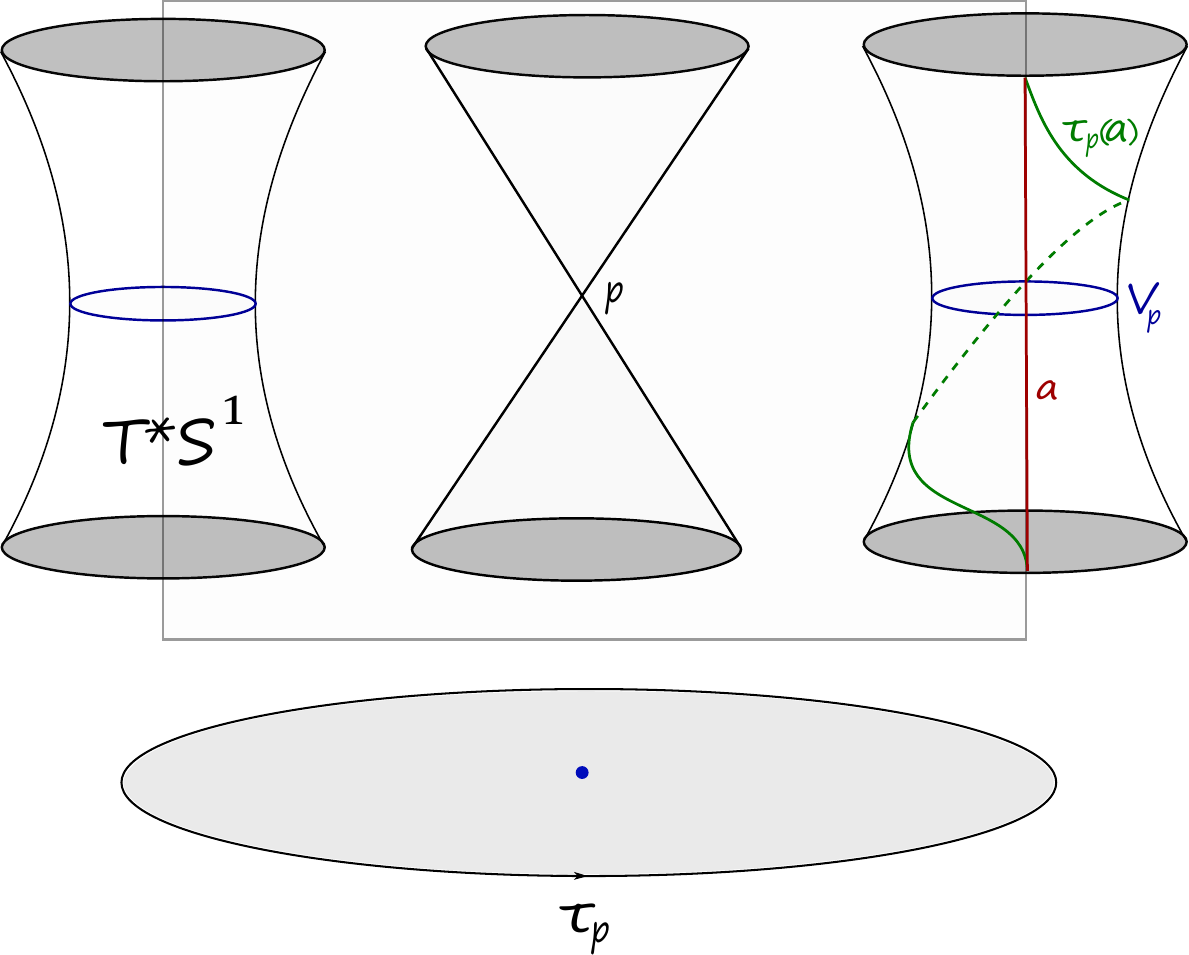}
    \caption{The local model for a Lefschetz singularity.}
    \label{fig:localmodel}
\end{figure}

Reciprocally, a $4$-dimensional Lefschetz fibration on $M$ over $\mathbb{D}^2$ is abstractly determined by the data of the regular fiber $P$ (a surface with non-empty boundary) and a collection of simple closed loops $V_1,\dots,V_n\subset P$. This determines a monodromy $\phi=\prod_{i=1}^n\tau_{V_i}$, a product of positive Dehn twists along the vanishing cycles $V_i$. The recipe to construct $M$ works as follows: decompose $P=\mathbb{D}^2\bigcup H_1\cup\dots \cup H_k$ into a handle decomposition with a single $0$-handle $\mathbb{D}^2$ and a collection of $2$-dimensional $1$-handles $H_1,\dots,H_k\cong \mathbb{D}^1\times \mathbb{D}^1$. One starts with the trivial Lefschetz fibration $M_0=\mathbb{D}^2\times \mathbb{D}^2\rightarrow \mathbb{D}^2$ with disk fiber; and then one attaches (thickened) $4$-dimensional $1$-handles $H_i\times \mathbb{D}^2$ to $M_0$ to obtain the trivial Lefschetz fibration $M_1=P\times \mathbb{D}^2\rightarrow \mathbb{D}^2$ with fiber $P$. In order to add the singularities, one attaches one $4$-dimensional $2$-handle $H=\mathbb{D}^2\times \mathbb{D}^2$ along $V_i\subset P\times \{1\}\subset \partial M_1$, viewed as the attaching sphere $V_i=S^1\times \{0\}\subset S^1\times \mathbb{D}^2\subset \partial H$. At each step of the $2$-handle attachments, we obtain a fibration with monodromy representation $\rho_i$ extending $\rho_{i-1}$ and satisfying $\rho_i(g_i)=\tau_{V_{i}}$, starting from the trivial representation $\rho_0=\mathds{1}: \pi_1(\mathbb{D}^2)=\{1\}\rightarrow \mbox{MCG}(P,\partial P)$. We denote the resulting manifold as $M=\mathbf{LF}(P,\phi)$, for which we have a handle description with handles of index $0,1,2$. 

\begin{remark}
The notation $\mathbf{LF}(P,\phi)$, although simple, is a bit misleading: we need to remember the factorization of $\phi$, since different factorizations lead in general to different smooth $4$-manifolds. One should perhaps use $\mathbf{LF}(P;V_1,\dots,V_n)$ instead, although we hope this will not lead to confusion.
\end{remark}

Having said that, we summarize this discussion in the following.
\begin{lemma}[\textbf{Relationship between Lefschetz fibrations and open books}] We have
$$
\partial \mathbf{LF}(P,\phi)=\mathbf{OB}(P,\phi),
$$
for $\phi=\prod_{i=1}^n\tau_{V_i}$ a product of positive Dehn twists along a collection of vanishing cycles $V_1,\dots,V_n$ in $P$. 
\end{lemma}

While so far this has been a discussion in the smooth category, one may upgrade this to the symplectic/contact category. While we have seen that open books support contact structures in the sense of Giroux, Lefschetz fibrations also support symplectic structures. This is encoded in the following.
\begin{definition}[\textbf{Symplectic Lefschetz fibrations}] An (exact) \emph{symplectic} Lefschetz fibration on an exact symplectic $4$-manifold $(M,\omega=d\lambda)$ is a Lefschetz fibration $\pi$ for which the vertical and horizontal boundary are convex, and the fibers $\pi^{-1}(b)$ are symplectic with respect to $\omega$, also with convex boundary.
\end{definition}

Here, convexity means that the Liouville vector field is outwards pointing. Note that, by Stokes's theorem and exactness of $\omega$, a symplectic Lefschetz fibration cannot have contractible vanishing cycles, since otherwise there would be a non-constant symplectic sphere in a fiber. The description of Lefschetz fibrations in terms of handle attachments can also be upgraded to the sympectic category via the notion of a \emph{Weinstein handle}. After smoothing out the corner $\partial_h M\cap \partial_v M$, the boundary $\partial M$ becomes contact-type via $\alpha=\lambda \vert_{\partial M}$, and the contact structure $\xi=\ker \alpha$ is supported by the open book at the boundary. The contact manifold $(\partial M,\xi)$ is said to be \emph{symplectically filled} by $(M,\omega)$ (see the discussion below on symplectic fillings of contact manifolds). 

Since the space of symplectic forms on a two-manifold is convex and hence contractible, one can show that, given the Lefschetz fibration $\mathbf{LF}(P,\phi)$, an \emph{adapted} symplectic form (i.e.\ as in the definition above) exists and is unique up to symplectic deformation. Therefore, similarly as in Giroux's correspondence, one can talk about $\mathbf{LF}(P,\phi)$ as a symplectomorphism class of symplectic manifolds, and use the short-hand notation $\mathbf{LF}(P,\phi)=(M,\omega)$.  

\begin{example}
An example which is relevant for the spatial CR3BP is that of $T^*S^2$. We consider the \emph{Brieskorn variety}
$$
V_\epsilon=\left\{(z_0,\dots,z_n)\in \mathbb{C}^{n+1}: \sum_{j=0}^{n}z_j^2=\epsilon\right\},
$$
and the associated \emph{Brieskorn manifold} $\Sigma_\epsilon=V_\epsilon \cap S^{2n+1}.$ If $\epsilon=0$, $V_0$ has an isolated singularity at the origin, and $\Sigma_0$ is called the \emph{link of the singularity}. For $\epsilon\neq 0$, the domain $V^{cpt}_\epsilon=V_\epsilon\cap B^{2n+2}$ is a smooth manifold, with boundary $\Sigma_\epsilon\cong \Sigma_0$; the manifold $V_\epsilon$ also inherits a symplectic form by restriction of $\omega_{std}$ on $\mathbb{C}^{n+1}$. Similarly, $\Sigma_\epsilon$ inherits a contact form by restriction of the standard contact form $\alpha_{std}=i\sum_j z_j d\overline{z}_j-\overline{z}_jdz_j$. In fact, $V_\epsilon$ is a Stein manifold, and $V_\epsilon^{cpt}$ is a Stein filling of $\Sigma_\epsilon$; see the discussion on Stein manifolds above, and fillings below.
\end{example} 

A standard fact is the following. 

\begin{proposition}\label{prop:Brieskorn_model}
The map $$(V_1,\omega_{std})\rightarrow (T^*S^n\subset T^*\mathbb{R}^{n+1},\omega_{can}),\;z=q+ip\mapsto (\Vert q \Vert^{-1}q,\Vert q\Vert p)$$ is a symplectomorphism, which restricts to a contactomorphism 
$$
(\Sigma_0,\alpha_{std})\rightarrow (S^*S^n\subset T^*\mathbb{R}^{n+1},\lambda_{can}).
$$

The standard Lefschetz fibration on $T^*S^n$ can be obtained from the Brieskorn variety model as
$$
V_1\rightarrow \mathbb{C},\; (z_0,\dots,z_n)\mapsto z_0.
$$
This induces the geodesic open book on $S^*S^n$ at the boundary, given by the same formula.     
\end{proposition}

\begin{figure}
    \centering
    \includegraphics[width=1 \linewidth]{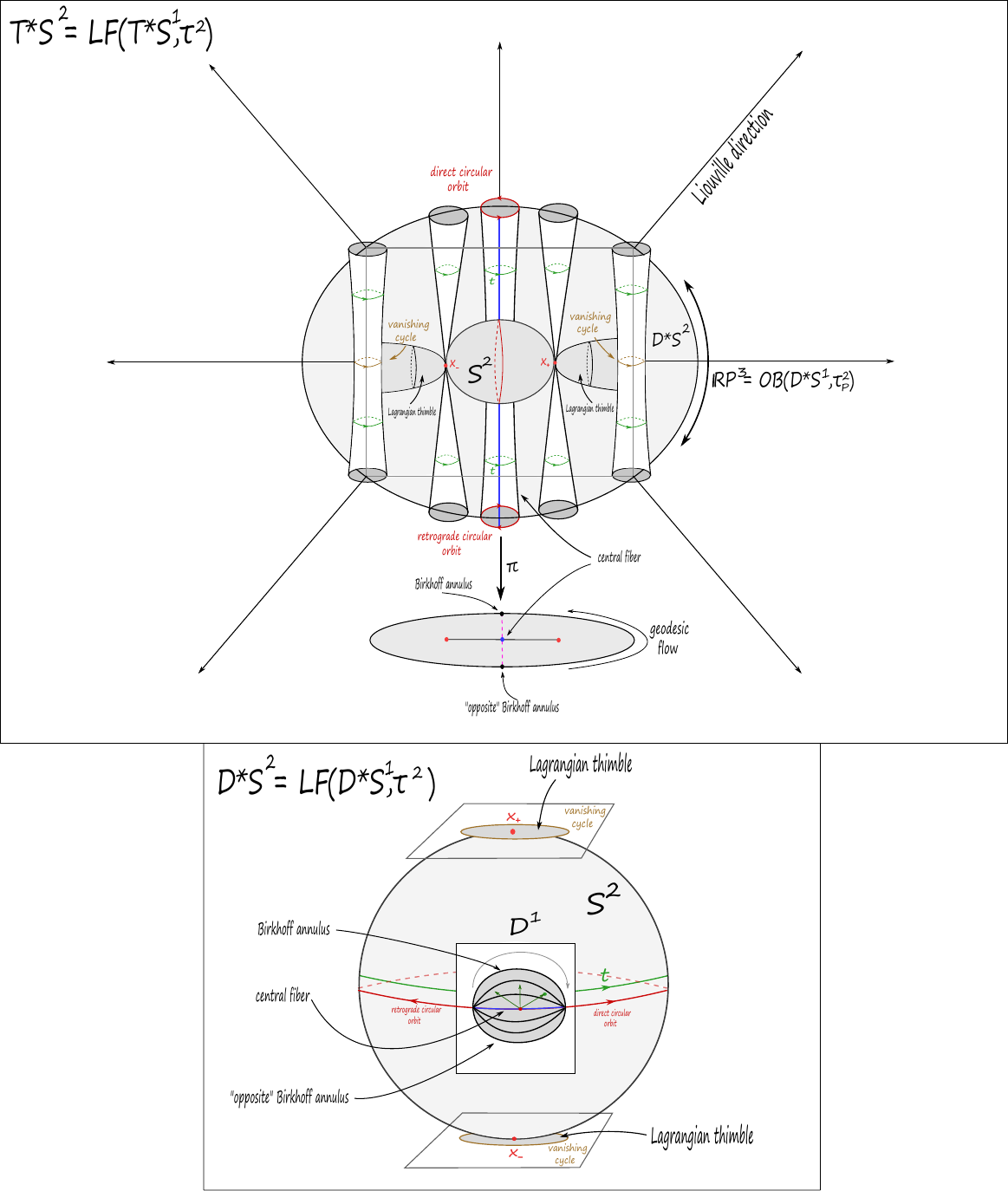}
    \caption{The anatomy of the standard Lefschetz fibration on $\mathbb{D}^*S^2=\mathbf{LF}(\mathbb{D}^*S^1,\tau^2)$, where $\tau$ is the Dehn twist along the zero section $S^1\subset \mathbb{D}^*S^1$. In the picture above, we draw $T^*S^2$ with its (non-compact) fibers $T^*S^1$, and the fibers on $\mathbb{D}^*S^2$ are obtained by projecting along the Liouville direction, and so become $\mathbb D^*S^1$ (compact). These are drawn in the picture below. The two critical points induce the monodromy $\tau^2$. We call the equators transversed in both directions the direct/retrograde (circular) orbits, for reasons that will become apparent. These are the asymptotics of the fibers $T^*S^1$, or equivalently the boundary components of the compact versions $\mathbb D^*S^1$.}
    \label{fig:LFDS2}
\end{figure}

The above map induces the Lefschetz fibration $T^*S^2=\mathbf{LF}(T^*S^1,\tau^2),$ where $\tau$ is the Dehn twist along the vanishing cycle $S^1\subset T^*S^1$, the zero section. We conclude again that $S^*S^2=\mathbb{R}P^3=\mathbf{OB}(\mathbb{D}^*S^1,\tau^2)$. See Figure \ref{fig:LFDS2}.

\medskip

To tie the above discussion with classical algebraic geometry, we introduce the following notion (in the closed case):
\begin{definition}[\textbf{Lefschetz pencil}] Let $M$ be a closed, connected, oriented, smooth $4$-manifold. A \emph{Lefschetz pencil} on $M$ is a Lefschetz fibration $\pi: M\backslash L \rightarrow \mathbb{C}P^1,$ where $L\subset M$ is a finite collection of points, such that near each base point $p \in L$ there exists a complex coordinate chart $(z_1, z_2)$ in which $\pi$ looks like the Hopf map $\pi(z_1, z_2)=[z_1 : z_2]$.
\end{definition}

\begin{figure}
    \centering
    \includegraphics[width=0.7 \linewidth]{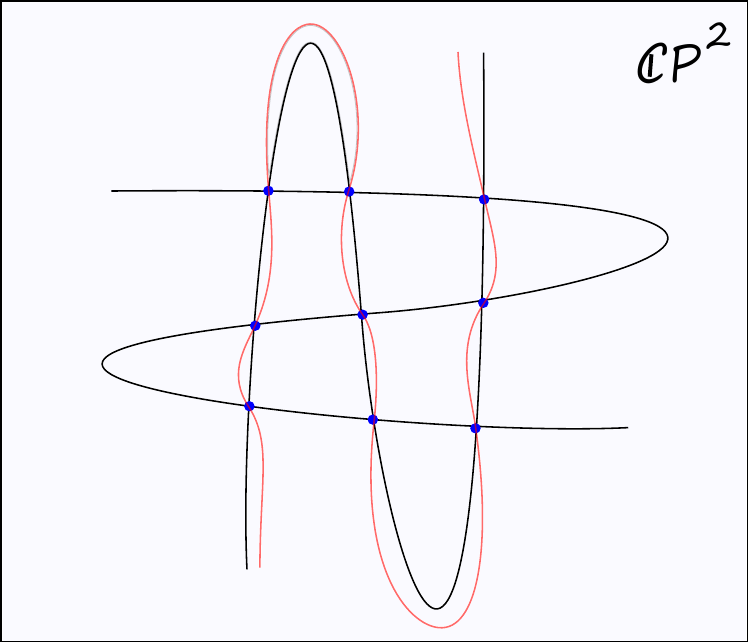}
    \caption{A cartoon of a pencil of cubics, where $L$ consists of $9$ points, and each fiber has genus $1$.}
    \label{fig:pencil}
\end{figure}

Lefschetz pencils arise naturally in the study of projective varieties, and linear systems of line bundles over them. The basic construction is the following. Consider two distinct homogeneous polynomials $F(x,y,z),G(x,y,z)$ of degree $d$ in projective coordinates $[x:y:z]\in \mathbb{C}P^2$ (i.e.\ sections of the holomorphic line bundle $\mathcal{O}(d)$), generic in the sense that $V(F)=\{F=0\}$ and $V(G)=\{G=0\}$ are smooth degree $d$ curves, of genus $g=\frac{(d-1)(d-2)}{2}$ by the genus-degree formula, and so that the base locus $V(F)\cap V(G)=L$ consists of a collection of $d^2$ distinct points (by B\'ezout's theorem). Consider the \emph{degree $d$ pencil} $\{C_{[\lambda:\mu]}\}_{[\lambda:\mu]\in \mathbb{C}P^1},$ where
$$
C_{[\lambda:\mu]}=V(\lambda F + \mu G)\subset \mathbb{C}P^2. 
$$
Through any point in $\mathbb{C}P^2\backslash L$, there is a unique $C_{[\lambda:\mu]}$ which contains it. We then have a Lefschetz pencil 
$$
\pi:\mathbb{C}P^2\backslash L \rightarrow \mathbb{C}P^1, 
$$
where $\pi([x:y:z])=[\lambda:\mu]$ if $C_{[\lambda:\mu]}$ is the unique degree $d$ curve in the family passing through $[x:y:z]$.

By construction, every curve in the pencil meets at the $d^2$ points in $L$. One can further perform a complex blow-up along each of these points, by adding an exceptional divisor (a copy of $\mathbb{C}P^1$) of all possible incoming directions at a given point, and the result is a Lefschetz fibration
$$
Bl_L\pi:Bl_L\mathbb{C}P^2\rightarrow \mathbb{C}P^1.
$$
By construction, this Lefschetz fibration has plenty of spheres, i.e.\ the exceptional divisors, which are sections of the fibration.

The above construction also extends to the case of closed $4$-dimensional projective varieties in some ambient projective space. Moreover, as we have already mentioned, projective varieties are K\"ahler, and in particular symplectic. It is a very deep fact that the above construction extends beyond the algebraic case to the general case of \emph{all} closed symplectic $4$-manifolds:

\begin{thm}[Donaldson \cite{D99}] Any closed symplectic $4$-manifold $(M,\omega)$ admits Lefschetz pencils with symplectic fibers. In fact, if $[\omega]\in H^2(M;\mathbb{Z})$ is integral, the fibers are Poincar\'e dual to $k[\omega]$ for some sufficiently large $k\gg 0$.
\end{thm} 

The above implies that techniques from algebraic geometry can also be applied in the symplectic category, and the interplay is very rich. From the above discussion, after blowing up a finite number of points on the given closed symplectic $4$-manifold $(M,\omega)$, we obtain a Lefschetz fibration.

\section{Digression: symplectic cobordisms and fillings}\label{sec:Liouville} We have already seen the fundamental relationship between contact and symplectic geometry. We now touch upon this a bit further.

\begin{definition}[\textbf{Symplectic cobordism}] A (strong) symplectic cobordism from a closed contact manifold $(X_-,\xi_-)$ to a closed contact manifold $(X_+,\xi_-)$ is a compact symplectic manifold $(M,\omega)$ satisfying:
\begin{itemize}
    \item $\partial M=X_+\bigsqcup X_-$;
    \item $\omega=d\lambda_\pm$ is exact near $X_\pm$, and the (local) Liouville vector field $V_\pm$ (defined via $i_{V_\pm}\omega=\lambda_\pm$) is inwards pointing along $X_-$ and outwards pointing along $X_+$;
    \item $\ker \lambda_\pm\vert_{X_\pm}=\xi_\pm$.
\end{itemize}
\end{definition}

If $\omega=d\lambda$ is globally exact and the Liouville vector field is outwards/inwards pointing along $X_\pm$, we say that $(M,\omega)$ is a \emph{Liouville} cobordism. The boundary component $X_+$ is called \emph{convex} or \emph{positive}, and $X_-$, \emph{concave} or \emph{negative}. Note that a symplectic cobordism is \emph{directed}; in general there might be such a cobordism from $X_-$ to $X_+$ but not viceversa. In fact, the relation $(X_-,\xi_-)\preceq (X_+,\xi_+)$ whenever there exists a symplectic cobordism as above, is reflexive, transitive, but \emph{not} symmetric. We remark that the opposite convention on the choice of \emph{to} and \emph{from} are also used in the literature.

\begin{definition}[\textbf{Symplectic filling/Liouville domain}] A (strong, Liouville) symplectic filling of a contact manifold $(X,\xi)$ is a (strong, Liouville) compact symplectic cobordism from the empty set to $(X,\xi)$. A Liouville filling is also called a \emph{Liouville domain}.
\end{definition}

The \emph{Liouville manifold} associated to a Liouville domain $(M,\omega)$ is its \emph{Liouville completion}, obtained by attaching a cylindrical end:
$$
(\widehat M,\widehat \omega=d\widehat \lambda)=(M,\omega=d\lambda)\bigcup_{\partial M} ([1,+\infty)\times \partial M, d(r\alpha)),
$$
where $\alpha=\lambda\vert_{\partial M}$ is the contact formm at the boundary. Liouville manifolds are therefore ``convex at infinity''. 

It is a fundamental question of contact topology whether a contact manifold is fillable or not, and, if so, how many fillings it admits (say, up to symplectomorphism, diffeomorphism, homeomorphism, homotopy equivalence, $s$-cobordism, $h$-cobordism,...). Note that, given a filling, one may choose to perform a symplectic blowup in the interior, which doesn't change the boundary but changes the symplectic manifold; in order to remove this trivial ambiguity one usually considers \emph{symplectically aspherical} fillings, i.e.\ symplectic manifolds $(M,\omega)$ for which $[\omega]\vert_{\pi_2(M)}=0$ (this holds if e.g.\ $\omega$ is exact, as the case of a Liouville filling). 

For example, the standard sphere $(S^{2n-1},\xi_{std})$ admits the unit ball $(B^{2n},\omega_{std})$ as a Liouville filling. A fundamental theorem of Gromov \cite[p.\ 311]{Gro85} says that this is unique (strong, symplectically aspherical=:ssa) filling up to symplectomorphism in dimension $4$; this is known up to diffeomorphism in higher dimensions by a result of Eliashberg--Floer--McDuff \cite{M91}, but unknown up to symplectomorphism. This was generalized to the case of \emph{subcritically} Stein fillable contact manifolds in \cite{BGZ19}. Another example is a unit cotangent bundle $(S^*Q,\xi_{std})$, which admits the standard Liouville filling $(\mathbb{D}^*Q,\omega_{std})$. There are known examples of manifolds $Q$ with $(S^*Q,\xi_{std})$ admitting only one ssa filling up to symplectomorphism (e.g.\ $Q=\mathbb{T}^2$, \cite{Wen10}; if $n\geq 3$ and $Q=\mathbb{T}^n$, this also holds up to diffeomorphism \cite{BGM19,GKZ19}), but there are other examples with non-unique ssa fillings which are not blowups of each other (e.g.\ $Q=S^n$, $n\geq 3$ \cite{Oba19}). See also \cite{SVHM, LMY17, LO18}. The literature on fillings is vast (especially in dimension $3$) and this list is by all means non-exhaustive.

\begin{remark}
There are also other notions of symplectic fillability: weak, Stein, Weinstein... which we will not touch upon. The set of contact manifolds admitting a filling of every such type is related via the following inclusions:
$$
\{\mbox{Stein}\}\subset \{\mbox{Weinstein}\} \subset \{\mbox{Liouville}\} \subset \{\mbox{strong}\}\subset \{\mbox{weak}\}.
$$
The first inclusion is an equality by a deep result of Eliashberg \cite{CE12}. All others are strict inclusions, something that has been in known in dimension $3$ for some time \cite{Bow12,Ghi05,Eli96}, but has been fully settled in higher-dimensions only very recently \cite{BGM19,BCS,ZZ20,MNW13}. 
\end{remark}

A very broad class for which very strong uniqueness results hold is the following. We say that a contact $3$-manifold $(X,\xi)$ is \emph{planar} if $\xi$ is supported (in the sense of Giroux) by an open book whose page has genus zero. 

\begin{thm}[Wendl \cite{Wen10}]\label{thm:Wendl} Assume that $(M,\omega)$ is a strong symplectic filling of a planar contact $3$-manifold $(X,\xi)$, and fix a supporting open book of genus zero pages, i.e.\ $M=\mathbf{OB}(P,\phi)$ with $g(P)=0$. Then $(M,\omega)$ is symplectomorphic to a (symplectic) blow-up of the symplectic Lefschetz fibration $\mathbf{LF}(P,\phi)$.
\end{thm}

If we assume that the strong filling is \emph{minimal}, in the sense that it doesn't have symplectic spheres of self-intersection $-1$ (i.e.\ exceptional divisors), such a filling is then uniquely determined. It follows as a corollary, that a planar contact manifold is strongly fillable if and only if every supporting planar open book
has monodromy isotopic to a product of positive Dehn twists. This reduces the study of strong fillings of a planar contact $3$-manifolds to the study of factorizations of a given monodromy into product of positive Dehn twists, a problem of geometric group theory in the mapping class group of a genus zero surface.

\medskip

\textbf{References.} A good introductory textbook to contact topology is Geiges' book \cite{G08}; see also \cite{G01} by the same author for a very nice survey on the history of contact geometry and topology, including connections to the work of Sophus Lie on differential equations (which gave rise to the contact condition), Huygens' principle on optics, and the formulation of classical thermodynamics in terms of contact geometry. For an introduction to symplectic topology, McDuff--Salamon \cite{MS17} is a must-read. Anna Cannas da Silva \cite{CdS01} is also a very good source, touching on K\"ahler geometry as well as toric geometry, relevant for the classical theory of integrable systems. For open books and Giroux's correspondence in dimension $3$, Etnyre's notes \cite{E06} is a good place to learn. For open books in complex singularity theory (i.e.\ Milnor fibrations), the classical book by Milnor \cite{M68} is a gem. For related reading on Brieskorn manifolds in contact topology, Lefschetz fibrations and further material, Kwon--van Koert \cite{KvK16} is a great survey. Another good source for symplectic geometry in dimension $4$, Lefschetz pencils, and its relationship to holomorphic curves and rational/ruled surfaces, is Wendl's recent book \cite{Wen18}.

\chapter{Contact geometry in the CR3BP}\label{ch:contact_geometry_in_the_CR3BP}

This chapter is devoted to the advent of the modern approach to the CR3BP, coming from the tools of contact and symplectic geometry. We will discuss how open books arise in the CR3BP, as well as the expected topological picture coming from iterating open books in a suitable way. We will further focus our attention on the return map for the CR3BP, and in particular give a very explicit study the return map in integrable case of the RKP. We will finish with a digression addressing the technicality that the symplectic form always degenerates at the boundary of a global hypersurface of section, a phenomenon which also arises in the context of billiards.

\section{Historical remarks}\label{sec:historical_remarks}

Before giving an overview of the new perspectives on the CR3BP coming from the modern approach to contact geometry, let us revisit some aspects pertaining to the history of the CR3BP. This section contains a historical account, from the Poincaré approach to finding closed orbits in the three-body problem, to some current developments in symplectic geometry. This is by all means non-exhaustive, and tilted towards the author's interests and biased understanding of the developments.

\medskip

\textbf{The perturbative philosophy.} One of the most basic approaches that underlies mathematics and physics is the perturbative approach. Basically, it means understanding a simplified situation first, where everything can be explicitly understood, and attempt to understand ``nearby'' situations by perturbing the parameters relevant to the problem in question. 

In the context of celestial/classical mechanics, this was precisely the approach of Poincar\'e. The idea is to start with a limit case, which is \emph{completely integrable} (i.e.\ an integrable system), perturb it, and study what remained. Integrable systems, roughly speaking, are those which allow enough symmetries so that the solutions to the equations of motion can be 
``explicitly'' solved for (however, quantitative questions need to allow sufficiently many functions, e.g.\ special functions such as elliptic integrals). The solutions tend to admit descriptions in terms of algebraic geometry. In the classical setting of celestial mechanics, if phase-space is $2n$-dimensional and the Hamiltonian $H$ Poisson-commutes with other $n-1$ Hamiltonians (which are therefore preserved under the Hamiltonian flow of $H$), the well-known Arnold--Liouville theorem provides action-angle coordinates in which the symplectic manifold is foliated by flow-invariant tori, along which the Hamiltonian flow is linear, with varying slopes (the \emph{frequencies}). In good situations, the generic tori are half-dimensional (and \emph{Lagrangian}, i.e.\ the symplectic form vanishes along them), whereas there might also be degenerate lower-dimensional tori. This is the natural realm of toric symplectic geometry, dealing with symplectic manifolds which admit a Hamiltonian action of the torus, and the study of the corresponding moment maps and their associated Delzant polytopes. There is also a related theory in algebraic geometry, where the polytope is replaced with a fan. However, in general (e.g.\ the Euler problem) we get only an $\mathbb{R}^n$-action, which is unfortunately beyond the scope of toric geometry. See \cite{HSW99} for more connections between the theory of integrable systems, and differential and algebraic geometry. 

The study of what remains after a small perturbation of an integrable system is the realm of KAM theory, as well as complementary weaker versions such as Aubry--Mather theory. Roughly speaking, the original version of the KAM theorem (due to Kolmogorov--Arnold--Moser) says that if one perturbs a ``sufficiently irrational'' Liouville torus, i.e.\ the vector of frequencies of the action is very badly approximated by rational numbers (it is \emph{diophantine}) and moreover the Hessian with respect to action variables is non-degenerate, then the Liouville tori survives to an invariant tori whose frequencies are close to the original one, and hence is foliated by orbits which are \emph{quasi-periodic}, in the sense that they are dense in the tori and never close up. Aubry--Mather theory is meant to deal with the rest of the tori, including resonant ones which are foliated by closed orbits and non-diophantine non-resonant ones, as well as large deformations (as opposed to sufficiently small perturbations). This theory provides invariant subsets which are usually Cantor-like, and obtained via measure-theoretical means (they are the supports of invariant measures minimizing certain action functionals).

\medskip

\textbf{The Poincar\'e--Birkhoff theorem, and the planar CR3BP.} The problem of finding closed orbits in the planar case of the CR3BP goes back to ground-breaking work in celestial mechanics of Poincar\'e \cite{P87,P12}, building on work of G.W.\ Hill on the lunar problem \cite{H77,H78}. The basic scheme for his approach may be reduced to:
\begin{itemize}
    \item[(1)] Finding a global surface of section for the dynamics;
    \item[(2)] Proving a fixed point theorem for the resulting first return map.
\end{itemize}
This is the setting for the celebrated Poincar\'e--Birkhoff theorem, proposed and confirmed in special cases by Poincar\'e and later proved in full generality by Birkhoff in \cite{Bi13}. The statement can be summarized as: if $f: A\rightarrow A$ is an area-preserving homeomorphism of the annulus $A=[-1,1]\times S^1$ that satisfies a \emph{twist} condition at the boundary (i.e.\ it rotates the two boundary components in opposite directions), then it admits infinitely many periodic points of arbitrary large period. The fact that the area is preserved is a consequence of Liouville's theorem for Hamiltonian systems; we have basically used this in our proof of Proposition \ref{prop:symplecto}. 

The whole point of a global surface of section is to reduce a \emph{continuous} flow on a $3$-manifold to the \emph{discrete} dynamics of a map on a $2$-manifold, thus reducing by one the degrees of freedom. It is perhaps fair to say, that this key (and beautiful) idea is responsible for motivating the well-studied area of dynamics on surfaces, a huge industry in its own right.

\medskip

\textbf{The direct and retrograde orbits.} The actual physical Moon is in \emph{direct} motion around the Earth (i.e.\ it rotates in the same direction around the Earth as the Earth around the Sun). The opposite situation is a \emph{retrograde} motion. In \cite{H77,H78}, while attempting to model the motion of the Moon, Hill indeed finds both direct and retrograde orbits. While still an idealized situation, such direct orbit is a reasonable approximation to the actual orbit of the Moon, and Hill even goes further to find better approximations via perturbation theory, something which deeply impressed Poincar\'e himself. Let us remark that direct orbits are usually the more interesting to astronomers, since most moons are in direct motion around their planet. Topologically, one may think of the retrograde/direct Hill orbits as obtained from a Hopf link in $S^3$, via the double cover to $\mathbb{R}P^3$. This is the binding of the open book $\mathbb{R}P^3= \mathbf{OB}(\mathbb{D}^*S^1,\tau^2)$, where $\tau$ is the positive Dehn twist along $S^1\subset \mathbb{D}^*S^1$.

\medskip

\textbf{Brouwer's and Frank's theorem.} In order to find the direct orbit away from the lunar problem, Birkhoff had in mind finding a disk-like surface of section whose boundary is precisely the retrograde orbit. The direct orbit would then be found via Brouwer's translation theorem: every area preserving homeomorphism of the open disk admits a fixed point. Removing the fixed point, we obtain an area preserving homeomorphism of the open annulus, which, via a theorem of Franks, admits either none or infinitely many periodic points. All this combined, one has: an area preserving homeomorphism of an open disk admits either one or infinitely many periodic points. Note that if the boundary is also an orbit, we obtain 2 or infinitely many. If furthermore we have twist, the Poincar\'e--Birkohff theorem provides infinitely many orbits. This is a classical heuristic for finding orbits that has survived to this day in several guises, as we will see below. See Figure \ref{fig:planarproblem}.

\medskip

\textbf{Perturbative results.} As we have seen, we have $\mathbb{R}P^3= \mathbf{OB}(\mathbb{D}^*S^1,\tau^2)$ as smooth manifolds, and one would hope that a concrete version of this open book is adapted to the (Moser-regularized) planar dynamics, and that the return map is a Birkhoff twist map. For $c<H(L_1)$ and $\mu \sim 0$ small, one can interpret from this perspective that Poincar\'e \cite{P12} proved this by perturbing the rotating Kepler problem (when $\mu=0$), which is an integrable system for which the return map is a twist map. Of course, he never stated it in these words. In the case where $c\ll H(L_1)$ is very negative and $\mu \in (0,1)$ is arbitrary, this was done by Conley \cite{C63} (also perturbatively), who checked the twist condition and used Poincar\'e--Birkhoff. In \cite{M69}, McGehee provides a disk-like global surface of section for the rotating Kepler problem problem for $c<H(L_1)$, and computes the return map.

\medskip

\textbf{Non-perturbative results.} More generally and \emph{non-perturbatively}, the existence of this adapted open book was obtained in \cite[Theorem 1.18]{HSW19} for the case where $(\mu,c)$ lies in the \emph{convexity range} via holomorphic curve methods due to Hofer--Wysocki--Zehnder \cite{HWZ98} (see also \cite{AFFHvK,AFFvK}). Here, the convexity range is the set of pairs $(\mu,c)$ of mass ratio and Jacobi constant for which the dynamics of the Levi--Civita regularization of the planar CR3BP is the dynamics at the boundary of a convex domain in $\mathbb R^4$ (induced by the standard Liouville form); it is a \emph{non-perturbative} set of parameters. This non-perturbative approach implies the use of modern techniques of symplectic and contact geometry, i.e.\ holomorphic curves.


\begin{figure}
    \centering
    \includegraphics[width=0.7 \linewidth]{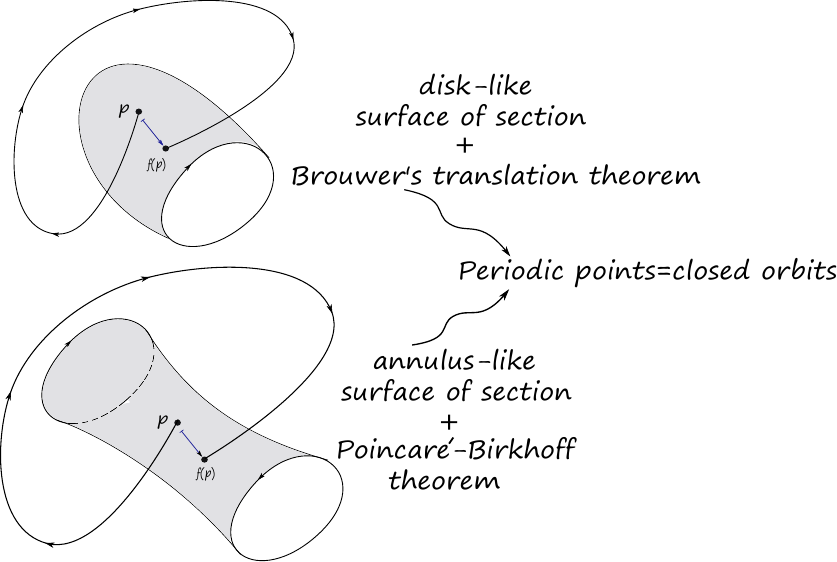}
    \caption{Obtaining closed orbits in the planar problem.}
    \label{fig:planarproblem}
\end{figure}

\medskip

\textbf{The search of closed geodesics: a very brief survey.} After suitable regularization, the round geodesic flow on $S^2$ appears as an integrable limit case in the planar CR3BP, when the Jacobi constant $c$ converges to $-\infty$. Poincaré was aware of this fact, which brought him, near the end of his life, to study the geodesic flow of ``near-integrable'' metrics on $S^2$, i.e.\ perturbations of the round one. One may well argue that this was one of the starting points of the very long and fruitful search of closed geodesics that ensued later throughout the 20th century.

A basic argument for finding closed geodesics, sometimes attributed to Birkhoff, was already present in work of Hadamard in 1898, who studied the case of surfaces with negative curvature. This is a variational argument on the loop space, in the sense that closed geodesics are viewed as loops which happen to be geodesics (as opposed to the dynamical point of view, where a closed geodesic is a geodesic path which happens to close up). It works as follows: on a compact manifold, one chooses a sequence of loops in a fixed homotopy class whose length converges to the infinimum in such class, and appeals to the Arzel\`a--Ascoli theorem. If the infimum is non-zero, this gives a non-trivial closed geodesic. This argument works if the fundamental group is non-trivial; it gives a geodesic in each non-trivial free homotopy class, and hence infinitely many if the genus is at least $1$. This leaves out the case of $S^2$, for which it gives nothing. The program of finding geodesics for general manifolds was picked up by Birkhoff in a more systematic way, who proved existence of at least one geodesic for the case of all surfaces and certain higher-dimensional manifolds including spheres. For the case where the infimum in the above variational argument is zero, Birkhoff introduced the famous minmax argument. For $S^2$, this works as follows: take the foliation of $S^2$ minus the north and south poles, whose leaves are the circles given by the parallels (think of the standard embedding, but where the metric is not the standard one). Choose a curve shortening procedure for each non-trivial leaf (there are several, the simplest one being replacing two nearby points on a loop by a geodesic arc; this is a tricky business, however, since the resulting loop might have self-intersections). This gives a sequence of foliations, and we may choose the loop with maximal length for each. These lengths are bounded from below for topological reasons. Again by Arzel\`a--Ascoli, the limit of such curves, being invariant under the shortening procedure, is a geodesic.  

Before Birkhoff, Poincaré himself \cite{P05} had the idea of obtaining a geodesic for the case of $S^2$ embedded in $\mathbb{R}^3$ as a convex surface $S$ (with the induced metric), by considering the shortest simple closed curve $\gamma$ dividing $S$ into two pieces of equal total Gaussian curvature. A simple argument using Gauss--Bonnet shows that $\gamma$ should be a geodesic. The full details of this beautiful argument were carried out by Croke in 1982 \cite{C82}, who considered the more general case of a convex hypersurface in $\mathbb{R}^n$.

Poincaré further proposed that, also in the case of a convex $S^2$ in $\mathbb{R}^3$, there should be at least $3$ closed geodesics with no self-intersections (i.e.\ simple). A short proof of this was published by Lusternik--Schnirelmann in 1929 \cite{LS29}. Their proof relied on two steps: firstly, to consider the space of all simple circles (great and small) and a continuous curve-shrinking procedure which keeps all such circles simple; and secondly, the fact that the space of non-oriented round geodesics is a copy of $\mathbb{R}P^2$ (it can be identified with the space of planes in $\mathbb{R}^3$ through the origin), together with the fact that every Morse function on $\mathbb{R}P^2$ has at least $3$ critical points. Unfortunately, there were gaps in both steps. These were filled in by Ballmann in 1978 \cite{B78}, who also considered the case of arbitrary genus; Gage--Hamilton and Grayson also developed the curvature flow (or curve-shortening flow), which may be viewed as the gradient flow of the length functional. It has the property that, if a smooth simple closed curve undergoes the curvature flow, it remains smoothly embedded without self-intersections.

Existence of at least one geodesic for arbitrary closed Riemannian manifolds was finally proved by Lusternik--Fet in 1951-1952 \cite{LF51,Fet53}. Their approach was based on Morse theory; and indeed the problem of finding geodesics was the initial motivation for Morse himself. Geodesics are the critical points of the energy functional on the loop space. Moreover, the space $\mathcal{L}M$ of parametrized closed curves on $M$ cannot be retracted into the subspace $\mathcal{L}^0M$  of homotopically trivial closed curves, and Lusternik--Schnirelmann theory applies to give a critical point outside of $\mathcal{L}^0M$. 

Even though the loop space of a manifold is infinite dimensional, if the manifold is compact then the energy functional satisfies the compactness condition of Palais--Smale, which in practice means that it behaves as a Morse function on a finite-dimensional manifold. However, the main difficulty in this approach is that each geodesic can be iterated, and this corresponds to distinct points in the loop space. Distinguishing two geometrically distinct geodesics is a subtle, hard problem. 

So far, all the above methods provide only finitely many geodesics, so how about infinitely many? In this direction, another beautiful idea due to Birkhoff, for a Riemannian $S^2$, is that of an annulus global surface of section; we have of course seen this in previous sections. One considers a closed geodesic $\gamma$ (which Birkhoff proved to exist via the minmax argument explained above), dividing the sphere in an upper and a lower hemishpere. One then considers vectors along $\gamma$ which point towards the upper hemisphere (this is an annulus) as initial values of geodesics, starts shooting orbits along these vectors, and considers the first return map. However, for this annulus to be a global surface of section, one needs that no geodesic gets ``trapped'' in the upper hemisphere (this will be satisfied for example when the Gaussian curvature is strictly positive). Moreover, one needs to further check the twist condition at the boundary in order to apply the Poincaré--Birkhoff theorem. Here, note that Birkhoff only stated the existence of at least two fixed points, but a simple argument which Birkhoff seems to have overlooked was provided by Neumann \cite{N77}, thus obtaining infinitely many periodic points (not related by iterations); this is the version of the Poincar\'e--Birkhoff theorem we stated above. In the case where we do have a well-defined Birkhoff map, what if the return map does not twist? This is where the theorem due to Franks from 1992 \cite{Fr92} that we mentioned above (which is a statement about the open annulus), comes into play; he obtained infinitely many geodesics on $S^2$ for this case. In the case where the Birkhoff annulus is not a global section and so there is no return map, an argument of Bangert from 1993 \cite{Ban93} shows that, if geodesics get trapped, they need to do so around a small ``waist'' (a ``short'' geodesic), or more formally, geodesics with no conjugate points. Moreover, he shows that the existence of a waist forces the existence of infinitely many geodesics. One key observation is that the Birkhoff return map sends a point on the boundary (lying on a geodesic) to its second conjugate point along this geodesic, and so some of the ideas where already present in Birkhoff's work. This filled in the general case, finally (after almost 90 years) obtaining the existence of infinitely many geodesics for an arbitrary metric on $S^2$. We further mention that in 1993 Nancy Hingston, building on work of other people (see \cite{Hi93} and references therein), also provided a full proof of a quantitative estimate on the growth of the number of geodesics with respect to length; if $N(l)$ is the number of geodesics with length at most $l$ then $N(l)\gsim l/log(l)$, i.e.\ the same growth rate as prime numbers.

One should further mention that Katok \cite{K73} (see also Ziller's account \cite{Z83}) has famously constructed examples of non-reversible Finsler metrics on $S^n, \mathbb{C}P^n$ with only finitely many closed geodesics. For instance, the case of $S^2$ can be described as the round geodesic flow, but on a frame rotating along the $z$-axis with irrational angle of rotation (and the metric is arbitrarily close to the round one); so that the only closed geodesics are the equator in both directions. This example shows that the general Finsler case is very different from the Riemannian case, an hence the $\mathbb{Z}_2$-action which allows to reverse geodesics should be used in a significant way in order to obtain infinitely many geodesics.

Another celebrated result in this story is that of Gromoll--Meyer 1969 \cite{GM69}: if the sequence of Betti numbers of the free loop space $\mathcal{L}M$ of $M$ is unbounded, then $M$ admits infinitely many geodesics (for \emph{any} metric). Morse had previously, in his 1932 book ``Calculus of variations in the large'' (although unfortunately with mistakes), computed the homology of $\mathcal{L}M$ in the non-degenerate case. For this, one may use a spectral sequence whose terms in the $E^1$-page consists of the homology of the base (constant loops) and the homology of each geodesic, endowed with a local coefficient system, and degree shifted by the Morse index.
Note that nondegeneracy is in the Morse--Bott sense, since we can always reparametrize loops (which we consider unoriented) via the action of $O(2)$ on $S^1$, and so we see one circle for each orientation in this homology group. Another ingredient is Bott's famous iteration formula for the index \cite{B56}, which implies that $\mu(\gamma^m)$ grows linearly with $m$. When combined with the homology computation via the above Morse--Bott spectral sequence, one sees that if the set of primitive geodesics is finite, then the Betti numbers of $\mathcal L M$ are bounded, and hence the result by Gromoll--Meyer follows in the non-degenerate case. The degenerate case, roughly speaking, is obtained by the fact that every degenerate orbit is the limit of a \emph{finite} number of non-degenerate ones, and contributes to the homology in a bounded index window. 

This leaves the question of when the sequence of Betti numbers of $\mathcal{L}M$ is unbounded. In \cite{VS76}, Vigué-Poirrier--Sullivan show, via the above result and algebraic calculations, that if $M$ has finite fundamental group, then the Betti numbers of $\mathcal{L}M$ are unbounded if and only if $H_*(M;\mathbb Q)$ requires at least $2$ generators as a ring. Ziller proves this holds for symmetric
spaces of rank $> 1$ \cite{Z77}. This covers many cases, but it leaves out many important ones e.g.\ $S^n, \mathbb{R}P^n, \mathbb{C}P^n, \mathbb{H}P^n, CaP^2$. 

On the other hand, one can consider the case of a \emph{generic} metric (or ``bumpy'', i.e.\ for which all geodesics are non-degenerate). For such a case, on any manifold with finite fundamental group, Gromov has also shown the following quantitative estimate: there exist constants $a,b$ such that $N(l)\geq \frac{a}{l}\sum_{i=1}^{bl} b_i(\mathcal{L}M)$. Rademacher \cite{R89} has shown the existence of infinitely many geodesics for bumpy metrics on manifolds with finite fundamental group. This result builds on work of Klingenberg--Takens \cite{KT72}, Klingenberg \cite{K78}, who reduced to the case where all orbits are hyperbolic; and Hingston \cite{Hi84}, who covered the bumpy case for $S^n, \mathbb{R}P^n, \mathbb{C}P^n, \mathbb{H}P^n, CaP^2$, under the hyperbolic-orbits-only assumption.

One therefore clearly sees that, while a ``simpler'' problem than finding closed orbits in the three-body problem, finding infinitely many closed geodesics is significantly complicated. This is a problem that has inspired enormous amounts of work, has spanned most of the 20th century, and still is not known in the general case. Indeed, it is still an open question whether any Riemannian metric on a given closed simply connected manifold admits
infinitely many closed geodesics. In particular, it is unknown for $S^n, n\geq 3$, for a general metric.

\medskip

\textbf{Remarks on Floer theory, and modern symplectic geometry.} As we have seen, symplectic geometry is the geometry of classical mechanics, dealing with Hamiltonians and their associated evolution equations, and in particular closed Hamiltonian orbits of period $1$. In this context, Arnold \cite{Ar65} proposed his famous conjecture on the minimal number of such orbits for a non-degenerate Hamiltonian on a closed symplectic manifold $M$: there should be at least as many as the sum of the Betti numbers of $M$. This is naturally related to the classical Morse inequalities. It is notable that Arnold proposed this conjecture \emph{as a version of the Poincaré--Birkhoff theorem} (here, note that the sum of Betti numbers of the annulus is $2$).

It was from this conjecture that one of the cornerstones of the modern methods of symplectic geometry was introduced; namely, Floer theory. Together with the introduction of holomorphic curves due to Gromov in 1985 \cite{Gro85}, these two developments form the building bricks of the symplectician's toolkit and daily musings. 

The approach of Floer to the Arnold conjecture \cite{F86,F88a,F88b,F89a,F89b,F89c} is again based on the ideas of Morse theory. Indeed, one can view Hamiltonian orbits as the critical points of a suitable action functional on the loop space, in such a way that flow-lines correspond to cylinders satisfying an elliptic PDE (the Floer equation). One defines a differential which counts these solutions, and the resulting homology theory is actually isomorphic to the Morse homology of the underlying manifold, so that the Arnold conjecture follows. Floer proved it under some technical assumptions, i.e.\ symplectic asphericity, and the symplectic Calabi--Yau condition; these have been lifted after work of several authors (Ono \cite{O95}, Hofer--Salamon \cite{HS95}, Liu--Tian \cite{LT98}, Fukaya--Ono \cite{FO99},...), at least for the case of rational coefficients. The technical details are very difficult, needing the introduction of virtual techniques. While Abouzaid--Blumberg \cite{AB21}, amongst other results, prove the Arnold conjecture with coefficients on every finite field by appealing to stable homotopy theory, very recent work of Bai--Xu derived the Arnold conjecture over the integers, by building global Kuranishi charts on the moduli spaces of Floer solutions.

As we have seen, a special case of closed Hamiltonian orbits is that of Reeb orbits in a contact-type level set. Since every contact manifold is contact-type in some symplectic manifold (i.e.\ its symplectization), one can view the problem of finding closed Reeb orbits as an odd-dimensional version of the Hamiltonian problem. In this setting, an important statement related to the Arnold conjecture is the \emph{Weinstein conjecture}, which claims the existence of at least one closed Reeb orbit for any contact form on a given contact manifold. Recalling that geodesic flows are particular cases of Reeb flows, this includes the statement that every Riemannian metric admits a closed geodesic (proved by Lusternik--Fet, as mentioned above). In dimension three, it was established by Taubes \cite{T07} (based on Seiberg--Witten theory), thus culminating a large body of work by several people extending over more than two decades. There are also further striking results in dimension $3$, e.g.\ Irie's results on equidistribution of closed orbits in the generic case \cite{I15, I18}, or the ``$2$ or infinitely many'' dichotomy for torsion contact structures \cite{CHP19}. This dichotomy uses the combination of Brouwer and Frank's theorem as discussed above as the fixed point theorem, and Hutching's embedded contact homology (ECH) to find the disk-like global surface of section; and so fits in well with the basic two-step approach by Poincar\'e. Irie's results rely on the relationship between volume and ECH capacities as proved by Cristofaro-Gardiner--Hutchings--Ramos \cite{CGHR15}. In higher dimensions, though
there are several partial results (e.g.\ \cite{AH,FHV90,HV89,HV92,V87}), the Weinstein conjecture is still open.

While the Arnold conjecture is stated for closed symplectic manifolds, a natural class of symplectic manifolds with non-empty boundary is that of Liouville domains. There is an associated Floer theory for such manifolds, which goes under the name of \emph{symplectic homology}. The first version of such theory was due to Floer--Hofer \cite{FH94}, and can be traced to the Ekeland--Hofer capacities and their relation to early versions of $S^1$-equivariant symplectic homology\footnote{This was discussed at the opening lectures by Hofer and Floer in Fall 1988 at the symplectic program at the MSRI Berkeley, although unfortunately is written nowhere. Hofer gave a lecture on capacities and the $S^1$-equivariant symplectic homology at a conference in Durham in 1989, whose proceedings are published in \cite{LMS89}, and contains the non-equivariant part of the story. I thank Hofer for these clarifications.}; see also section~5 in \cite{H89} for an even previous and non-equivariant version, called symplectology. There is also a version due to Viterbo \cite{V99, V18} (see also \cite{CO18,BO09} for more recent versions), who showed that symplectic homology of a cotangent bundle is the homology of the free loop space of the base, a bridge between the classical story of finding geodesics, and the modern Floer-theoretic approach (see also \cite{AS,SW06}). 

In the Liouville setting, as opposed to the closed setting, the difference is that the associated Floer theory recovers not only the homology of the manifold, but also dynamical data at the contact-type boundary (i.e.\ closed Reeb orbits). Of course, one of the motivations for such a theory is the Weinstein conjecture, at least for those contact manifolds which bound a Liouville domain (i.e.\ Liouville fillable ones). Heuristically, if the symplectic homology is infinite-dimensional or zero, then there is at least one orbit at the boundary (since the homology of the manifold is finite dimensional and non-zero, although, strictly speaking, here we need consider the case of finite-type Liouville domains; see e.g.\ \cite{O04} for a nice survey, containing these and related ideas). 

The Arnold conjecture is a statement about \emph{fixed} points (or $1$-periodic orbits) of Hamiltonian maps, and predicts a finite number of such. On the other hand, one could want to estimate the number of \emph{periodic} points (recall the same situation for the Poincaré--Birkhoff theorem, whose original version predicted $2$ fixed points, although one can also obtain infinitely many periodic points, as was observed after Birkhoff). The analogous statement for Hamiltonian or Reeb flows is the \emph{Conley conjecture}. Roughly speaking, for a ``vast'' collection of closed symplectic manifolds, every Hamiltonian map has infinitely many simple periodic orbits and, moreover, simple periodic orbits of unbounded minimal period whenever the fixed points are isolated. This was proved by Ginzburg for closed symplectically aspherical symplectic manifolds in \cite{Gi10} (see \cite{GG15} and references therein, for a survey and history of the problem; and \cite{GG19} for what the author understands is the current state of the art). One of the key inputs is a special class of critical points introduced by Hingston, and
later called symplectic degenerate maxima/minima (SDM) by Ginzburg. The presence of an SDM forces the
existence of infinitely many closed orbits (cf.\ \cite{Hi93,Hi97} for the case of geodesics on $S^2$).

\smallskip

We conclude this section with the following (clearly debatable but rather convincing from the above story) meta-mathematical claim: \emph{the three-body problem inspired large portions of modern symplectic geometry}. In all probability, it would also be fair to make the same claim for most of the modern theory of dynamical systems.

\medskip

\textbf{Final remark on different approaches.} Amongst the approaches that we have discussed (by all means non-exhaustive) we point out that the advantage of KAM theory (in the pertubative case), when compared to more abstract approaches via general fixed point theorems, is that in favourable situations one can localize periodic (or quasi-periodic) orbits in bounded regions of phase-space, and obtain better qualitative information on these. This is, of course, much more complicated in non-perturbative situations, where rigorous numerics is usually the preferred approach. See \cite{FvK18} for examples of return maps on a disk-like global surface of section, obtained numerically, for the planar problem. 

\medskip

\textbf{More references.} A nice basic introduction to the classical KAM theorem is e.g.\ \cite{W08}. Another very nice exposition on the basics behind Mather theory is e.g.\ \cite{S15}. A beautiful and very detailed account on the three-body problem and Poincar\'e's work are the notes by Chenciner \cite{Ch15}. A  very recent and detailed survey on open questions on geodesics, illustrating the vastness and richness of their search, is that of Burns and Matveev \cite{BM13}. I also based parts of the above brief survey on very nice lectures by Nancy Hingston given at the summer school ``Current Trends in Symplectic Topology'', July 2019, at the Centre de recherches mathématiques, Université de Montréal, Canada; where I happened to be in the audience. Of course, this is a classical story and there are plenty of other sources; see e.g.\ Oancea's much more detailed account \cite{O14} and references therein (as well as the appendix due to Hrynewicz on the story for $S^2$), with a view towards symplectic geometry.

\section{The advent of contact topology in the CR3BP}\label{sec:advent} The next result opens up the possibility of using modern techniques from contact and symplectic geometry on the CR3BP (holomorphic curves, Floer theory,...). 

Denote by $\overline{\Sigma}_c^E$ and $\overline{\Sigma}_c^M$ the bounded components of the Moser-regularized energy hypersurfaces for the spatial problem and $c<H(L_1)$, and let $\overline{\Sigma}_c^{E,M}$ be the connected sum bounded component, for $c\in (H(L_1),H(L_2))$. Similarly, use $\overline{\Sigma}_{P,c}^E$, $\overline{\Sigma}_{P,c}^M$ and $\overline{\Sigma}_{P,c}^{E,M}$ for the case of the planar problem. 

\begin{thm}[\cite{AFvKP12} (planar problem), \cite{CJK18} (spatial problem)]\label{thm:contact_type} If $c<H(L_1)$, the Moser-regularized energy hypersurfaces $\overline{\Sigma}_c^E, \overline{\Sigma}_c^M,\overline{\Sigma}_{P,c}^E, \overline{\Sigma}_{P,c}^M$ are all contact-type. The same holds for $\overline{\Sigma}_c^{E,M},\overline{\Sigma}_{P,c}^{E,M}$, if $c\in (H(L_1),H(L_1)+\epsilon)$ for sufficiently small $\epsilon>0$. As contact manifolds, we have
$$
\overline{\Sigma}_c^E\cong\overline{\Sigma}_c^M\cong (S^*S^3,\xi_{std}), \mbox{ if } c<H(L_1),
$$
$$
\overline{\Sigma}_{P,c}^E\cong\overline{\Sigma}_{P,c}^M\cong (S^*S^2,\xi_{std}), \mbox{ if } c<H(L_1),
$$
and 
$$
\overline{\Sigma}_c^{E,M}\cong (S^*S^3,\xi_{std})\#(S^*S^3,\xi_{std}), \mbox{ if } c\in (H(L_1),H(L_1)+\epsilon).
$$
$$
\overline{\Sigma}_{P,c}^{E,M}\cong (S^*S^2,\xi_{std})\#(S^*S^2,\xi_{std}), \mbox{ if } c\in (H(L_1),H(L_1)+\epsilon).
$$
In all above cases, the planar problem is a codimension-$2$ contact submanifold of the spatial problem. In particular, the dynamics of the CR3BP, in the low-energy range (see paragraph below) and near the primaries, is given by a Reeb flow.
$\hfill \square$
\end{thm}
Recall that the above just means that there exists a Liouville vector field which is transverse to the regularized level sets; in fact, this is just the fiber-wise Liouville vector field in regularized coordinates. The regularized level sets, as contact manifolds, are standard and well-known, so not very interesting from a geometrical perspective. However, their interest lies in the given non-standard \emph{dynamics} for the underlying standard geometry. The Hamiltonian dynamics for the problem now becomes the Reeb dynamics, and the planar problem (from a dynamical perspective rather than a geometric one) is actually invariant under the Reeb flow. We will refer as the \emph{low-energy range} to the interval $(-\infty,H(L_1)+\epsilon)$ of energies $c$ for which the above result holds. 

Conjecturally, the contact condition should extend all the way to $H(L_2)$. We will take this up this point later in Section \ref{sec:open_books}, when we discuss open books for the CR3BP.

\begin{remark}[\textbf{Weinstein handles}]
In the above statement, the connected sum is to be interpreted in the contact category; this amounts to attaching a \emph{Weinstein} $1$-handle to the disjoint union of two copies of $(S^*S^3,\xi_{std})$. Roughly speaking, this means removing two Darboux balls and identifying their boundaries via attaching a $1$-handle, which is endowed with the extra structure of a symplectic form which glues well to the symplectization form of the standard contact form at the boundary of each ball. The result is a \emph{Liouville/Weinstein cobordism} having $(S^*S^3,\xi_{std})\bigsqcup(S^*S^3,\xi_{std})$ at the negative end, and $(S^*S^3,\xi_{std})\#(S^*S^3,\xi_{std})$ at the positive one. Note that here the terms positive/negative are relevant: the Liouville vector field is outwards/inwards pointing at the corresponding boundary components, respectively, and so these cobordisms are oriented. This is always the local Morse-theoretical picture for a non-degenerate index $1$ critical point of a Hamiltonian (as is the case of $L_1$). To learn about Weinstein manifolds, see e.g.\ \cite{CE12}; this source also provides deep connections between this notion and that of Stein manifolds.
\end{remark}

\textbf{References.} For a very detailed and well-exposed overview of contact geometry and holomorphic curves in the planar case of the CR3BP, we refer to Frauenfelder--van Koert \cite{FvK18}.

\section{Open books in the spatial CR3BP}\label{sec:open_books}

In this section, we present results of the author, in co-authorship with Otto van Koert \cite{MvKa20}. The main direction is to generalize the approach of Poincar\'e in the planar problem (i.e.\ Steps (1) and (2) outlined above) to the \emph{spatial} problem. In this section, we will discuss a complete generalization of Step (1), while a possible approach to Step (2) is discussed in Chapter \ref{ch:Ham_twist_maps}, where a generalization of the classical Poincar\'e--Birkhoff theorem is presented.

\subsection{Step (1): Global hypersurfaces of section.} 

We first state a structural result, which provides the basic architecture and scaffolding for the problem.

\medskip

In the spatial CR3BP, fix a mass ratio $\mu\in [0,1]$. Denote a connected, bounded component of the regularized energy level set with energy $c$, by $\overline{\Sigma}_c$, which by Theorem \ref{thm:contact_type} is contact-type if $c$ is in the low energy range. We also denote by $\overline{\Sigma}_c^P\subset \overline{\Sigma}_c$ the (codimension-$2$) regularized planar problem. 

By Giroux correspondence, we know there exist supporting open book decompositions on $\overline{\Sigma}_c,\overline{\Sigma}_c^P$ for such $c$ in the low energy range. However, as we emphasized already, this correspondence does not give adapted open books whenever the dynamics is fixed. The content of the following result is that the \emph{given} dynamics of the spatial CR3BP in the low-energy range, and near the primaries, is given by a contact form which is a Giroux form for some concrete open book.

\begin{thm}[Moreno--van Koert \cite{MvK20a}]\label{thm:openbooks}

If $c$ lies in the low-energy range, $\overline{\Sigma}_c$ is of contact-type and admits a supporting open book decomposition for energies $c<H(L_1)$ that is adapted to the dynamics. 
Furthermore, if $\mu<1$, then there is $\epsilon>0$ such that the same holds for $c\in (H(L_1),H(L_1)+\epsilon)$.
The open books have the following abstract form:
\[
\overline{\Sigma}_c
\cong
\begin{cases}
(S^*S^3,\xi_{std})=\textbf{OB}(\mathbb{D}^*S^2,\tau^2), & \mbox{if }c<H(L_1) \\
(S^*S^3,\xi_{std})\#(S^*S^3,\xi_{std})=\textbf{OB}(\mathbb{D}^*S^2 \natural \mathbb{D}^*S^2,\tau_1^2 \circ \tau_2^2), & \mbox{if }c\in (H(L_1),H(L_1)+\epsilon), \mu<1.
\end{cases}
\]
In all cases, the binding is the planar problem 

\[
B=\overline{\Sigma}^P_c=
\begin{cases}
(S^*S^2,\xi_{std}), & \mbox{if }c<H(L_1) \\
(S^*S^2,\xi_{std})\#(S^*S^2,\xi_{std}), & \mbox{if }c\in (H(L_1),H(L_1)+\epsilon), \mu<1.
\end{cases}
\]
\end{thm}

\begin{figure}
    \centering
    \includegraphics[width=0.4 \linewidth]{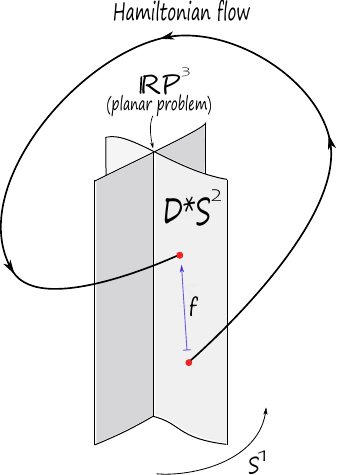}
    \caption{The open book for $\overline{\Sigma}_c$, with $c<H(L_1)$, and the return map $f$.}
    \label{fig:openbook}
\end{figure}

Here, $\mathbb{D}^*S^2$ is the unit cotangent bundle of the $2$-sphere, $\tau$ is the positive Dehn--Seidel twist along the Lagrangian zero section $S^2\subset \mathbb{D}^*S^2$, and $\mathbb{D}^*S^2 \natural \mathbb{D}^*S^2$ denotes the boundary connected sum of two copies of $\mathbb{D}^*S^2$. 
The monodromy of the second open book is the composition of the square of the positive Dehn--Seidel twists along both zero sections (they commute). See Figure \ref{fig:openbook} for an abstract representation. 

\smallskip

\begin{figure}
    \centering
    \includegraphics{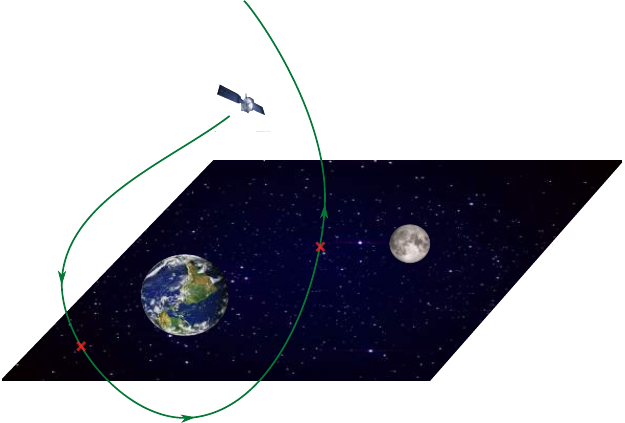}
    \caption{Theorem \ref{thm:openbooks} admits a physical interpretation: away from collisions, the orbits of the negligible mass point intersect the plane containing the primaries transversely. This is intuitively clear from a physical perspective (and it follows from the fact that the planar problem is invariant), and translates (after regularization) to the fact that the ``pages'' $\{q_3=0,p_3>0\}$, $\{q_3=0,p_3<0\}$ of the ``physical'' open book are global hypersurfaces of section outside of the collision locus. Unfortunately this does not extend continuously to the latter, as explained in Figure \ref{fig:verticaloribts}. The binding is the planar problem.}
    \label{fig:transversalitypic}
\end{figure}

\begin{remark}
The following remarks are in order.
\begin{itemize}
    \item (\textbf{Reduction}) The above result reduces the study of a continuous dynamics on a $5$-dimensional manifold, to the study of the discrete dynamics of the return map on a $4$-manifold. Moreover, the topology of this section is completely understood, and which is in a concrete sense the simplest possible topology, as the open book may not be destabilized (recall that stabilization increases the complexity of the page of the open book). 

\medskip
    
    \item (\textbf{Non-perturbative}) Theorem \ref{thm:openbooks} holds for $c$ in the whole low-energy range, and therefore it is an inherently \emph{non-perturbative} result. A heuristical reason is the following. while in the planar case finding the invariant subset, i.e.\ the binding, is non-trivial (the search for the direct and retrograde orbits indeed has a long history), the invariant subset in the spatial case is immediately obvious; it is the planar problem. Moreover, one can work with global coordinates. This is elementary but quite involved, and will be discussed, in some detail, in Section \ref{sub:basic_idea} below.

\medskip

    \item (\textbf{Explicitness}) The open book is explicit, i.e.\ it can be written down in coordinates. This makes it amenable for numerical work. The difficulty is that the hypersurface of section is $4$-dimensional, which makes visualization of the return map harder.

    \item (\textbf{Energy range validity}) The open book in the above result actually exists for energy below $H(L_2)$ (and near the primaries). As mentioned in Section \ref{sec:advent}, the contact condition should also extend to $H(L_2)$. This would mean that the associated Poincar\'e return map is a symplectomorphism also in this energy range.
\end{itemize}
\end{remark}
The technique of proof does not rely on holomorphic curves, since one can directly write down the open book explicitly; it is rather elementary, but the calculations are very involved.

\subsection{The basic idea}\label{sub:basic_idea} Theorem \ref{thm:openbooks} is motivated by the following observation. We consider a Stark--Zeemaan system satisfying Assumptions (A1) and (A2). In unregularized (or physical) coordinates, we put 
$$
B_u:=\{(\vec q,\vec p)\in H^{-1}(c)~|~q_3=p_3=0 \},
$$
the planar problem. Its normal bundle is trivial, and we have the following map to $S^1$: 
\begin{equation}\label{pinonreg}
\pi_u: H^{-1}(c) \setminus B_u \longrightarrow S^1,
(\vec q,\vec p) \longmapsto \frac{q_3+ip_3}{\Vert q_3+ip_3\Vert}.
\end{equation}
We will refer to this map as the \emph{physical} open book. We consider the angular $1$-form
$$
d\pi_u:=\frac{\Omega^u_p}{p_3^2+q_3^2},
$$
where
\begin{equation}\label{omegapunreg}
\Omega^u_p=p_3 dq_3 -q_3 dp_3,
\end{equation}
is the unregularized numerator. We need to see whether $d\pi_u(X_H)$ is non-negative, and vanishes only along the planar problem.

From Equation \eqref{Hamvf}, we have
\begin{equation}\label{omegaXH}
d\pi_u(X_H)=\frac{p_3^2+q_3^2\left( \frac{g}{\Vert \vec q\Vert^3}+\frac{1}{q_3} \frac{\partial V_1}{\partial q_3}(\vec q)\right)}{p_3^2+q_3^2}.
\end{equation}

Note that Assumption (A2) implies that $\frac{\partial V_1}{\partial q_3}(\vec q)=a q_3+O(q_3^2)$ near $q_3=0$, and so $\frac{1}{q_3} \frac{\partial V_1}{\partial q_3}(\vec q)$ is well-defined at $q_3=0$. In order for the above expression to satisfy the required non-negativity condition, we impose the following.

\begin{assumptions}

(A3) We assume that the function 
    $$
    F(\vec q)=\frac{g}{\Vert \vec q\Vert^3}+\frac{1}{q_3} \frac{\partial V_1}{\partial q_3}(\vec q)
    $$ is everywhere positive.
\end{assumptions}
Note that it suffices that the second summand be non-negative.

\begin{remark}
In the CR3BP, from Equation (\ref{potential}), we obtain
$$
\frac{\partial V_1}{\partial q_3}(\vec q)=q_3\frac{1-\mu}{\Vert \vec q - \vec e\Vert^3},
$$
and therefore the corresponding expression in Equation (\ref{omegaXH}) is non-negative, vanishing if and only if $p_3=q_3=0$.
\end{remark}

\begin{figure}
    \centering
    \includegraphics{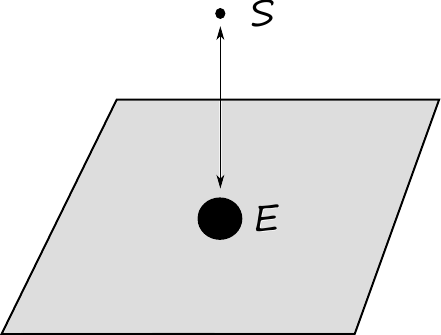}
    \caption{For the RKP, there exist (regularized) collision orbits which are periodic and ``bounce'' vertically over a primary, always staying on the region $q_3>0$ (or $q_3<0$). We call them the \emph{polar} orbits. This means that the ``pages'' $\{q_3=0,p_3>0\}$, $\{q_3=0,p_3<0\}$ are \emph{not} transverse to the regularized dynamics, as these orbits always stay on the same side of the planar problem.}
    \label{fig:verticaloribts}
\end{figure}

The obvious problem of the above computation is that it a priori does not extend to the collision locus, and indeed it cannot (see Figure \ref{fig:verticaloribts}). In fact, one needs to interpolate with the \emph{geodesic} open book described in Section \ref{sec:geodesic}, which is well-behaved near the collision locus. This creates an interpolation region where fine estimates are needed, and this is the main difficulty in the proof; we refer to \cite{MvK20a} for the details. 

\subsection{Symmetries.} Recall that the natural symmetry group of the CR3BP is given by $\mathbb Z_2\oplus \mathbb Z_2$. There is a symplectic involution of $(\mathbb{R}^6,dp\wedge dq)$ given by
$$
r:\mathbb R^6\rightarrow \mathbb R^6,
$$
$$
r(q_1,q_2,q_3,p_1,p_2,p_3)= (q_1,q_2,-q_3,p_1,p_2,-p_3).
$$
We also have the anti-symplectic involutions
$$
\rho_1,\rho_2:\mathbb R^6\rightarrow \mathbb R^6, 
$$
$$
\rho_1(q_1,q_2,q_3,p_1,p_2,p_3) =(q_1,-q_2,-q_3,-p_1,p_2,p_3),
$$
$$
\rho_2(q_1,q_2,q_3,p_1,p_2,p_3) =(q_1,-q_2,q_3,-p_1,p_2,-p_3),
$$
satisfying the relations $\rho_1\circ \rho_2=\rho_2 \circ \rho_1=r$, and so generating the abelian group $\{1,r,\rho_1,\rho_2\}\cong \mathbb{Z}_2 \oplus \mathbb{Z}_2$.

After regularization, the symplectic involution admits the following intrinsic description. Consider the smooth reflection $R:S^3 \rightarrow S^3$ along the equatorial sphere $S^2\subset S^3$. Then $r$ is the physical transformation it induces on $T^*S^3$, given by
$$
r: T^*S^3\rightarrow T^*S^3
$$
$$
r(q,p)=(R(q),(d^*_qR)^{-1}(p)).
$$
This map preserves the unit cotangent bundle $S^*S^3$. The maps $\rho_1,\rho_2$ also have regularized versions. The following emphasizes the symmetries present in our setup.

\begin{proposition}[\cite{MvK20a}]\label{prop:symmetries}
Let $c<H(L_1)$, and consider the symplectic involution $r:S^*S^3\rightarrow S^*S^3$. The open book decomposition on $\overline{\Sigma}_c=\mathbf{OB}(\mathbb{D}^*S^2,\tau^2)$ is symmetric with respect to $r$, in the sense that
$$
r(P_\theta)=P_{\theta+\pi},\;\;\mbox{Fix}(r)=B=\overline{\Sigma}_c^P.
$$
Moreover, the anti-symplectic involutions preserve $B$ and satisfy
$$
\rho_1(P_\theta)=P_{-\theta},\;\rho_2(P_\theta)=P_{-\theta + \pi}.
$$
In particular, $\rho_1$ preserves $P_0$ and $P_\pi$, whereas $\rho_2$ preserves $P_{\pi/2}$ and $P_{-\pi/2}$.  
\end{proposition}

In other words, the open book is compatible with all the symmetry group $\mathbb{Z}_2 \oplus \mathbb{Z}_2$.

\subsection{The return map.}\label{sec:return_map} First, we recall a standard definition. We say that a symplectomorphism $f:(M,\omega)\rightarrow (M,\omega)$ is \emph{Hamiltonian} if $f=\phi_K^1$, where $K: \mathbb{R}\times M\rightarrow \mathbb{R}$ is a smooth (time-dependent) Hamiltonian, and $\phi_K^t$ is the Hamiltonian isotopy it generates. This is defined by $\phi_K^0=id$, $\frac{d}{dt}\phi_K^t=X_{K_t}\circ \phi_K^t$, and $X_{H_t}$ is the Hamiltonian vector field of $H_t$ defined via $i_{X_{H_t}}\omega=dH_t$. Here we write $K_t=K(t,\cdot)$.

In the spatial CR3BP, for $c<H(L_1)$, and after fixing a page $P=\pi^{-1}(1)$ of the corresponding open book, Theorem \ref{thm:openbooks} implies the existence of a Poincar\'e return map $f:\mbox{int}(P)\rightarrow \mbox{int}(P)$. Moreover, as in Proposition \ref{prop:symplecto}, we can consider the $2$-form $\omega$ obtained by restriction to $P$ of $d\alpha$, where $\alpha$ is the contact form on $\overline{\Sigma}_c$ for the spatial problem, whose restriction to the binding $\alpha_P$ is the contact form for the planar problem. Recall that $\omega$ is symplectic only along the \emph{interior} of $P$. Moreover, we have a smooth identification int$(P)\cong \mbox{int}(\mathbb{D}^*S^2)$, giving a symplectomorphism $G:\mbox{int}(P)\rightarrow \mbox{int}(\mathbb{D}^*S^2)$ on the interior which extends smoothly to the boundary $B$, but its inverse $G^{-1}$, although continuous at $B$, is \emph{not} differentiable along $B$ since $\omega$ becomes degenerate there. After conjugating $f$ with $G$ and considering $\widetilde\omega=G_*\omega$, we obtain a symplectomorphism  $\widetilde{f}:=G \circ f \circ G^{-1}:(\mbox{int}(\mathbb{D}^*S^2),\widetilde\omega)\rightarrow (\mbox{int}(\mathbb{D}^*S^2),\widetilde \omega)$, where $\widetilde\omega$ is a Liouville filling of $(B,\alpha_P)$. In particular, $\widetilde\omega$ is non-degenerate at $B$. This phenomenon will be discussed in detail in Section \ref{sec:degenerate}, i.e.\ the page of the open book is a \emph{degenerate} Liouville domain, as defined there.

\begin{thm}[Moreno--van Koert \cite{MvK20a}]\label{thm:returnmap}
For every $\mu \in (0,1]$, $c<H(L_1)$, the associated Poincar\'e return map $f$ extends smoothly to the boundary $\partial P$, and in the interior it is an exact symplectomorphism 
$$
f=f_{c,\mu}:(\mbox{int}(P),\omega)\rightarrow (\mbox{int}(P),\omega),
$$
where $\omega=d\alpha$ (depending on $c,\mu$). Moreover, $f$ is Hamiltonian in the interior. 

After conjugating with $G$, $\widetilde{f}$ extends continuously to the boundary, is Hamiltonian in the interior, and the Liouville completion of $\widetilde\omega$ is symplectomorphic to the standard symplectic form $\omega_{std}$ on $T^*S^2$. 
\end{thm}

\begin{figure}
    \centering
    \includegraphics{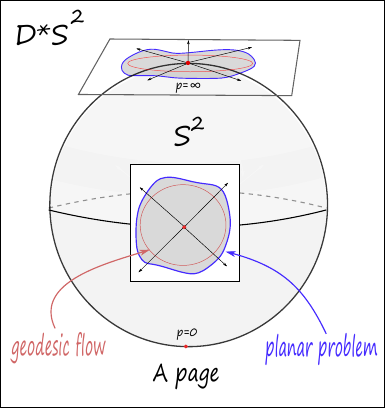}
    \caption{A page of the open book as a symplectic filling of the planar problem, viewed as a fiber-wise star-shaped domain in $T^*S^2$. The geodesic flow corresponds to the unit cotangent bundle.}
    \label{fig:pageasfilling}
\end{figure}

The fact that $f$ is an exact symplectomorphism follows from Proposition \ref{prop:symplecto}. The fact that $f$ extends to the boundary is non-trivial, and relies on second order estimates near the binding: it suffices to show that the Hamiltonian giving the spatial problem is positive definite on the symplectic normal bundle to the binding. This nondegeneracy condition can be interpreted as a convexity condition that plays the role, in this setup, of the notion of dynamical convexity due to Hofer--Wysocki--Zehnder \cite{HWZ98}. Note that if a continuous extension exists, then by continuity it is unique. 

The fact that $f$ is Hamiltonian in the interior follows from:
\begin{enumerate}
    \item The monodromy of the open book is Hamiltonian (here, the Hamiltonian is allowed to move the boundary);
    \item The fact that the return map $f$ is, under the above convexity assumption, always symplectically isotopic to a representative of the monodromy (after correcting the symplectic form), via a boundary-preserving isotopy;
    \item $H^1(P;\mathbb{R})=0$, so that every symplectic isotopy is Hamiltonian.
\end{enumerate}

\section{Iterated picture}\label{sec:iterated_picture}

We now describe a richer topological picture for the regularized low energy level sets in the CR3BP, in terms of \emph{iterated planar} open books. We focus on the case $c<H(L_1)$. Theorem \ref{thm:openbooks} provides an open book for the spatial problem, which exists for all values of $\mu$, of the abstract form $(S^*S^3,\xi_{std})=\mathbf{OB}(\mathbb D^*S^2,\tau^2)$. We have also already mentioned the open book from \cite{HSW19} for the planar problem, which exists for all values of the convexity range, of abstract form $(S^*S^2,\xi_{std})=\mathbf{OB}(\mathbb D^*S^1,\tau_P^2)$ where $\tau_P$ is the Dehn twist on the annulus (it is an open question whether the latter open book extends over the whole low-energy range). While these open books are in some sense independent of each other, they are also inter-related. The aim for this section is to describe precisely how they interact with each other, and what is the expected topological picture (conjecturally for the whole low energy range). Roughly speaking, under the presence of the spatial open book, the presence of the planar open book forces an \emph{iterated} structure, described as follows.

\subsection{Iterated planar open books} Topologically and abstractly, the situation may be understood as follows. The Stein manifold $\mathbb{D}^*S^2$ carries a Lefschetz fibration structure, whose smooth fibers are the annuli $\mathbb{D}^*S^1$, and its monodromy is precisely $\tau^2_P$ along the vanishing cycle $S^1 \subset \mathbb{D}^*S^1$ (the square of the standard Dehn twist $\tau_P$). 
We write $\mathbb{D}^*S^2=\mathbf{LF}(\mathbb{D}^*S^1,\tau_P^2)$. 
By restricting this Lefschetz fibration to the boundary, we obtain the above open book for $\mathbb{R}P^3=\mathbf{OB}(\mathbb D^*S^1,\tau_P^2)$. The Lefschetz fibration on the pages $\mathbb{D}^*S^2$ gives $(\mathbb S^*S^3,\xi_{std})$ the structure of an \emph{iterated planar} contact $5$-manifold, which has been studied in \cite{Acu, AEO20, AM18}. To motivate its definition, we note the following result of Giroux--Pardon.

\begin{thm}[\cite{GP}]
 Up to deformation equivalence, every Weinstein domain $W$ is obtained by attaching critical Weinstein handles to a stabilization $W_0\times \mathbb{D}^2$ of another Weinstein domain $W_0$ (the \emph{fiber}) along a collection of Lagrangian spheres.
\end{thm}
In other words, every Weinstein domain admits a Weinstein Lefschetz fibration over $\mathbb{D}^2$ with Weinstein fibers. Once we are given an exact symplectic Lefschetz fibration on a Weinstein domain, one can further construct an exact symplectic Lefschetz fibration on the codimension two Weinstein fiber, and iterate this process until the fiber is $4$-dimensional. This idea motivates the following. 

\begin{definition}[\cite{Acu07}] \label{IPLF}
An \textit{iterated planar Lefschetz fibration} $f: (W^{2n}, \omega) \to \mathbb{D}^2$ on a $2n$-dimensional Weinstein domain $(W^{2n}, \omega)$ is an exact symplectic Lefschetz fibration satisfying the following properties: 

\begin{enumerate}

\item There exists a sequence of exact symplectic Lefschetz fibrations $f_i:(W^{2i}, \omega_i) \to \mathbb{D}^2$ for $i=2, \dots, n$ with $f=f_n$.

\item The total space $(W^{2i}, \omega_{i})$ of $f_{i}$ is a regular fiber of $f_{i+1}$, for $i=2, \dots, n-1$.

\item $f_2: (W^4, \omega_2) \to \mathbb{D}^2$ is a planar Lefschetz fibration, i.e. the regular fiber of $f_2$ is a genus zero surface with nonempty boundary, which we denote by $W^2$.
\end{enumerate}

\end{definition} 

If $f: W \rightarrow \mathbb{D}^2$ is an iterated planar Lefschetz fibration, then the boundary of $W$ inherits an open book decomposition whose pages are diffeomorphic to the regular fibers of $f$. The following definition is motivated by looking at the open book decomposition induced by the boundary restriction of an iterated planar Lefschetz fibration. 

\begin{definition}[\cite{Acu07}] \label{IPOB}
An \textit{iterated planar (IP) open book decomposition} of a contact manifold $(M, \xi)$ is an open book decomposition $(M, \xi)=\textbf{OB}(W, \varphi)$ whose page $W$ admits an iterated planar Lefschetz fibration.
\end{definition}

We can then give the main definition.

\begin{definition}[\cite{Acu07}] \label{IPCM}
For any $n>1$, an \textit{iterated planar (IP) contact manifold} $(M, \xi)$ is a $(2n+1)$-dimensional contact manifold supported by an open book decomposition whose Weinstein page admits an iterated planar Lefschetz fibration.
\end{definition}

The relevant examples for the CR3BP are the contact manifolds $(S^*S^n,\xi_{std})$ and their connected sums. More generally, we have the following.

\subsection{Examples: $A_k$-singularities} For $n\geq 2$, the co-disk bundle $\mathbb{D}^{*}S^n$ admits an iterated planar Lefschetz fibration where each regular fiber is $\mathbb{D}^{*}S^{n-1}$, and the Lefschetz fibration on $\mathbb{D}^{*}S^2$ is planar with fibers $\mathbb{D}^{*}S^1=[0,1] \times S^1$; see Proposition \ref{prop:Brieskorn_model}. This gives an IP structure to the contact manifolds $(S^*S^n,\xi_{std})=\mathbf{OB}(\mathbb D^*S^{n-1},\tau^2)$, which appear in the CR3BP.

More generally, consider the $A_k$-singularity given by $$A_k= \{(z_1, \dots, z_n)\in \mathbb{C}^n \mid z_1^2+\dots+z_{n-1}^2+z_{n}^{k+1}=1\} \subset (\mathbb{C}^n, \omega_{std})$$ for $n\geq 3$ and $k\geq2$. The $A_k$-singularity can be expressed as a \emph{plumbing} of $k$ copies of $\mathbb{D}^{*}S^{n-1}$. Moreover, the Lefschetz fibration on the $A_k$-singularity, defined by the projection onto the last coordinate $z_n$, is $\mathbb{D}^{*}S^{n-1}$. This observation together with the existence of an iterated planar Lefschetz fibration on $\mathbb{D}^{*}S^{n-1}$ imply that the $A_k$-singularity admits an iterated planar Lefschetz fibration.\smallskip

\begin{remark}
    We should remark that in \cite{Acu07} (see also \cite{AM}) the Weinstein conjecture is proven for all IP contact manifolds, i.e.\ every contact form admits a periodic orbit. In particular, it follows from Theorem \ref{thm:contact_type} that the CR3BP (in the low energy range and near the primaries) always admits a periodic orbit. This, however, follows from an older result of Hofer--Viterbo \cite{HV89}, who prove the Weinstein conjecture for contact-type hypersurfaces in contangent bundles.
\end{remark}

We finish by stating the following result, which follows directly from the above discussion, and Theorem \ref{thm:openbooks}.

\begin{thm}[\textbf{IP open books on the spatial CR3BP}]\label{thm:IP} For $c$ in the low energy range and near the primaries, the Moser-regularized energy level set $\overline{\Sigma}_c$ has the structure of an IP contact $5$-manifold, endowed with a concrete supporting IP open book decomposition adapted to the dynamics. 
\end{thm}

\subsection{Iterating beyond the CR3BP} The iterated planar picture which arises in the CR3BP can be of course be iterated indefinitely to higher dimensions, well beyond the CR3BP, at least at the topological level. The full (abstract) iterated picture that becomes apparent is summarized in the following diagram.

\begin{center}\label{diag:reg_pic}
    \begin{tikzcd}[scale cd=0.75]
        S^1=\mathbf{OB}(D^*S^0,\tau_0)\arrow[r, hook] \arrow[d,hook]& S^3=\mathbf{OB}(D^*S^1,\tau_1) \arrow[r,hook]\arrow{d}{2:1}[swap]{\mathcal{L}} & S^5=\mathbf{OB}(D^*S^2,\tau_2) \arrow[r,hook]& S^7=\mathbf{OB}(D^*S^3,\tau_3) \; \cdots\\
        S^*S^1=\mathbf{OB}(D^*S^0,\tau_0^2) \arrow[ur,hook]\arrow[r, hook]\arrow[d,"2:1"] & S^*S^2=\mathbf{OB}(D^*S^1,\tau_1^2) \arrow[ur,hook]\arrow[d]\arrow[r,hook]\arrow[bend right=60,swap,rightarrowtail]{u}{2:1} 
        & S^*S^3=\mathbf{OB}(D^*S^2,\tau_2^2)\arrow[ur,hook]\arrow[d]\arrow[r,hook]\arrow[u,rightarrowtail, "2:1"]
        &S^*S^4=\mathbf{OB}(D^*S^3,\tau_3^2)\arrow[d]\arrow[u,rightarrowtail, "2:1"]\;  \cdots
        \\
        S^1 \arrow[r,hook]& S^2=\mathbb{C}P^1 \arrow[r,hook] & S^3 \arrow[r,hook] & S^4=\mathbb{H}P^1 \; \cdots\\
    \end{tikzcd}
\end{center}   

We now explain the diagram in more detail. The bottom horizontal arrows are inclusions of equators, and the bottom vertical arrows are the natural projection maps $S^*S^n\rightarrow S^n$. The bottom squares clearly commute. The map $\tau_n: D^*S^n\rightarrow D^*S^n$ is the Dehn-Seidel twist, and $\mathcal{L}$ is the 2:1 cover induced by the Levi-Civita regularization. Recall the Brieskorn variety 
$$
V_n=\left\{z=(z_0,\dots,z_n) \in \mathbb
C^{n+1}:\sum_{i=0}^n z_i^2=0\right\}\cong T^*S^n
$$
and the corresponding Brieskorn manifold
$$
\Sigma_n=V_n\cap S^{2n+1}=\left\{z=(z_0,\dots,z_n) \in \mathbb
C^{n+1}:\sum_{i=0}^n z_i^2=0, \sum_{i=0}^n\vert z_i\vert ^2=1\right\}\cong S^*S^n.
$$

We have the open book $\Sigma_n=\mathbf{OB}(D^*S^{n-1},\tau_{n-1}^2)$ given by
$$
\pi_n: \Sigma_n\backslash \Sigma_{n-1}\rightarrow{S^1}
$$
$$
z=(z_0,\dots,z_n)\mapsto \frac{z_n}{\vert z_n\vert},
$$
with binding $\Sigma_{n-1}=\Sigma_n\cap \{z_n=0\}\cong S^*S^{n-1}$. The middle horizontal arrows in the diagram are then inclusions of bindings of these open books into the ambient manifold.

Denoting $z=(\underline z,z_n)\in \mathbb C^n\times \mathbb C$, $\xi=(\underline \xi,\xi_n)\in \mathbb{R}^n\times \mathbb R$, $\eta=(\underline \eta,\eta_n)\in \mathbb{R}\times \mathbb R^n$, following \cite{Or69}, we have a 2:1 branched cover
$$
\Phi: \Sigma_n\rightarrow S^{2n-1}\subset \mathbb C^n,
$$
given by 
$$
z=(\underline z,z_n)\mapsto \frac{\underline z}{\sqrt{1-\vert z_n\vert^2}},
$$
branched along $\Sigma_{n-1}=\Sigma_n\cap \{z_n=0\}=S^*S^{n-1}$, the binding of $\pi_n$. These 2:1 branched covers are the top vertical arrows in the diagram, and the middle horizontal inclusions are also the inclusions of the branching loci.

Moreover, the sphere $S^{2n-1}=\mathbf{OB}(D^*S^{n-1},\tau_{n-1})=\left\{(w_1,\dots,w_n)\in \mathbb C^n: \sum_{i=1}^n \vert w_i \vert^2=1 \right\}$ comes with the Milnor open book
$$
\Pi_n: S^{2n-1}\backslash \Sigma_{n-1}\rightarrow S^1
$$
$$
(w_1,\dots,w_n)\mapsto \frac{w_1^2+\dots+w_n^2}{\vert w_1^2+\dots+w_n^2 \vert},
$$
with binding $\Sigma_{n-1}=S^{2n-1}\cap\left\{w_1^2+\dots+w_n^2=0\right\}$ and monodromy $\tau_{n-1}$. Note that the inclusion 
$$
i_n: S^{2n-1}\hookrightarrow S^{2n+1}
$$
$$
(w_1,\dots,w_n)\mapsto (w_1,\dots,w_n,0)
$$
is compatible with the Milnor open books, i.e.
$$
\Pi_{n}=\Pi_{n+1} \circ i_n,
$$
and so in particular maps the binding $\Sigma_{n-1}\cong S^*S^{n-1}$ of $\Pi_n$ to the binding $\Sigma_n\cong S^*S^n$ of $\Pi_{n+1}$ via the natural inclusion. Moreover, the binding of $\pi_n$ is mapped by $\Phi$ to the binding of $\Pi_n$. This discussion implies the commutativity of the diagram (after removing the arrow $\mathcal{L}$).

All of the open books in the diagram are IP open books, supporting the corresponding standard contact structures. Note moreover that the top horizontal arrows are also inclusions of bindings of open books, namely the trivial open books
$$
S^{2n+1}\backslash S^{2n-1} \rightarrow S^1
$$
$$
(w_0,\dots,w_n)\mapsto \frac{w_n}{\vert w_n\vert},
$$
of abstract type $S^{2n+1}=\mathbf{OB}(\mathbb D^{2n},\mathrm{Id})$.

It would be interesting to investigate whether there are regularization schemes for Hamiltonian systems in arbitrary degrees of freedom which fit into the above iterated picture, with explicit systems admitting adapted open books of the type described. A natural candidate is the classical KS-regularization \cite{KS}, which is obtained from the Levi--Civita regularization by replacing the complex numbers with the quaternions.


\section{Return map for the RKP}\label{sec:RKP_returnmap} The aim of this section is to explicitly study the return map for integrable limit case of the CR3BP given by the RKP, corresponding to $\mu=0$. We will give a sketch of the proof of the following result, following the exposition in  Appendix A in \cite{MvK20a}.

\begin{thm}[\cite{MvK20a}, Integrable case]\label{thm:integrablecase} In the RKP, for the low energy range and near the primary, the return map preserves the annuli fibers of a concrete symplectic Lefschetz fibration of abstract type $\mathbb{D}^*S^2=\mathbf{LF}(\mathbb{D}^*S^1,\tau_P^2)$, where it acts as a classical integrable twist map on regular fibers, and fixes the two (unique) nodal singularities on the singular fibers. The boundary of each of the symplectic fibers coincides with the direct/retrograde planar circular orbits (a Hopf link in $\mathbb{R}P^3$).
\end{thm}

The two fixed points are the north and south poles of the zero section $S^2$, and correspond to the two periodic collision orbits bouncing on the primary (one for each of the half-planes $q_3>0$, $q_3<0$), which we call the \emph{polar} orbits; see Figure \ref{fig:verticaloribts}. In order to prove this, we will derive an explicit formula for the return map, as follows.

Recall that, in unregularized coordinates, the RKP is described by the Hamiltonian $H=K+L$, where
$$
K= \frac{1}{2} \Vert p \Vert^2  -\frac{1}{\Vert q\Vert},
\quad
L= q_1 p_2 - q_2 p_1.
$$
The regularized Hamiltonian is $Q(\xi,\eta)=\frac{1}{2}f^2(\xi,\eta)\Vert \eta \Vert^2$ where $$f(\xi,\eta)=1+(1-\xi_0)(-c-1/2+\xi_2\eta_1-\xi_1\eta_2).$$ An explicit adapted open book for the low energy RKP is given by the geodesic open book 
$$
(
\xi, \eta) \longmapsto \frac{\xi_3 +i\eta_3}{\Vert \xi_3+i\eta_3 \Vert}.
$$
Indeed, we have the following.
\begin{lemma}\label{app:lemmarotKepler}
The geodesic open book is a supporting open book for the rotating Kepler problem for $c<-3/2$.
\end{lemma}

\begin{proof}
Note that the pages of the geodesic open book, which is adapted to the Kepler problem $K$, are also invariant under the Hamiltonian flow of $L$ (which acts by rotation along the $(\xi_0,\eta_0)$-axis inside a given page). Recall that $\phi^t_H=\phi^t_K \circ \phi^t_L$. While $\phi^t_L$ leaves the pages invariant, the flow $\phi^t_K$ is transverse to them. This implies the claim. \end{proof}

In order to study the return map, we consider the page $$
P= \left\{ (\xi;\eta) \in T^*S^3:\;Q(\xi,\eta)=\frac{1}{2},\; \xi_3=0,\; \eta_3 \geq 0 \right\},
$$ which is the easiest to visualize (it is analogous to the $0$-page in Figure \ref{fig:birkhoffannulus}). 

The return time in unregularized coordinates is simply the minimal Kepler period for the corresponding Kepler energy $K$ (which is preserved under the flow of $H$). By Kepler's third law, this return time depends only $K$, and is given by 
$$
T = T(K)=\frac{\pi}{2(-K)^{3/2}}.
$$
From $\phi_K^T=id$, we obtain $\phi^T_H=\phi^T_K \circ \phi^T_L=\phi^T_L$, and therefore the return mo is given by
$$
R:P \rightarrow P
$$
$$
R(q_1,q_2,q_3;p_1,p_2,0)=
\left( \mathrm{Rot}_{\frac{\pi}{2(-K)^{3/2}}}(q_1,q_2),q_3; \mathrm{Rot}_{\frac{\pi}{2(-K)^{3/2}}}(p_1,p_2),0\right),
$$
where $\mathrm{Rot}_\phi$ is the rotation by angle $\phi$. 
This is generated by the Hamiltonian $L$ restricted to the global hypersurface of section (as the time $T$-map). To obtain an explicit formula for a Hamiltonian generating $R$ in time-$1$, using the relation $K+L=c$, we write 
$$
R=\phi^{X_L}_{T(K)}=\phi^{T(K)X_L}_{1}=\phi^{T(c-L)X_L}_1.
$$
We see that there is a function $g(L)$ such that $R=\phi^{X_{g(L)}}$, by noting that $X_{g(L)}=g'(L)X_{L}$. With $g'(L)=\frac{\pi}{2(L-c)^{3/2}}$, we can compute $g(L)$ as
$$
g(L)=-\pi\left( 2(L-c)\right)^{-1/2}.
$$
We may now describe the return map in the Moser regularized coordinates, which is given by
$$
R:P\rightarrow P
$$
$$
R(\xi_0,\xi_1,\xi_2,0;\eta_0,\eta_1,\eta_2,\eta_3)=
\left(\xi_0,\mathrm{Rot}_{T(c-L)}(\xi_1,\xi_2),0;\eta_0,\mathrm{Rot}_{T(c-L)}(\eta_1,\eta_2),\eta_3\right).
$$
The Hamiltonian $L$ is given in these coordinates by $L=\xi_2 \eta_1 -\xi_1 \eta_2$.

Now we can see that the two polar orbits are correspond to the fixed points $x_\pm=(\pm 1,0,0,0;0,0,0,\eta_3)$ of the return map, where $\eta_3>0$ is such that they lie on the level set $Q^{-1}(\frac{1}{2})$. The fixed point $x_-$ corresponds to the southern polar orbit; in unregularized coordinates it is the point in the $q_3$-axis that is maximally far from the origin.
The fixed point $x_+$ corresponds to the northern polar orbit, and this fixed point corresponds to the periodic collision point. This orbit and nearby periodic orbits for $\mu\gtrsim 0$ were studied by Belbruno in \cite{B81}.

We can now appeal to the fact from Proposition \ref{prop:Brieskorn_model}, i.e.\ that the Brieskorn model 
$$
V=\left\{(z_0,z_1,z_2)\in \mathbb C^{3}: z_0^2+z_1^2+z_2^2=1\right\}
$$
for $T^*S^3$ admits the concrete Lefschetz fibration
$$
\pi: V\rightarrow \mathbb C, 
$$
$$
\pi(z_0,z_1,z_2)=z_0,
$$
of abstract type $\mathbf{LF}(T^*S^1,\tau_P^2)$. Pulling back to $T^*S^2$ via the symplectomorphism $\psi$ of Proposition \ref{prop:Brieskorn_model}, we obtain the Lefschetz fibration
$$
\Theta=\pi \circ \psi^{-1}: T^*S^2\rightarrow \mathbb C
$$
$$
\Theta(\xi;\eta)=\sqrt{\frac{1+\sqrt{1+4\Vert\eta\Vert^2}}{2}} \xi_0+i \sqrt{ \frac{2}{1+\sqrt{1+4\Vert\eta\Vert^2}}} \eta_0.
$$
Note that the return map preserves the expression defining the Lefschetz fibration, i.e.\ 
\begin{equation}\label{eq:invariance}
 \Theta \circ R =\Theta.   
\end{equation}
This is not the end of the story, however, as the page $P$, while diffeomorphic to $\mathbb D^*S^2$, is not symplectomorphic to $T^*S^2$. The symplectic form in its \emph{interior} needs to be modified to make it symplectomorphic to $T^*S^2$, which can be done by a deformation in the $\eta_3$ coordinate near the boundary; only after this modification does Equation \ref{eq:invariance} makes sense on the interior of $P$. Moreover, the fibers of the resulting fibration are asymptotic to the equator $S^1\subset S^2$ transversed in both directions, and these orbits do not necessarily agree with the direct/retrograde orbits (this holds only for the Kepler problem), as they lie over circles in the upper hemisphere of $S^2$ which are translates of the equator. A further modification of the fibers of the Lefschetz fibration near their boundary is needed, to make them asymptotic to these orbits, while keeping the fibers symplectic. We refer the reader to Appendix B in \cite{Mvk20a} for details.

The return map acts on the cylinder fibers by a rotation on each circle that makes them, with angle $T(c-L)$, which varies from circle to circle. The direct orbit is rotated with an angle $T(c-L_{dir})$, while the retrograde, with an angle $T(c-L_{ret})$, which point in opposite directions. This is the staple of a classical integrable twist map, and finishes the sketch of the proof.

\section{Digression: degenerate Liouville domains and billiards}\label{sec:degenerate} The discussion in Section \ref{sec:return_map} touches on a fundamental \emph{feature} of return maps arising on adapted open books, which is also a feature of \emph{billiards} and their associated billiard maps. Namely, the degeneration of the symplectic form at the boundary of a global hypersurface of section. This can be fixed by a conjugation which is only continuous at the boundary, as in Theorem \ref{thm:returnmap}. However, even if the return map extends smoothly to the boundary, its conjugation will only extend continuously. In other words, either we choose $\omega$ to be symplectic everywhere and the return map to be only $C^0$ at the boundary, or we choose the return map to extend smoothly, paying the price that $\omega$ degenerates at the boundary. Getting around this fact is one of the main technical difficulties when trying to use methods from Floer theory to bear into the setup. This situation can be captured by the following definitions.

We now introduce the notion of a \emph{degenerate} Liouville domain, which roughly speaking is a Liouville domain away from the boundary, but the symplectic form degenerates at the latter.

\begin{definition}\label{def:deg_Liouville}
A \emph{degenerate} Liouville domain is a pair $(W,\lambda,\alpha_B)$ where:
\begin{itemize}
    \item $W$ is a smooth manifold with boundary $B=\partial W$;
    \item $\alpha_B$ is a contact form on $B$; and
    \item $\lambda$ is a $1$-form such that $\omega=d\lambda$ is symplectic in the interior int$(W)$, but degenerates at the boundary $B$ along a normal direction.
\end{itemize}
The last condition means that there is a collar neighbourhood $C=[0,1]\times B\subset W$ of the boundary, with collar coordinate $r\in [0,1]$ such that $B=\{r=1\}$, along which $\lambda=A\cdot\alpha_B$, where $A=A(b,r)$ for $b\in B$, $r \in [0,1]$, is a smooth function satisfying
\begin{itemize}
    \item $\partial_rA>0$ for $r<1$;
    \item $\partial_rA\vert_{r=1}\equiv 0$;
    \item $A\vert_{r=1}\equiv 1$.
\end{itemize}
The $1$-form $\lambda$ is a \emph{degenerate} Liouville form.

We will say that $\lambda$ is \emph{convex} at the boundary if
\begin{itemize}
\item $\partial_r^2A\vert_{r=1}>0$. 
\end{itemize}

\end{definition}

Note that the notion of \emph{ideal} Liouville domains due to Giroux \cite{Gir20} is different from the above. Indeed, in the ideal case, the Liouville form has a pole at the boundary, which is not the case here. Moreover, while degenerate Liouville domains have finite symplectic volume, ideal ones have infinite volume. 

Sources of degenerate Liouville domains are pages of open book decompositions, and domains of billiard maps, as follows.

\subsection{Example: open books} Let $(M,\xi=\ker \alpha)$ be a contact manifold where $\alpha$ is adapted to an open book decomposition $(B,\theta)$ of abstract type $M=\mathbf{OB}(W,\phi)$. Let $\omega=d\lambda$ with $\lambda=\alpha\vert_{\mbox{int}(W)}$. Then $\omega\vert_B$ is degenerate, as the Reeb vector field of $\alpha$ is tangent to $B$. Moreover, \cite[Proposition 3.1]{DGZ} provides a neighbourhood $B\times \mathbb{D}^2\subset M$ of the binding in which $$\alpha=A(\alpha_B+s^2 d\theta),$$ where $(s,\theta)$ are polar coordinates, for a smooth positive function $A$ satisfying $A\equiv 1$ along $s=0$ and $\partial_sA<0$ for $s>0$. The contact condition also implies that $A=A(s)$ only depends on $s$, and the fact that $R_\alpha$ is tangent to $B$ implies $\partial_sA\vert_{s=0}=0$. Then $\lambda=A(s)\alpha_B$, and by changing coordinates to $r=1-s$, we see explicitly that $(W,\omega)$ is a degenerate Liouville domain. The return map $$f: (\mbox{int}(W),\omega)\rightarrow (\mbox{int}(W),\omega)$$ may not extend to the boundary, but it does whenever $\lambda$ is convex at the boundary; see \cite{MvK20a}.  

\subsection{Example: billiards} Let $D\subset \mathbb R^2,$ $D\cong \mathbb D^2$ be a planar and strictly convex billiard table. Recall that the associated billiard is the dynamical system on $S^*D=D\times S^1$ which follows the straight line defined by the initial condition until it hits the boundary of the table, then reflects according to the standard law from optics (i.e.\ that the angle of incidence equals the angle of reflection), and continues as before. The associated Birkhoff annulus $$W\cong [0,\pi]\times S^1\cong\mathbb D^*S^1$$ is the set of all unit vectors pointing inwards along $\partial D$, which can be parametrized by two variables $\theta \in [0,\pi]$ and $\varphi \in S^1$. $W$ admits the two-form $$\omega=d(-\cos(\theta)d\varphi)=\sin(\theta)d\theta\wedge d\varphi,$$ so that $(W,\omega)$ is a degenerate Liouville domain. The associated symplectic billiard map $$f:(\mbox{int}(W),\omega)\rightarrow (\mbox{int}
(W),\omega),$$ mapping a vector to the next vector along a billiard trajectory, extends smoothly to the boundary as the identity, and preserves $\omega$. In a different description where the $2$-form is non-degenerate at the boundary, the billiard map is only continuous at the boundary; we will present an intrinsic description below.

\subsection{Degeneration and non-degeneration.} We now make the trade-off we alluded to more precise.

A degenerate Liouville domain $(W,\omega=d\lambda)$ can be turned into a Liouville domain $(W,\omega_Q=d\lambda_Q)$, which is unique up to Liouville isotopy, i.e.\ deformation of the Liouville form. We will call $(W,\omega_Q)$, the \emph{non-degeneration} of $(W,\omega)$. Similarly, a Liouville domain has an associated \emph{degeneration}. More precisely:

\begin{lemma}\textbf{(Degeneration and non-degeneration)}\label{lemma:degeneration} We have the following.
\begin{itemize}
    \item[(1)]\textbf{(From degenerate to non-degenerate)} Let $(W,\lambda,\alpha_B)$ be a degenerate Liouville domain. Then there exists a boundary-preserving homeomorphism $Q:W\rightarrow W$ which is smooth in $\mbox{int}(W)$ and only continuous at $B=\partial W$, such that:

\medskip

    \begin{itemize}
        \item[$\bullet$]  $(W,\lambda_Q=Q^*\lambda)$ is a Liouville domain with strict contact-type boundary $(B,\alpha_B)$. This Liouville domain is unique up to Liouville isotopy. 

\medskip
        
        \item[$\bullet$]   Given a smooth map $f:(W,\lambda,\alpha_B)\rightarrow (W,\lambda, \alpha_B)$, then $f_Q=Q \circ f \circ Q^{-1}$ is smooth in the interior, where it preserves $\omega_Q$, but only continuous at $B$.
    \end{itemize}

\medskip
    
\item[(2)] \textbf{(From non-degenerate to degenerate)} Let $(W,\lambda)$ be a Liouville domain with strict contact boundary $(B,\alpha_B)$. Then, there exists a boundary-preserving homeomorphism $S=Q^{-1}:W\rightarrow W$ which is smooth but its inverse $Q$ is as above, such that $(W,\lambda_S=S^*\lambda, \alpha_B)$ is a degenerate Liouville domain. 
   
\end{itemize}

\end{lemma}

Here, $S$ is called a \emph{squaring} map, and $Q$ is called a \emph{square root} map. Both will simply be reparametrizations in directions normal to the boundary. It follows that if $f$ is generated by $H_t$ on the interior, $f_Q$ is generated by $H_t^Q=H_t\circ Q$ on the interior, but the isotopy is not-well defined at the boundary. 

To put it simply, the trade-off expressed in Lemma \ref{lemma:degeneration} is that in general we can either choose a coordinate description of the setup in which the $2$-form is symplectic but the map only continuous at the boundary, or alternatively where the $2$-form degenerates at the boundary but the map is smooth (e.g.\ as in the CR3BP, see below). 

\begin{remark}\label{rk:Hams} Below, we will assume that Hamiltonians in the non-degenerate picture are $C^1$, whereas Hamiltonians in the degenerate picture are only $C^0$, and are related by the formula $E_t = H_t \circ Q$.
\end{remark}

\begin{proof}[Proof of Lemma \ref{lemma:degeneration}] Take a collar $[0,1]\times B$ where $\lambda=A(s)\alpha_B$ is as in Definition \ref{def:deg_Liouville}. We consider a map $Q: W \rightarrow W$, which is the identity away from $[0,1]\times B\subset W$, and on $(0,1]\times B$ is of the form $Q(s,b)=(\varphi(s,b),b)$ where $\varphi$ solves the ODE
$$
\varphi^\prime(b,s)=-\frac{1}{\partial_sA(\varphi(s,b))}>0, \mbox{ for } s>0,
$$
which then has a smooth solution for $s>0$. Integrating this equation with respect to $s$, we see that
$$
A(\varphi(b,s),s)=1-s,
$$
for $s>0$, by choosing the integration constant to be $1$. By continuity, we see that this equation has a unique solution satisfying $\varphi(0)=0$, which is only continuous at $s=0$. Then $Q$ is a diffeomorphism in the interior of $W$, but it extends only continuously to the boundary. 

If we define $\lambda_Q=Q^*\lambda, \omega_Q=Q^*\omega=d\lambda_Q$, then $\omega_Q$ is a symplectic form satisfying $\omega_Q=d((1-s)\alpha_B)$. If we make the change of coordinates $r=1-s$, we see that $\omega_Q=d(r\alpha_B)$ near $B=\{r=1\}=\{s=0\}$, and therefore $(W,\omega_Q)$ is a Liouville domain with strict contact-type boundary $(B,\alpha_B)$. Note that the maps $Q,Q^{-1}$ are given in the $r$-coordinate as
$$
Q(r,b)=\left(F(r,b),b\right),\;Q^{-1}(r,b)=(A(r,b),b),
$$
where $F(r,b)=1-\varphi(1-r,b)$. 

The only choice we made above consisted of the coordinates in which $\lambda=A\alpha_B$. Any other choice of such coordinates differs by a reparametrization in the $s$-direction, which doesn't change the Liouville isotopy class of the resulting domain. The statement about self-maps is obvious.

For the converse, write $\lambda=(1-s)\alpha_B$ where $s=1-r$ on $[0,1]\times B$. Let $S:W\rightarrow W$ which is the identity away from $[0,1]\times B\subset W$, and on $(0,1]\times B$ is of the form $S(s,b)=(s^2,b)$. Then $\lambda_S=S^*\lambda=(1-s^2)\alpha_B$ defines a degenerate Liouville form. Note that to remove it, the corresponding map $Q$ is $Q(s)=\sqrt{s}$. This finishes the proof.
\end{proof}

\subsection{Billiards as fiber-wise star-shaped domains} What follows is an intrinsic way of thinking about the non-degeneration of the domain of billiard maps. 

We consider the geodesic flow $\phi_t$ on $T\mathbb{R}^n$, given by $\phi_t(x,v)=(x+tv,v)$, which is the Hamiltonian flow of $H:T\mathbb{R}^n\rightarrow \mathbb R$, $H(x,v)=\frac{1}{2}\Vert v \Vert^2$. We can view $H$ as a moment map of the $\mathbb{R}$-action on $T\mathbb{R}^n$ given by $t\cdot (x,v)=\phi_t(x,v)$. Therefore the symplectic quotient
$$
Q=H^{-1}(1/2)\backslash\mathbb{R}=S\mathbb R^n\backslash\mathbb R
$$
inherits a symplectic structure $\omega_{red}$, defined by $\pi^*\omega_{red}=\omega\vert_{H^{-1}(1/2)}$, where $\pi: H^{-1}(1/2)\rightarrow Q$ is the quotient map, and $\omega=dx\wedge dv=d(x\cdot dv)$ is the standard symplectic form on $T\mathbb{R}^n$. Note that $Q$ is naturally identified with the set of oriented lines in $\mathbb{R}^n$. Moreover, we have the following.
\begin{lemma}
$Q$ is exact symplectomorphic to $T^*S^{n-1}$ with its standard Liouville form.
\end{lemma}
\begin{proof}
We may write 
$$
T^*S^{n-1}=\{(q,p)\in \mathbb{R}^n\oplus \mathbb{R}^n: \vert q \vert =1, p\cdot q=0\},
$$
with Liouville form given by $p\cdot dq$. Consider the map
$$
F: Q \rightarrow T^*S^{n-1}
$$
$$
l\mapsto (q,p),
$$
which maps an oriented line $l$ to $(q,p)$, where $q$ is the unit vector which directs it, and $p$ is the unique vector through the origin which lies in $l$ and is orthogonal to $q$, i.e.\ $p$ is uniquely determined by $q\cdot p=0$ and $l=\{t.q+p: t \in \mathbb{R}\}$. The map $F$ is clearly a diffeomorphism. Note that in terms of the $(x,v)$ coordinates, it is given by $p=x$ and $q=v$, i.e.\ the roles of the coordinates are swapped. The Liouville form $x\cdot dv$ then corresponds to $p\cdot dq$ under $F$. This finishes the proof.
\end{proof}

Now, consider a billiard table $B\subset \mathbb{R}^n,$ $B\cong \mathbb{D}^n$, with smooth convex boundary $C=\partial B\cong S^{n-1}$. We assume without loss of generality that the origin lies in the interior of $B$. Consider further the open subset $U_B\subset Q$ consisting of those oriented lines in $\mathbb{R}^n$ which point inwards to $B$ along $C$. Then $U_B$ inherits a symplectic from given by $\omega_B=\omega_{red}\vert_{U_B}$. 

Moreover, $U_B$ is naturally identified with a fiber-wise convex domain $D\subset T^*S^{n-1}$, $D\cong \mathbb{D}^*S^{n-1}$, under the diffeomorphism $F$ in the proof of the above lemma. Indeed, if $q$ is the unit vector that directs an oriented line $l$, we consider the hyperplane $H_q\subset \mathbb{R}^n$ orthononal to $q$, $H_q\cong \mathbb{R}^{n-1}$, which intersects $B$ in a $n-1$-dimensional disk $D_q=H_q\cap B\cong \mathbb D^{n-1}$. Then $H_q$ is identified with the cotangent fiber $T^*_qS^{n-1}$ under $F$, and $D_q\subset H_q$ is identified with a convex disk by convexity of $B$; we may then define $D$ as $D=\bigsqcup_{q\in S^{n-1}} D_q$, which we can view as a fiber-wise convex domain in $T^*S^{n-1}$. We further have that $U_B=\bigsqcup_{q\in S^{n-1}} U_{B,q}$, where $U_{B,q}$ consists of those lines in $U_B$ directed by $q$, and so the identification $U_B\cong D$ follows by identifying $U_{B,q}$ with $D_q$. 

Note that $B$ is uniquely determined by $D$, and viceversa (as long as $D$ is obtained from a billiard table), and so a convex billiard table can be indeed thought of as a fiber-wise convex domain in $T^*S^{n-1}$, and in particular, as a Liouville domain $(D,\omega_D)$. With this identification, the associated billiard map, mapping a unit vector pointing inwards to that corresponding to the next billiard segment, is then a map $f: (D,\omega_D)\rightarrow (D,\omega_D),$ which extends continuously to the boundary as the identity. We see that this is then the non-degenerate description of the setup, which is moreover intrinsic.

\subsection{A simple example: the billiard on a ball}\label{sec:billiard_ball} We consider the toy case of the billiard map on a ball in $\mathbb R^{n}$. As noted above, the billiard map is a map $T:\mathbb D^*S^{n-1}\rightarrow \mathbb D^*S^{n-1}$, and we shall give explicit expressions. The following computations are due to Otto van Koert. 

Let us focus first in the case $n=2$, i.e.\ a disk in the plane. The associated billiard map is simply
$$
f: \mathbb{D}^*S^1 \rightarrow \mathbb D^*S^1
$$
$$
f(\varphi,\theta)=(\varphi +2\theta, \theta)
$$
preserving $\omega = \sin (\theta)d\varphi \wedge d\theta$, and extending smoothly to the boundary. Note that this is a Dehn twist, which is a completely integrable twist map. In order to show it is Hamiltonian (in the interior), noting that it is generated by the vector field $2\theta \partial_\varphi$, we need to solve the equation $i_{2\theta\partial_\varphi}\omega=2\theta \sin(\theta)d\theta=dH,$ which is solved by the Hamiltonian
$$
H(\theta)= 2(-\theta \cos(\theta) + \sin(\theta)),
$$
i.e.\ $f=\phi_1^H$. Note that the Hamiltonian extends smoothly to the boundary. If we now conjugate to the standard symplectic form (i.e.\ we perform the non-degeneration) $\omega_0 = d\varphi\wedge d\delta$, this is achieved with the change of coordinates $-\cos(\theta)=\delta$, i.e.\ $\theta=\mathrm{arccos}(-\delta)=\pi - \mathrm{arccos}(\delta)\in [-1,1]$. The billiard map in these new coordinates is then
$$
T : \mathbb D^*S^1\rightarrow \mathbb D^*S^1
$$
$$
T(\varphi, \delta)=  (\varphi - 2 \mathrm{arccos}(\delta), \delta).
$$
Recalling that $\mathrm{arccos}'(\delta)=-\frac{1}{\sqrt{1-\delta^2}},$ which explodes at $\delta=\pm 1$, we clearly see that $T$ is only continuous at the boundary. The Hamiltonian generating $T$ is 
$$
\widetilde H(\delta)=H(\pi - \mathrm{arccos}(\delta))=2(\mathrm{arccos}(-\delta)\delta+\sin(\mathrm{arccos}(\delta))),
$$
whose derivative is $\partial_\delta \widetilde H=2\mathrm{arccos}(-\delta)$, i.e.\ it is $C^1$ a the boundary, but not $C^2$.

For the higher-dimensional case, we can use the geodesic flow on the round sphere. The billiard map $T$ itself actually acts non-trivially on homology (for even $n$), as so it is
not Hamiltonian. However, the square $\tau = T^2$ \emph{is} Hamiltonian. Indeed, writing
$$
\mathbb D^*S^{n-1}=\left\{(q,p)\in \mathbb R^n \oplus \mathbb R^2: \vert q \vert =1, p\cdot q=0, \vert p \vert \leq 1\right\},
$$
we have
$$
\tau: \mathbb D^*S^{n-1} \rightarrow \mathbb D^*S^{n-1},
$$
$$
\tau (q,p) = \phi^R_{
f(\vert p \vert)}(q, p),
$$
where $R$ is the vector field generating the geodesic flow, and
$$
f(r) = 2\pi - 4 \mathrm{arccos}(r).
$$
The Hamiltonian generating $\tau$ is
$$
\widetilde{H}(q,p) = \widetilde{H}(\vert p \vert)=(2\pi - 4 \mathrm{arccos}(\vert p \vert))\vert p \vert + 4\sqrt{1-\vert p \vert ^2},
$$
which is smooth for $\vert p \vert < 1$, and $C^1$ but not $C^2$ on the boundary (its derivative is $\widetilde H'(\vert p \vert)=2 (\pi - 2\mathrm{arccos}(\vert p \vert))$).

\chapter{Floer homology}\label{ch:Floer_homology}

This chapter is intended as a rough overview of the different flavors of Floer homology which arise in the study of the CR3BP. The treatment will err on the side of conciseness, and we will usually treat the simplest cases in order to avoid unnecessary technicalities. In particular, we will not provide proofs of the statements we will make, but refer to the literature where appropriate. 

\section{The CZ-index}

The aim of this section is to give a number of different ways to define the CZ-index, which is part of the index theory of the symplectic group. This will not only be used in the construction of Floer homology, but will be important later when discussing practical aspects. In practical terms, the CZ-index helps understand which families of orbits connect to which (as it stays constant if no bifurcation occurs, and jumps under bifurcation in a prescribed way). We will give two definitions, the second one being more amenable for numerical implementation, and provide useful formulas for low dimensional cases, which are relevant for planar orbits in the CR3BP.

\subsection{The Maslov index} Now we will define the Robbin--Salamon index, as well as the CZ-index index, as a suitable intersection number with the \emph{Maslov cycle}. We will follow the exposition in \cite{FvK}, where we refer the reader to details here omitted.

We consider the \emph{Lagrangian Grassmannian} $\Lambda(n)$, the manifold consisting of all Lagrangian subspaces of $\mathbb C^n$. It comes with an obvious transitive action of $U(n)$ with stabilizer $O(n)$, which gives it the structure of a homogeneous space
$$
\Lambda(n)=U(n)/O(n).
$$
We then consider the map
$$
\rho: U(n)/O(n)\rightarrow S^1,
$$
$$
[A]\mapsto \det A^2.
$$
Then $\rho$ is a fibre bundle with fiber $SU(n)/SO(n)$, which is simply connected. The long exact sequence of homotopy groups then implies that $\rho$ induces an isomorphism of fundamental groups, and therefore
$$
\pi_1(\Lambda(n))\cong \mathbb Z.
$$
If $\lambda:S^1\rightarrow \Lambda(n)$ is a loop of Lagrangian subspaces, then the \emph{Maslov index} of $\lambda$ is
$$
\mu(\lambda)=\deg(\rho \circ \lambda)\in \mathbb Z,
$$
i.e.\ $\mu(\lambda)=[\lambda]\in \pi_1(\Lambda(n))=\mathbb Z$.

We now extend this notion to \emph{paths} of Lagrangians. Given a basepoint $L_0 \in \Lambda(n)$, let
$$
\Lambda^k(n)=\{L\in \Lambda(n): \dim(L\cap L_0)=k\},
$$
for $0\leq k\leq n$. This gives a stratification
$$
\Lambda(n)=\bigcup_{k=0}^n \Lambda^k(n),
$$
and the codimension of $\Lambda^k(n)$ inside $\Lambda(n)$ is 
$$
\text{codim}(\Lambda^k(n),\Lambda(n))=\frac{k(k+1)}{2}.
$$
See Proposition 10.2.1 in \cite{FvK}. In particular, $\Lambda^1(n)$ is a codimension $1$ submanifold of $\Lambda(n)$, with boundary $$\partial \Lambda^1(n)=\bigcup_{k=2}^n\Lambda^k(n)$$ having codimension at least $2$ inside $\Lambda^1(n)$. 

Now, one can identify the tangent space $T_L\Lambda(n)$ with
$$
T_L\Lambda(n)=S^2(L),
$$
the space of symmetric bilinear forms on $L$. Indeed, we can take coordinates $\mathbb R^{2n}=\mathbb R^n\oplus i \mathbb R^n$ such that $L=\mathbb R^n$, and where the symplectic form is $\omega(v,w)=\langle v_1,w_2\rangle-\langle v_2,w_1\rangle$ with respect to this splitting. If $L'$ is another Lagrangian with $L\cap L'=\{0\}$, then $L'$ is the graph of a linear map $S:L\rightarrow \mathbb R^n$, and the Lagrangian condition on $L'$ implies that $S$ is symmetric:
$$
0=\omega(x+iSx,y+iSy)=\langle x, Sy\rangle - \langle y, Sx\rangle.
$$
Then the Lagrangian $L'$, viewed as a nearby Lagrangian to $L$, is identified with the quadratic form $Q_{L'}$ given by $S$. This identification is independent of choices. Moreover, one can show that the tangent space to $\Lambda^k(n)$ is
$$
T_L\Lambda^k(n)=\{L'\in T_L\Lambda(n): Q_{L'}\vert_{L_0\cap L}=0\}.
$$
See Lemma 10.2.4 in \cite{FvK}. This allows to define a co-orientation of $\Lambda^1(n)$ inside $\Lambda(n)$, by defining $L' \in T_L\Lambda(n)/T_L\Lambda^1(n)$ to be positive if $Q_{L'}\vert_{L_0\cap L}$ is positive definite.

Given a smooth loop $\lambda:S^1\rightarrow \Lambda(n)$, by dimensional reasons, we can assume up to perturbation that it intersects $\Lambda^1$ along its interior, and transversely. Then we can define its intersection number with $\Lambda^1$ as
$$
\mu(\lambda)=\sum_{t\in \lambda^{-1}(\Lambda^1)} \nu(t),
$$
where $\nu(t)=\pm 1$ according to whether $\partial_t\lambda$ is positive in $T_{\lambda(t)}\Lambda(n)/T_{\lambda(t)}\Lambda^1(n)$ or not. This is an equivalent definition of the Maslov index, see Theorem 10.2.6 in \cite{FvK}. 

\subsection{The RS-index}

We now define the Maslov index for \emph{paths}, which is called the RS-index, and is due to Robbin--Salamon \cite{RS93}.

Fix a basepoint $L_0\in \Lambda(n)$. Given a smooth path $\lambda:[0,1]\rightarrow \Lambda(n)$ of Lagrangians, the \emph{crossing form} at $t$ is defined as
$$
C(\lambda,L_0,t):=Q_{\dot \lambda(t)}\vert_{\lambda(t)\cap L_0}.
$$
A \emph{crossing} $t\in [0,1]$ is a point in $\lambda^{-1}\left(\bigcup_{k=1}^n\Lambda^k(n)\right)$, and it is \emph{regular} if $C(\lambda,L_0,t)$ is non-degenerate. Up to perturbation, all interior crossings $t\in (0,1)$ can be assumed regular, and we can moreover assume $\lambda(t)\notin \Lambda^k(n)$ for $k\geq 2$ (although the endpoints $\lambda(0),\lambda(1)$ might lie in higher codimension stratum $\Lambda^k(n)$ with $k\geq 2$). Then its \emph{RS-index} is
$$
\mu_{L_0}(\lambda)=\frac{1}{2}\text{sign}\;C(\lambda,L_0,0)+\sum_{t\in (0,1)} \text{sign}\; C(\lambda,L_0,t)+\frac{1}{2}\text{sign}\;C(\lambda,L_0,1)\in \frac{1}{2}\mathbb Z.
$$

The RS-index satisfies the following properties.

\medskip

\begin{itemize}
    \item \textbf{(Invariance)} If two paths are homotopic relative endpoints, then they have the same RS-index.

\medskip
    
    \item \textbf{(Concatenation)} The RS-index of the concatenation of two paths is the sum of their RS-indices.

\medskip
    
\item \textbf{(Loop)} The RS-index of a loop coincides with its Maslov index.
    
\end{itemize}

\subsection{The CZ-index} We now define the CZ-index as a RS-index. Given a symplectic vector space $(V,\omega)$, let $Sp(V)$ denote the space of linear maps $A:V\rightarrow V$ satisfying $A^*\omega=\omega$. Then for $A\in Sp(V)$ the graph
$$
\Gamma_A=\{(x,Ax);x\in V\}\subset V\oplus V
$$
is a Lagrangian subspace of $(V\oplus V,\omega\oplus -\omega)$. In particular, the diagonal
$$
\Delta=\Gamma_{id}=\{(x,x):x\in V\}
$$
is Lagrangian.

Given a smooth path $\Phi:[0,1]\rightarrow Sp(V)$, we can consider the associated path of graphs
$$
\Gamma_\Phi: [0,1]\rightarrow \Lambda(V\oplus V)
$$
$$
\Gamma_\Phi(t)=\Gamma_{\Phi(t)}.
$$
The path $\Phi$ is \emph{non-degenerate} if 
$$
\det(\Phi(1)-id)\neq 0,
$$
which is equivalent to $\Gamma_{\Phi(1)}\in \Lambda_{\Delta}^0$, i.e.\ it intersects the diagonal in a zero dimensional subspace. For a non-degenerate path, we can then define the \emph{Conley--Zehnder} index as 
$$
\mu_{CZ}(\Phi)=\mu_{\Delta}(\Gamma_\Phi)=\frac{1}{2}\text{sign}\; C(\Gamma_\Phi,\Delta,0)+\sum_{t\in (0,1)}\text{sign}\; C(\Gamma_\Phi,\Delta,t) \in \mathbb Z.
$$
This is an integer since sign $C(\Gamma_\Phi,\Delta,1)=0$ by nondegeneracy, and sign $C(\Gamma_\Phi,\Delta,0)\in 2\mathbb Z$ as $C(\Gamma_\Phi,\Delta,1)$ is a non-degenerate quadratic form in the even dimensional space $\Delta$.

\medskip

We also have the following alternative definition of the CZ-index, which is amenable for numerical implementation. We consider the \emph{Maslov cycle}
$$
\Sigma:=\{A\in Sp(2n): \det(A-\mathds 1)=0\},
$$
consisting of all symplectic matrices with $1$ as an eigenvalue. Then the complement of the Maslov cycle inside $Sp(2n)$ consists of two connected components
$$
Sp(2n)^+:=\{A\in Sp(2n): \det(A-\mathds 1)>0 \},
$$
$$
Sp(2n)^-:=\{A\in Sp(2n): \det(A-\mathds 1)<0 \}.
$$
We then fix basepoints in each of them, given by
$$
B_+=-\mathds 1 \in Sp(2n)^+,\; B_-=\text{diag}(2,1/2,-1,\dots,-1) \in Sp(2n)^-. 
$$
From the polar decomposition, every symplectic matrix can be written in the form
$$
A=UP,
$$
where $U$ is unitary, and $P$ is symmetric and positive definite. In terms of $A$, the matrix $U$ can be written as
$$
U=(AA^t)^{-1/2}A.
$$
This defines a retraction
$$
\rho: Sp(2n)\rightarrow U(n)\subset Sp(2n)
$$
$$
\rho(A)=U=(AA^t)^{-1/2}A=\left(\begin{array}{cc}
  X   & -Y \\
  Y   & X
\end{array} \right)=X+iY.
$$
Note that $\rho(B_-)=\mathds 1$, $\rho(B_+)=-\mathds 1$.

Now, given $\Phi:[0,1]\rightarrow Sp(2n)$ with $\Phi(0)=\mathds 1$ and $\Phi(1)$ non-degenerate, we connect $\Phi(1)$ with $B_+$ if $\Phi(1)\in Sp(2n)^+$ via a path inside $B_+$, or with $B_-$ if $\Phi(1)\in Sp(2n)^-$ via a path inside $B_-$ (this uses that $B_\pm$ are path-connected). Denote by
$$
\widetilde \Phi: [0,2]\rightarrow Sp(2n)
$$
the resulting concatenation. We then obtain a map
$$
\gamma_\Phi:[0,2]\rightarrow S^1,
$$
$$
\gamma_\Phi(t)=\text{det}_{\mathbb C}(\rho \circ \widetilde \Phi(t)).
$$
In order to turn $\gamma_\Phi$ into a loop in $S^1$, we square it. The CZ-index is then defined as
$$
\mu_{CZ}(\Phi)=\deg \gamma_\Phi^2 \in \mathbb Z.
$$
This is independent of choices, as $B_\pm$ are in fact simply connected.


\subsection{Useful formulas} 

\begin{figure}
    \centering
    \includegraphics[width=0.8\linewidth]{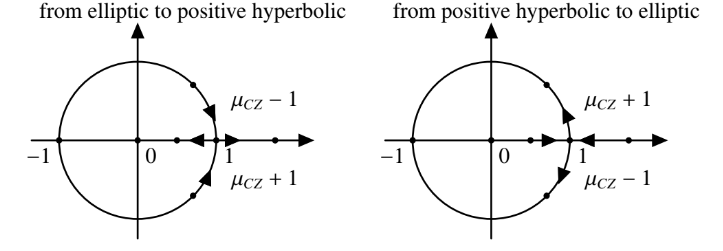}
    \caption{$\mu_{CZ}$ jumps by $\pm 1$ when crossing $1$, according to direction of bifurcation, as shown. If it stays elliptic, the jump is by $\pm 2$. This is determined by the $B$-sign.}
    \label{fig:CZ-jumps}
\end{figure}

We now provide simple formulas to compute the CZ-index in practice.

\medskip

\textbf{Planar case.} The following are formulas for planar orbits, i.e.\ periodic orbits of Hamiltonian systems in two degrees of freedom, e.g.\ the planar CR3BP. Consider $x$ a planar orbit with (reduced) monodromy $M^{red}_x$, and $x^k$ its $k$-fold cover.

\medskip

\begin{itemize}
    \item \textbf{Elliptic case:} $M^{red}_x$ is conjugated to a rotation,
    \begin{equation}\label{reduced_planar}
        M^{red}_x \sim \left(\begin{array}{cc}
       \cos \varphi  & -\sin \varphi \\
    \sin \varphi &  \cos \varphi
    \end{array}\right),
    \end{equation}
    with Floquet multipliers $e^{\pm 2\pi i\varphi}$. Then

    $$
    \mu_{CZ}(x^k)=1+2\cdot \lfloor k\cdot\varphi/2\pi\rfloor
    $$     
    
    In particular, it is odd, and jumps by $\pm$ 2 if the eigenvalue $1$ is crossed in a family. 
    
    For a symmetric periodic orbit, the monodromy matrix at a symmetric point is a Wonenburger matrix $M^{red}_x = \left(\begin{array}{cc}
       a  & b \\
    c &  a
    \end{array}\right)$. In the case (\ref{reduced_planar}), we see that if $b=-\sin \varphi < 0$ then the rotation is determined by $\varphi$ and if $b=-\sin \varphi>0$ then the rotation is determined by $-\varphi$; this determines the CZ-index jump, see Figure \ref{fig:CZ-jumps}.

\medskip
    
    \item \textbf{Hyperbolic case:} $M_x^{red}$ is diagonal up to conjugation, $$
    M^{red}_x \sim \left(\begin{array}{cc}
\lambda  & 0 \\
    0 &  1/\lambda
    \end{array}\right),
    $$
    with Floquet multipliers $\lambda,1/\lambda$. Then

    $$
    \mu_{CZ}(x^k)=k\cdot n,
    $$    
    where the linearized flow $d\phi_t^H$ rotates the eigenspaces by angle $\frac{\pi n t}{T}$, with $n$ even/odd if $x$ positive/negative hyperbolic. 
\end{itemize}    

\medskip
    
\textbf{Spatial case.} We consider now periodic orbits which are planar, but we view them inside a spatial problem. In other words, we assume that the Hamiltonian admits a planar problem in the $xy$-plane (the \emph{ecliptic}), and that the reflection along such $xy$-plane gives rise to a symplectic symmetry of $H$ (e.g.\ the CR3BP). If $x\subset \mathbb R^2$ is a planar orbit, recall that we have a symplectic splitting into planar and spatial blocks
$$
M^{red}_x\sim\left(\begin{array}{cc}
    M_p^{red} &  0\\
    0 & M_s
\end{array}\right)\in Sp(4),\quad M_p^{red}, M_s \in Sp(2).
$$
Then $$\mu_{CZ}(x)=\mu^p_{CZ}(x)+\mu^s_{CZ}(x),$$ where each summand corresponds to $M_p^{red}$ and $M_s$ respectively, and they can be computed by the formulas provided for the planar case. 

The following will be useful when discussing bifurcations in Section \ref{sec:bifurcations}, which correspond to the eigenvalue $1$ being crossed in a family of periodic orbits. The CZ-index jumps whenever this occurs, in a way described in Figure \ref{fig:CZ-jumps}. Whenever the family is planar, it can undergo planar to planar bifurcation, by which new planar families arise; but it may also undergo planar to spatial bifurcations, by which new spatial families arise. We have that:

\medskip
    
\begin{itemize}
    \item Planar to planar bifurcations correspond to jumps in $\mu_{CZ}^p$.

\medskip
    
    \item Planar to spatial bifurcations correspond to jumps of $\mu_{CZ}^s$.
\end{itemize}

\section{Hamiltonian Floer homology} Hamiltonian Floer homology is the version of Floer homology introduced by Floer, for closed symplectic manifolds, in the context of proving the Arnold conjecture. An excellent book with all details is Audin--Damian's \cite{AD10}. In what follows, we will stick to the symplectically aspherical case, in order to avoid having to discuss virtual techniques, needed for the general case, which are well beyond the scope of this book.

Let $(M,\omega)$ be a closed symplectic manifold, which is \emph{symplectically aspherical}, i.e.\ $\omega\vert_{\pi_2(M)}=0$, and which is also \emph{symplectically Calabi-Yau}, i.e.\ $c_1(M,\omega)\vert_{\pi_2(M)}=0$. Consider 
$$
H:S^1\times M\rightarrow \mathbb R, H_t=H(t,\cdot),
$$
a time-periodic $C^2$ Hamiltonian, with Hamiltonian vector field $X_{H_t}$ defined via $i_{X_H}\omega=dH.$ Consider the space $\mathcal{L}(M)$ of $C^1$ contractible loops in $M$, and the action functional from classical mechanics, given by
$$
\mathcal{A}_H:\mathcal{L}(M)\rightarrow \mathbb R,
$$
$$
\mathcal{A}_H(\gamma)=-\int_{\mathbb D^2} \overline{\gamma}^*\omega - \int_{S^1} H_t(\gamma(t))dt,
$$
where $\overline{\gamma}:\mathbb D^2\rightarrow M$ is a $C^1$-map satisfying $\overline{\gamma}\vert_{\partial \mathbb D^2}=\gamma$ (which exists as $\gamma$ is contractible). By symplectic asphericity, the above functional is well-defined.

Choose $J=\{J_t\}_{t\in S^1}$ a family of $\omega$-compatible almost complex structures on $M$, and denote the corresponding Riemannian metric by $g_t=\omega(\cdot,J_t\cdot).$ The associated gradient is $\nabla H_t=-J_tX_{H_t}$. A short computation yields that the differential of $\mathcal{A}_H$ at $\gamma$ is
$$
d_\gamma\mathcal{A}_H:T_\gamma\mathcal{L}(M)=\Gamma(\gamma^*TM)\rightarrow \mathbb R,
$$
$$
d_\gamma\mathcal{A}_H(\eta)=\int_{S^1}\omega(\dot \gamma(t) -X_{H_t}(\gamma(t)),\eta(t))dt,
$$
and therefore we see that the set of critical points of $\mathcal{A}_H$ is
$$
\text{crit}(\mathcal{A}_H)=\mathcal{P}(H),
$$
where we denote by $\mathcal{P}(H)$ the contractible $1$-periodic orbits of $H$, i.e.\ solutions to $\dot \gamma=X_{H_t}(\gamma)$. We assume that $H$ is \emph{generic}, in the sense that all its $1$-periodic orbits are non-degenerate (which is equivalent to nondegeneracy as critical points of the action functional).

We also define an $L^2$-metric on $\mathcal{L}(M)$ given by
$$
g(\eta_1,\eta_2)=\int_{S^1} g_t(\eta_1(t),\eta_2(t))dt.
$$
From the computation of the differential, we see that
$$
d_\gamma\mathcal{A}_H(\eta)=g(\nabla_\gamma \mathcal{A}_H,\eta),
$$
where 
$$
\nabla_\gamma \mathcal{A}_H=J_t(\dot \gamma - X_{H_t}(\gamma)).
$$
The negative gradient equation for a path $s\mapsto u_s\in \mathcal{L}(M)$, i.e.\ a map $u:\mathbb R\times S^1\rightarrow M$, given by
$\partial_s u =-\nabla_{u(s)}\mathcal{A}_H$, is then equivalent to
\begin{equation}\label{eq:Floer_equation}
\partial_s u +J_t\partial_tu=\nabla X_{H_t}.
\end{equation}
This is the \emph{Floer equation}.

Every element $\gamma\in \mathcal{P}(H)$ can be assigned a grading, as follows. Consider a spanning disk $\overline{\gamma}:\mathbb D^2\rightarrow M,$ and take a symplectic trivialization of the bundle $\overline{\gamma}^*TM$ over $\mathbb D^2$. This induces a symplectic trivialization $$\tau:\gamma^*TW\stackrel{\cong}{\longrightarrow} S^1\times \mathbb R^{2n}.$$ If $T$ is the period of $\gamma$, the trivialization induces a path of symplectic matrices $\Phi:[0,T]\rightarrow Sp(2n)$, given by the linearization of the Hamiltonian flow along $\gamma$ viewed in the trivialization, i.e.\
$$
\Phi(t)= \tau \circ d_{\gamma(t)}\phi_H^t \circ \tau^{-1}.
$$
We can then assign the grading
$$
\vert x \vert=-\mu_{CZ}(\Phi),
$$
given by the CZ-index of the path $\Phi$. This is well-defined, as different trivializations $\tau,\tau'$ arising from spanning disks $\overline{\gamma},\overline{\gamma}'$ induce CZ-indices which differ by
$$
\mu_{CZ}(\Phi)-\mu_{CZ}^{'}(\Phi')=2c_1(A)=0,
$$
where $A$ is the sphere obtained by gluing together the two disks (here we use the symplectic Calabi-Yau assumption).

Given critical points $x,y\in \mathcal{P}(H)$, we can then consider the moduli space of Floer solutions which are asymptotic to these, i.e.\
$$
\widetilde{\mathcal{M}}(x,y)=\left\{u:\mathbb R\times S^1\rightarrow M: u\;\text{ satisfies (\ref{eq:Floer_equation}) and } \lim_{s\rightarrow +\infty}u(s)=x,\lim_{s\rightarrow -\infty}u(s)=y \right\}.
$$
This comes with a natural $\mathbb R$-action given by reparametrization $t\cdot u(s,\cdot)=u(s+t,\cdot)$, and we denote the quotient by
$$
\mathcal{M}(x,y)=\widetilde{\mathcal{M}}(x,y)/\mathbb R.
$$
If the path $J_t$ is \emph{generic}, then this quotient is a manifold of dimension
$$
\dim \mathcal{M}(x,y)=\vert y \vert - \vert x \vert -1.
$$
We now consider the \emph{Floer chain complex}
$$
FC_k(H)=\bigoplus_{x \in \mathcal{P}(H), \vert x \vert=k} \mathbb Z \cdot \langle x\rangle
$$
Assigning coherent orientations to the moduli spaces $\mathcal{M}(x,y)$ as in \cite{FH93}, one may count the $0$-dimensional ones with signs, and therefore define a degree $-1$ differential given by
$$
\partial: FC_*(H)\rightarrow FC_{*-1}(H)
$$
$$
\partial \langle x \rangle =\sum_{y: \dim \mathcal{M}(x,y)=0}\#\mathcal{M}(x,y)\langle y \rangle.
$$
In order to show that this is well-defined, one needs a compactness theorem, ensuring that the number of $y$ such that $\mathcal{M}(x,y)$ is non-empty is finite, for which a suitable notion of energy is relevant. This compactness theorem, combined with a gluing theorem, then implies that $\partial^2=0$. One can then define the \emph{Floer homology} as
$$
FH_*(H,J)=H_*(FC_*(H),\partial).
$$
The fundamental property of Floer homology is that it does \emph{not} depend on $(H,J)$, i.e.\ a (regular) homotopy $(H_s,J_s)$ induces, via a \emph{continuation map}, an isomorphism 
$$
FH_*(H_0,J_0)\cong FH_*(H_1,J_1).
$$
Therefore one can assume that $H$ is autonomous, Morse and $C^2$-small. In such case, one can show that the only $1$-periodic orbits are critical points of $H$, and that Floer solutions are $S^1$-invariant and hence Morse flow lines. Moreovoer, the Morse index and the CZ-index of a critical point are related by
$$
\vert x\vert = \dim(M)/2-\text{ind}(x)
$$
The conclusion is then:

\begin{thm}[\textbf{Floer}] 
The Floer homology of a closed is (up to a constant grading shift) isomorphic to Morse homology of the underlying smooth manifold.
\end{thm}

Using that Morse homology is in turn isomorphic to singular homology, one obtains the Arnold conjecture as a corollary:

\begin{thm}[Arnold conjecture, \textbf{Floer}] 
The number of $1$-periodic orbits of a non-degenerate Hamiltonian on a closed symplectic manifold is bounded from below by the sum of the Betti numbers of the manifold.
\end{thm}

\section{Symplectic homology} Symplectic homology is the adaptation of Hamiltonian Floer homology, but for the non-closed case, i.e.\ for Liouville manifolds. Excellent surveys are Oancea's \cite{O04}, Seidel's \cite{Sei06}, Salamon's \cite{Sal97} (see also Wendl's notes \cite{Wen}). We will follow the sign conventions of \cite{BO09}. 

Let $(M,\omega=d\lambda)$ be a Liouville domain, with contact-type boundary $(\partial M,\xi=\ker \alpha)$. We consider its Liouville completion
$$
(\widehat M,\widehat \omega)=(M,\omega)\cup([1,+\infty)\times \partial M,d(r\alpha)).
$$
The construction of symplectic homology for $(M,\omega)$ will follow the same strategy as for the closed case, although it will differ in a few crucial points. 

Given a Hamiltonian $H_t:M\rightarrow \mathbb R$, since $\omega$ is exact, the extension of a periodic orbit to a disk is unnecessary, and one can define the action functional on the (not necessarily contractible part of) the loop space as
$$
\mathcal{A}_H(x)=-\int_{S^1}x^*\lambda-\int_{S^1} H_t(x(t))dt.
$$
The critical points are again the $1$-periodic orbits of $H$, and the Floer equation is derived as before.

However, the main issue in defining a Floer homology is that Floer solutions may a priori escape to infinity along the cylindrical end, which would impose issues on compactness. More generally, the required compactness theorem relies on three ingredients:
\begin{itemize}
    \item Floer cylinders must satisfy a priori $C^0$-bounds.
    \item Floer cylinders must satisfy uniform bounds on the \emph{energy}
    $$
    E(u):=\frac{1}{2}\int_{\mathbb R\times s^1}\left(\vert \partial_su\vert^2+\vert \partial_tu-X_{H_t}(u)\vert^2\right)ds\wedge dt.
    $$
    \item The moduli space of holomorphic spheres that could ``bubble off'' from a sequence of Floer solutions must have codimension at least $2$.
\end{itemize}
The last condition is irrelevant in the case of an exact symplectic form, as there are non-constant holomorphic spheres. For a Floer solution $u:\mathbb R\times S^1\rightarrow M$ asymptotic to $x,y\in \mathcal{P}(H)$ at $\pm \infty$, the energy is indeed bounded, as we have
\begin{equation}
\begin{split}
\mathcal{A}_H(x)-\mathcal{A}_H(y)&=-\int_{-\infty}^{\+\infty}\frac{d}{ds}\mathcal{A}_H(u(s))ds\\
&=-\int_{-\infty}^{\+\infty}g(\nabla \mathcal{A}_H(u(s)),\partial_s u(s))ds\\&=\int_{-\infty}^{\+\infty}\vert \partial_s u\vert^2ds \wedge dt=E(u).
\end{split}
\end{equation}

The crucial condition is then the first condition, which is no longer automatic in the non-closed case. The required bounds are usually obtained via a \emph{maximum principle}. The key observation is that, under suitable assumptions on the Hamiltonian and almost complex structure (yet to be specified), given a Floer solution $u$ which intersects the cylindrical end $[1,+\infty)\times \partial M$, the function
$$
r\circ u: u^{-1}([1,+\infty)\times \partial M)\rightarrow \mathbb R,
$$
where $r$ is the $\mathbb R$-coordinate, is, up to lower order terms, \emph{subharmonic}. That is, it satisfies
$$
\Delta (r\circ u)+l.o.t.\geq 0,
$$
where $\Delta$ is the Laplacian (see e.g.\ Lemma 1.5 in \cite{O04}). This in turn implies that $r\circ u$ cannot have interior maxima. In particular, if the asymptotic orbits of a Floer solution lie away from the cylindrical end, the whole cylinder also does, and therefore it cannot escape to infinity.

The conditions we need to impose on the Hamiltonians and almost complex structures, in order to have a maximum principle, are as follows. First, if $\mathcal P(\alpha)$ denotes the set of all periodic Reeb orbits of $\alpha$ at the boundary, then the \emph{spectrum} of $\alpha$ 
$$
\mbox{spec}(\alpha)
=\{ a \in \R ~|~\text{there exists }\gamma \in \mathcal P(\alpha) \text{ such that } a= \mathcal A(\gamma) \},
$$
where the \emph{action} is defined as $\mathcal{A}(\gamma)=\int_\gamma \alpha$. A Hamiltonian $H:\widehat W \rightarrow \mathbb R$ is \emph{non-degenerate} if all its $1$-periodic orbits are non-degenerate, and it is \emph{linear at infinity} if at the cylindrical end $ H$ has the form $ H(r,b,t) =cr+d$ for some constants $c>0$ and $d$, where its \emph{slope} is slope$( H):=c$. A Hamiltonian $ H$ that is non-degenerate and linear at infinity with slope$(H) \notin \mbox{spec}(\alpha)$ will be called \emph{admissible}. We call an $S^1$-family of almost complex structures $J=J_t$ on $\widehat W$ \emph{admissible} if it is $\omega$-compatible, and on the cylindrical end it is translation invariant, $J \xi =\xi$, and $J\partial_r =R_\alpha$. 

If $H$ is admissible, then its Hamiltonian vector field on the cylindrical end is $X_H=-cR_\alpha,$ i.e.\ it is collinear with the Reeb field. By choice of $c$, it follows that $H$ has no $1$-periodic orbits along the cylindrical end. Moreover, if $(H,J)$ is an admissible pair, then the maximum principle holds, and therefore compactness also holds. If its moreover generic, then we have a well-defined Floer homology $FH_*(H,J)$.

The main difference with the closed case is that this Floer homology is \emph{not} independent of $H$, and it only detects Reeb orbits at the boundary whose action is at most $c$. In order to address this, we take a limit with respect to an increasing sequence of slopes. For this, we have to give more details on how continuation maps work. For a $1$-parameter family of admissible Hamiltonians $H^s$, $s\in [0,1]$, with slopes $c_s\rightarrow +\infty$ monotonically, together with accompanying admissible $J^s$, we consider the \emph{parametric} Floer equation
$$
\partial_su+J_t^s\partial_tu=\nabla X_{H_t}^s.
$$
The monotonic behavior of the slopes implies that the maximum principle holds for the parametric equation as well. Counting the solutions to this equation in the $0$-dimensional moduli spaces $\mathcal{M}(x,y,H^s,J^s)$ (with fixed asymptotics $x,y$) yields the continuation map
$$
\Phi_{H^s,J^s}: FC_*(H^0)\rightarrow FC_*(H^1)
$$
$$
\Phi_{H^s,J^s}(\langle x \rangle)= \sum_{y: \dim \mathcal{M}(x,y,H^s,J^s)=0} \# \mathcal{M}(x,y,H^s,J^s) \langle y \rangle.
$$
Suitable compactness and gluing theorems imply that this is a chain map, and so induces a map in homology
$$
\Phi_{H^s,J^s}: FC_*(H^0,J^0)\rightarrow FC_*(H^1,J^1).
$$
Moreover, these maps are independent on the homotopy (among homotopies with monotonically increasing slope), and concatenation of homotopies yield compositions of the corresponding maps. Therefore we can define \emph{symplectic homology} as the direct limit
$$
SH_*(M,\lambda)= \lim_{\stackrel{\longrightarrow}{c}}HF_*(H_c,J_c),
$$
where the slope of $H_c$ is $c\rightarrow +\infty$ monotonically, and the direct limit is defined with respect to the above continuation maps. This definition is in fact independent of the contact form and the almost complex structure, as can be seen via continuation maps for homotopies of $J$, see \cite{Sei06}.

\medskip

One fundamental computation of symplectic homology is for the case of cotangent bundles, for which we have the following result.

\begin{thm}[\cite{V18,AS06,SW06}] The symplectic homology of a cotangent bundle $(T^*Q,\omega_{std})$ is isomorphic to the singular homology of the free loop space $\mathcal{L}Q$ of $Q$, i.e.\ 
$$
SH_*(T^*Q,\lambda_{std})=H_*(\mathcal{L}Q).
$$
    
\end{thm}

\section{Lagrangian Floer homology} Lagrangian Floer homology is a version of Floer homology for Lagrangians in a closed symplectic manifold (introduced by Floer), whereas its wrapped version is its further adaptation to the case of Liouville manifolds. A very nice introduction to the subject, in the context of Fukaya categories, is Auroux's \cite{Au14}, from where we borrow some parts of the exposition.

The first observation to make is that Hamiltonian Floer homology can be recast as a theory whose generators correspond to intersection points of two suitable Lagrangians. Indeed, given a closed symplectic manifold $(M,\omega)$ and a Hamiltonian $H_t:(M,\omega)\rightarrow \mathbb R$, consider the symplectic manifold 
$$
(M\times M,\omega\oplus - \omega).
$$
Since the time-$1$ map $\phi$ of the Hamiltonian flow is a symplectomorphism, the graph 
$$
\Gamma_\phi=\{(x,\phi(x))\in M\times M\}
$$
is a Lagrangian submanifold. In particular, the diagonal
$$
\Delta=\{(x,x)\in M\times M\}
$$
is also. Then by definition, $1$-periodic orbits of $H_t$ correspond to intersection points in $\Gamma_\phi\cap \Delta,$ and this intersection is transverse if and only if $H_t$ is non-degenerate. Moreover, note that $\Gamma_\phi=(id \times \phi)(\Delta)$ is the image of the diagonal under the Hamiltonian diffeomorphism $id\times \phi$. A bit more work shows that for suitable almost complex structure, Floer cylinders in $M$ correspond to Floer \emph{strips} in $M\times M$ with boundary in $\Gamma_\phi$ and $\Delta$, i.e.\ maps
$$
u:\mathbb R\times [0,1]\rightarrow M\times M,
$$
satisfying 
$$
\lim_{s\rightarrow \pm \infty} u(s,\cdot)=x_\pm \in \Gamma_\phi\cap \Delta,
$$
$$
u(\mathbb R\times \{0\})\subset \Delta, u(\mathbb R\times \{1\})\subset \Gamma_\phi,
$$
and the Cauchy--Riemann equation
$$
\partial_s u + J_t\partial_tu=0,
$$
for suitable $J_t$. 

In general, consider a closed Lagrangian submanifold $L\subset (M,\omega)$ of a closed symplectic manifold, and a Hamiltonian diffeomorphism $\phi:M\rightarrow M$. The following theorem for the case of Lagrangians then recovers the Arnold conjecture.

\begin{thm}[\textbf{Floer} \cite{F88a}] Assume that the symplectic area of any disc in $M$ with boundary in $L$ vanishes (i.e.\ $\omega\vert_{\pi_2(M,L)}=0$), and that $\phi(L)$ and $L$ intersect transversely. Then the number of intersection points in $\phi(L)\cap L$ is bounded from below by the sum of the Betti numbers of $L$.
\end{thm}

In order to prove the above theorem, Floer introduced a version of Floer homology for Lagrangian intersections, with the following formal properties. For a pair of Lagrangians $(L_0,L_1)$ intersecting transversely, we have an associated chain complex 
$$
FC_*(L_0,L_1)=\bigoplus_{p \in L_0\cap L_1} \Lambda\cdot \langle p \rangle,
$$
generated by intersection points over the \emph{Novikov field} $$\Lambda=\left\{ \sum_{i=0}^\infty a_iT^{\lambda_i}: a_i \in \mathbb K,\lambda_i\in \mathbb R^+,\lim_{i\rightarrow +\infty}\lambda_i=+\infty\right\},$$ where $\mathbb K$ is some base field, equipped with a differential $\partial$ which satisfies:
\begin{itemize}
    \item $\partial^2=0$, so that there is a well-defined \emph{Lagrangian} Floer homology $FH_*(L_0,L_1)$;
    \item If $L_1,L_2$ are Hamiltonian isotopic, then $FH_*(L_0,L_1)\cong FH_*(L_0,L_2)$;
    \item If $L_1=\phi(L_0)$ for a Hamiltonian diffeomorphism $\phi,$ then $FH_*(L_0,\phi(L_0))\cong H_*(L_0)$.
\end{itemize}
The above theorem follows directly from these formal properties of Lagrangian Floer homology, as the rank of $H_*(L,\phi(L))=H_*(L)$ is then bounded from above by $\#(L\cap \phi(L))$.

Similarly as Hamiltonian Floer homology, which can be viewed as an analogue of Morse homology for the action functional from mechanics, Lagrangian Floer homology can be thought of as an infinite-dimensional analogue of Morse homology for another action functional. This is defined on the universal cover of the path space 
$$
\mathcal{P}(L_0,L_1) = \{\gamma: [0,1] \rightarrow M :\gamma(0) \in L_0, \gamma(1) \in L_1\},
$$
and is given by
$$
\mathcal{A}(\gamma,[\Gamma]) = -\int_{\Gamma}\omega,
$$
where $(\gamma, [\Gamma]) \in \widetilde{\mathcal{P}}(L_0, L_1)$ consists of a path $\gamma \in \mathcal{P}(L_0, L_1)$ and an equivalence class $[\Gamma ]$ of a homotopy $\Gamma : [0,1] \times [0,1] \rightarrow M$ between $\gamma$ and a fixed base point in the connected component of $\mathcal{P}(L_0,L_1)$ containing $\gamma$. The critical points of $\mathcal{A}$ correspond to constant paths at intersection points, and its gradient flow lines are pseudo-holomorphic strips bounded by $L_0$ and $L_1$.

In fact, we will actually have to perturb the setup with a time-dependent Hamiltonian $H_t$, and rather than consider intersection points $p,q  \in L_0\cap L_1$, we will need to consider intersection points $p,q\in L_0\cap (\phi_1^H)^{-1}(L_1)$, where $\phi_1^H$ is the time-$1$ map of the Hamiltonian flow. In other words, generators will correspond to Hamiltonian \emph{chords} from $L_0$ to $L_1$, i.e.\ paths $\gamma:[0,1]\rightarrow M$ satisfying $\dot \gamma(t)=X_{H_t}(\gamma(t))$, with $\gamma(0)\in L_0$ and $\gamma(1)\in L_1$. Note that this also makes sense in the case $L_0=L_1$. Then a \emph{Floer strip} from $q$ to $p$ is a map $u:\mathbb R\times [0,1]\rightarrow M$ solving the associated Floer equation (for some choice of time-dependent and compatible $J_t$), satisfying the boundary conditions
$$
u(s,0)\in L_0,\;u(s,1)\in L_1, \text{ for all } s,\;\lim_{s\rightarrow +\infty}u(s,t)=p, \lim_{s\rightarrow -\infty}u(s,t)=q,
$$
and the finite energy condition
$$
E(u)=\int u^*\omega=\int_{\mathbb R\times S^1}\vert \partial_su \vert^2 ds\wedge dt <\infty.
$$
For a relative homotopy class $A\in \pi_2(M,L_0\cup L_1)$, denote by $\widetilde{\mathcal{M}}(p,q;A,J)$ the space of Floer strips in class $A$, and by $\mathcal{M}(p,q;A,J)$ its quotient by the natural $\mathbb R$-action $(a.\cdot u)(s,t)=(s+a,t)$. In case where \emph{regularity/transversality} holds (i.e.\ the linearization of the Floer equation is a surjective operator), then $\mathcal{M}(p,q;A,J)$ is a smooth manifold; this is the main technical point of the construction. Two other crucial ingredients are compactness of these moduli spaces, for which a suitable compactness theorem is required,  and orientability, by which one should be able to orient these moduli spaces. For the latter condition one usually assumes that $L_0,L_1$ are oriented and spin. The choice of a spin structure then canonically determines an orientation of the moduli spaces, see \cite{FOOO,Sei08}. There is also the issue of gradings, by which one should be able to assign gradings to the generators of the theory.

Assuming all these technical issues are addressed, the \emph{Floer differential} is then defined as
$$
\partial:CF_*(L_0,L_1)\rightarrow CF_{*-1}(L_0,L_1)
$$
$$
\partial\langle q \rangle=\sum_{\substack{p\in L_0\cap L_1\\ \dim \mathcal{M}(p,q;A,J)=0}} \#\mathcal{M}(p,q;A,J) T^{\omega(A)}p,
$$
where $\omega(A)=\langle \omega,A\rangle\in \mathbb R$ is the natural evaluation of the symplectic form over the homotopy class $A$. 

In our setup, the Hamiltonian perturbation to the CR-equation we have considered will be sufficient to achieve transversality, so that the moduli spaces are manifolds with the expected dimension. In order to make sure that the above sum is well-defined over the Novikov field, one needs to appeal to a version of Gromov compactness, which ensures that for any given energy bound $E$, there are only finitely many homotopy classes $A$ for which $\mathcal{M}(p,q;A,J)$ is non-empty. 

Moreover, in the case where the Lagrangians are \emph{exact}, i.e.\ whenever $\omega=d\lambda$ is itself exact, and $\lambda\vert_{L_i}=df_i$ for some smooth functions $f_i:L_i\rightarrow \mathbb R$, we can get rid of Novikov coefficients. Indeed, noting that $\omega(A)=f_1(q)-f_0(q)-(f_1(p)-f_0(p))$, we can rescale each generator by $p\mapsto T^{f_1(p)-f_0(p)}p$ and thus eliminate the weights from the differential, and work directly over $\mathbb K$. 

We now say a few words about how the grading and indices of Floer strips are defined. Consider a Floer strip $u: \mathbb R\times [0,1]\rightarrow M$. Since its domain is contractible, the symplectic bundle $u^*TM$ is trivial. Choosing a symplectic trivialization, then the paths
$$
l_0=u^*\vert_{\mathbb R\times \{0\}}TL_0,\;l_1=u^*\vert_{\mathbb R\times \{1\}}TL_1
$$
can be viewed as paths in the Lagrangian Grassmannian $\Lambda(n)$, one joining $T_pL_0$ to $T_qL_0$, and the other joining $T_pL_1$ to $T_pL_0$. In order to turn these paths in to a loop, one can concatenate with the \emph{canonical short paths} $\lambda_p,\lambda_q$ respectively joining $T_pL_0$ to $T_pL_1$, and $T_qL_0$ to $T_qL_1$. Here, the canonical short path between two transverse Lagrangian subspaces $\lambda_0,\lambda_1\in \Lambda(n)$ is defined by first choosing coordinates such that $\lambda_0=\mathbb R^n$ and $\lambda_1=i\mathbb R^n$, and then considering the path $\lambda_t=e^{-i\pi t/2}\mathbb R^n$, $t\in [0,1]$. Then, the index of the Floer strip $u$ is defined as the Maslov index of the resulting loop in $\Lambda(n)$. This is then the expected dimension of the moduli space $\widetilde{\mathcal{M}}(p,q;A,J)$ in which $u$ lies.

In order to define a $\mathbb Z$-grading for the Floer complex, we make the following assumptions:
\begin{itemize}
    \item $2c_1(M,\omega)=0$; and
    \item The \emph{Maslov class} $\mu_L\in H^1(L;\mathbb Z)$ vanishes.
\end{itemize}

The first condition allows to find a nowhere vanishing section $\Theta$ of $(\Lambda_\mathbb C^nT^*M)^{\otimes 2}$, which can be used to associate to a Lagrangian $L$ its \emph{phase function} $\varphi_L:L\rightarrow S^1$ via $\varphi_L(p)=\text{arg}(\Theta\vert_{T_pL})\in S^1$. The second condition then allows to coherently lifts this phase function to $\mathbb R$, as the Maslov class is by definition the homotopy class $\mu_L=[\varphi_L]\in [L,S^1]=H^1(L;\mathbb Z)$. If this vanishes, then one can lift this map to $\widetilde \varphi_L:L\rightarrow \mathbb R$. A \emph{graded} Lagrangian $\widetilde L$ is then a Lagrangian $L$ endowed with a choice of such lift $\widetilde \varphi_L$. Given two graded Lagrangians $\widetilde L_0,\widetilde L_1$ and $p \in L_0\cap L_1$, one obtains a preferred homotopy class of paths connecting $T_pL_0$ to $T_pL_1$ by connecting their lifts by a path. Combining this with the opposite of the canonical short path joining them we obtain a loop. The \emph{degree} of $p$ is then the (opposite of the) Maslov index of this loop, denoted $\vert p \vert$. Then we have the formula
$$
\mbox{ind}(u)=\vert q\vert - \vert p \vert, 
$$
and so the Floer differential indeed has degree $-1$. 

In general, if we do not impose the above assumptions, but only impose that the Lagrangians be oriented, we only get a $\mathbb Z_2$-grading, where an intersection point is graded by the sign of the intersection, according to the orientations.

We now briefly discuss compactness. Gromov compactness says that a sequence $u_n$ of Floer strips with bounded energy can converge, up to reparametrization, to a \emph{nodal tree} of Floer strips. This means that the following types of degenerations can occur:
\begin{itemize}
    \item (Strip breaking) energy concentrates at $s \rightarrow \pm \infty$, and suitable translations of $u_n$ converge to a non-constant strip; 
    \item (disk bubbling) energy concentrates at a point of the boundary of the strip, and suitable rescalings of $u_n$ converge to a holomorphic disk with boundary in $L_0$ or $L_1$;
    \item (sphere bubbling) energy concentrates at an interior point of the strip, and suitable rescalings of $u_n$ converge to a holomorphic sphere. 
\end{itemize}

Strip breaking (and a gluing theorem) is what implies $\partial^2=0$, provided one can exclude disk or sphere bubbling. The standard way to do so is to assume that $[\omega]\cdot \pi_2(M,L_i)=0$ for $i=0,1$. Under this assumption, the following was obtained by Floer, putting together what we have discussed.

\begin{thm}[\textbf{Floer}] Assume that $[\omega]\cdot \pi_2(M,L_i)=0$ for $i=0,1$, and when $char(\mathbf K) = 2$ assume that $L_0, L_1$ are oriented and spin.
Then the Floer differential $\partial$ is well-defined, it squares to zero, and the associated Floer homology is, up to isomorphism, independent of the chosen
almost-complex structure and invariant under Hamiltonian isotopies of each Lagrangian.
\end{thm}

The independence of $J$ is again obtained via suitable continuation maps associated to homotopies of $J$. We refer e.g.\ to \cite{Au14} for more details on this result.

An important case is when $L_0=L_1=L$ (oriented and spin). By Weinstein's neighbourhood theorem, a neighbourhood of $L$ in $M$ is symplectomorphic to $T^*L$, so we can restrict to the latter. If we consider a Morse function $f$ on $L$, and we denote by $L_1$ the graph of $\epsilon df$ for some small $\epsilon >0$, then $L_1$ is an exact Lagrangian Hamiltonian isotopic to the zero section $L$ (the isotopy being generated by $\epsilon \pi^*f$ where $\pi:T^*L\rightarrow L$ is the natural projection). The intersections of $L$ with $L_1$ correspond to critical points of $f$, which are transverse as these are non-degenerate. A choice of grading on $L$ then induces a grading on $L_1$ via the isotopy, so that the Floer complex is then graded. Moreover, the Morse index of a critical point $p$ satisfies
$$
\text{ind}(p)=\dim(L)-\vert p\vert,
$$
so that up to a shift the grading in the Floer complex coincides with the Morse grading. A result of Floer then implies that Floer solutions $u$ (for a suitable time-dependent $J$) are in correspondence with Morse flow lines $\gamma$ of $f$, via the correspondence $\gamma(s)=u(s,0)$. It follows that the Floer complex is isomorphic to the Morse complex of $f$, after rescaling the generators via $p\mapsto T^{\epsilon f(p)}p$. Therefore we obtain 
$$
HF_*(L,L)\cong HF_*(L,L_1)\cong H_{\dim(L)-*}(L;\Lambda).
$$
In general, using energy bounds, one can moreover show that, under the assumption that $L$ does not bound non-constant symplectic disks, all Floer solutions between $L$ and a small Hamiltonian perturbation of $L$ will lie in a collar neighbourhood of $L$. Therefore the calculation of $HF_*(L,L)$ reduces to the above calculation, and we obtain the following.

\begin{thm}[\textbf{Floer}] If $[\omega]\cdot \pi_2(M,L)=0$, then $HF_*(L,L)\cong H_{\dim(L)-*}(L;\Lambda)$.

\end{thm}

As explained above, the Arnold conjecture follows from this computation. 

\section{Wrapped Floer homology} We now briefly discuss the adaptation of Lagrangian Floer homology to the case of Liouville manifolds, introduced by Abouzaid--Seidel \cite{AS10} (see also Abbondandolo--Schwarz \cite{AS10}, and \cite{Au14,KKK22}).

Given a Liouville manifold $(M,\omega=d\lambda)$, with contact-type boundary $(\partial M,\xi=\ker \alpha)$, we consider exact Lagrangians which are \emph{cylindrical} at infinity, i.e.\ they are non-compact and coincide with the cone $(1,+\infty)\times\partial L$ along the cylindrical end, where $\partial L$ is Legendrian in $\partial M$. We also call such Lagrangians \emph{admissible}.

We again consider Hamiltonians $H:M\rightarrow \mathbb R$ which are admissible in the sense of symplectic homology, i.e.\ non-degenerate and linear at infinity. This implies that the Hamiltonian vector field is collinear with the Reeb field along the cylindrical end. Given two exact and cylindrical Lagrangians $L_0,L_1$, the wrapped Floer chain complex $CW_*(L_0,L_1)$ is then generated over the base field $\mathbb K$ by elements in $\phi^H_1(L_0)\cap L_1,$ (which we assume transverse) i.e.\ Hamiltonian chords from $L_0$ to $L_1$. By choice of the Hamiltonian, the chords come in two types: those in the interior, and those which lie in the cylindrical end. The latter are Reeb chords from the Legendrians $\partial L_0$ to $\partial L_1$. Transversaility of the intersection is equivalent to \emph{nondegeneracy} of the chords, i.e.\ $d\phi_H^1(T_{c(0)}L_0)\cap T_{c(1)}L_1=\{0\}$ for all chords $c$. We will further impose in the definition of admissibility that the slope of the Hamiltonian at the cylindrical end does lie in the \emph{spectrum}
$$
\text{spec}(\alpha,\partial L_0,\partial L_1)=\left\{T=\int_c\alpha\;\Big\vert\; c:[0,T]\rightarrow M \text{ chord from } \partial L_0 \text{ to } \partial L_1\right\}.
$$

The differential again counts Floer strips with boundary in $L_0$ and $L_1$, with respect to an auxiliary cylindrical $J$. Because of the assumptions on the Lagrangians and Hamiltonian, one can show via a priori energy estimates that all Floer strips which converge to a given generator remain within a bounded subset of $M$, and therefore the differential is well-defined (a finite linear combination). The resulting Floer homology is independent of $J$ but depends on the Hamiltonian. Taking a direct limit over such Hamiltonians with monotonically increasing slopes as in the definition of symplectic homology, the resulting \emph{wrapped} Floer homology is denoted $HW_*(L_0,L_1)$. 

We refer to e.g.\ \cite{Au14} for more details (including the product structure, and more generally the $A_\infty$ operations on wrapped Floer homology).

\subsection{Cotangent bundles} A fundamental computation is the case of cotangent bundles. Consider a compact spin Riemannian manifold $Q$, and consider $M=T^*Q$ with its standard Liouville structure, and with the Hamiltonian $H=\Vert p \Vert^2$.

\begin{thm}[Abbondandolo--Schwarz \cite{AS06}]\label{thm:AS} Let $L=T_q^*Q$, the cotangent fiber over $q\in Q$. Then there is an isomorphism
$$
HW_*(L,L)\cong H_{*}(\Omega_qQ),
$$
between the wrapped Floer homology of $L$ and the homology of the loop space $\Omega_qQ$ of $Q$ based at $q$.    
\end{thm}

The above isomorphism also preserves suitable products \cite{AS10}, and has also been refined in \cite{Ab12}, by upgrading the statement to the chain level, including the $A_\infty$ structure.

\section{Local Floer homology} While non-degeneracy of the Hamiltonians has been so far necessary in order to introduce Floer homology, the \emph{local} versions of all Floer homology theories are meant to address possible degeneracies of the generators (which is certainly to be expected for Hamiltonians arising in celestial mechanics). These ideas go back to Floer, see e.g.\ \cite{F89b,F89c}. More recently, the local theory has been instrumental in the proof of the Conley conjecture, see \cite{Gi10, GG10}.

The general scheme for defining such local versions is the following. We consider first the case of Hamiltonian Floer homology. Given a possibly degenerate Hamiltonian $H$ on a closed symplectic manifold $(M,\omega)$, we consider a possibly degenerate $1$-periodic orbit $\gamma$, which we assume to be isolated. We then fix a small collar neighbourhood $U$ of $\gamma$ which intersects no other $1$-periodic orbit, and perturb the Hamiltonian, via a perturbation supported in $U$, so that now $\gamma$ bifurcates into a finite collection of nearby non-degenerate $1$-periodic orbits, all of them lying inside the a priori fixed isolating neighbourhood $U$. One then proves a localizing theorem, which ensures that all Floer trajectories between the locally arising orbits stay within the localizing neighbourhood (this follows from the analysis in e.g.\ \cite{FHS}). By compactness and gluing, every broken Floer trajectory between such orbits will also lie in $U$.

Then one defines the \emph{local} Floer homology $HF^{loc}_*(H,\gamma)$ of $H$ at $\gamma$ as the homology of the complex generated by the local non-degenerate generators near $\gamma$, with differential counting Floer solutions between them. This is actually independent on the perturbation and the almost complex structure, which can be shown via the standard continuation argument. 

\medskip

In the Liouville case, one can still carry out the above definition for an admissible Hamiltonian $H$, as it is completely localized around the given periodic orbit. Taking a limit to get rid of the dependence on $H$, we obtain \emph{local} symplectic homology $SH_*^{loc}(\gamma,M,\lambda)$. Moreover, we also gave a Lagrangian version: given two cylindrical Lagrangians $L_0,L_1$, one can define \emph{local wrapped} Floer homology at a possibly degenerate \emph{chord} $\gamma$, denoted $HW_*^{loc}(\gamma,L_0,L_1)$, following the same scheme as outlined above. This was done in detail in Limoge's PhD thesis \cite{L25}, and partially in \cite{ML24}.

\medskip

As an example, if $\gamma$ is actually non-degenerate and with CZ-index $k$, then $HF_*^{loc}(H,\gamma)$ coincides with the base field, supported in degree $k$. As another interesting example, if a degenerate orbit $z$ undergoes period-doubling bifurcation (see Section \ref{sec:bifurcations} below), then before bifurcation there are no orbits, and after bifurcation there are two orbits $x,y$ whose indices differ by $1$ and such that $\partial x = y$. We then see that this orbit $z$ is \emph{homologically invisible}, i.e.\ its local Floer homology vanishes, and so one can think of period-doubling bifurcations as homologically irrelevant. Moreover, an important property is that local Floer homology is invariant under homotopies of the Hamiltonian such that $\gamma$ is uniformly isolated for the homotopy.

The \emph{support} of $HF_*^{loc}(H, \gamma)$, denoted supp $HF_*^{loc}(H, \gamma)$, is the collection of integers $k$ such that $HF^{loc}_k(H, \gamma) \neq 0$. Since $HF_*^{loc}(H, \gamma)$ is finitely generated, the support is finite. An important property is that the support satisfies
$$
supp \;HF_*^{loc}(H, \gamma)\subset [\Delta(\gamma) - n, \Delta(\gamma) + n],
$$
i.e.\ it lies in a window of length at most $2n$, where $\Delta(\gamma)$ is the \emph{mean index} of $\gamma$ (see e.g.\ section~ 3.1.1 of  \cite{GG15}), and $n=\dim(M)/2$.

\chapter{Fixed point theory of Hamiltonian twist maps}\label{ch:Ham_twist_maps}

Recall from our discussion on Poincar\'e's work that the second step in order to prove the existence of periodic orbits in the CR3BP consisted in proving a fixed-point theorem, i.e.\ the Poincar\'e--Birkhoff theorem. This chapter is devoted to introducing two versions of a higher-dimensional generalization of this result, obtained first by the author and Otto van Koert in \cite{MvK20b}, together with an improved version, obtained by the author in collaboration with Arthur Limoge in \cite{ML25}. Both versions requires a generalization of the classical twist condition to higher-dimensions. The first result imposes a very rigid version of the twist condition, which we will therefore call the \emph{strong} twist condition. For the second result, we will impose a \emph{weak} twist condition, which is a $C^1$-open condition, as opposed to the strong version. We will present a sketch of the main proof, following \cite{MvK20b}, which is the same for both results, i.e.\ the difference lies in the setups. We will also discuss a relative version for finding chords on Lagrangians. 

We should emphasize that, while the notion we introduced is rather natural (especially from the perspective of Floer theory as discussed in Chapter \ref{ch:Floer_homology}), and moreover we have general abstract theorems, the twist condition in either of its forms is at this stage unsatisfactory, as it is unknown as to whether they apply to the CR3BP, for reasons discussed below. While the full generalization of the first step is provided by Theorem \ref{thm:openbooks}, the second step, i.e.\ an application of an abstract fixed point theorem to the CR3BP, will require further work. We will discuss this in more detail in what follows.

\section{Hamiltonian twist maps}

The periodic points of the return map for the CR3BP are either boundary periodic points, which give planar orbits, or interior periodic points which are in 1:1 correspondence with spatial orbits. We are interested in finding \emph{interior} periodic points, and we follow Poincaré's philosophy to try to find them. First we need to discuss a subtle but important point.

\subsection{A trade-off}\label{sec:trade-off} In the setup which arises in the CR3BP, when we consider the global hypersurfaces of section provided by the open book as in Theorem \ref{thm:openbooks}, we encounter the following heuristic trade-off, which is a feature of the setup which may not be avoided. We have one of the following two situations:

\begin{itemize}
    \item[\textbf{(A)}] The return map extends smoothly to the boundary, while the symplectic form degenerates at the boundary; or
    \item[\textbf{(B)}] The return map extends continuously to the boundary (but admits no $C^1$ extension), while the symplectic form extends to the boundary also as a symplectic form.
\end{itemize}

Here, recall that the relevant symplectic form is $\omega=d\alpha\vert_P$, where $\alpha$ is the ambient contact form for the spatial CR3BP (recall Theorem \ref{thm:contact_type}), and $P$ is a page of the open book provided by Theorem \ref{thm:openbooks}. The boundary degeneracy arises because the Reeb vector field is tangent to the boundary, and by definition, its spans the kernel of $d\alpha$. Both setups are equivalent to each other (in the case of the CR3BP), and are related by a change of coordinates which is smooth in the interior but only continuous at the boundary. In general, one needs to prove the continuous extension of the map to the boundary, which is not at all an obvious fact. For the CR3BP, Otto van Koert, together with the author, proved this in \cite{MvK20a}. Of course, for the 2-body problem, where the return map is the identity, the boundary extension is smooth. We will explain this in detail below. 

First, we will assume the ideal model where both the return map admits a $C^1$-extension to the boundary, and symplectic form does not degenerate, and we will later explain how to adapt the setup.

\subsection{The strong Hamiltonian twist condition.} 
We first propose a generalization of the twist condition introduced by Poincar\'e, for the Hamiltonian case and for arbitrary Liouville domains, which we will call the \emph{strong} twist condition. This is because we will also introduce a weaker notion, which is an open condition, whereas the strong version is not. 

Let $(W,\omega=d\lambda)$ be a $2n$-dimensional Liouville domain, and consider a Hamiltonian symplectomorphism $\tau$. Let $(B,\xi)=(\partial W,\ker \alpha)$ be the contact manifold at the boundary where $\alpha=\lambda\vert_B$, and $R_\alpha$ the Reeb vector field of $\alpha$. The Liouville vector field $V_\lambda$ is defined via $i_{V_\lambda}\omega=\lambda$. 


\begin{definition}(\textbf{Hamiltonian twist map})\label{def:twistmap}
We say that $\tau$ is a \emph{Hamiltonian twist map} (with respect to $\alpha$), if $\tau$ is generated by a $C^2$ Hamiltonian $H:\mathbb{R}\times W \rightarrow \mathbb{R}$ which satisfies $X_{H_t}\vert_B=h_tR_\alpha$ for some \emph{positive} and smooth function $h:\mathbb{R}\times B \rightarrow \mathbb{R}^+$.
\end{definition}

We say that the Hamiltonian $H_t$ above is an \emph{positively wrapping} generating Hamiltonian for $\tau$. In particular, $H_t\vert_B \equiv const$ on $B$, and $\tau(B)\subset B$, so that $X_{H_t}\vert_{B}\in TB$. We have $h_t=dH_t(V_\lambda)\vert_B$ is the derivative of $H_t$ in the Liouville direction $V_\lambda$ along $B$, which we assume strictly positive. Also, $\tau\vert_B$ is the time-$1$ map of a positive reparametrization of the Reeb flow on $B$. But note that, while the latter condition is only localized at $B$, the twist condition is of a \emph{global} nature, as it requires global smoothness of the generating Hamiltonian.

\begin{remark}[\textbf{Shortcomings of the twist condition}]\label{rk:twist_condition}
    As we have discussed in Sections \ref{sec:return_map} and \ref{sec:degenerate}, the global hypersurface of section in the CR3BP is not directly a Liouville domain, as it is a \emph{degenerate} one. This can be fixed, but paying the price that the return map becomes only continuous at the boundary, as explained in Section \ref{sec:trade-off}. In Definition \ref{def:twistmap}, we have consciously ignored this, as the standard setup of Floer theory is not immediately adapted to a degenerate situation. This is, in part, what makes the above definition rather unsatisfactory. In Section \ref{sec:twist_degenerate}, we will provide alternative definitions to address this issue, and adapt the theorem; see Remark \ref{rk:adapting_thm} below. Another shortcoming of the Definition \ref{def:twistmap} is that it is \emph{not} open as soon as the Liouville domain has dimension at least $4$. Indeed, as opposed to the $2$-dimensional case where there is only one positive direction at the boundary (given by the orientation), there are a priori many directions which are ``positive'' at the boundary, i.e.\ the natural notion is being positively transverse to the contact structure. This \emph{weak} version of the twist condition is the best adapted to a fixed-point theorem (a version for degenerate Liouville domains, see Section \ref{sec:twist_degenerate}), since the methods of our proof adapt, after smoothing the map near the boundary. 
    
    Moreover, a further shortcoming of the proof under the strong case, relies on controlling index growth of periodic orbits, and therefore requires the Hamiltonian to be at least $C^2$ in order to consider the linearization of its Hamiltonian vector field. As the examples on billiards of Section \ref{sec:billiard_ball} show, this condition is not always satisfied at the boundary, as for billiard maps one expects $C^1$-extensions there (this is not an issue in dimension $2$, however, as control on the index growth can be replaced by a filtration by homotopy classes as in \cite{MvK20b}).
 \end{remark}

Here is a simple example illustrating why sufficient smoothness of the Hamiltonian is relevant for the purposes of fixed points.

\begin{example}[Integrable twist maps]\label{ex:inttwistmap} Let $M=S^n$ for $n\geq 1$ with the round metric, and $H: T^*M \rightarrow \mathbb{R}$, $H(q,p)=2\pi|p|$ (\emph{not} smooth at the zero section); $\phi_H^1$ extends to all of $\mathbb{D}^*M$ as the identity. It is a positive reparametrization of the Reeb flow at $S^*M$, a full turn of the geodesic flow, and all orbits are fixed points with fixed period. If we smoothen $H$ near $|p|=0$ to $K(q,p)=2\pi g(|p|)$, with $g(0)=g^\prime(0)=0$, then $\tau=\phi_K^1: \mathbb{D}^*M \rightarrow \mathbb{D}^*M$, $\tau(q,p)=\phi_H^{2\pi g^\prime(|p|)}(q,p)$, is now a Hamiltonian twist map. If $g^\prime(|p|)=l/k\in \mathbb{Q}$ with $l,k$ coprime, then $\tau$ has a simple $k$-periodic orbit; therefore $\tau$ has simple interior orbits of arbitrary large period (cf.\ \cite[p.\ 350]{KH95, M86}, for the case $M=S^1$).
\end{example}


\subsection{The weak Hamiltonian twist condition.} 

We now discuss the weak twist condition, for both the degenerate and non-degenerate case.

\medskip

\textbf{Setup (A): Twist condition for degenerate Liouville domains.}\label{sec:twist_degenerate} The following is an adaptation of the definition of Hamiltonian twist map, to the degenerate case.  

\begin{definition}\textbf{(Hamiltonian twist maps on degenerate Liouville domains)} Let $f:(W,\lambda,\alpha_B)\rightarrow(W,\lambda,\alpha_B)$ be a smooth map on a degenerate Liouville domain. We say that $f$ is a \emph{Hamiltonian twist map}, if the following hold:

\medskip

\begin{itemize}
    \item $f\vert_{\mbox{int}(W)}=\phi_H^1$ is generated by a $C^2$ Hamiltonian $H_t:\mbox{int}(W)\rightarrow \mathbb{R}$;

\medskip
    
    \item The Hamiltonian $H_t$ admits a $C^0$-extension (but not necessarily a $C^1$-extension) to the boundary, but the maps $f$ extends smoothly; and

\medskip
    
    \item Near the boundary $B$, the generating Hamiltonian vector field satisfies $h_t=\alpha_B(X_{H_t})>0$ near $B$, and satisfies $h_t\uparrow +\infty$ as one approaches $B$.
\end{itemize}
 
We say that the generating Hamiltonian isotopy $H_t$ is \emph{strictly twisting} or \emph{positively wrapping}. We say that the map satisfies the \emph{strong} weak condition.
\end{definition}

\subsection{Setup (B): Twist condition for Liouville domains.} Recalling from section \ref{sec:degenerate} that to a degenerate Liouville domain there corresponds an honest one, we can give a definition of $C^0$-Hamiltonian twist maps. Namely, these are the self-maps $f_Q:(W,\omega_Q)\rightarrow (W,\omega_Q)$ that correspond under degeneration to Hamiltonian twist maps on degenerate Liouville domains.

\begin{definition}\label{def:C0Hamtwistmaps} \textbf{($C^0$-Hamiltonian twist maps on Liouville domains)} Let $f: (W,\omega)\rightarrow (W,\omega)$ be a map on a Liouville domain. We say that it is a \emph{$C^0$-Hamiltonian twist map}, if it corresponds under degeneration to a Hamiltonian twist map on a degenerate Liouville domain. More precisely:

\begin{itemize}
    \item $f\vert_{\mbox{int}(W)}=\phi_H^1$ is generated by a $C^2$-Hamiltonian $H_t: \mbox{int}(W)\rightarrow \mathbb R$;

\medskip
    
    \item Both $f$ and the Hamiltonian $H_t$ admit $C^0$-extensions to the boundary, but not necessarily $C^1$-extensions; and

\medskip
    
    \item Near the boundary $B$, the generating Hamiltonian vector field satisfies $h_t=\alpha(X_{H_t})$, where $h_t>0$, and $h_t\uparrow +\infty$ as we approach $B$.
\end{itemize}

We say that the isotopy $H_t$ is \emph{infinitely strictly wrapping} or \emph{infinitely positively wrapping}, or simply \emph{infinitely wrapping}. We say that the map satisfies the \emph{weak} twist condition.

\end{definition}

The clear advantage of the weak twist condition with respect to the strong one is that it is open.

\subsection{Local models near the boundary} The following computations are key in order to understand the above definitions. If $(W,\omega)$ is a Liouville domain with $\omega=d(r\alpha_B)$ near $B$, and $E_t:(W,\omega)\rightarrow \mathbb{R}$ is an arbitrary $C^1$-Hamiltonian, its Hamiltonian vector field near $B$ is given by
\begin{equation}\label{eq:general}
X_{E_t}=(\partial_rE_t)R_\alpha +\frac{1}{r}\left(X_{E_t}^\xi-dE_t(R_\alpha)Y\right),
\end{equation}
where $R_\alpha$ is the Reeb vector field of $\alpha$, $Y=r\partial_r$ is the Liouville vector field in $W$, and $X_{E_t}^\xi$ is defined implicitly via $i_{X_{E_t}^\xi}d\alpha_B=-dE_t\vert_{\xi}$. 

Now, recall Remark \ref{rk:Hams}, by which we assume Hamiltonians on the degenerate side of the picture are $C^1$. In this vein, if $E_t=H_t\circ Q$ where $Q:W\rightarrow W$ is a square root map (see Section \ref{sec:degenerate}), $Q=(F,id)$ near $B$, and $H_t:(W,\lambda_S,\alpha_B)\rightarrow \mathbb R$ is $C^1$ on the non-degeneration of $(W,\omega)$, then the above formula implies, noting that the derivatives of $E_t$ along the $B$-direction coincide with those of $H_t$, that
\begin{equation}\label{eq:XQ}
X_{E_t}=\left[((\partial_rH_t)\circ Q)\cdot F^\prime(r)\right]R_\alpha +\frac{1}{F(r)}\left[\left(X_{H_t}^\xi-dH_t(R_\alpha)Y\right)\circ Q\right].
\end{equation}
This is clearly not defined at $r=1$ since $F^\prime(r)$ gives a pole for the first summand, $X_1$, at $r=1$. The second summand $X_2$ is irrelevant for our purposes (as it is bounded), since the first summand gives the infinite wrapping. 

Indeed, recall from the proof of Lemma \ref{lemma:degeneration} that $$F'(r)=\frac{1}{\partial_r A \circ Q'},$$ and at $r=1$ we have have $\partial_r A\vert_{r=1}=0$, since $F(1)=1$, so this gives the pole, as expected. Therefore, the coefficient of $X_1$ is
$$
C^Q_t:=\frac{\partial_r H_t \circ Q}{\partial_r A \circ Q}.
$$
Noting that the denominator of $C_t^Q$ has a zero, we have then the following cases:

\begin{itemize}
    \item[(1)] $\partial_rH_t > 0$ near $r=1$, i.e.\ the numerator of $C_t^Q$ is strictly positive, and therefore $C_t^Q$ has a (positive) pole, since $\partial_r A > 0$. This is the twist map case.

    \medskip
    
    \item[(2)] $\partial_r H_t\vert_{r=1} \equiv 0$, with the same order as $\partial_rA\vert_{r=1}\equiv 0$. In this case, $C_t^Q$ extends continuously to the boundary, and could have mixed sign.

    \medskip

    \item[(3)] $\partial_rH_t\vert_{r=1}<0$ near $r=1$, in which case the numerator of $C_t^Q$ is strictly positive, and therefore $C_t^Q$ has a (negative) pole, since $\partial_r A > 0$. 

    \item[(4)] $\partial_rH_t\vert_{r=1}$ changes sign, in which case $C_t^Q$ has both positive and negative poles at $r=1$.
\end{itemize}

\subsection{Index and Action growth}

The Hamiltonian twist condition will be used to extend the Hamiltonian to a Hamiltonian that is admissible for computing symplectic homology. The extended Hamiltonian can have additional $1$-periodic orbits and these, as well as $1$-periodic orbits on the boundary, need be distinguished from the interior periodic points of $\tau$. There are two mechanisms to separate the boundary orbits from the interior ones: \emph{index growth}, and \emph{action growth}. We will mostly focus on the latter, and leave index growth to a digression.

\section{A generalized Poincar\'e--Birkhoff theorem}

Within the context of Case (B) of the trade-off of Section \ref{sec:trade-off}, we propose the following generalization of the Poincar\'e--Birkhoff theorem.

\begin{theorem}[Limoge--Moreno \cite{LM25}, based on Moreno--van Koert \cite{MvK20b}, Generalized Poincar\'e--Birkhoff theorem, \textbf{Long interior orbits}]\label{thm:main_thm_LM}
Let $f:(W,\omega)\rightarrow (W,\omega)$ be a $C^0$-Hamiltonian twist map on a Liouville domain. Assume the following:
\begin{itemize}
\item\textbf{(fixed points)} All fixed points of $f$ are isolated;
\item \textbf{(First Chern class)} $c_1(W)=0$ if $\dim W \geq 4$;
\item\textbf{(Symplectic homology)} $SH_\bullet(W)$ is non-zero in infinitely many degrees.
\end{itemize}
Then $f$ has simple interior periodic points of arbitrarily large (integer, minimal) period.
\end{theorem}

\begin{remark}\label{rk:thmA} Let us discuss some aspects of the theorem:

\begin{enumerate}

    \item\textbf{(Twist condition and regularity)} We should emphasize that the above result, as opposed to the one in \cite{MvK20b}, addresses the natural boundary degeneracy explained in the trade-off of Section \ref{sec:trade-off}. Moreover, it drops the assumption on index-definiteness, which in particular means that one no longer needs to assume that the contact structure at the boundary is symplectically trivial. What is more, the mechanism used for separating boundary from interior orbits only needs that the Hamiltonian is $C^1$, but not $C^2$ (as it uses action and not index). However, the main proof needs $C^2$ regularity, as it uses the grading.

    \medskip

    \item\textbf{(Fixed points)} The assumption on fixed points is not really an assumption, as failure would imply that there are already infinitely many periodic orbits (of bounded period, however), as e.g.\ in the 2-body problem, in which case the map is the identity. One could rephrase the conclusion by saying: either infinitely many orbits of unbounded period, or infinitely many orbits of bounded period. ``Generically'', one expects finitely many fixed points.

    \medskip

    \item\textbf{(Grading)} We need impose the assumptions $c_1(W)\vert_{\pi_2(W)}=0$ (i.e.\ $W$ is symplectic Calabi--Yau) to have a well-defined integer grading on symplectic homology. This is a mild condition, as e.g.\ it holds for cotangent bundles over orientable basis.
    
    \medskip

    \item\textbf{(Symplectic homology)} One can relax the assumption on symplectic homology to the case $SH_\bullet(W)\neq 0$, by adaptation of Ginzburg's arguments in the proof of the Conley conjecture, see \cite{Gi10,H11}. Conjecturally, this condition should be equivalent to being infinite-dimensional. The orthogonal case of vanishing symplectic homology is related to recent work of Masci \cite{Mas25}, who obtained a Poincar\'e--Birkhoff theorem for asymptotically linear Hamiltonian maps in Euclidean space. Perhaps similar ideas can be used for Liouville domains with vanishing symplectic homology.

    \medskip

    \item\textbf{(Non-degeneracy)} We should emphasize that there is no assumption on non-degeneracy on the Hamiltonian twist map; we deal with this situation similarly as in Ginzburg's proof of the Conley conjecture \cite{Gi10}, i.e.\ by appealing to local symplectic homology.

    \medskip

    \item\textbf{(Surfaces)} If $W$ is a surface, then the condition that $SH_\bullet(W)$ is infinite dimensional just means that $W\neq D^2$; for $D^2$ we have $SH_\bullet(D^2)=0$, and a rotation on $D^2$ gives an obvious counterexample to the conclusion. 
    In the surface case, the argument simplifies, and one can simply work with homotopy classes of loops rather than the grading on symplectic homology. The Hamiltonian twist condition recovers the classical twist condition for $W=\mathbb{D}^*S^1$, due to orientations, and hence the above is clearly a version of the classical Poincaré--Birkhoff theorem.
  
   \medskip
    
    \item\textbf{(Cotangent bundles)} The symplectic homology of the cotangent bundle of a closed manifold is infinite dimensional, due to a result of Viterbo \cite{V99, V18} (see also \cite{AS,SW06}), combined e.g.\ with a theorem of Gromov \cite[Sec.\ 1.4]{G78}.

     \medskip
    
    \item\textbf{(Long orbits)} If $W$ is a global hypersurface of section for some Reeb dynamics, with return map $\tau$, interior periodic points with long (integer) period for $\tau$ translates into spatial Reeb orbits with long (real) period. See Appendix C in \cite{MvK20b}.
    
     \medskip
    
    \item\textbf{(Katok examples)} There are well-known examples due to Katok \cite{K73} of Finsler metrics on spheres with only finitely many simple geodesics, which are arbitrarily close to the round metric. Moreover, they admit global hypersurfaces of section with Hamiltonian return maps, for which the index-definiteness and the condition on symplectic homology both hold. It follows that the return map does not satisfy the twist condition for any choice of Hamiltonians.
    
     \medskip
    
    \item\textbf{(Spatial CR3BP)} From the above discussion and \cite{MvK}, we gather: the only standing obstruction for applying the above result to the spatial CR3BP, is checking the weak Hamiltonian twist condition. This would give a proof of existence of \emph{spatial} long orbits in the spirit of Conley \cite{C63}, which could in principle be collision orbits (these may be excluded, at least perturbatively, by different methods). Since the geodesic flow on $S^2$ arises as a limit case (i.e.\ the Kepler problem), it should be clear from the discussion on Katok examples that this is a subtle condition, i.e.\ one cannot simply argue perturbatively to obtain the twist condition. In \cite{MvK}, we have computed a generating Hamiltonian for the integrable case of the rotating Kepler problem (see Section \ref{sec:RKP_returnmap}); it does \emph{not} satisfy the twist condition in the spatial case (in the planar case, a Hamiltonian twist map was essentially found by Poincar\'e). This does not mean a priori that there is not \emph{another} generating Hamiltonian which does, but this seems rather difficult to check.

\end{enumerate}

\end{remark}

As a particular case of Thm.\ \ref{thm:GBP}, we state the above result for star-shaped domains in cotangent bundles, as a case of independent interest (cf.\ \cite{H11}), and with the simplest statement:

\begin{thm}[Limoge--Moreno \cite{LM25}, based on Moreno--van Koert \cite{MvK20b}]
\label{thmB}
Suppose that $W$ is a fiber-wise star-shaped domain in the Liouville manifold $(T^*M,\lambda_{can})$, where $M$ is simply connected, orientable and closed, and assume that $\tau:W\to W$ is a $C^0$-Hamiltonian twist map. Then $\tau$ has infinitely many periodic points.
\end{thm} 

\subsection{Sketch of proof of Theorem \ref{thm:GBP}} We now provide a brief sketch of the proof of Theorem \ref{thm:GBP}, and refer to \cite{Mvk20b,LM25} for details. The main proofs of \cite{Mvk20b} and \cite{LM25} agree, the difference relying on the setup, and the mechanism for differentiating boundary from interior orbits. In what follows, we will solely focus on the setup for \cite{LM25}.

\medskip
We first address a $C^1$-version, which is needed for the proof. For this, we have the following definition.

\begin{definition}\textbf{($C^1$-Hamiltonian twist map)} Let $f: (W,\omega)\rightarrow (W,\omega)$ be a map on a Liouville domain. We say that it is a \emph{$C^1$-Hamiltonian twist map}, if
\begin{itemize}
    \item\textbf{(Hamiltonian)} $f=\phi_H^1$ is generated by a $C^1$ Hamiltonian $H_t:W\rightarrow \mathbb R$;
    \item\textbf{(Weakened Twist Condition)} At the boundary $B$, the generating Hamiltonian vector field satisfies $h_t:=\alpha(X_{H_t})>0$.
\end{itemize}

\end{definition}

Let $(W,\omega = \mathrm{d}\lambda)$ be a Liouville domain with boundary $B$, and $f : W \to W$ a $C^1$-Hamiltonian twist map. In what follows, we will further need to make the following assumptions.

\medskip

\textbf{Assumptions.} \textbf{(Quantitative twist condition)}\label{cond:QTC}
    $f$ (or $H_t$) is said to satisfy the \emph{quantitative} twist condition if it can be generated by a $C^1$-Hamiltonian $H_t : W\to\mathbb{R}$ such that

    \smallskip

    \begin{enumerate}
        \item $H_t|_{B} > 0$,

        \item $\min\limits_{B} h_t > \max\limits_{B} H_t$.
    \end{enumerate}

    \smallskip
    Here, $$h_t = \partial_r H_t=\alpha(X_{H_t})=\langle X_{H_t}, R_\alpha\rangle,$$ with $r$ the Liouville coordinate near the boundary, so that $Y = r\partial_r$ is the Liouville vector near $B$. Condition (2), when combined with (1), gives a quantitative version of the weak twist condition.

\smallskip

Note that the quantitative twist condition cleary implies the weakened twist condition. As it turns out, the quantitative twist condition will be for free in the case that the $C^1$-Hamiltonian twist map comes from smoothing a $C^0$-Hamiltonian twist map, as in Section \ref{sec:smoothing} (which is the only case of interest). The key estimate is the following.

\begin{proposition}[Action growth, \cite{LM25}]\label{prop:actiongrowth}
    Let $(W, \omega = \mathrm{d}\lambda)$ be a Liouville domain with boundary $B$, and $H_t : W\to \mathbb{R}$ a Hamiltonian satisfying the quantitative twist condition. Then we can construct an extension $\widehat{H}$ of $H$ to the completion $\widehat W$ such that $\widehat H$ is linear at infinity, i.e 

    \vspace{-0.8em}

    \begin{equation*}
        \widehat{H} = ar - \varepsilon
    \end{equation*}

    \smallskip
    
    \hspace{-1.1em} and there exist constants $c > 0, d \in \mathbb{R}$ such that for every trajectory $x : [0,T] \to [1,+\infty)\times B$ we have

    \vspace{-0.6em}

    \begin{equation*}
        \mathcal{A}_{\widehat H}(x) < -c\cdot T + d.
    \end{equation*}
    Moreover, the constant $c$ grows to infinity if $\partial_r H$ grows to infinity in $C^0$-norm.
\end{proposition}

In particular, the action of the trajectories in the cylindrical end goes to $-\infty$ linearly with $T$.

\section{Smoothing a $C^0$-Hamiltonian twist map}\label{sec:smoothing} The key idea for the proof of Theorem \ref{thm:main_thm_LM} is to approximate a given $C^0$-Hamiltonian twist map $f$ by a family of $C^1$-Hamiltonian twist maps $f_\epsilon$, where $f_\epsilon$ is generated by a $C^1$ Hamiltonian $H_\epsilon$ and converges to $f$ as $\epsilon \rightarrow 0$, in such a way that the slopes of the extensions $\widehat H_\epsilon$ grow to infinity. Then taking a limit, we compute the symplectic homology of $W$, and then the arguments of \cite{MvK20b} apply.

\medskip

Let $f:(W,\omega)\rightarrow (W,\omega)$ be a $C^0$-Hamiltonian twist map on a Liouville domain, as in the statement, with infinitely wrapping generating Hamiltonian $H_t: W \rightarrow \mathbb R$. We will prove the following.

\begin{theorem}\label{thm:smoothing}
    For $\epsilon \geq 0$, there exists a family of $C^1$-Hamiltonian twist maps $f_\epsilon$ on a Liouville domain $(W,\omega_\epsilon)$, such that:
    \begin{itemize}
        \item $f_\epsilon$ is generated by a $C^1$ Hamiltonian $H_\epsilon = H_{t,\epsilon}$ which converges in $C^0$ to $H_t$ as $\epsilon \rightarrow 0$.
        \item Along $B$, the function $h_{t,\epsilon}=\alpha(X_{H_{t,\epsilon}})=\partial_r H_{t,\epsilon}$ diverges uniformly and monotonically as $\epsilon \rightarrow 0$, but all derivatives of $H_{t,\epsilon}$ in directions tangent to $B$ remain uniformly bounded. 
    \item As $\epsilon \rightarrow 0$, $\omega_\epsilon$ converges to $\omega$ in $C^\infty$.
    \item The completion $\widehat \omega_\epsilon$ on $\widehat W$ is independent of $\epsilon>0$, i.e.\ it agrees with $\widehat \omega$.
    \item The slope of the extension $\widehat H_\epsilon$ on $\widehat W$ is bounded below by $C/\epsilon$ with $C>0$, and so monotonically diverges as $\epsilon \rightarrow 0$.
    \end{itemize}
\end{theorem}

A direct corollary of the above is the following.

\begin{corollary}\label{cor:limit} Let $H=H_t: W\rightarrow \mathbb R$ be an infinitely wrapping Hamiltonian on a Liouville domain $W$ generating a $C^0$-Hamiltonian twist map, and let $H_\epsilon=H_{t,\epsilon}$ be as in Theorem \ref{thm:smoothing}. Then:
\begin{itemize}
    \item We have $$
\lim_{\epsilon \rightarrow 0}HF(\widehat W,\widehat{H}_\epsilon)=SH(W).
$$
\item For $\epsilon$ sufficiently small, $H_\epsilon$ satisfies the quantitative twist condition.
\end{itemize}

\end{corollary}

\begin{proof}[Proof of Theorem \ref{thm:smoothing}] We view $(W,\omega=\omega_Q=Q^*\omega)$ as the non-degeneration of its degeneration $(W,\omega_S)$ via a square root map $Q: (W,\omega_Q)\rightarrow (W,\omega_S)$ as in Lemma \ref{lemma:degeneration}, and we now take a smooth approximation of $Q$. Assume $$Q(r,b)=(F(r),b)=(1-\varphi(1-r),b)$$ on a collar $(1-\epsilon,1]\times B$ as in Lemma \ref{lemma:degeneration}, where $B=\{r=1\}$. In particular, $F'(r)=\varphi'(1-r)$ has a pole at $r=1$. We write $H_t=E_t \circ Q=E_t^Q$ where $E_t$ is a $C^1$ Hamiltonian on $(W,\omega_S)$.

We let $g_\epsilon:[0,1]\rightarrow (0,\infty)$, for $\epsilon>0$, be a family of positive $C^\infty$ functions such that $g_\epsilon(s)=g_0(s)$ for $\epsilon\leq s\leq 1$, where $g_0=\varphi^\prime$, and such that $g_\epsilon(0)=1/\epsilon$. Then the $g_\epsilon$ are smooth truncations of $g_0$ near $s=0$, and $g_\epsilon\rightarrow g_\infty$ in the $C^0$ topology. We define the functions $$\varphi_\epsilon(s)=\int_0^sg_\epsilon(x)dx,$$ which are smooth, $\varphi_\epsilon(0)=0$ for all $\epsilon$, and $\varphi_\epsilon\rightarrow \varphi_0:=\varphi$ in $C^0([0,1])$. We consider a smooth map
$$
Q_\epsilon:W\rightarrow W
$$
given by the identity away from the collar $(1-\epsilon,1]\times B$, and which coincides with
$$
Q_\epsilon(s,b)=(\varphi_\epsilon(s),b)
$$
near $B$, so that $Q_\epsilon\rightarrow Q_0:=Q$ in $C^0$. Let $$F_\epsilon(r)=1-\varphi_\epsilon(1-r)$$ denote the change of coordinates from the $s$-coordinate to the $r$-coordinate corresponding to $Q_\epsilon$ (so that $Q_\epsilon(r,b) = (F_\epsilon(r),b)$), satisfying $F_\epsilon\vert_{r=1}\equiv 1$ for all $\epsilon$, and converging to $$F_0(r):=F(r)=1-\varphi(1-r).$$ Finally, we define $H_{t,\epsilon}:(W,\omega_\epsilon)\rightarrow \mathbb R$ by $H_{t,\epsilon}=E_t\circ Q_\epsilon$, so that $H_{t,\epsilon}$ coincides with $H_t$ away from the collar, converges to $H_t=E_t^Q$ in $C^0$, and
$$
H_{t,\epsilon}(r,b)=E_t(F_\epsilon(r),b).
$$ Note that their Hamiltonian vector fields, when computed with respect to the symplectic form $\omega_\epsilon=Q_\epsilon^* \omega$ (which looks like $\mathrm{d}(F_\epsilon(r)\alpha)$ near the boundary), are given near $B$ by
\begin{equation}\label{eq:truncation}
X_{H_{t,\epsilon}}=\left[((\partial_rE_t)\circ Q_\epsilon)\cdot g_\epsilon(1-r)\right]R_\alpha +\frac{1}{F_\epsilon(r)}\left[\left(X_{H_t}^\xi-\mathrm{d}H_t(R_\alpha)Y\right)\circ Q_\epsilon\right].
\end{equation}
The effect is that the first summand in this vector field no longer has a pole at $r=1$, and hence $X_{H_{t,\epsilon}}$ is indeed smooth. Moreover, note that as $\epsilon$ goes to zero, the vector field becomes more and more collinear with the Reeb vector field. And while the symplectic form a priori depends on $\epsilon$, since $F_\epsilon$ is positive, the effect is to change the contact form at the boundary, and their completions $\widehat \omega_\epsilon$ on $\widehat W$ are all symplectomorphic. That is, their completion is an $\epsilon$-independent $2$-form $\widehat \omega$ on $\widehat W$, so that the Liouville structure is fixed.

By assumption, we have
$$
X_{H_t}=h_t R_\alpha + Z,
$$
where $h_t>0$, $h_t\vert_{r=1}\equiv \infty$, and $Z$ is linearly independent of $R_\alpha$. If we look at the truncated version, Equation (\ref{eq:truncation}) implies $X_{H_{t,\epsilon}}$ is of the form
$$
X_{H_{t,\epsilon}}=h_{t,\epsilon}R_\alpha + Z_\epsilon,
$$
where $h_{t,\epsilon}>0$, and $$h_{t,\epsilon}\vert_{r=1}=\frac{1}{\epsilon}\cdot (\partial_rH_t)\vert_{r=1}\rightarrow \infty$$ uniformly and monotonically as $\epsilon\rightarrow 0$. In other words, $H_{t,\epsilon}$ generates a $C^1$-Hamiltonian twist map $f_\epsilon$. This admits the following interpretation: if $H_t$ strictly twists at the boundary with respect to $\omega$, then $H_t^Q$ gives a strict ``infinite twist'' at the boundary when computed with respect to $\omega$, and $H_{t,\epsilon}$ is a truncation that strictly twists up to $1/\epsilon$ times the original twisting of $H_t$. Moreover, making $\epsilon$ smaller makes $h_{t,\epsilon}$ uniformly larger in $C^0$ norm, while keeping all $B$ derivatives of $H_{t,\epsilon}$ uniformly bounded in $C^0$ norm. 

We immediately see that the extension $\widehat H_\epsilon$ of $H_\epsilon$ to $\widehat W$ is of the form
$$
\widehat H_\epsilon = C^\epsilon_1 (r-1) + C_0^\epsilon,
$$
with $C_1^\epsilon \geq \max_{B,t}(h_{t,\epsilon})\geq \frac{1}{\epsilon}C$, with $C:=\min_{B,t} \partial_r H_t>0$. This finishes the proof.
\end{proof}

Given a $C^0$-Hamiltonian twist map $\tau: W \rightarrow W$ under the assumptions as in the statement of the theorem, consider the Liouville completion $(\widehat W,\widehat \omega)$ of $(W,\omega=d\lambda)$, and let $H_t$ be a generating and positively wrapping Hamiltonian. The first step is to use the twist condition at the boundary to extend $\tau$ to a map $\widehat \tau$ on $\widehat W$, generated by a Hamiltonian $\widehat H_t$ which is admissible for computing symplectic homology, i.e.\ it is linear at the cylindrical end. For this, one extends $W$ inside $\widehat W$ by adding a small collar neighbourhood to the boundary, and interpolates $H_t$ along it using a careful Taylor expansion in the Liouville direction, to make it linear with slope larger than $\max h_t$ near the new extended boundary. Then symplectic homology can be computed by taking a limit of the Floer homologies of powers of the extension, i.e.\

\begin{equation}\label{eq:limit}
SH_\bullet(W) \cong \varinjlim_{k} HF_\bullet(\widehat{H}^{\natural k}),   
\end{equation}
where $\widehat{H}^{\natural k}$ is the Hamiltonian that generates $\widehat{\tau}^{k}$ in time $1$. Here we use the version of symplectic homology which allows for degeneracies in the Hamiltonians used to define it, which needs that fixed points be isolated as it appeals to local symplectic homology (this can be achieved for the extension by a small perturbation on the extension collar). Moreover, there is a spectral sequence constructed from the local symplectic homologies of all orbits which converges to the global symplectic homology, whose $E^1$-page is given by
\[
E^1_{pq}(\tau)=\bigoplus_{\substack{\gamma \in \mathcal P(H)\\f(p-1)<\mathcal A_H(\gamma)<f(p)}} HF^{loc}_{p+q}(H,\gamma),
\]
where we order the action values of the finitely many (isolated) fixed points in a strictly increasing sequence $\{ a_i \}_{i=1}^k$ and choose a strictly increasing function $f:\N_0\to \R$ such that $f(i) < a_{i+1}<f(i+1)$. The key point is that, while a degenerate orbit may bifurcate into a collection of periodic orbits, the number of such is finite, and their indices lies in a bounded action window (i.e.\ a degenerate orbit is responsible only for finitely many contributions to symplectic homology).

By assumption symplectic homology is non-zero in infinitely many degrees (here we implicitly use that $c_1(W)\vert_{\pi_2(W)}=0$ in order to have a well-defined $\mathbb Z$-grading). Then the isomorphism (\ref{eq:limit}) and the above discussion on degenerate orbits suggests that there should be plenty of generators in the right hand side of (\ref{eq:limit}), which is generated by periodic points of $\widehat \tau$ (i.e.\ fixed points of powers of $\widehat \tau$), all of which should be geometrically different. However, the interpolation procedure potentially creates new orbits in the extension collar, which need to be ruled out. This is where the action growth condition comes in, as it implies that every Reeb orbit at the boundary has action which grows with period. Once this is obtained, one then observes that as there are no contributions from the spurious orbits, as they escape any finite action window and are therefore thrown out in the limit in (\ref{eq:limit}) (here we need to use the natural action filtration in symplectic homology, as used e.g.\ in the above spectral sequence). In other words, all non-trivial contributions to symplectic homology come from interior periodic orbits, which are then orbits of the original map $\tau$. This is the main intuitive idea for the proof. What follows are more details on how this is implemented.

The main argument in the proof, which formalizes what we have already sketched above, and simplifies Ginzburg's proof of the Conley Conjecture \cite{Gi10}. Namely, assuming the interior fixed points of $\tau$ are isolated, further assume by contradiction that the minimal periods of all interior periodic points of $\tau$ are, in increasing order, given by $m_0=1,m_1,\ldots, m_\ell$.
Pick an increasing sequence $\{ p_i\}_{i=1}^\infty$ going to infinity, such that each $p_i$ is indivisible by the $m_1,\ldots, m_\ell$. Due to the choice of $p_i$'s, all fixed points of $\widehat \tau^{p_i}$ are isolated, and we can compute symplectic homology as $SH_*(W)=\varinjlim_i HF_*(\widehat H^{\#p_i}).$ 
Since symplectic homology is non-zero in infinitely many degrees, for all $N>2n k$, where $\dim(W)=2n$, we find distinct degrees $i_1,\ldots,i_N$ such that $SH_{i_j}(W)\neq 0$, ordered by increasing absolute value. By the action growth property, we can choose $p_i$ sufficiently large such that the following hold:
\begin{enumerate}
\item Each fixed point of $\widehat{\tau}^{\# p_i}$ that is contained in $\widehat W \setminus \mbox{int}(W)$ has action whose absolute value is larger than $|i_N|+2n$;
\item the Floer homology groups $HF_{i_j}(\widehat H^{\# p_i})$ are non-trivial for $j=1,\ldots, N$. 
\end{enumerate}

Since we know from (1) that no $1$-periodic orbit in $\widehat W \setminus \mbox{int}(W)$ can contribute to local Floer homology of degree $i_j$, and from (2) that there
must be non-trivial contributions in this degree, we conclude that (using the aforementioned spectral sequence) contributions to symplectic homology must come from the local Floer homology of an orbit $\gamma$ in int$(W)$. Because we have assumed that the $p_i$'s are indivisible by $m_1,\ldots, m_\ell$ we conclude that each such orbit $\gamma$ must be an iterate of one of the orbits $\gamma_1,\ldots,\gamma_k$, say $\gamma=\gamma_j^{p_i}$. Moreover, we may also use the key property of the mean index $\Delta(x)$ of a Hamiltonian trajectory $x$: $$\mbox{supp}\;HF^{loc}_*(\widehat H^{\# p_i},\gamma_j^{p_i})\subset[p_i \Delta(\gamma_j)-n, p_i \Delta(\gamma_j)+n].$$ This is saying what we explained above, i.e.\ that degenerate orbits $\gamma$ are responsible for a finite index window in local symplectic homology (and therefore in symplectic homology via the spectral sequence). This covers at most $2n k$ different degrees, leaving some of the degree $i_j$ uncovered as we had chosen $N>2n k$. This is a contradiction, and finishes the sketch of the proof.

\section{A relative Poincar\'e--Birkhoff theorem}

The relative version of a periodic orbit is a \emph{chord} between two Lagrangian submanifolds, i.e.\ a Hamiltonian path from one to the other. For instance, a chord between two cotangent fibres over two points corresponds to a particle travelling from one point to the other (i.e.\ the velocities are not fixed). In the context of the CR3BP, there are several Lagrangians of geometric and dynamical interest. An example is the collision locus, which corresponds to a collision of the small mass with one of the large masses. Its chords are \emph{consecutive} collisions, i.e.\ the small mass collides once, bounces back, and collides again. As we have seen, this is allowed only as a mathematical artifact through collision regularization. Of course, collision orbits are \emph{not} in and of themselves of interest from a practical perspective. Indeed, navigational engineers would very much prefer to avoid crashing a multi-million dollar spacecraft against the surface of Mars. However, collision orbits can (if non-degenerate as chords) be perturbed to actual orbits that pass close to the large celestial body. These may in turn be used for e.g.\ gravitational slingshots to reach another target, thus saving large amounts of (very expensive!) fuel. See Bolotin--Mackay \cite{BM06} for existence of infinitely many such spatial flyby orbits in the CR3BP (for the perturbative regime where $\mu$ is small), obtained from perturbing non-degenerate spatial consecutive collision orbits in the Kepler problem. Other Lagrangians of interest and their chords detect solar or lunar eclipses (see e.g.\ \cite{Ru22}), or trajectories which are normal to the $xz$-plane (e.g.\ the well-known halo orbits, see e.g.\ \cite{DHRK23} and references therein). 

Theorem \ref{thm:main_thm_LM} also has a relative version addressed at proving the existence of chords, in a very general setting, but inspired by trying to prove the existence of spatial consecutive collision orbits in the CR3BP, as well as orbits normal to the $xz$-plane. In order to give its statement, we need to introduce certain basic notions.

Given a Hamiltonian map $\tau: W\rightarrow W$ of a symplectic manifold $W$, and a Lagrangian $L\subset W$, a chord of $\tau$ of order $m\geq 1$ with respect to $L$ is a pair $(x,m)$, where $x\in L$ is such that $\tau^m(x)\in L$. The minimal order of a chord $(x,m)$ is the minimal $m'$ such that $(x,m')$ is a chord of order $m'$. A $k$-periodic chord is a chord $(x,k)$ of order $k$ such that $\tau^k(x)=x$; its minimal period is the minimal such $k$ (which might differ from its minimal order). 
A $1$-periodic chord is a \emph{fixed} chord. An iterate of a periodic chord $(x,k)$ of minimal period $k$ is a (periodic) chord of the form $(x,n\cdot k)$ for $n\geq 1$. A \emph{sub-chord} of a chord $(x,k)$ of (not necessarily minimal) order $k$ is a chord $(x',l)$ where  $\tau^n(x)=x'$ for some $n\geq 0$ and $l+n\leq k$. 

\begin{remark} For a generic Hamiltonian, end points of chords are never starting points of chords (Lemma 8.2 of \cite{AS10}). Hence, the minimal order is the same thing as the order, and periodic chords of period > 1 simply do not exist. However, since we are interested in systems which are not necessarily generic (e.g.\ the CR3BP), we do not have the luxury of making such genericity assumptions.
\end{remark}

With these notions in place, the statement of the theorem is the following.

\begin{thm}[Moreno--Limoge \cite{ML24, ML25}, \textbf{Long interior chords}]\label{thm:chords}

Suppose that $\tau$ is a $C^0$-Hamiltonian twist map of a connected Liouville domain $(W,\lambda)$. Let $\alpha=\lambda\vert_B$, and let $L\subset (W,\lambda)$ be an exact, spin, Lagrangian with Legendrian boundary. Assume the following.

\medskip

\begin{itemize}

\item \textbf{(Periodic chords)} there are finitely many periodic chords;

\medskip

\item\textbf{(Index-definiteness)} If $\dim \partial W>1$, then assume $c_1(W)\vert_{\pi_2(W)}=0$;

\medskip

\item\textbf{(Wrapped Floer homology)} $HW_*(L,L)$ is supported in infinitely many degrees.
\end{itemize}
Then $\tau$ admits interior chords with respect to $L$, of arbitrary large order, which are not sub-chords of any periodic chord.
\end{thm}

The above is, roughly speaking, obtained from Theorem \ref{thm:GBP} by replacing symplectic homology by wrapped Floer homology. In the case of a cotangent fibre in a cotangent bundle $T^*Q$, $HW_*(L,L)$ is well-known to be infinite dimensional, as it is isomorphic to the homology of the based loop space of $Q$ as in Theorem \ref{thm:AS}. The proof of Theorem \ref{thm:chords} is also analogous to that of \ref{thm:GBP}, with the added difficulty that the introduction of a \emph{local} version of wrapped Floer homology was needed, in order to address possible degeneracies. While at this point these are standard methods, the construction of local wrapped Floer homology, plus the expected local-to-global spectral sequence made from local wrapped Floer homologies but converging to the full wrapped Floer homology, should be considered the main technical contributions of \cite{ML24}.

This result is relevant for the collision locus in the CR3BP, which is the cotangent fibre over the North pole in $W=\mathbb D^*S^2$. If the twist condition were to hold (cf.\ Remark \ref{rk:twist_condition} however), then the above theorem would give infinitely many spatial consecutive collision orbits for the (low-energy, near-primary) dynamics of the CR3BP.

This is illustrated in the following example.

\begin{example}[\textbf{RKP}] We now study collisions in the integrable case of the RKP. Recall from Section \ref{sec:RKP_returnmap} that the return map is
$$
R: P\rightarrow P,
$$
$$
R(\xi_0,\xi_1,\xi_2,0;\eta_0,\eta_1,\eta_2,\eta_3)=
\left(\xi_0,\mathrm{Rot}_{T(c-L)}(\xi_1,\xi_2),0;\eta_0,\mathrm{Rot}_{T(c-L)}(\eta_1,\eta_2),\eta_3\right),
$$
where $$P= \left\{ (\xi;\eta) \in T^*S^3:\;Q(\xi,\eta)=\frac{1}{2},\; \xi_3=0,\; \eta_3 \geq 0 \right\}
$$
is the global hypersurface of section. The collision locus is $\mathcal C=\{(\xi,\eta)\in P: \xi_0=1\}\cong \mathbb D^2$, i.e.\ the $2$-disk cotangent fiber over the north pole $N=(1,0,0,0)\in S^3$ (also the north pole in the base of $P\cong \mathbb{D}^*S^2$). Using that $L\vert_{\mathcal C}=0$, we clearly see that $\mathcal C$ is invariant under $R$, and $R:\mathcal C\rightarrow \mathcal C$ is a rotation by angle $T(c)$, i.e.\ $R\vert_\mathcal{C}=\mathrm{Rot}_{T(c)}$. 

We conclude that there are infinitely many chords of every order $k\geq 1$, although all of them have minimal order $1$. These are non-isolated, since they come in a family parametrized by the disk $\mathcal{C}$, whose boundary circle corresponds to planar chords. The origin is always fixed, and corresponds to the northern polar collision orbit. Moreover, if $T(c)/2\pi$ is irrational, there are no periodic chords except for the origin. If $T(c)=2\pi p/q$ is a rational multiple of $2\pi$, every point in $\mathcal{C}$ different from the origin is periodic of the same minimal period $q$, and so gives a periodic chord of minimal period $q$. The conclusion of the theorem then holds, although the situation described here is of course very non-generic.
    
\end{example}

\section{Digression: index growth, and a further Poincar\'e--Birkhoff theorem} We consider a suitable index growth condition on the dynamics on the boundary of a Liouville domain, which is satisfied in the three-body problem whenever the \emph{planar} dynamics is strictly convex. This assumption allows us to separate boundary and extension orbits from interior ones via the index. This is an alternative way, when compared to the action growth used above, which is weaker, as it imposes more conditions. This is how the first version of a generalized Poincar\'e--Birkhoff theorem was proved in \cite{MvK20b}.

\begin{definition}

We call a strict contact manifold $(Y,\xi=\ker \alpha)$ \emph{strongly index-definite} if the contact structure $(\xi,d\alpha)$ admits a symplectic trivialization $\epsilon$ with the property that there are constants $c>0$ and $d\in \mathbb{R}$ such that for every Reeb chord $\gamma:[0,T]\rightarrow Y$ of Reeb action $T=\int_0^T \gamma^*\alpha$ we have
    $$
    \vert\mu_{RS}(\gamma;\epsilon)\vert\geq c T+d,
    $$
    where $\mu_{RS}$ is the Robbin--Salamon index \cite{RS93}.
    
\end{definition}

Index-positivity is defined similarly, where we drop the absolute value. A variation of this notion was explored in Ustilovsky's thesis \cite{U99}. 
He imposed the additional condition $\pi_1(Y) = 0$, so that index positivity becomes independent of the choice of trivialization, although the exact constants $c$ and $d$ still depend on the trivialization $\epsilon$. The global trivialization is important when considering extensions of our Hamiltonians, as it allows us to measure the index growth of potential new orbits. The point in the above definition is that the index of boundary orbits grows to infinity under iterations of our return map, and so these do not contribute to symplectic homology.

\medskip

A general condition for index-positivity to hold, which is also relevant for the CR3BP, is the following.

\begin{lemma}[\cite{MvK20b}]\label{lemma:indexposapp}
Suppose that $(\Sigma,\alpha)$ is a strictly convex hypersurface in $\mathbb{R}^4$. Then $(\Sigma,\alpha)$ is strongly index-positive.
\end{lemma}

The generalized Poincar\'e--Birkhoff proved in \cite{MvK20b} is then:

\begin{thm}[Moreno--van Koert \cite{MvK20b}. Generalized Poincar\'e--Birkhoff theorem]
\label{thm:GBP}
Suppose that $\tau$ is an exact symplectomorphism of a connected Liouville domain $(W,\lambda)$, and let $\alpha=\lambda\vert_B$. Assume the following.

\begin{itemize}
\item\textbf{(Hamiltonian twist map)} $\tau$ is a Hamiltonian twist map;

\item\textbf{(Fixed points)} all fixed points of $\tau$ are isolated;

\item\textbf{(Index-definiteness)} If $\dim W\geq 4$, then assume $c_1(W)\vert_{\pi_2(W)}=0$, and $(\partial W, \alpha)$ is strongly index-definite; 

\item\textbf{(Symplectic homology)} $SH_*(W)$ is infinite dimensional.
\end{itemize}
Then $\tau$ has simple interior periodic points of arbitrarily large (integer) period.
\end{thm}

\section{Digression: From Hamiltonian twist maps to return maps} While our original motivation to introduce the notion of a Hamiltonian twist was to apply the generalized Poincar\'e--Birkhoff theorem to the return map of the CR3BP, we now focus on the converse. Namely, starting from a Hamiltonian twist map, we realize it as the return map for some adapted open book. 

\medskip

Let $(W,\omega)$ be a $2n$-dimensional Liouville domain with strict contact boundary $(B,\alpha_B)$. Recall that an exact symplectomorphism is a map $\varphi: W\rightarrow W$ satsfying $\lambda-\varphi^*\lambda=dS$ for some smooth function $S:W\rightarrow \mathbb R$ which vanishes near the boundary. The aim of this section is to prove the following.

\begin{thm}
    Let $f:W\rightarrow W$ be a $C^1$-Hamiltonian twist map, and $\varphi: W\rightarrow W$ a compactly supported exact symplectomorphism. Then there is an adapted contact form on $\mathbf{OB}(W,\varphi)$ whose Poincaré return map is $f$.
\end{thm}

We follow a similar strategy as in \cite{AGZ22} (who deal with the case of the disk). We first focus on the case $\varphi=$ Id is trivial. Consider a positively wrapping generating Hamiltonian $H_t$ for $f$, which we may assume to be periodic, i.e.\ $t\in S^1$. By adding a sufficiently large constant to $H_t$, we assume that $H_t\vert_B=C_t>0$. On $S^1\times W$, we then define
$$
\alpha:=H_t dt + \lambda.
$$
To check the contact condition, we compute
\begin{equation}
\begin{split}
\alpha\wedge d\alpha^n&=(H_tdt+\lambda)\wedge (dH_t\wedge dt+\omega)^n\\
&=(H_tdt+\lambda)\wedge (\omega^n+n\cdot dH_t\wedge dt\wedge d\lambda^{n-1})\\
&=H_t dt \wedge \omega^n + n\cdot dH_t\wedge dt \wedge \lambda \wedge d\lambda^{n-1}\\
&=dt\wedge (H_t\omega^n +n\cdot \lambda \wedge dH_t \wedge d\lambda^{n-1}). 
\end{split}
\end{equation}

We can add a sufficiently large constant to $H_t$ which makes the first summand large, without changing the second one. Then $\alpha$ becomes contact along $S^1\times W$.

We note that $n\cdot dH_t\wedge \lambda \wedge d\lambda^{n-1}=\lambda(X_{H_t})\cdot \omega^n$, which one verifies by contracting with $X_{H_t}$ on either side (and using the convention $i_{X_{H_t}}\omega=dH_t$), so that the contact condition is equivalent to
$$
\lambda(X_{H_t})<H_t.
$$
As $\lambda=i_V\omega$ where $V$ is the Liouville vector field, we have $\lambda(X_{H_t})=\omega(V,X_{H_t})=dH_t(V)$, so that the contact condition becomes
\begin{equation}\label{eq:contact_condition}
dH_t(V)<H_t.   
\end{equation}

In a collar neighbourhood $[0,1]\times B$, where $\lambda=r\alpha_B$ and $V=r\partial_r$, this is the condition
\begin{equation}\label{eq:contact_cond2}
r\partial_r H_t<H_t,
\end{equation}
which at $r=1$ is 
$$
h_t<C_t,
$$
where $h_t=\partial_r H_t\vert_{r=1}$. Note that the twist condition implies in particular that $h_t>0$, and that the contact condition implies the quantitative twist condition \ref{cond:QTC}.

The Reeb vector field of $\alpha$ is then 
$$
R_\alpha=\frac{X_{H_t}+\partial_t}{H_t+\lambda(X_{H_t})},
$$
whose first return map is $f$.

In case where $\varphi$ is not trivial, we take a smooth function $F: [0,1]\rightarrow [0,1]$ which vanishes near $0$ and equals $1$ near $1$. We define on $W\times [0,1]$ 
$$
\lambda^\varphi=F(t) \varphi^*\lambda+(1-F(t))\lambda,
$$
$$
H_t^\varphi=F(t)\varphi^*H_t+(1-F(t))H_t,
$$
which agree respectively with $\lambda$ and $H_t$ near $B=\partial W$, and
$$
\alpha:=\alpha^\varphi:=H_t^\varphi dt + \lambda^\varphi.
$$
This gives a well-defined $1$-form on the mapping torus
$$
W_\varphi=W\times \mathbb R/(x,t+1)\sim (\varphi(x),t).
$$
The contact condition is checked as above. Note that $d\lambda^\varphi=\omega$, when $\lambda^\varphi$ is viewed as a $t$-dependent form on $W$. Viewing $H_t^\varphi$ as a $t$-dependent family of Hamiltonians on $W$, we can then define $X_{H_t^\varphi}$ via $i_{X_{H_t^\varphi}}\omega=-dH_t^\varphi$ as before. We then obtain that the Reeb vector field of $\alpha$ is
$$
R_\alpha=\frac{X_{H^\varphi_t}+\partial_t}{H^\varphi_t+\lambda(X_{H^\varphi_t})}.
$$
Note that
$$
dH_t^\varphi=F(t)d(H_t \circ \varphi)+(1-F(t))dH_t,
$$
where we do not differentiate in the $t$-direction. It follows that
$$
X_{H_t^\varphi}=F(t)\varphi^*X_{H_t}+(1-F(t))X_{H_t},
$$
and so indeed $H_t^\varphi$ generates $f$ in time $1$.

\medskip

It remains to show that we can collapse the boundary of $W_\varphi$ and get a well-defined contact form in the quotient. We will make the following mild assumption, to be in a similar setup as in \cite{AGZ22}. We assume that:
\begin{assumptions}
\begin{enumerate}
    \item\label{cond:jet} the $1$-jet of $H_t$ at $B=\partial W$ is $t$-independent.
    \end{enumerate}  
\end{assumptions}

The assumption implies that $\alpha\vert_{S^1\times B}$ is invariant under the $S^1$-action generated by $Y=\partial_t$. The associated moment map is 
$$
\mu=\alpha(Y)=H_t,
$$
which is constant along $B=\{\mu=C\}$. Moreover, $C$ is a regular value, as $\partial_r \mu\vert_{S^1\times B}=h_t=\partial_r H_t\vert_{r=1}>0$. Then, as in \cite{AGZ22}, one can apply Lerman's contact cut construction to obtain a well-defined contact form on $\mathbf{OB}(W,\varphi)=W_\varphi/\sim$, obtained by collapsing the orbits of $Y$ along the boundary $\partial W_\varphi$ to points.

\section{Digression: Morrison's example and outlook} The notions discussed above touch on what is a vastly unexplored line of research: the study of (non-compactly supported) Hamiltonian maps on Liouville domains. This is a natural higher-dimensional version of the study of dynamics on surfaces, which is a huge industry in its own right, and arguably inspired by the planar CR3BP and Poincaré's work. An interesting and concrete starting question is naturally the following.

\medskip

\textbf{Question.} Given a Hamiltonian map $f:W\rightarrow W$ on a Liouville domain, does it have interior periodic points? How many? Are there obstructions on $f$ and/or $W$?

\medskip

The most basic example of a Liouville domain is the ball $W=B^{2n}$. It turns out that in this case, the situation is already different than from the $2$-dimensional disk, where Brouwer's translation theorem applies (and is purely topological). Indeed, the following example due to Morrison \cite{Morr82} gives a Hamiltonian map on the ball of dimension at least $4$ which has no interior fixed point. This means that in the above question, for $W=B^{2n}$, one should indeed impose constraints on the map and/or the dynamics on the boundary.

Consider the Hopf fibration 
$$
\pi: S^{2n-1}\rightarrow \mathbb CP^{n-1}
$$
$$
(z_1,\dots,z_{n})\mapsto [z_1:\dots:z_{n}],
$$

and further consider a map that collapses the closed ball $B^{2n}$ by identifying boundary points if they lie in the same fiber of the Hopf fibration, and coinciding with the identity in the interior. The quotient space is $\mathbb CP^{n+1}$ and the quotient map is
$$
h: B^{2n}\rightarrow \mathbb CP^{n} 
$$
$$
z=(z_1,\dots,z_{n})\mapsto (1-\vert z \vert^2: z_1 : \dots :z_{n}). 
$$

This is a diffeomorphism in the interior, and on the boundary it coincides with the Hopf fibration $\pi$ onto $\mathbb CP^{n-1}=\{z_0=0\}\subset \mathbb CP^{n}$. In other words, $h$ is the map attaching the top cell to $\mathbb CP^{n}$, corresponding to the cell-decomposition $\mathbb CP^{n}=\mathbb CP^{n-1}\cup B^{2n}$. Consider the perfect Morse function
$$
\varphi: \mathbb CP^{n}\rightarrow \mathbb R,
$$
$$
\varphi([z_0:\dots:z_{n}])=\sum_{j=1}^{n+1} j\vert z_{j-1}\vert^2,
$$
with $n+1$ critical points $p_i=[0:\dots: 1 : \dots: 0]$, $i=0,\dots,n$. Pick a point $q \in \mathbb CP^{n-1}=\pi(S^{2n-1})$ different from the $p_i$, e.g.\ $q=[0:1:\dots:1]$. Then one can construct an isotopy $F_t: \mathbb CP^n \rightarrow \mathbb CP^n$ such that $F_0=$Id, $F_t$ fixes $p_i$ for $i=1,\dots,n$, and $F_1(p_0)=q$, say by choosing a narrow tube around a path connecting $p_0$ to $q$, and pushing with the finger. This can be done in a smooth way. Define the function
$$
H: B^{2n}\rightarrow \mathbb R,
$$
$$
H=\varphi \circ F_1^{-1} \circ h. 
$$
By construction, all critical points of $\varphi \circ F_1^{-1}$ lie in $\mathbb CP^{n-1}$, therefore all the critical points of $H$ lie on the boundary (where they come in critical \emph{circles}, i.e.\ fibers of the Hopf fibration). The time-$\epsilon$ flow of the Hamiltonian vector field $X_H$ (which is tangent to the boundary) is then a Hamiltonian map whose fixed points are in correspondence with critical points of $H$ for sufficiently small $\epsilon$, and hence has no interior fixed points. This finishes the construction. Note that this of course does not work for $n=1$.

\medskip

We finish this section with a version of the Arnold conjecture for Hamiltonian twist maps.

\begin{conjecture}
If $f:B^{2n}\rightarrow B^{2n}$ is a $C^0$-Hamiltonian twist map with $n\geq 2$, with no boundary fixed points, then it has a fixed point in the interior. More generally, if $f:W\rightarrow W$ is a non-degenerate $C^0$-Hamiltonian twist map on a Liouville domain with no boundary fixed points, then
$$
\# \mathrm{Fix}(f)\geq \sum_{i=0}^{\dim(W)} b_i(W),
$$
where $b_i(W)=\dim_\mathbb Q H^i(W;\mathbb Q)$ is the $i$-th Betti number of $W$.
\end{conjecture}

Note that in the Morrison construction all of the fixed points have been pushed to the boundary, where they form critical circles, so that the assumptions of the conjecture rule out this example. The Hamiltonian that generates the Morrison map is moreover not one for which one can define any reasonable Floer homology.

\chapter{Symplectic dynamics}\label{ch:symp_dynamics}

This chapter is devoted to outlining the basic notions of \emph{symplectic dynamics}. This is a framework introduced by Helmut Hofer, which integrates notions from symplectic geometry and dynamics, and in particular builds on Floer theory and the theory of pseudo-holomorphic curves. As opposed to the perturbative methods described in the historical remarks of Section \ref{sec:historical_remarks} of Chapter \ref{ch:contact_geometry_in_the_CR3BP}, this is is an inherently \emph{non-perturbative} approach to the study of dynamics. This will pave the way to explain the main construction of \cite{M20}, which we do at the end of the chapter.

\section{Pseudo-holomorphic curves in symplectizations}

\subsection{Finite energy foliations} Finite energy foliations of contact $3$–manifolds were introduced by Hofer, Wysocki and Zehnder in \cite{HWZ98,HWZ03}. These have had several uses in the study of Reeb dynamics, symplectic fillings, and pseudo-rotations of the disk, see e.g.\ \cite{BH,Br15a,Br15b,FH18,Wen08,Wen10}, and references therein. We will now discuss the basic definitions. 

Consider the symplectization $(M,\omega)=(\mathbb{R}\times \Sigma,d(e^t\alpha))$ of a contact $3$-manifold $(\Sigma,\alpha)$. Its tangent space splits as $TM=\xi \oplus \langle \partial_t,R_\alpha\rangle$. A (cylindrical, $\alpha$-compatible) almost complex structure is an endomorphism $J\in \mbox{End}(TM)$ satisfying:
\begin{itemize}
    \item $J^2=-\mathds{1}$;
    \item $J(\xi)=\xi$, $J(\partial_t)=R_\alpha$;
    \item $J$ is $\mathbb{R}$-invariant;
    \item $g=d\alpha(\cdot, J\cdot)$ defines a $J$-invariant Riemannian metric on $\xi$.
\end{itemize}

A \emph{$J$-holomorphic curve} is then a map 
$$
u=(a,v):(\dot S,j)\rightarrow (M,J),
$$
where $\dot S=S\backslash \Gamma$ is a Riemann surface $(S,j)$ with a finite collection $\Gamma$ of points removed, such that $u$ intertwines the complex structures, i.e.\ it satisfies the non-linear \emph{Cauchy--Riemann} equation
$$
J \circ d u = d u \circ i.
$$
The \emph{Hofer-energy} of such a curve is the quantity
$$
\mathbf{E}(u)=\sup_{\varphi \in \mathcal{P}} \int_{\dot{S}} u^*\omega_\varphi,
$$
where $\mathcal{P}=\{\varphi:\mathbb{R}\rightarrow (0,1):\varphi^\prime \geq 0\}$ is the set of orientation preserving diffeomorphisms between $\mathbb{R}$ and $(0,1)$, and $\omega_\varphi=d(e^{\varphi(t)}\alpha)$ is a symplectic form. The choice of $J$ implies that the integrand is point-wise non-negative and so $\mathbf{E}(u)\geq 0$. We say that the holomorphic curve has \emph{finite energy} if $\mathbf{E}(u)<+\infty$. A fundamental property is that non-constant finite energy $J$-holomorphic curves are asymptotic to closed Reeb orbits at their punctures (originally noted by Hofer in his proof of the Weinstein conjecture for overtwisted contact $3$-manifolds).

\begin{proposition}[\cite{H93}]\label{Hofer_result}
Let $u=(a,v):\dot{\mathbb D}=\mathbb D\backslash \{0\}\rightarrow \mathbb R\times \Sigma$ be a non-constant $J$-holomorphic curve of finite energy. If it is bounded, then it extends to a $J$-holomorphic map on $\mathbb D$. Otherwise it is proper, $-\lim_{R\rightarrow 0}\int_{\partial \mathbb D_R} v^*d\alpha:=T\neq 0$ and is finite, and there exists a sequence $R_k\rightarrow +\infty$ such that $$\lim_{k\rightarrow +\infty} u(e^{-2\pi(R_k+ it)})=\gamma(tT),$$ for a closed Reeb orbit $\gamma$ of period $\vert T \vert$.
\end{proposition}
Moreover, if $\gamma$ is non-degenerate, the above convergence is exponential and we have $\lim_{R\rightarrow +\infty} u(e^{-2\pi(R+ it)})=\gamma(tT)$, $\lim_{R\rightarrow +\infty} a(e^{-2\pi(R+ i t)})/R=T$. The puncture is \emph{positive} or \emph{negative} depending on the sign of $T$. Therefore we have a decomposition $\Gamma=\Gamma^+\cup \Gamma^-$ into positive and negative punctures.

A \emph{trivial cylinder} over a $T$-periodic orbit $\gamma$ is the finite energy cylinder
$$
u:\mathbb R\times S^1\rightarrow \mathbb R\times \Sigma,
$$
$$
u(s,t)=(Ts,\gamma(tT)).
$$
The above proposition then says that, roughly speaking, finite energy holomorphic curves are asymptotic to trivial cylinders at  the punctures.

\begin{definition} [\textbf{Finite energy foliation}] A \emph{finite energy foliation} for $(\Sigma, \alpha, J)$ is a smooth two-dimensional foliation $\mathcal F$ of $\mathbb R \times \Sigma$ such that
\begin{itemize}
    \item Each leaf $F \in \mathcal F$ is the image of an embedded
finite energy $J$–holomorphic curve, and there exists a constant
that bounds the energy of every leaf uniformly.
\item  The foliation is $\mathbb R$-invariant, i.e.\ for every leaf $F \in \mathcal F$ and $a \in \mathbb R$, then $F+a\in \mathcal{F}$. 
\end{itemize}
\end{definition}

If $\gamma \subset \Sigma$ is a periodic orbit which is an asymptotic limit for some leaf $u \in \mathcal F$, then the orbit cylinder $\mathbb R \times \gamma$ is also a leaf of $\mathcal{F}$.

Now, a fundamental property for holomorphic planes is \emph{positivity of intersections}; since $M$ is $4$-dimensional, generically two curves intersect at a finite number of points, and if they are holomorphic the intersection numbers are positive. However, there is an an obvious drawback: curves are non-compact and so the classical intersection pairing is not homotopy invariant, since intersections can disappear to infinity. The solution to this issue was provided by Siefring \cite{Sie11}, who, using the very explicit asymptotic behaviour of finite energy curves, defined an intersection pairing with all the desired properties. In particular, it is homotopy invariant, takes into consideration interior intersections as well as those ``coming from infinity'', and two holomorphic curves have vanishing Siefring intersection if and only if their images do not intersect at all. Moreover, in such a case, their projections to $\Sigma$ do not intersect unless their images coincide. 

This has a number of consequences. Denote by $B \subset \Sigma$ the union of all the closed orbits that occur as asymptotic limits for leaves of $\mathcal F$, which we call the \emph{binding} of the foliation. This is the projection down to $\Sigma$ of all the orbit cylinders in $\mathcal F$. Then if $\mathcal F$ is a finite energy foliation, the projections of its
leaves from $\mathbb R \times \Sigma$ to $\Sigma$ yield a smooth foliation of $\Sigma \backslash B$, and each leaf is transverse to the Reeb flow. 

\begin{figure}
    \centering
    \includegraphics{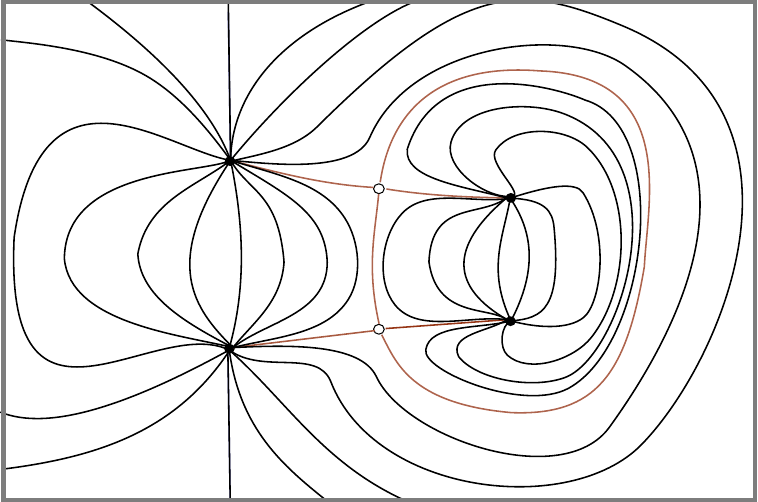}
    \caption{A cross section of a finite energy foliation in $S^3=\mathbb R^3\cup \{\infty\}$, with binding consisting of three orbits, one hyperbolic and two elliptic.}
    \label{fig:enter-label}
\end{figure}

A finite energy foliation is \emph{stable} if its binding consists of non-degenerate orbits, and all the leaves of the foliation are the images of \emph{regular} holomorphic curves (i.e.\ the operator arising by linearizing the CR-equation is surjective, see e.g.\ \cite{Wen08,Wen10b}). In practice, this means that the foliation persists under perturbations of $J$ and $\alpha$.

\section{Reeb dynamics on convex $3$-spheres}\label{sec:HWZ}

This section is devoted to reviewing the groundbreaking work of Hofer--Wysocki--Zehnder for the study of Hamiltonian dynamics in dimension four \cite{HWZ98}. 

We begin with a definition. A connected compact hypersurface $\Sigma \subset \mathbb{R}^{4}$ is said to \emph{be strictly convex} if there exists a domain $W \subset \mathbb{R}^{4}$ and a smooth function $\phi: \mathbb{R}^{4}\rightarrow \mathbb{R}$ satisfying:
\begin{itemize}
    \item[(i)](Regularity) $\Sigma=\{\phi=0\}$ is a regular level set;
    \item[(ii)](Bounded domain) $W = \{z\in \mathbb{R}^{4}: \phi(z) \leq 0\}$ is bounded and contains the origin; and
    \item[(iii)](Positive-definite Hessian) $\nabla^2\phi_{z}(h, h) > 0$ for $z\in W$ and for each non-zero tangent vector $h\in T\Sigma$.
\end{itemize}
In this case, the radial vector field is transverse to $\Sigma$, and so $\Sigma$ is a contact-type $3$-sphere, inheriting a contact form $\alpha$ induced by the standard Liouville form in $\mathbb{R}^{4}$. 

In \cite{HWZ98}, Hofer--Wysocki--Zehnder prove the following.

\begin{thm}[\cite{HWZ98}]\label{thm:HWZ}
A strictly convex hypersurface $(\Sigma,\alpha)\subset \mathbb{R}^4$ has either $2$ or infinitely many periodic orbits.
\end{thm}

The strategy of the proof is finding a disk-like global surface of section, and use the Brouwer--Franks combination mentioned as a heuristic in Chapter \ref{ch:contact_geometry_in_the_CR3BP}. The difficulty is precisely finding the section. These are to be thought of as the (holomorphic) pages of a trivial open book on $\Sigma\cong S^3=\mathbf{OB}(\mathbb{D}^2,\mathds{1})$, which is adapted to the given Reeb dynamics. The rough idea is as follows.

One assumes the existence of a special Reeb orbit $\gamma$, in the sense that is unknotted and linked to every other Reeb orbit (necessary conditions to be the binding of a trivial open book for $S^3$), non-degenerate, has minimal period, and satisfies $\mu_{CZ}(\gamma)=3$. Here, while in general the CZ-index of an orbit depends on a trivialization of the tangent bundle along a choice of disk bounded by $\gamma$, in the case of $S^3$, where $\pi_2(S^3)=0$, this is independent of choices. One then considers the moduli space $\mathcal{M}$ of finite energy $J$-holomorphic planes asymptotic to this Reeb orbit $\gamma$, and having vanishing Siefring self-intersection, modulo the action of $\mathbb{R}$-translation in the image (recall $J$ is $\mathbb{R}$-invariant) and conformal reparametrizations of the domain $\mathbb{C}$. Its expected dimension is $\dim \mathcal{M}=\mu_{CZ}(\gamma)-2=1$, by the Riemann--Roch formula for the Fredholm index. Moreover, the miraculous $4$-dimensional phenomenon of automatic transversality shows that $\mathcal{M}$ is a manifold for any cylindrical $J$. The properties of the Siefring pairing implies that the projections of planes in $\mathcal{M}$ are immersed, do not intersect, and provide a local foliation of $\Sigma$. A further step needed in order to get a global foliation is a way to compactify $\mathcal{M}$. This is provided by Gromov's compactification (or the SFT compactification), obtained by adding strata of nodal curves and ``holomorphic buildings'' with potentially several ``floors''; strictly speaking, these a priori are no longer planes. However, the fact that $\gamma$ is linked to every other orbit can be used to show that no extra strata needs to be added to $\mathcal{M}$, and is in fact \emph{a priori} compact. The result is that $\mathcal{M}\cong S^1$, and projecting the planes in $\mathcal{M}$ to $\Sigma$ provides a global foliation of $\Sigma$. The leaves of this foliation are the $S^1$-family of pages of an open book with binding $\gamma$, and are in fact global surfaces of section for the Reeb dynamics.

While the assumption on the existence of $\gamma$ above might seem far-fetched, it is implied by \emph{dynamical convexity} \cite[Thm.\ 1.3]{HWZ98}. 

\begin{definition}
    We say that $(\Sigma,\alpha)$ is \emph{dynamically convex} if $\mu_{CZ}(\gamma)\geq 3$ for Reeb every orbit $\gamma$. 
\end{definition}

This condition is implied by strict convexity \cite[Thm.\ 3.4]{HWZ98}; intuitively, this implies that there is ``enough winding'' of the linearized Reeb flow along each orbit, so that nearby orbits intersect the leaves of the foliation (and so, at the end of the day when the open book is obtained, we also obtain a return map). The special Reeb orbit is found by first considering the case of an ellipsoid, in which it is explicitly found, then interpolating to the dynamically convex case by considering a symplectic cobordism, and finally using properties of finite energy planes in cobordisms; see Section 4 in \cite{HWZ98}.

\smallskip

\textbf{Conclusion.} The main message to take away from this discussion is that the global surfaces of section are the (holomorphic) pages of a trivial open book on $\Sigma\cong S^3=\mathbf{OB}(\mathbb{D}^2,\mathds{1})$, which is \emph{a posteriori} adapted to the given Reeb dynamics, as they form the leaves of the projection of a stable finite energy foliation of the symplectization. The way that this result ties up with the planar CR3BP is via the Levi-Civita regularization; one says that $(\mu,c)$ lies in the convexity range whenever the Levi-Civita regularization is dynamically convex (cf.\ Prop.\ \ref{prop:doublecover}). The holomorphic open book provided by Hofer--Wysocki--Zehnder, given suitable symmetries, descends to a \emph{rational} open book on the Moser-regularized hypersurface $\mathbb{R}P^3$ (i.e.\ the pages are disks, but their boundary is doubly covered). Alternatively, \cite[Thm.\ 1.18]{HSW} provides an honest open book with annuli fibers for $\mathbb{R}P^3=\mathbf{OB}(\mathbb{D}^*S^1,\tau^2)$, adapted to the planar dynamics. This circle of ideas has also been fruitfully exploited in e.g.\ \cite{H12, H14}; see \cite{HS20} for a very nice survey and references therein, especially for the applications on the planar CR3BP.

\section{Siefring higher-dimensional intersection theory}

This section is devoted to outlining the basic workings of Richard Siefring's (still work in progress) higher-dimensional intersection theory. This is based on the exposition article \cite{MS19}.

The main goal of this section is to understand the intersection properties of punctured
pseudoholomorphic curves with
asymptotically cylindrical pseudoholomorphic hypersurfaces.  The main difficulty 
arises from the noncompactness of the manifolds in question.
Indeed, a punctured pseudoholomorphic curve whose image is not contained in the
holomorphic hypersurface may have punctures limiting to a periodic orbits lying in
the hypersurface.
In this case, it's not a priori clear that the intersection number between the curve and the hypersurface is finite.
Even assuming this intersection number is finite, it is not homotopy invariant as intersections can be lost or created
at infinity. These difficulties can be dealt with via higher-dimensional analogs of the
techniques developed by Siefring in \cite{Sie11}. We explain the basic definitions and results in what follows.

\medskip

Let $M^{2n+1}$ be a closed, orientable manifold.  A pair
$\Ha=(\lambda, \omega)\in\Omega^{1}(M)\times\Omega^{2}(M)$ is called a
\emph{stable Hamiltonian structure} on $M$ if
\begin{itemize}
\item $\lambda\wedge \omega^{n}$ is a volume form on $M$,
\item $d\omega=0$, and
\item $d\lambda \subset \ker \omega$.
\end{itemize}
A stable Hamiltonian structure on $M$ determines a splitting
\[
TM=\R\X\oplus (\xi, \omega|_{\xi})
\]
of the tangent space of $M$ into a symplectic hyperplane distribution
$(\xi=\ker\lambda, \omega|_{\xi})$ and a line bundle given by the span of the 
\emph{Reeb vector field} $\X$, i.e.\ the unique vector field satisfying
\[
\lambda(\X)=1 \qquad\text{ and }\qquad i_{\X}\omega=0.
\]
We will refer to the triple $(M, \lambda, \omega)$ as a \emph{stable Hamiltonian manifold}.
A stable Hamiltonian structure $(\lambda, \omega)$ on $M$ is said to be \emph{nondegenerate} if
all periodic orbits of the Reeb vector field are nondegenerate.

A codimension-$2$ submanifold $V\subset M$ is a
\emph{stable Hamiltonian hypersurface} of $M$ if
the pair $\Ha':=(\lambda', \omega')$ defined by
\[
\lambda':=i^{*}\lambda \qquad \omega':=i^{*}\omega,
\]
where $i:V\hookrightarrow M$ is the inclusion map, is a stable Hamiltonian structure on $V$.
In this case, the hyperplane distribution $\xi':=\ker\lambda'$ is identified via
$i_{*}$ with $TV\cap\xi$.
We say a stable Hamiltonian hypersurface $V\subset M$ is a
\emph{strong stable Hamiltonian hypersurface}
if, in addition, $V$ is invariant under the flow of $\X$.
Along a strong stable Hamiltonian hypersurface $V\subset M$, we obtain a splitting of the tangent space of $M$
\begin{equation}\label{e:splitting}
TM|_{V}=\R \X\oplus(\xi', \omega')\oplus (\xi_{V}^{\perp}, \omega)
\approx\R X_{\Ha'}\oplus (\xi\cap TV, \omega)\oplus (\xi_{V}^{\perp}, \omega)
\end{equation}
into a line bundle spanned by $\X$ and two symplectic vector bundles,
where
\begin{equation}\label{e:symp-comp-def}
\xi_{V}^{\perp}=\{v\in \xi|_{V}:\omega(v, i_{*}w)=0, \forall w\in \xi'\},
\end{equation}
is the symplectic complement of $\xi'\approx TV\cap\xi$ in $\xi|_{V}$. We note that the first two summands give $TV$.
The linearized flow of $\X$ along $V$ preserves this splitting along with the symplectic structure on the
second two summands.

Let $\gamma:S^{1}\approx\R/\Z\to M$ be a $T$-periodic orbit of $\X$, i.e.\ $\gamma$ satisfies the equation
\[
\dot\gamma(t)=T\cdot \X(\gamma(t))
\]
for all $t\in S^{1}$.
Assuming that 
$\gamma(S^{1})\subset V$, 
we can choose a symplectic trivialization of the
hyperplane distribution
$\xi=\xi'\oplus \xi_{V}^{\perp}|_{\gamma(S^{1})}$ along $\gamma$ which respects the splitting
\eqref{e:splitting}, i.e.\ one of the form
\[
\Phi=\Phi_{T}\oplus\Phi_{N}:\xi'\oplus \xi_{V}^{\perp}|_{\gamma(S^{1})} \to
S^{1}\times (\R^{2n-2}, \omega_{0})\oplus(\R^{2}, \omega_{0})
\]
with $\Phi_{T}$ and $\Phi_{N}$ symplectic trivialization of $\xi'|_{\gamma(S^{1})}$ and
$\xi_{V}^{\perp}|_{\gamma(S^{1})}$ respectively.
Given such a trivialization we can define the Conley--Zehnder index of the orbit
$\gamma$ viewed as an orbit in $M$ as usual by
\[
\mu^{\Phi}(\gamma)
:=\mu_{CZ}(\Phi(\gamma(t))\circ d\psi_{Tt}(\gamma(0))\circ\Phi(\gamma(0))^{-1})
\]
where $\psi:\R\times M\to M$ is the flow generated by $\X$ and where
$\mu_{CZ}$ is the CZ-index. But since
$d\psi_{t}$ preserves the splitting \eqref{e:splitting}, we can also consider the
Conley--Zenhder indices that arise from the restrictions of $d\psi_{t}$ to $\xi'|_{\gamma}$ and
$\xi_{V}^{\perp}|_{\gamma}$.
In particular we define
\[
\mu^{\Phi_{T}}_{V}(\gamma)
:=\mu_{CZ}(\Phi_{T}(\gamma(t))\circ d\psi_{Tt}(\gamma(0))|_{\xi'}\circ\Phi_{T}(\gamma(0))^{-1})
\]
which is the Conley--Zehnder index of $\gamma$ viewed as a periodic orbit lying in $V$, and
\[
\mu_{N}^{\Phi_{N}}(\gamma)
:=\mu_{CZ}(\Phi_{N}(\gamma(t))\circ d\psi_{Tt}(\gamma(0))|_{\xi_{V}^{\perp}}\circ\Phi_{N}(\gamma(0))^{-1})
\]
which we will call the \emph{normal Conley--Zehnder index of $\gamma$ relative to $\Phi_{N}$}. These quantites are related by
\[
\mu^{\Phi}(\gamma)=\mu_{V}^{\Phi_{T}}(\gamma)+\mu_{N}^{\Phi_{N}}(\gamma).
\]

An almost complex structure $J$ on $\xi$ is said to be \emph{compatible} with the
stable Hamiltonian structure $(\lambda, \omega)$
if the bilinear form on $\xi$ defined by
$\omega(\cdot, J\cdot)|_{\xi\times\xi}$
is a metric.
We will denote the set of compatible complex structures on $\xi$ by $\J(M, \xi)$.
Given a
strong stable Hamiltonian hypersurface
$V\subset M$
and a choice of $J\in\J(M, \xi)$, $J$ is \emph{$V$-compatible} if $J$ fixes the
hyperplane distribution
$\xi'=TV\cap\xi$ along $V$.
We will denote the set of such complex structures by
$\J(M, V, \xi)$.
Such a $J\in\J(M, V, \xi)$ necessarily also fixes
$\xi_{V}^{\perp}$, as is easily seen. We note that since both the linearized flow $d\psi_{t}$ of $\X$
and a compatible $J\in\J(M, V, \xi)$ preserve the splitting
\eqref{e:splitting}, the asymptotic operator
\[
\A_\gamma h(t):=-J\left.\frac{d}{ds}\right|_{s=0}d\psi_{-Ts}h(t+s)
\]
of a periodic orbit $\gamma$ lying in $V$ also preserves the splitting.  We will write
\[
\A_\gamma=\A_\gamma^{T}\oplus\A_\gamma^{N}:W^{1,2}(\xi')\oplus W^{1,2}(\xi_{V}^{\perp})
\to L^{2}(\xi')\oplus L^{2}(\xi_{V}^{\perp})
\]
to indicate the resulting splitting of the operator.

We then extend 
$J$ to an $\R$-invariant almost complex structure $\tilde J$ on $\R\times M$, so that it maps the $\mathbb R$-direction to the Reeb vector field. Then the submanifold
$\R\times V$ of $\R\times M$ is $\tilde J$-holomorphic.
One can then consider $\tilde J$-holomorphic hypersurfaces which are asymptotic to cylindrical $\tilde J$-holomorphic hypersurfaces of the form
$\R\times V$ with $V$ a
strong stable Hamiltonian hypersurface.
Before giving a more precise definition,
we introduce some more geometric data on our manifold.

Given a $J\in\J(M, \xi)$ we can define a Riemannian metric
\begin{equation}\label{e:metric-M}
g_{J}(v, w)=\lambda(v)\lambda(w)+\omega(\pi_{\xi}v, J\pi_{\xi}w)
\end{equation}
where
$\pi_{\xi}:TM\approx\R \X\oplus\xi\to\xi$
is the projection onto $\xi$ along $\X$.
We can extend $g_{J}$ to a metric $\tilde g_{J}$ on $\R\times M$
by 
\begin{equation}\label{e:metric-RxM}
\tilde g_{J}:=da\otimes da+\pi^{*}g_{J}.
\end{equation}
We will denote the exponential maps of
$g_{J}$ and $\tilde g_{J}$ by $\exp$ and
$\widetilde\exp$ respectively, which are related by
\[
\widetilde\exp_{(a, p)}(b, v)=(a+b, \exp_{p}v).
\]
If $J$ is $V$-compatible for some strong
stable Hamiltonian hypersurface $V\subset M$,
then the symplectic normal bundle $\xi_{V}^{\perp}$ is the
$g_{J}$-orthogonal complement of $TV$ in $TM|_{V}$, and that
$\pi^{*}\xi_{V}^{\perp}$ is the $\tilde g_{J}$-orthogonal complement of $T(\R\times V)$ in
$T(\R\times M)|_{\R\times V}$.
Since $V$ is assumed to be compact, the restrictions of
$\exp$ and $\widetilde\exp$ to $\xi_{V}^{\perp}$ and $\pi^{*}\xi_{V}^{\perp}$ respectively
are embeddings on some neighborhood of the zero sections.

Now consider a pair $V_{+}$, $V_{-}$ of
strong stable Hamiltonian hypersurfaces
and assume that
$V:=V_{+}\cup V_{-}$ is also a
strong stable Hamiltonian hypersurface,
i.e.\ that all components of $V_{+}$ and $V_{-}$ are either disjoint or identical.
We let $J\in\J(M, V, \xi)$ be a $V$-compatible $J$ with associated $\R$-invariant almost complex structure
$\tilde J$ on $\R\times M$.
We will consider $\tilde J$-holomorphic submanifolds which outside of a compact set can be described
by exponentially decaying sections of the normal bundles to $V_{+}$ and $V_{-}$.  More precisely, a $\tilde J$-holomorphic submanifold $\tilde V\subset\R\times M$ is
\emph{positively asymptotically cylindrical over $V_{+}$} and
\emph{negatively asymptotically cylindrical over $V_{-}$} if there exists an $R>0$ and sections
\begin{gather*}
\eta_{+}:[R, +\infty)\to C^\infty(\xi_{V_{+}}^{\perp}) \\
\eta_{-}:(-\infty, -R] \to C^\infty(\xi_{V_{-}}^{\perp}) \\
\end{gather*}
so that
\begin{gather*}
\tilde V\cap\left([R, +\infty)\times M\right)=\bigcup_{(a, p)\in [R, +\infty)\times V_{+}} \widetilde\exp_{(a, p)}\eta_{+}(a, p) \\
\tilde V\cap\left((-\infty, -R]\times M\right)=\bigcup_{(a, p)\in (-\infty, -R]\times V_{-}} \widetilde\exp_{(a, p)}\eta_{-}(a, p)
\end{gather*}
and so that there exist constants $M_{i}>0$, $d>0$ satisfying
\[
\vert\widetilde\nabla^{i}\eta_{\pm}(a, p)\vert\le M_{i}e^{-d \vert a\vert}
\]
for all $i\in\N$ and $\pm a\in [R, +\infty)$, where
$\tilde\nabla$ is the extension of a connection $\nabla$ on  $\xi_{V}^{\perp}$
to a connection $\tilde \nabla$ on $\pi^{*}\xi_{V}^{\perp}$ defined by requiring
$\tilde\nabla_{\partial_{a}}\eta(a, p)=\partial_{a}\eta(a, p)$.
We will refer to the sections $\eta_{+}$ and $\eta_{-}$ respectively as \emph{positive}
and \emph{negative asymptotic representatives of $\tilde V.$}

Before presenting the relevant results it will be convenient
to establish some standard assumptions and notations for the next several
definitions and results.

\begin{assumptions}\label{standing-assumptions}
We assume that:
\begin{enumerate}\renewcommand{\theenumi}{\alph{enumi}}
\item $(M, \lambda, \omega)$ is a closed, manifold with nondegenerate stable Hamiltonian structure
$(\lambda, \omega)$ with $\xi=\ker\lambda$ and $\X$ the associated Reeb vector field,

\item $V_{+}\subset M$, $V_{-}\subset M$, and $V=V_{+}\cup V_{-}$ are
strong stable Hamiltonian hypersurfaces
of $M$,

\item $J\in\J(M, V, \xi)$ is a $V$-compatible complex structure on $\xi$ and
$\tilde J$ is the $\R$-invariant almost complex structure on $\R\times M$ associated to $J$,

\item $\tilde g_{J}$ is the Riemannian metric on $\R\times M$ defined by \eqref{e:metric-RxM} and
$\widetilde\exp$ is the associated exponential map,

\item $\tilde V\subset\R\times M$ is a $\tilde J$-holomorphic hypersurface which is positively asymptotically cylindrical over $V_{+}$ and negatively asymptotically cylindrical over $V_{-}$,

\item $\eta_{+}:[R, +\infty)\to C^\infty (\xi_{V_{+}}^{\perp})$ and
$\eta_{-}:(-\infty, -R]\to C^{\infty}(\xi_{V_{-}}^{\perp})$ are, respectively, positive and negative asymptotic representatives of $\tilde V$,

\item $C=[S, j, \Gamma=\Gamma^{+}\cup\Gamma^{-}, \widetilde{u}=(a, u)]$ is a finite-energy $\tilde J$-holomorphic curve, and
at $z\in\Gamma$, $C$ is asymptotic to $\gamma_{z}^{m_{z}}$
(with $\gamma_{z}^{m_{z}}$ indicating the $m_{z}$-fold covering of a simple periodic orbit $\gamma_{z}$),

\item $\Phi$ is a trivialization of $\xi_{V}^{\perp}$ along every periodic orbit lying in $V$
which occurs as an asymptotic limit of $C$. 

\end{enumerate}
\end{assumptions}

The following theorem can be seen as a generalization of
Theorem 2.2 in \cite{Sie11}.  

\begin{thm}[Siefring, to appear]\label{t:normal-asymptotics}
Under the above assumptions, assume that at $z\in\Gamma^{+}$, $C$ is asymptotic to a periodic orbit
$\gamma_{z}^{m_{z}}\subset V_{+}$.
Then there exists an $R'\in\R$, a smooth map
\[
u_{T}:[R', \infty)\times S^{1}\to [R, \infty)\times V_{+}
\]
and a smooth section
\[
u_{N}:[R', \infty)\times S^{1}\to u_{T}^{*}\pi^{*}\xi_{V_{+}}^{\perp}
\]
so that the map
\begin{equation}
(s, t)\mapsto\widetilde\exp_{u_{T}(s, t)}u_{N}(s, t)
\end{equation}
parametrizes $C$ near $z$.
Moreover,
if we assume that the image of $C$ is not a subset of the asymptotically cylindrical hypersurface $\tilde V$,
then
\begin{equation}\label{e:normal-asymptotics}
u_{N}(s, t)-\eta_{+}(u_{T}(s, t))=e^{\mu s}[e(t)+r(s, t)]
\end{equation}
for all $(s, t)\in[R', \infty)\times S^{1}$
where:
\begin{itemize}
\item $\mu<0$ is a negative eigenvalue of the normal asymptotic operator
$\A_{\gamma_{z}^{m_{z}}}^{N}$,
\item $e\in\ker(\mathbf{A}_{\gamma_{z}^{m_{z}}}^{N}-\mu)\setminus\{0\}$ is an eigenvector with eigenvalue $\mu$, and
\item $r:[R', \infty)\times S^1 \to u_{T}^{*}\pi^{*}\xi_{V}^{\perp}$ is a smooth section satisfying
exponential decay estimates of the form
\begin{equation}\label{e:normal-asymptotics-remainder}
\vert\tilde\nabla^{i}_{s}\tilde\nabla^{j}_{t} r(s, t)\vert\le M_{ij}e^{-d\vert s\vert }
\end{equation}
for some positive constants $M_{ij}$, $d$ and all $(i, j)\in\N^{2}$
\end{itemize}

Similarly, if we 
assume that at $z\in\Gamma^{-}$, $C$ is asymptotic to a periodic orbit
$\gamma_{z}^{m_{z}}\subset V_{-}$., then there exists an $R'\in\R$, a smooth map
\[
u_{T}:(-\infty, R'] \times S^{1}\to (-\infty, -R]\times V_{-}
\]
and a smooth section
\[
u_{N}:(-\infty, R'] \times S^{1}\to u_{T}^{*}\pi^{*}\xi_{V_{-}}^{\perp}
\]
so that the map
\[
(s, t)\mapsto\widetilde\exp_{u_{T}(s, t)}u_{N}(s, t)
\]
parametrizes $C$ near $z$.
Moreover, if the image of $C$ is not contained in $\tilde V$, then 
$u_{N}(s, t)-\eta_{-}(u_{T}(s, t))$ satisfies a formula of the form \eqref{e:normal-asymptotics}
for all $(s, t)\in(-\infty, R']\times S^{1}$, where now:
\begin{itemize}
\item $\mu>0$ is a positive eigenvalue of the normal asymptotic operator
$\A_{\gamma_{z}^{m_{z}}}^{N}$,
\item $e\in\ker(\A_{\gamma_{z}^{m_{z}}}^{N}-\mu)\setminus\{0\}$, as before,
is an eigenvector with eigenvalue $\mu$, and
\item $r:(-\infty, R']\to u_{T}^{*}\pi^{*}\xi_{V}^{\perp}$ is a smooth section satisfying
exponential decay estimates of the form \eqref{e:normal-asymptotics-remainder}
for some positive contants $M_{ij}$, $d$.
\end{itemize}
\end{thm}

The bundles of the form
$u_{T}^{*}\pi^{*}\xi_{V}^{\perp}$
occurring in the statement of this theorem are trivializable
since they are complex line bundles over a space which retracts onto $S^{1}$.
In any trivialization the eigenvector $e$ from formula \eqref{e:normal-asymptotics}
satisfies a linear, nonsingular ODE,
and thus is nowhere vanishing since we assume it is not identically zero.
Since the ``remainder term'' $r$ in the formula \eqref{e:normal-asymptotics}
converges to zero, the functions
$u_{N}(s, t)-\eta_{\pm}(u_{T}(s, t))$ are nonvanishing for sufficiently large $\vert s\vert$.
However, since zeroes of this function can be seen to correspond to intersections between
the curve $C$ and the hypersurface $\tilde V$ occuring sufficiently close to the punctures of $C$, we conclude
that all intersections between $C$ and $\tilde V$ are contained in a compact set.
Moreover, since intersections between $C$ and $\tilde V$ can be shown to be isolated and of positive local order,
we conclude that the algebraic intersection number between $C$ and $\tilde V$ is finite:

\begin{corollary}\label{c:finite-intersections}
Under the above assumptions, assume also that
no component of the curve $C$
has image contained in 
the $\tilde J$-holomorphic hypersurface $\tilde V$.  Then the algebraic intersection number
$C\cdot \tilde V$, defined by summing local intersection indices, is finite and nonnegative, and
$C\cdot \tilde V=0$ precisely when $C$ and $\tilde V$ do not intersect.
\end{corollary}

This corollary deals with the first difficulty in understanding intersections between punctured curves and asymptotically cylindrical hypersurfaces described above, namely the finiteness of the intersection number. A second consequence of the asymptotic formula from Theorem \ref{t:normal-asymptotics}, again stemming
from the fact that the quantities $u_{N}(s, t)-\eta_{\pm}(u_{T}(s, t))$ are nonzero for sufficiently large
$\vert s\vert$, is that the normal approach of the curve $C$ has a well-defined winding number
relative to a trivialization $\Phi$ of $\xi_{V}^{\perp}|_{\gamma_{z}}$.  This winding will be given by
the winding of the eigenvector from formula \eqref{e:normal-asymptotics} relative to $\Phi$, and for a given 
puncture $z$ of $C$, we will denote this quantity by
\[
\wind_{rel}^\Phi((C; z), \tilde V)=\wind(e).
\]
Combining this observation with the characterization of the Conley--Zehnder
index in terms of the asymptotic operator from Definition 3.9/Theorem 3.10 in \cite{HWZ95} leads
to the following corollary.

\begin{corollary}\label{c:normal-rel-winding}
Under the above assumptions, also assume that
no component of the curve
$C=[S, j, \Gamma_{+}\cup\Gamma_{-}, \tilde u=(a, u)]$
has image contained in
the holomorphic hypersurface $\tilde V$.  Then:
\begin{itemize}
\item If $z\in\Gamma$ is a positive puncture at which $\tilde u$ limits to $\gamma_z^{m_z}\subset V_{+}$ then
\begin{equation}\label{e:wind-rel-pos}
\wind_{rel}^{\Phi}((C; z), \tilde V)\le \fl{\mu_{N}^{\Phi}(\gamma_z^{m_z})/2}=:\alpha_N^{\Phi;-}(\gamma_z^{m_z}).
\end{equation}

\item If $z\in\Gamma$ is a negative puncture at which $\tilde u$ limits to $\gamma_z^{m_z}\subset V_{-}$ then
\begin{equation}\label{e:wind-rel-neg}
\wind_{rel}^{\Phi}((C; z), \tilde V)\ge \ceil{\mu_{N}^{\Phi}(\gamma_z^{m_z})/2}=:\alpha_N^{\Phi;+}(\gamma_z^{m_z}).
\end{equation}
\end{itemize}
\end{corollary}

The numbers $\alpha_N^{\Phi;-}(\gamma)$ and $\alpha_N^{\Phi;+}(\gamma)$ are, respectively, the biggest/smallest winding number achieved by an eigenfunction of the normal asymptotic operator of any orbit $\gamma$ corollary responding to a negative/positive eigenvalue. Observe that we have the formulas 
\begin{equation}\label{alphas}
\mu_{N}^{\Phi}(\gamma)= 2\alpha_N^{\Phi;-}(\gamma)+p_N(\gamma)=2\alpha_N^{\Phi;+}(\gamma)-p_N(\gamma),
\end{equation}
where $p_N(\gamma) \in \{0,1\}$ is the \emph{normal parity} of the orbit $\gamma$ (which is independent of the trivialization $\Phi$).

We will see in a moment that this corollary
can be used to 
deal with the second difficulty in understanding intersections between 
punctured curves and asymptotically cylindrical hypersurfaces described above, namely, the fact
that the algebraic intersection number may not be invariant under homotopies.
We first introduce some terminology.
Assuming that
no component of
$C$ is a subset of $\tilde V$, we define the asymptotic
intersection number at the punctures of $C$ in the following way:
\begin{itemize}
\item
If for the positive puncture $z\in\Gamma_{+}$, $\gamma^{m_z}_{z}\subset V_{+}$, we define the 
\emph{asymptotic intersection number $\delta_{\infty}((C; z); \tilde V)$ of $C$ at $z$ with $\tilde V$} by
\begin{equation}\label{e:local-asymp-inum-pos}
\delta_{\infty}((C; z), \tilde V)=\fl{\mu_{N}^{\Phi}(\gamma_z^{m_z})/2}-\wind_{rel}^{\Phi}((C; z), \tilde V).
\end{equation}

\item 
If for the negative puncture $z\in\Gamma_{-}$, $\gamma^{m_z}_{z}\subset V_{-}$, we define the 
\emph{asymptotic intersection number $\delta_{\infty}((C; z); \tilde V)$ of $C$ at $z$ with $\tilde V$} by
\begin{equation}\label{e:local-asymp-inum-neg}
\delta_{\infty}((C; z), \tilde V)=\wind_{rel}^{\Phi}((C; z), \tilde V)-\ceil{\mu_{N}^{\Phi}(\gamma^{m_z}_{z})/2}.
\end{equation}

\item For all other punctures $z\in\Gamma_{\pm}$
(i.e. those for which $\gamma_{z}$ is not contained in $V_{\pm}$),
we define
\begin{equation}\label{e:local-asymp-inum-zero}
\delta_\infty((C; z), \tilde V)=0.
\end{equation}
\end{itemize}
We then define the \emph{total asymptotic intersection number of $C$ with $\tilde V$} by
\begin{equation}\label{e:global-asymp-inum}
\delta_{\infty}(C, \tilde V)=\sum_{z\in\Gamma}\delta_{\infty}((C; z), \tilde V).
\end{equation}

We observe that as a result of Corollary \ref{c:normal-rel-winding} the local and total asymptotic intersection numbers are always nonnegative.

\vspace{0.5cm}

We can use the trivialization $\Phi$ of $\xi_{V}^{\perp}$ along the asymptotic periodic orbits of $C$ lying
in $V$ to construct a perturbation $C_{\Phi}$ of $C$ in the following way.
For each puncture $z\in\Gamma$ for which the asymptotic limit $\gamma_{z}^{m_{z}}$ lies in $V$,
we first extend $\Phi$ to a trivialization
$\Phi:\xi_{V}^{\perp}|_{U_{z}}\to U_{z}\times\R^{2}$ on some open neighborhood 
$U_{z}\subset V$ of the asymptotic limit $\gamma_{z}$.  Then we consider the asymptotic parametrization
\[
(s, t)\mapsto \widetilde\exp_{u_{T}(s, t)}u_{N}(s, t)
\]
from Theorem \ref{t:normal-asymptotics} above
for $(s, t)\in [R,+\infty)\times S^{1}$ or $(-\infty, -R]\times S^{1}$ as appropriate,
where $R>0$ is chosen large enough so that $u_{T}$ has image contained in the neighborhood $U_{z}$ of
$\gamma_{z}$ on which the trivialization $\Phi$ has been extended.
We then perturb the map by replacing the above parametrization of $C$ near $z$ by the map
\[
(s, t)\mapsto \widetilde\exp_{u_{T}(s, t)}\bp{u_{N}(s, t)+\beta(\abs{s})\Phi(u_{T}(s, t))^{-1}\varepsilon}
\]
where $\beta:[0, \infty)\to[0, 1]$ is a smooth cut-off function equal to $0$ for $s<\abs{R}+1$ and equal to
$1$ for $\abs{s}>\abs{R}+2$, and $\varepsilon \ne 0$ is thought of as a number in $\C\approx\R^{2}$.
Given this, we can then define the 
\emph{relative intersection number $i^{\Phi}(C, \tilde V)$
of $C$ and $\tilde V$ relative to the the trivialization $\Phi$}
by
\[
i^{\Phi}(C, \tilde V):=C_{\Phi}\cdot \tilde V.
\]
It can be shown that this number is independent of choices made in the construction of $C_{\Phi}$ provided
the perturbations are sufficiently small.

We now define the
\emph{holomorphic intersection product of $C$ and $\tilde V$} by
\[
C*\tilde V:=i^{\Phi}(C,\tilde V)
+\sum_{\stackrel{z\in\Gamma_{+}}{\gamma_{z}\subset V_{+}}}\fl{\mu_{N}^{\Phi}(\gamma_{z}^{m_{z}})/2}
-\sum_{\stackrel{z\in\Gamma_{-}}{\gamma_{z}\subset V_{-}}}\ceil{\mu_{N}^{\Phi}(\gamma_{z}^{m_{z}})/2}
\]
The key facts about the holomorphic intersection product are now given in the following theorem which generalizes Theorem 2.2/4.4 in
\cite{Sie11}.
\begin{theorem}[Generalized positivity of intersections, Siefring, to appear]\label{t:intersection-positivity}
With $\tilde V$ and $C$ as in the above assumptions, assume that $C$ is not contained in $\tilde V$,
the holomorphic intersection product $C*\tilde V$ depends only on the relative homotopy classes of $C$ and $\tilde V$.
Moreover, if the image of $C$ is not contained in $\tilde V$, then
\[
C*\tilde V=C\cdot \tilde V+\delta_{\infty}(C, \tilde V)\ge 0
\]
where $C\cdot \tilde V$ is the algebraic intersection number, defined by summing local intersection indices,
and $\delta_{\infty}(C, \tilde V)$ is the total asymptotic intersection number, defined by
\eqref{e:local-asymp-inum-pos}--\eqref{e:global-asymp-inum}.
In particular, $C*\tilde V\ge 0$ and equals zero if and only if $C$ and $\tilde V$ don't intersect and
all asymptotic intersection numbers are zero.
\end{theorem}
The proof of this theorem follows very similar lines to Theorem 2.2/4.4 in \cite{Sie11}.
The essential point is that the relative intersection number can be shown to be given by the formula
\[
i^{\Phi}(C, \tilde V)=
C\cdot \tilde V
-\sum_{\stackrel{z\in\Gamma_{+}}{\gamma_{z}\subset V_{+}}}\wind^\Phi_{rel}((C; z), \tilde V)
+\sum_{\stackrel{z\in\Gamma_{-}}{\gamma_{z}\subset V_{-}}}\wind^\Phi_{rel}((C; z), \tilde V).
\]
The result will then follow from Corollary \ref{c:normal-rel-winding} above. 

Analogous to the case in four dimensions studied in \cite{Sie11}, the
$\R$-invariance of the almost complex structure allows one to compute
the holomorphic intersection number and (in some cases) the algebraic intersection
number of a holomorphic curve and a holomorphic hypersurface with respect
to asymptotic winding numbers and intersections of each object with the asymptotic limits
of the other.

Before stating the relevant results we will first make some additional assumptions. We will henceforth assume that:
\begin{assumptions}\label{assumptions-2}$\;$
\begin{enumerate}\renewcommand{\theenumi}{\alph{enumi}}
\item $V_{+}$ and $V_{-}$ are disjoint,
\item $\xi_{V}^{\perp}$ of $V=V_{+}\cup V_{-}$ is trivializable, 
\item $\Phi:\xi_{V}^{\perp}\to V\times \R^{2}$ is a global trivialization,
\item $\tilde V$ is connected, and
\item the projection $\pi(\tilde V)$ of $\tilde V$ to $M$ is an embedded codimension-$1$ submanifold of
$M\setminus V$.
\end{enumerate}
\end{assumptions}
Under these assumptions, $\tilde V$ has a well-defined winding
$\wind_{\infty}^{\Phi}(\tilde V, \gamma)$ relative to $\Phi$
around any orbit $\gamma\subset V=V_{+}\cup V_{-}$
which can be defined by considering the asymptotic representatives $\eta_{+}$ or $\eta_{-}$ as appropriate
and computing
\begin{equation}\label{e:V-gamma-wind-1}
\wind_{\infty}^{\Phi}(\tilde V, \gamma)=\lim_{\abs{s}\to\infty}\wind \Phi^{-1}\eta_{\pm}(s, \gamma(\cdot)),
\end{equation}
or, equivalently by
\begin{equation}\label{e:V-gamma-wind-2}
\wind_{\infty}^{\Phi}(\tilde V, \gamma)=\wind_{rel}^{\Phi}((\R\times\gamma; \pm\infty), \tilde V).
\end{equation}
As in Corollary \ref{c:normal-rel-winding} above, it follows from the
asymptotic formula from Theorem \ref{t:normal-asymptotics} above that
\[
\wind_{\infty}^{\Phi}(\tilde V, \gamma)\le \fl{\mu_{N}^{\Phi}(\gamma)/2}=\alpha_N^{\Phi;-}(\gamma)
\]
if $\gamma\subset V_{+}$ and
\[
\wind_{\infty}^{\Phi}(\tilde V, \gamma)\ge \ceil{\mu_{N}^{\Phi}(\gamma)/2}=\alpha_N^{\Phi;+}(\gamma)
\]
if $\gamma\subset V_{-}$.

The following theorem, which can be seen as the higher dimensional version of Corollary 5.11 in
\cite{Sie11}, gives a computation of the algebraic intersection number of
$\tilde V$ and $\R$-shifts of the curve $C$ in terms of the asymptotic data, and the intersections of
each object with the asymptotic limits of the other.

\begin{thm}[Siefring, to appear]\label{t:intersection-computation}
Under the above general assumptions, also assume that the curve $C$ is connected and
not equal to an orbit cylinder and not contained in $\R\times V$.
For $c\in\R$, denote by $C_{c}$ the curve obtained from translating $C$ in the $\R$-coordinate by $c$.
Then for all but a finite number of value of $c\in\R$, the algebraic intersection number
$C_{c}\cdot \tilde V$ is given by the formulas:

\begin{equation*}
\begin{split}
C_{c}\cdot \tilde V
&=C\cdot (\R\times V_{+}) \\
&\hskip.25in
+\sum_{\stackrel{z\in\Gamma_{+}}{\gamma_{z}\in V_{+}}} 
\bp{\max\br{m_{z}\wind_{\infty}^{\Phi}(\tilde V, \gamma_{z}), \wind^{\Phi}_{rel}((C;z), \R\times V_{+})}-\wind^{\Phi}_{rel}((C;z), \R\times V_{+})} \\
&\hskip.25in
+ \sum_{z\in\Gamma_{-}}m_{z}(\R\times \gamma_{z})\cdot \tilde V  \\
&\hskip.25in
+\sum_{\stackrel{z\in\Gamma_{-}}{\gamma_{z}\in V_{-}}}
\bp{m_{z}\wind_{\infty}^{\Phi}(\tilde V, \gamma_{z})
-
\min\br{m_{z}\wind_{\infty}^{\Phi}(\tilde V, \gamma_{z}), \wind^{\Phi}_{rel}((C;z), \R\times V_{-})}} \\
&\hskip.25in
+\sum_{\stackrel{z\in\Gamma_{-}}{\gamma_{z}\in V_{+}}}
\bp{\wind_{rel}^{\Phi}((C;z), \R\times V_{+})-m_{z}\wind_{rel}^{\Phi}(\tilde V, \gamma_{z})} \\
\end{split}
\end{equation*}
\begin{equation*}
\begin{split}
&=\sum_{z\in\Gamma_{+}}m_{z}(\R\times \gamma_{z})\cdot \tilde V  \\
&\hskip.25in
+\sum_{\stackrel{z\in\Gamma_{+}}{\gamma_{z}\in V_{+}}}
\bp{\max\br{m_{z}\wind_{\infty}^{\Phi}(\tilde V, \gamma_{z}), \wind_{rel}^{\Phi}((C; z), \R\times V_{+})}-m_{z}\wind_{\infty}^{\Phi}(\tilde V, \gamma_{z})} \\
&\hskip.25in
+C\cdot (\R\times V_{-}) \\
&\hskip.25in
+\sum_{\stackrel{z\in\Gamma_{-}}{\gamma_{z}\in V_{-}}}
\bp{\wind_{rel}^{\Phi}((C;z), \R\times V_{-})-\min\br{m_{z}\wind_{\infty}^{\Phi}(\tilde V, \gamma_{z}),\wind_{rel}^{\Phi}((C;z), \R\times V_{-})}} \\
&\hskip.25in
+\sum_{\stackrel{z\in\Gamma_{+}}{\gamma_{z}\in V_{-}}}
\bp{m_{z}\wind_{\infty}^{\Phi}(\tilde V, \gamma_{z})-\wind_{rel}^{\Phi}((C;z), \R\times V_{-})}
\end{split}
\end{equation*}

with each of the grouped terms always nonnegative.
\end{thm}

The nonnegativity of the terms in the above formulas allows us to establish a
convenient set of conditions which will guarantee that the  projections 
$\pi(C)$ and $\pi(\tilde V)$ of the curve and hypersurface to $M$ do not intersect.
Indeed, if $\pi(C)$ and $\pi(\tilde V)$ are disjoint then
$C_{c}$ and $\tilde V$ are disjoint for all values of $c\in\R$ and hence
the algebraic intersection number of $C_{c}$ and $\tilde V$ is zero for all values of $c\in\R$.
Since the formulas from Theorem \ref{t:intersection-computation} compute this number (for all but a finite
number of values of $c\in\R$) in terms of nonnegative quantities, we can conclude that all terms in the above formulas vanish.
We thus obtain the following corollary, which generalizes Theorem 2.4/5.12 in \cite{Sie11}.

\begin{corollary}[Siefring, to appear]\label{c:nonintersection-conditions}
Assume that all of the hypotheses of Theorem \ref{t:intersection-computation}
hold and that $\pi(C)$ is not contained in $\pi(\tilde V)$.
Then the following are equivalent:

\begin{enumerate}
\item $\pi(C)$ and $\pi(\tilde V)$ are disjoint.

\item  All of the following hold:
\begin{enumerate}
\item
None of the asymptotic limits of $C$ intersect $\pi(\tilde V)$.

\item
$\pi(C)$ does not intersect $V=V_{+}\cup V_{-}$.

\item
For any puncture $z$ at which $C$ has asymptotic limit $\gamma^{m}$ lying in $V=V_{+}\cup V_{-}$,
$\wind^{\Phi}_{rel}((C; z), V)=m \wind_{\infty}^{\Phi}(\tilde V, \gamma)$.
\end{enumerate}
\end{enumerate}
\end{corollary}

\section{Pseudo-holomorphic dynamics}

\begin{figure}
    \centering
    \includegraphics[width=0.65\linewidth]{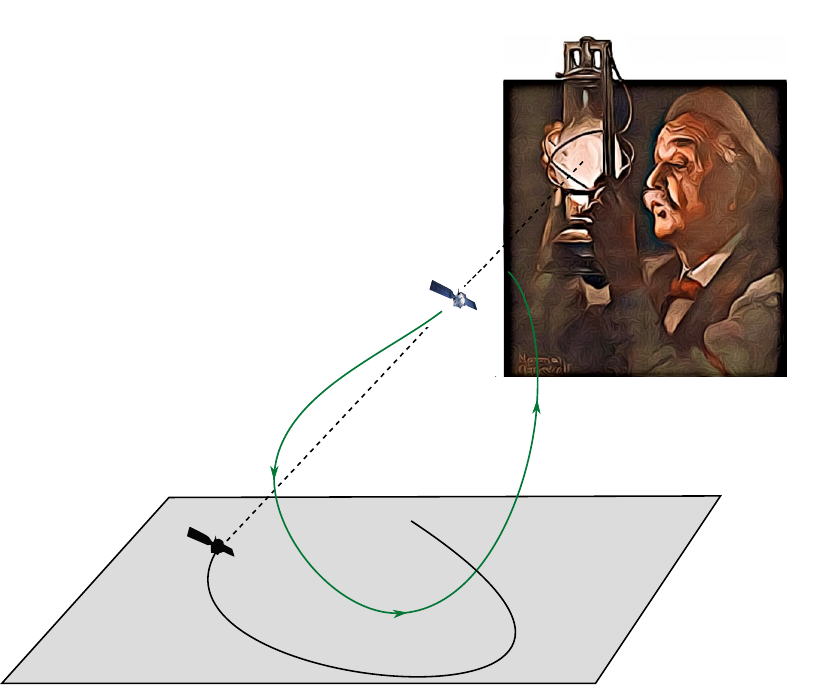}
    \caption{\textbf{Philosophy:} To shed some light on a complicated higher-dimensional problem, try first to look at the shadow that your lantern is producing!}
    \label{fig:lantern}
\end{figure}

This section is based on \cite{M20}. We will associate to the (low-energy, near primary) dynamics of the spatial CR3BP, a Reeb flow on $S^3$, which is in some sense a ``shadow'' of the original dynamics, by using a finite energy foliation on the symplectization of the $5$-dimensional level sets. The idea is to look at a lower-dimensional avatar of the original $5$-dimensional dynamics, in the hope to learn something about the latter from knowledge on the former. See Figure \ref{fig:lantern}. 

The aim to foliate the whole $5$-dimensional space, with a foliation whose leaf space is $S^3$, and define a geometric structure and a dynamical system on this leaf space. The construction of this foliation uses Siefring's intersection theory in a crucial way, as one first finds a foliation by hypersurfaces, and the intersection theory implies that the holomorphic curves in the foliation all lie in the leaves of the holomorphic hypersurface foliation, where the four dimensional techniques apply. This construction will \emph{concretely} realize the iterated picture explained in Section \ref{sec:iterated_picture}, but now with an adapted \emph{IP foliation}, i.e.\ the pages of the open book in the spatial CR3BP of Theorem \ref{thm:IP} are now foliated by Lefschetz fibrations which are compatible with the given dynamics, in the sense that the fibers are symplectic for the symplectic form preserved by the return map. We make this precise in the following definition.

\begin{definition}[\textbf{IP foliation}] Let $(M,\xi)=\mathbf{OB}(P,\varphi)$ be an IP contact $5$-fold, with $P=\mathbf{LF}(F,\varphi_F)$ and $L=\partial F$, and let $\alpha$ be a Giroux form, adapted to a concrete open book $\theta: M\backslash B\rightarrow S^1$ of this abstract type. Denote $P_\varphi=\overline{\theta^{-1}}(\varphi)$. An adapted \emph{IP foliation} consists of the following data.
\begin{itemize}
    \item A concrete open book $\theta_B: B\backslash L\rightarrow S^1$ adapted to $\alpha_B=\alpha\vert_B$.
    \item A concrete Lefschetz fibration $f_\varphi: P_\varphi \rightarrow \mathbb D^2$, which is symplectic with respect to $\omega_\varphi=d\alpha\vert_{P_\varphi}$, and which induces the open book $\pi_B$ at the boundary, for every $\varphi\in S^1$.
\end{itemize}
    \end{definition}

In particular, the foliation is adapted to the given dynamics in the sense that every leaf is symplectic with respect to $d\alpha$, where $\alpha$ is the contact form giving the ambient dynamics. This is a $5$-dimensional version of Hofer--Wysocki--Zehnder's notion of a finite energy foliation. The following is proven in \cite{M20}. 

\begin{thm}[\cite{M20}, \textbf{IP foliations for the spatial CR3BP}]\label{thm:IPfoliation} For $(c,\mu)$ in the convexity range, and near the primaries, the Moser-regularized energy level set $\overline{\Sigma}_c$ admits an adapted IP foliation.
\end{thm}

Here, we need to recall that the page $\mathbb{D}^*S^2=\mathbf{LF}(\mathbb{D}^*S^1,\tau_P^2)$ of the open book of Theorem \ref{thm:openbooks} has a Lefschetz fibration with genus zero fibers over the $2$-disk, with monodromy the Dehn twist $\tau_P$ ($P$ here is for ``planar'', to differentiate from the monodromy $\tau$ used for the spatial case; recall Figure \ref{fig:LFDS2}). This gives an iterated planar structure on the $5$-dimensional contact manifold $(S^*S^3,\xi_{std})$, and therefore on the Moser-regularized low-energy level sets in the CR3BP (see Theorem \ref{thm:IP}). The content of Theorem \ref{thm:IPfoliation} is then to realize the abstract data with concrete data, adapted to the given dynamics.

Note that for the integrable case, i.e.\ the rotating Kepler problem, we have already constructed such an IP foliation in Section \ref{sec:RKP_returnmap}, by hand. Such a foliation has the property that it is perserved under the dynamics. However, in the non-integrable cases, that will certainly \emph{not} be the case in general. In the general case, the leaves of the IP foliation will be projections to the $5$-dimensional level set of asymptotically cylindrical holomorphic curves in the symplectization, all of which are annuli and asymptotic to the direct/retrograde planar orbits.

\subsection{Construction} The strategy for the construction of this foliation is sketched as follows. First, we view the planar CR3BP as a strong stable hypersurface $B$ of the spatial CR3BP, in the sense of the previous section, i.e.\ as the binding of the open book from Theorem \ref{thm:openbooks}. The symplectization $\mathbb R\times B$ of the planar CR3BP is then a holomorphic hypersurface in the symplectization $\mathbb R\times M$ of the spatial CR3BP, if we choose a $B$-compatible almost complex structure; we can moreover choose the restriction of this almost complex structure to $\mathbb R\times B$ as we please. One can then find a foliation $\mathcal F$ of the $6$-dimensional symplectization $\mathbb R\times M$ by holomorphic hypersurfaces which are asymptotically cylindrical to $\mathbb R\times B$, and which project to $M$ as the pages of the open book from Theorem \ref{thm:openbooks}; as in Section \ref{sec:HWZ}, this is called a \emph{holomorphic} open book decomposition (cf.\ \cite{Wen10d}). If we moreover assume that we have a holomorphic open book on $\mathbb R\times B$ with genus zero pages, for some almost complex structure, we then choose the ambient one to coincide with the given one along the $\mathbb R\times B$. Now, the moduli space in $\mathbb R\times B$ needs to extend to a moduli space $\mathcal M$ in $\mathbb R\times M$, which is foliated by holomorphic hypersurfaces. The claim is then that the holomorphic curves in $\mathcal M$ need to lie in the leaves of the holomorphic hypersurface foliation $\mathcal{F}$. This is proved via Siefring interesection theory. 

Indeed, on the one hand, if a curve $u$ in $\mathcal M$ does not lie in a leaf of $\mathcal{F}$, the Siefring intersection pairing between $u$ and some hypersurface $H$ in the foliation would be strictly positive, by the fact that interior intersections contribute positively. But on the other hand, one could compute that this intersection pairing is zero, which gives a contradiction.\footnote{We can also give an alternative argument for the setting of the CR3BP: one can construct an explicit foliation (by hand) for the integrable limit case (RKP), see Appendix A in \cite{MvK20a}. A continuation argument as in Theorem 3.9 in \cite{M18}, appealing to the implicit function theorem, yields the claim.} Indeed, using that the hypersurfaces have no negative asymptotics, by homotoping the hypersurfaces towards the negative ends and using homotopy invariance (or alternatively appealing to Theorem \ref{t:intersection-computation}, noting that the winding numbers vanish), we have
$$
u*H=u*(\mathbb{R}\times B).
$$
From homotopy invariance with respect to $u$, it suffices to compute the above in the case where $u\in \mathcal{M}_B$ lies completely in $\mathbb{R}\times B$. Then we can appeal to the following intersection formula:
$$
u*(\mathbb{R}\times B)=\frac{1}{2}(\mu_N^\tau(u)-\#\Gamma(u)_{odd}),
$$
where $\mu_N^\tau(u)$ is the total normal Conley-Zehnder index of $u$ with respect to a trivialization $\tau$ of the symplectic normal bundle to $\mathbb{R}\times B$, and $\#\Gamma(u)_{odd}$ is the number of asymptotics of $u$ which have odd normal Conley-Zehnder index. 

In order to obtain that the curves in $\mathcal{M}$ yield a foliation of the pages by Lefschetz fibrations, we appeal to Wendl's theorem for fillings (Theorem \ref{thm:Wendl}). We apply this theorem (in a parametric way) to the hypersurfaces in $\mathcal{F}$, all of which can be seen as fillings of the planar problem $B$, with its genus zero open book. This finishes the construction of the foliation.\footnote{Note also that, in the context of the CR3BP, this foliation exists in the integrable limit case $\mu=0$, and is constructed by hand, see \cite{MvK20a}. It may also be constructed as in the current way. Uniqueness of gluing, i.e.\ the implicit function theorem, gives the desired foliation, agreeing with the one constructed by hand for the limit case.}

In order to conclude the existence of such a foliation in the spatial CR3BP, the required holomorphic open book on $\mathbb R \times B$ is provided by \cite{HSW}, but only for the convexity range. This finishes the sketch of the construction.

\subsection{Dynamics on moduli spaces} 

\begin{figure}
        \centering
        \includegraphics[width=0.5 \linewidth]{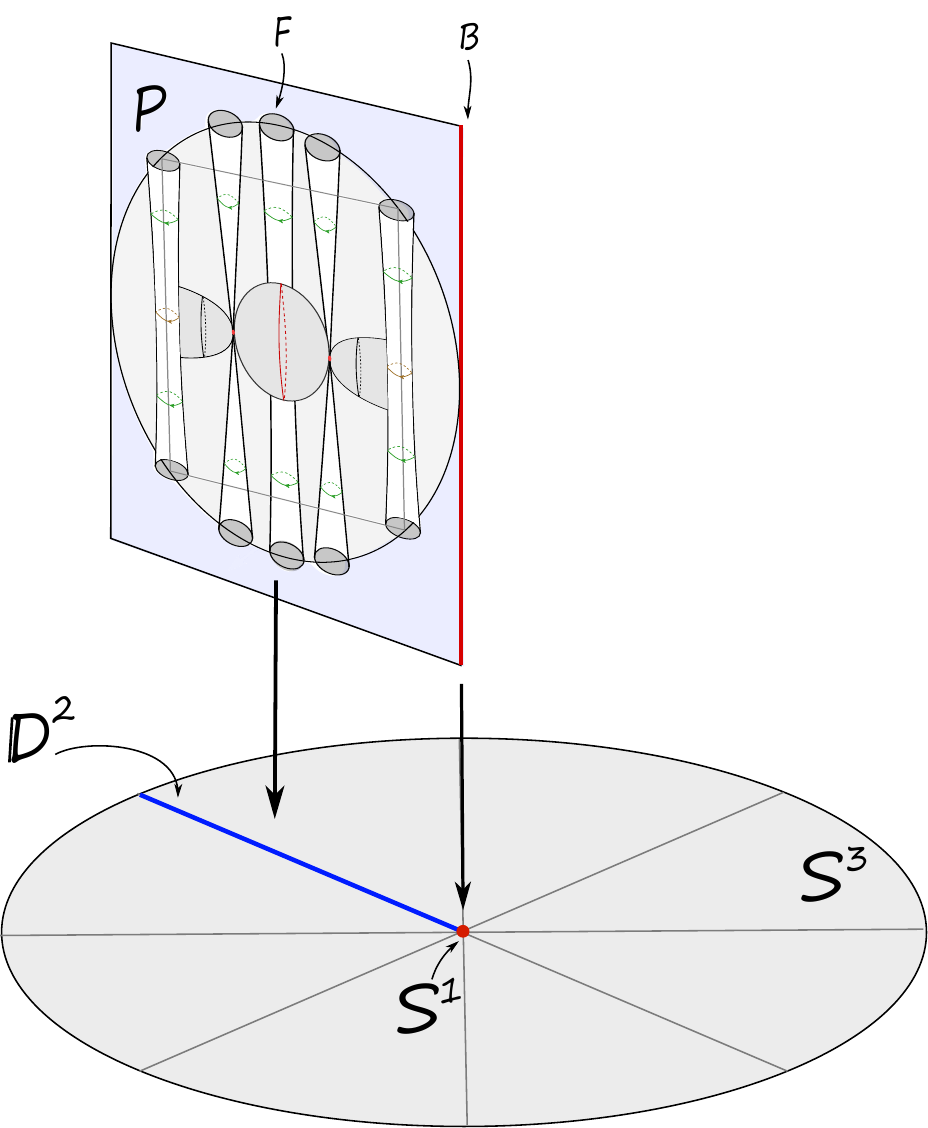}
        \caption{The moduli space of curves is a copy of $S^3=\mathbf{OB}(\mathbb{D}^2,\mathds{1})$.}
        \label{fig:MODULI}
    \end{figure}

An IP foliation can be used to induce a dynamical system on $S^3$. The main topological observation is the following. The leaf space $\mathcal{M}$ of the leaves of an IP foliation (i.e.\ the moduli space parametrizing them) is a copy of $S^3$. Indeed, each page $P$ of the open book is a $2$-disk worth of fibers; we moreover have an $S^1$-family of such pages, all of them sharing the boundary $B$ (the binding), and such that their Lefschetz fibration all induce the $S^1$-family of pages of the open book $B=\mathbf{OB}(F,\varphi_F)$. It follows that the leaf space carries the trivial open book $\mathcal{M}=\mathbf{OB}(\mathbb{D}^2,\mathds{1})\cong S^3$, whose disk-like page corresponds to the base of the Lefschetz fibration in the page $P$, and whose binding $\mathcal{M}_B$ is the $S^1$-family of pages for $B$. See Figure \ref{fig:MODULI}. 

But the story does not end here, as we can also endow the leaf space itself with a geometric structure, and a dynamical system.

\begin{thm}[\cite{M20}, \textbf{contact and symplectic structures on moduli}]\label{thm:contactstructure} Given an IP foliation on $(M,\xi)$, the leaf space $\mathcal{M}$ carries a natural contact structure $\xi_\mathcal{M}$ which is supported by the trivial open book on $S^3$. Moreover, the symplectization form on $\mathbb{R}\times M$ associated to any Giroux form $\alpha_M$ on $M$ induces a tautological symplectic form on $\mathbb{R}\times \mathcal{M}$, which is naturally the symplectization of a contact form $\alpha_\mathcal{M}$ for $\xi_{\mathcal{M}}$, whose Reeb flow is adapted to the trivial open book on $\mathcal{M}$.
\end{thm}

\begin{figure}
    \centering
    \includegraphics[width=0.5 \linewidth]{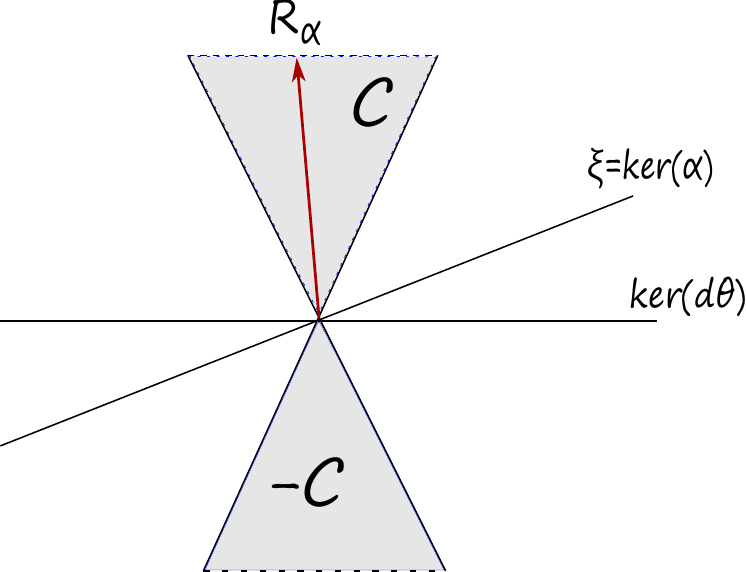}
    \caption{A cone structure adapted to an open book.}
    \label{fig:adapted}
\end{figure}

\begin{figure}
    \centering
    \includegraphics[width=0.8\linewidth]{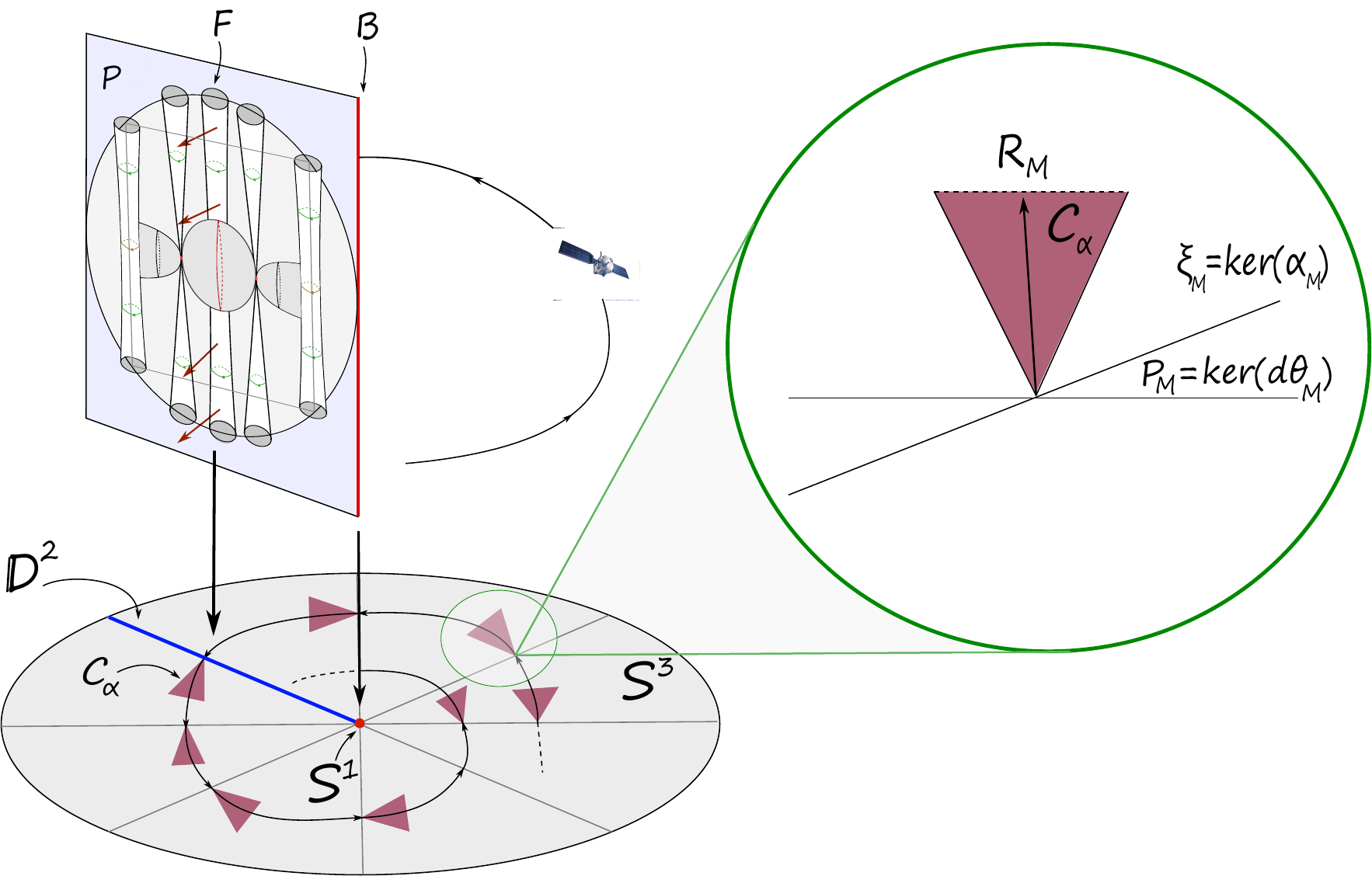}
    \caption{The shadowing cone.}
    \label{fig:shadowing cone}
\end{figure}

As the standard contact structure on $S^3$ is the unique contact structure supported by the trivial open book, we see that $\xi_\mathcal{M}$ is isotopic to the standard contact structure. The contact form can be written down via the following tautological formula, via leaf-wise integration:
$$
(\alpha_\mathcal{M})_u(v)=\int_{z\in u}\alpha_z(v(z))dz,
$$
where $u\in \mathcal{M}$, $v\in T_u\mathcal{M}=\ker \mathbf D_u$ for $\mathbf{D}_u$ the linearized CR-operator of $u$, and $dz=d\alpha\vert_u$ is an area form along $u$. The contact structure $\xi_\mathcal{M}=\ker \alpha_\mathcal{M}$ and the $1$-dimensional distribution $\ker d\alpha_\mathcal{M}$ can then be thought of as the average of the contact planes $\xi_z$, and respectively of $\ker d\alpha_z$, for $z \in u$, i.e.\
$$
\xi_\mathcal{M}=\int_{z\in u}\pi_*(\xi_z)dz,
$$
$$
\ker d\alpha_\mathcal{M}=\int_{z\in u} \pi_*(\ker d\alpha_z)dz,
$$
where $\pi: M\backslash L\rightarrow S^3$ is the quotient map to the leaf space. This means that the Reeb vector field $R_\mathcal{M}$ of $\alpha_\mathcal{M}$ spans the average direction in the ``shadowing cone'' $C_\alpha=\pi_*(\ker d\alpha)\subset TS^3$; see Figure \ref{fig:shadowing cone}.

One can therefore think of the holomorphic shadow as the ``guiding direction'' of the cone. We encode the properties of this cone in the following general definition. In what follows, a cone structure on a manifold $M$ is a smooth choice $C_x\subset T_xM$ of a subset which is closed under multiplication by a positive scalar at each fiber. The cone is everywhere non-trivial if $C_x\neq 0$ for every $x$.

\begin{definition}\label{def:stronglyadapted}
Consider an everywhere non-trivial cone structure $C$ on a manifold $M$, where $M$ is endowed with an open book $\theta: M\backslash B\rightarrow S^1$. We say that $C$ is strongly adapted (or simply adapted) to $\theta$ if
\begin{enumerate}
    \item[(1)] $C\vert_B\subset TB$;
    \item[(2)] $d\theta$ is a section of $C\vert_{M\backslash B}$; 
\end{enumerate}
and if there exists a Giroux form $\alpha$ for the open book such that
\begin{enumerate}
    \item [(3)] $\alpha$ is a section for $C$;
    \item [(4)] The Reeb vector field $R_\alpha$ is interior to $C$.
\end{enumerate}
Here, a section for $C$ is a $1$-form which is strictly positive on non-zero vectors of $C$. See Figure \ref{fig:adapted}. Moreover, a vector is interior to a conical set if the smallest disk which covers some base\footnote{A \emph{base} for a cone inside a tangent space is a section of the cone, i.e.\ a continuous choice of vector in each direction tangent to the cone; a base for a cone structure is a continously varying choice of bases for each tangent space, see e.g.\ \cite{S76}.} of the cone contains the vector. We note that the cone structure arising in the CR3BP is 2-dimensional (i.e.\ the image under a map with a 2-dimensional source), and therefore is not a convex cone, i.e.\ does not contain its interior.

\end{definition}

\smallskip

\textbf{The holomorphic shadow.} We define the \emph{holomorphic shadow} of the Reeb dynamics of $\alpha_M$ on $M$ to be the Reeb dynamics of the associated contact form $\alpha_\mathcal{M}$ on $S^3$, provided by Theorem \ref{thm:contactstructure}. The flow of $\alpha_{\mathcal{M}}$ can be viewed as a flow $\phi_t^{M;\mathcal{M}}$ on $M\backslash L$ which leaves the holomorphic foliation $\mathcal{M}$ invariant (i.e.\ it maps holomorphic curves to holomorphic curves). It is the ``best approximation'' of the Reeb flow of $\alpha_M$ with this property, as its generating vector field is obtained by reparametrizing the projection of the original Reeb vector field to the tangent space of $\mathcal{M}$, via a suitable $L^2$-orthogonal projection. Concretely, we have
$$
R_\mathcal{M}(u)=\frac{P_u(R_\alpha\vert_u)}{(\alpha_\mathcal{M})_u(P_u(R_\alpha\vert_u))} \in T_u\mathcal{M},
$$
where $P_u: W^{1,2}(N_u)\rightarrow \ker \mathbf{D}_u$ denotes the $L^2$-orthogonal projection with respect to the metric $$
g_u(v,w)=\int_{z\in u} g_z(v(z),w(z))dz,
$$
with $g_z=d\alpha_z(\cdot,J \cdot)+\alpha_z\otimes \alpha_z+dt \otimes dt$, and $v,w \in W^{1,2}(N_u)$ sections of the normal bundle $N_u$ to $u$. It may also be viewed as a Reeb flow $\phi_t^{S^3;\mathcal{M}}$ on $S^3$, related to the one on $M$ via a semi-conjugation
\begin{center}\label{diag:semiconj}
    \begin{tikzcd}
        M\backslash L\arrow[r, "\phi_t^{M;\mathcal{M}}"] \arrow[d, "\pi"]&
        M\backslash L \arrow[d, "\pi"]\\
        S^3 \arrow[r, "\phi_t^{S^3;\mathcal{M}}"]& S^3 \\
    \end{tikzcd}
\end{center}
where $\pi$ is the projection to the leaf-space $\mathcal{M}\cong S^3$. We will now focus on the global properties of the correspondence $\alpha_M \mapsto \alpha_\mathcal{M}$.

For $F$ a genus zero surface, let $\mathbf{Reeb}(F,\phi_F)$ denote the collection of contact forms whose flow is adapted to some concrete planar open book $\pi_B: B\backslash L \rightarrow S^1$ on a given $3$-manifold $B$, of abstract form $B=\mathbf{OB}(F,\phi_F)$. Iteratively, we define $\mathbf{Reeb}(\mathbf{LF}(F,\phi_F),\phi)$ to be the collection of contact forms with flow adapted to some concrete IP open book $\pi_M: M\backslash B \rightarrow S^1$ on a $5$-manifold $M$, of abstract form $M=\mathbf{OB}(\mathbf{LF}(F,\phi_F),\phi)$, whose restriction to the binding $B=\mathbf{OB}(F,\phi_F)$ belongs to $\mathbf{Reeb}(F,\phi_F)$. We call elements in $\mathbf{Reeb}(\mathbf{LF}(F,\phi_F),\phi)$ \emph{IP contact forms}, or \emph{IP Giroux forms}.

We then have a map
$$
\mathbf{HS}:\mathbf{Reeb}(\mathbf{LF}(F,\phi_F),\phi)\rightarrow \mathbf{Reeb}(\mathbb{D}^2,\mathds{1}),
$$
given by taking the holomorphic shadow with respect to an auxiliary almost complex structure $J$ associated to $\alpha_M$. We refer to $\mathbf{HS}^{-1}(\alpha_{std})$ as the \emph{integrable fiber}, where $\alpha_{std}$ denotes the standard contact form in $S^3$.

\begin{thm}[\cite{M20}, Reeb flow lifting theorem]\label{thm:lifting} $\mathbf{HS}$ is surjective.
\end{thm}

In other words, fixing an auxiliary some $J$, we may lift any Reeb flow on $S^3$ adapted to the trivial open book, as the holomorphic shadow of the Reeb flow of an IP Giroux form adapted to \emph{any} choice of concrete IP contact $5$-fold. The map $\mathbf{HS}$ is clearly not in general injective, as it forgets dynamical information in the fibers. The above theorem says that Reeb dynamics on an IP contact $5$-fold is at least as complex as Reeb dynamics on the standard contact $3$-sphere. Recalling that the Levi-Civita regularization of the planar CR3BP (for subcritical energy) gives a Reeb flow on $S^3$, this gives a rough ``measure'' of the complexity of the spatial CR3BP. Namely, the spatial problem belongs to a space of dynamics which is at least as complex as the space of dynamics in which the planar problem lives. The (informal) motto is the following. 

\medskip

\textbf{Motto.} ``\emph{The \textbf{spatial} C3RBP is dynamically at least as complex as the \textbf{planar} CR3BP}''.

\medskip

The above is of course not a formal statement, as the two systems may not be related by the shadow map, but it is at least suggestive. 

\medskip

\textbf{Dynamical Applications.} We wish to apply the above results to the spatial CR3BP. We first introduce the following general notion. Consider an IP $5$-fold $M$ with an IP Reeb dynamics, endowed with an IP foliation $\mathcal{M}$. Fix a page $P$ in the IP open book of $M$, and consider the associated Poincar\'e return map $f: \mbox{int}(P) \rightarrow \mbox{int}(P)$. A (spatial) point $x\in \mbox{int}(P)$ is said to be \emph{leaf-wise} (or \emph{fiber-wise}) \emph{$k$-recurrent} with respect to $\mathcal{M}$ if $f^k(x) \in \mathcal{M}_x$, where $\mathcal{M}_x$ is the leaf of $\mathcal{M}$ containing $x$, and $k\geq 1$. This means that $f^k(\mbox{int}(\mathcal{M}_x))\cap \mbox{int}(\mathcal{M}_x)\neq \emptyset$. This is, roughly speaking, a symplectic version of the notion of \emph{leaf-wise intersection} introduced by Moser \cite{M78} for the case of the isotropic foliation of a coisotropic submanifold, as the leaves are symplectic rather than isotropic. 

\begin{figure}
    \centering
    \includegraphics[width=1 \linewidth] {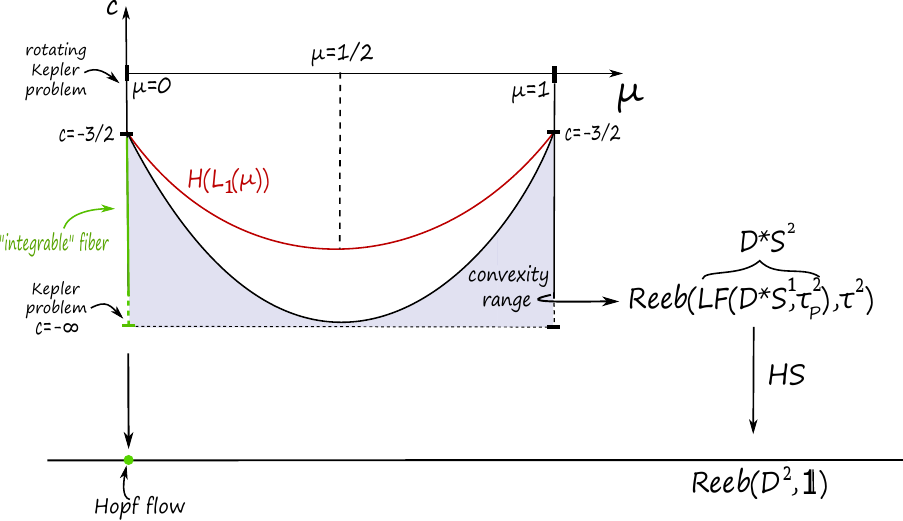}
    \caption{An abstract sketch of the convexity range in the SCR3BP (shaded), for which the holomorphic shadow is well-defined. We should disclaim that the above is not a plot; the convexity range is not yet fully understood, although it contains (perhaps strictly) a region which qualitatively looks like the above, cf.\ \cite{AFFHvK,AFFvK}.}
    \label{fig:my_label}
\end{figure}

In the integrable case of the RKP, the holomorphic foliation provided by Theorem \ref{thm:IPfoliation} can be obtained directly as in Section \ref{sec:RKP_returnmap}. Denote this ``integrable'' holomorphic foliation on $S^*S^3$ by $\mathcal{M}_{int}$. Since the return map for $\mu=0$ preserves fibers as explained in Theorem \ref{thm:integrablecase}, every point is leaf-wise $1$-recurrent with respect to $\mathcal{M}_{int}$. If the mass ratio is sufficiently small, then the leaves of $\mathcal{M}_{int}$ will still be symplectic with respect to $d\alpha$, where $\alpha$ is the corresponding perturbed contact form on the unit cotangent bundle $S^*S^3$. 

We have the following perturbative result.

\begin{thm}[\cite{M20}]\label{thm:application}
In the spatial CR3BP, for any choice of page $P$ in the open book of Theorem \ref{thm:openbooks}, for any fixed choice of $k\geq 1$, for sufficiently small $\mu$ (depending on $k$), for energy $c$ below the first critical value $H(L_1(\mu))$, along the bounded components of the Hill region, and for every $l\leq k$, there exist infinitely many points in $\mbox{int}(P)$ which are leaf-wise $l$-recurrent with respect to $\mathcal{M}_{int}$.
\end{thm}

In simpler words, the spatial CR3BP admits an abundance of leaf-wise recurrent points, at least in the perturbative regime.

\begin{remark}
The same conclusion holds for arbitrary $\mu \in [0,1]$, but sufficiently negative $c\ll 0$ (depending on $\mu$ and $k$).
\end{remark}

In fact, the conclusion of the Theorem \ref{thm:application} holds whenever the relevant return map is sufficiently close to a return map which preserves the leaves of the holomorphic foliation of Theorem \ref{thm:IPfoliation} (i.e.\ which coincides with its holomorphic shadow on $M$). We now give the main ideas for this fact.

\subsection{Proof of Theorem \ref{thm:application}} The idea for the proof of Theorem \ref{thm:application} is the following. Consider an IP foliation on an IP contact $5$-fold $(M,\xi)=\mathbf{OB}(P,\varphi)$, with binding $B=\mathbf{OB}(F,\varphi_F)$, $L=\partial F$, and adapted IP Giroux form $\alpha$. 

\begin{figure}
    \centering
    \includegraphics[]{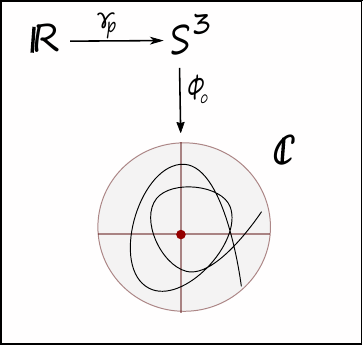}
    \caption{The qualitative image of a transverse path $\gamma_p$ under the defining map $\Phi_0:S^3\rightarrow \mathbb{C}$ for the trivial open book $\pi_0=\Phi_0/\vert \Phi_0\vert$.}
    \label{fig:transversepath}
\end{figure}

We will keep track of which leaves are intersected by each Reeb orbit in $M$, without changing the original dynamics. Namely, for $p \in M\backslash L$, we consider the path
$$
\gamma_p(t)=\pi(\phi_t^M(p))\in \mathcal{M},
$$
where $\phi_t^M$ is the flow of $R_\alpha$, and $\pi: M\backslash L\rightarrow \mathcal{M}$ is the projection to the leaf space. If $p \in B\backslash L$, this is a parametrization of $\mathcal{M}_B \cong S^1$, the binding of the open book in $\mathcal{M}$; if $p \in M\backslash B$, this is a path in $S^3$ which is tangent to the shadowing cone $C_\alpha$, an so is positively transverse to each disk-like page $P_\varphi^\mathcal{M}$ as well as to the contact structure $\xi_\mathcal{M}$. Note that different choices of $p$ might induce paths which intersect each other (corresponding to their orbits intersecting the same leaf), and even self-intersect (corresponding to an orbit intersecting the same leaf multiple times), so these paths are not orbits of an autonomous flow. See Figure \ref{fig:transversepath}. We will refer to the collection
$$
\mathbf{TS}(\alpha,J)=\{\gamma_p: p \in M\backslash L\}
$$
as the \emph{transverse shadow} of the Reeb flow of $\alpha$ on $M$, with respect to $J$, which is by definition the collection of those orbits of $C_\alpha$ coming from orbits of $\alpha$ on $M$.

One may further choose to ``package'' these transverse paths in different ways, e.g.\ by considering those paths induced by points on a section of the Lefschetz fibration at a given page, as follows. 

Fix the $0$-page $P_0$ in $M$, with its concrete Lefschetz fibration $f_0: P_0\rightarrow \mathbb D_0^2$, and the associated Poincar\'e return map $f:\text{int}\;P_0 \rightarrow \text{int}\;P_0$, which we assume for simplicity extends to the boundary (as in the CR3BP). Consider a two disk $D \cong \mathbb{D}^2$ satisfying:
\begin{itemize}
    \item $D\subset P_{0}$;
    \item $\partial D \subset \partial P_{0}=B$ is a loop which is disjoint from the binding $L$ of the concrete open book in $B$, and transverse to the interior of each of its pages and to the contact structure $\xi_B$;
    \item $D_0$ is a symplectic section of the Lefschetz fibration $f_0$, i.e.\ $D$ intersects each fiber of $f_0$ precisely once, and hence $D=im(s)$ for $s: \mathbb{D}_0^2\rightarrow P_0$ satisfying $f_0 \circ s=id$.
\end{itemize}

We refer to such a disk $D$ as a \emph{(horizontal) symplectic tomography} for the Reeb dynamics on $M$. Note that, if $\partial D$ is a Reeb orbit of $\alpha_B$ which is linked with $L$, then $f(\partial D)=\partial D$ is invariant under the return map. 

\begin{figure}
    \centering
    \includegraphics[width=0.6 \linewidth]{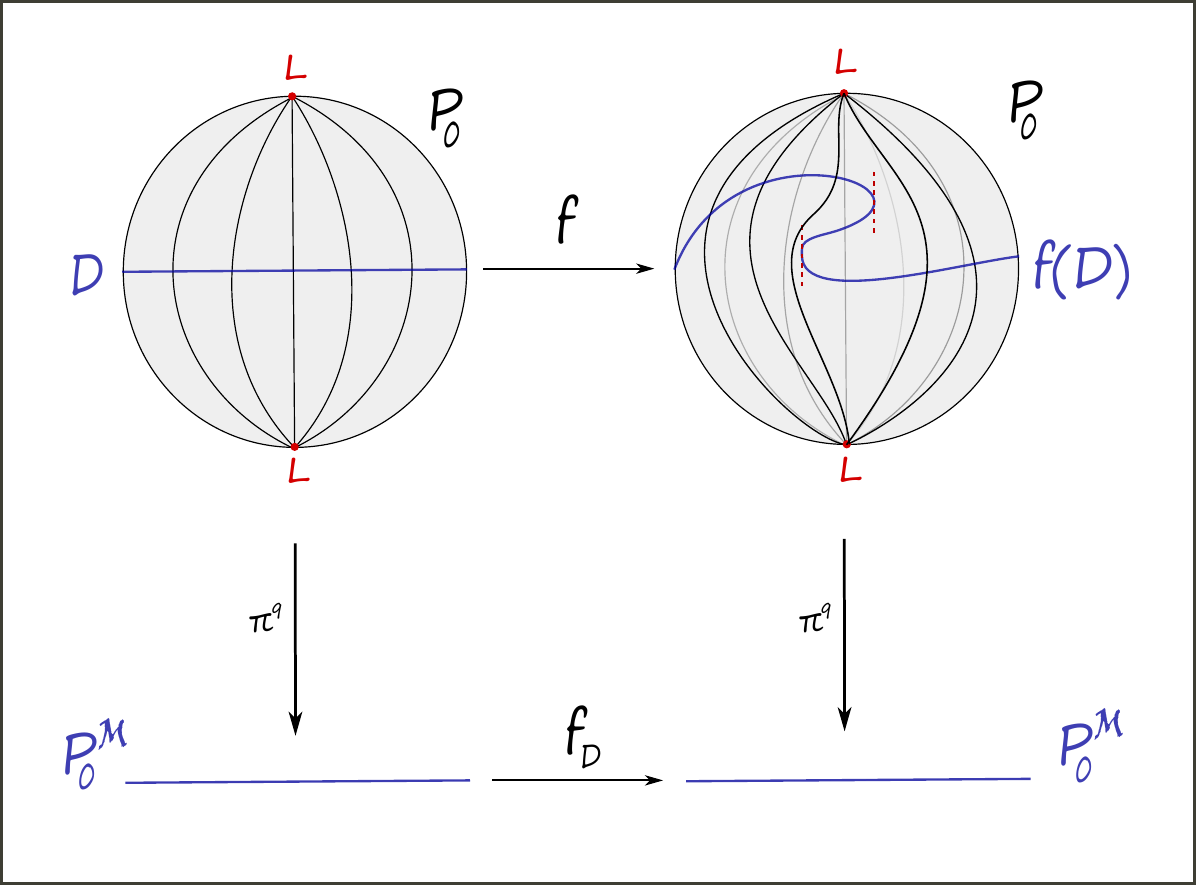}
    \caption{The return map $f_D$ associated to the tomography $D$. Open tangencies of $f(D)$ with the vertical foliation (which could, a priori, theoretically arise from ``foldings'' of the disk $f(D)$) might prevent in general that $f_D$ preserves area. This does \emph{not} happen perturbatively, i.e.\ when we perturb a foliation-preserving map, however. Note that $f(L)=L$.}
    \label{fig:verticaltang}
\end{figure}

For each such symplectic tomography $D$, we have an associated return map 
$$f_D: P_0^\mathcal{M}\rightarrow P_0^\mathcal{M}$$ on the $0$-page of the moduli space, as follows. We identify $P_0^\mathcal{M}$ with $\mathbb{D}_0^2$, and define $f_D$ by
$$
f_D(u)=\gamma_{s(u)}(\tau(u,D)) \in P_0^\mathcal{M},
$$
where $\tau(u,D)=\min\{t>0: \gamma_{s(u)}(t)\in P_0^\mathcal{M}\}$ is the first return time of the transverse path $\gamma_{s(u)}$ to the $0$-page $P_0^\mathcal{M}$. See Figure \ref{fig:verticaltang}.

The symplectic disk $(D,d\alpha\vert_D)$ is symplectomorphic to $(P_0^\mathcal{M},d\alpha_{\mathcal{M}}\vert_{P_0^\mathcal{M}})$, and both have finite symplectic area. In general, $f_D$ might a priori decrease area. Indeed, $f(D)$ is a symplectic disk in $P_0$ with the same symplectic area as $D$, but it might have an open set of vertical tangencies, i.e.\ intersecting a fiber along an open subset of positive area (as opposed to $D$, which intersects fibers at a single point). Nevertheless, this is not possible for perturbative situations where one perturbs a fiber-wise preserving map, in which case the perturbed $f_D$ still preserves area. 

On the other hand, if $f_D$ is easily seen to be surjective, i.e.\ every holomorphic fiber of $f_0$ is intersected by the symplectic disk $f(D)$. Indeed, since $f$ is homotopic to the identity by a smooth homotopy which preserves the boundary, $f(D)$ is homologous to $D$ relative boundary. Moreover, $\partial f(D)$ is a push-off in the Reeb direction of $\partial D$ (they agree if $\partial D$ is an orbit, as observed above), hence they can be homotoped to each other via the Reeb flow of $R_B$, and in particular away from $L$ (the boundary of the holomorphic fibers). It follows that the homological intersection number of $f(D)$ with the fibers agrees with that of $D$, i.e.\ it is $1$. In general, $f_D$ is not necessarily injective. However, this is certainly true in the case where $f$ is close to a fiber-wise preserving map, since otherwise $f(D)$ would have vertical tangencies.  

if $f$ is sufficiently close to a fiber-wise preserving map, then $f_D$ is an area-preserving homeomorphism of the $2$-disk for every tomography $D$. By Brouwer's translation theorem, we find an interior fixed point for $f_D$; by construction this corresponds to an (interior) fiber-wise $1$-recurrent point in the fixed page $P_0$. Varying vertically the tomography $D$ along $P_0$, we obtain infinitely many such points. If $k\geq 1$, fiber-wise $k$-recurrent points correspond to interior fixed points of the return map $$f_{k,D}(u)=\gamma_{s(u)}(\tau_k(u,D)),$$ where $\tau_k(u,D)$ is the $k$-th return time of the transverse path $\gamma_{s(u)}$ to $P_0^\mathcal{M}$. Note that this map is in general \emph{different} from $f_D^k$; recall that $f_D$ is not the return map of an autonomous flow. Having fixed $k$, using that $D$ and $P_0$ vary in compact families, we can take a sufficiently small perturbation so that $f_{l,D}$ is still an area-preserving homeomorphism for every $l\leq k$, $D$ and any choice of page, and apply the same argument. This finishes the proof of Theorem \ref{thm:application}.

\section{Feral curves}\label{sec:feral} This section provides a novel application of the theory of (feral) pseudoholomorphic curves to the \emph{planar} CR3BP. As far as the author knows, there is no analogous result in the literature pertaining to the classical problem which resembles our result (Theorem \ref{thm:goal}). It is very much non-perturbative, and the proof is based on the methods of feral curves recently introduced by Fish--Hofer. One should disclaim that the heavylifting behind the result was all carried out by the author's colleague Rohil Prasad, while the author's contribution consisted simply in making the connection to the CR3BP (and understanding the main ideas in Prasad's paper \cite{Pr24}).

As we have discussed previously, a fundamental theorem of Hofer (Proposition \ref{Hofer_result}) yields existence of periodic orbits for contact-type level sets of Hamiltonian systems. These orbits arise as limit sets of punctured pseudoholomorphic curves in symplectic cobordisms, under the assumption of having finite Hofer energy. If the assumption of having the contact condition on the level set is dropped, and so the dynamics is no longer given by a Reeb flow, or the finiteness of the energy is not assumed, then the standard machinery does not immediately apply, and one needs new technology to address the problem of existence of dynamical objects (orbits, or more generally, invariant subsets). In this direction, Fish--Hofer \cite{FH23} introduced the notion of a \emph{feral curve}, which is a generalization of a (finite energy) pseudoholomorphic curve, which allows for wilder asymptotic behavior (and hence more interesting limit sets than simply periodic orbits), while still keeping some control on the behavior of the curves in question. This allowed them to prove a conjecture of Herman from his 1998 ICM address, namely:

\begin{thm}[\cite{FH23}]\label{herman}
Let $H:\mathbb R^4\rightarrow \mathbb R$ be a Hamiltonian defined on $\mathbb R^4$, with its standard symplectic structure. Assume that $M:=H^{-1}(0)$ is a non-empty compact and regular level set. Then the Hamiltonian flow of $H$ on $M$ is not minimal, i.e.\ there exists a closed nonempty proper invariant subset.
\end{thm}

The sketch proof is fairly straightforward to give, as it uses ideas that are by now standard. The real work lies in introducing the notion of feral curve, and proving all the necessary technical results (e.g.\ compactness, area bounds, etc.) needed to carry out the proof. Indeed, the sketch goes roughly as follows: Given the level set $M\subset \mathbb R^4$, one symplectically embeds a neighbourhood of $M$ inside $\mathbb CP^2$, stretches the neck along this hypersurface, use Gromov's existence result for degree one curves in $\mathbb CP^2$, shows that these curves stretch towards the negative end of the resulting symplectic cap, and establishes a compactness result yielding a non-compact feral pseudoholomorphic curve whose limit set is the desired closed invariant subset. The assumption on the space being $4$-dimensional as opposed to higher-dimensional is used when proving that the subset is proper, and basically boils down to positivity of intersections.

For completeness, we give a precise definition of a feral curve in what follows, and leave understanding the details of the proof of the above result to the (very) interested reader. First, a \emph{framed Hamiltonian manifold} is a $(2n+1)$-dimensional manifold $M$ together with a $1$-form $\lambda$ and a $2$-form $\omega$ satisfying $d\omega=0$ and $\lambda\wedge \omega^n$ is a volume form. Note that this generalizes the definition of a (strict) contact manifold, i.e.\ when $\omega=d\lambda$. The key point is that the level set of a Hamiltonian, while perhaps not contact-type, always carries the structure of a framed Hamiltonian manifold. Indeed, if $H:(W,\Omega)\rightarrow \mathbb R$ is a Hamiltonian in a symplectic manifold, and $0$ is a regular value, then on $M:=H^{-1}(0)$ we can define $\omega=\Omega\vert_M$, and $\lambda$ via $\ker(\lambda)=TM\cap J(TM)$ where $J$ is a choice of $\Omega$-compatible almost complex structure, and setting $\lambda(X_H)=1$. The Reeb field $X_\eta$ of a framed Hamiltonian structure $\eta=(\lambda,\omega)$ is defined implicitly via:
$$\lambda(X_\eta)=1,\;\omega(X_\eta,\cdot)=0.$$ The notion of an almost complex structure $J$ on $\mathbb R_t \times M$ compatible with a framed Hamiltonian manifold is analogous as to the contact case, i.e.\ $J$ is $\mathbb{R}$-invariant, maps $\partial_t$ to $X_\eta$, and preserves $\ker dt \cap \ker \lambda$ on each $t$-slice, where it is $\omega$-compatible.

Now we can define a feral curve in a symplectization (the definition for other target manifolds is completeley analogous).

\begin{definition}[feral curve]
Let $(M,\eta)$ be a framed Hamiltonian manifold, and let $J$ be an almost complex structure compatible with $\eta$ on $\mathbb R \times M$. A feral curve is then a tuple $\mathbf{u}=(u,\Sigma,j,\mu,D)$ where $(\Sigma,j,\mu,D)$ is a Riemann surface with a finite set of marked points $\mu$ and a finite set $D$ of nodal singularities, with finitely many connected components and finite genus, $u:(\Sigma,j)\rightarrow \mathbb (R\times M,J)$ is a map satisfying the Cauchy--Riemann equation, having finite $\omega$-energy (i.e.\ $\int_\Sigma u^*\omega<\infty$), and finally $u$ is assumed to have finitely many \emph{generalized punctures} Punct$(u)$.
\end{definition}

The set of generalized punctures Punct$(u)$ is defined as follows. Let $u:\Sigma \rightarrow \mathbb R \times M=:W$ be a smooth proper map. Choose an exhaustion of $W$ by open sets $W_k\subset W$, i.e.\ sastisfying $W_k\subset W_{k+1}$ for all $k$, $W=\bigcup_k W_k$. Let Punct$_k(u)$ be the number of non-compact path-connected components of $\Sigma\backslash u^{-1}(W_k)$. Then
$$\mbox{Punct}(u):=\lim_k \mbox{Punct}_k(u),$$ in the Hausdorff topology.

\medskip

More recently, Prasad was able to simplify the arguments of Fish--Hofer and strengthen their results. In particular, he obtained the following density result for closed invariant subsets, which is a significant improvement of Theorem \ref{herman}:

\begin{thm}[Prasad \cite{Pr24}, 2024]\label{thm:Prasad} Let $H : \mathbb R^4 \rightarrow \mathbb R$ be a smooth Hamiltonian and let $M$ be a nonempty compact regular level set of $H$. Then $M$ contains an infinite family of pairwise distinct, proper, compact invariant subsets whose union is dense in $M$.
\end{thm}

For clarity, we should emphasize that these invariant subsets are not necessarily disjoint. One should also emphasize that there are examples due to Katok of Hamiltonian dynamical systems on the $3$-sphere (as a star-shaped domain in $\mathbb R^4$) with exactly two periodic orbits, which are also the only minimal subsets, and exactly three invariant ergodic measures (the total measure and the ones supported by the two orbits). In other words, minimal sets, periodic orbits and ergodic measures can be simple, but there are always more complicated closed invariant subsets, as per the above theorem. Moreover, in the Katok example, no collection of three closed invariant subsets can be disjoint. This implies that Theorem \ref{thm:Prasad} cannot be upgraded to disjoint invariant sets, i.e.\ it is sharp in this sense.

The strategy of the proof follows similar lines as to Fish--Hofer's, with the difference that not only degree 1 curves are used, but rather curves of arbitrary large degree $d \gg 1$ in $\mathbb CP^2$ (in other words, the fact that $\mathbb CP^2$ has a very rich Gromov--Witten theory is used), and the procedure for obtaining limiting sets is completely reworked. The latter can be said to be most original contribution of Prasad's paper.

Let us now come to the main observation concerning how these ideas bear on the classical, \emph{planar}, CR3BP (so that phase space is $4$-dimensional). One can adapt the proof of Prasad to the case where $\mathbb R^4$ is replaced by $T^*S^2$, which is precisely what is needed to apply the conclusion of Prasad's theorem to the planar CR3BP, as the Hamiltonian for the latter is autonomous, and therefore the dynamics at a fixed Jacobi constant takes places in a level set inside $T^*S^2$ after Moser regularization. For this, we appeal to Lemma 5.7.1 in \cite{FvK18}. Now, we need to impose that the energy value lies below $H(L_2)$, so that the regularized level set is actually compact.

With this in mind, the adaptation to $T^*S^2$ works as follows. Following \cite[Section 3.6]{Pr24}, using Weinstein's neighbourhood theorem, one symplectically compactifies $T^*S^2$ to $S^2\times S^2$, where the zero section $S^2\subset T^*S^2$ is identified with the Lagrangian anti-diagonal $\overline{\Delta}$. Moreover, the diagonal $\Delta \subset S^2\times S^2$ is a symplectic divisor which can be assumed to be disjoint from the Weinstein neighbourhood. One then needs to obtain a large collection of holomorphic curves, and apply an analogous neck stretching procedure. The relevant curves are provided by Proposition A.5 in \cite{Pr24}, which is a rather technical input which appeals to a version of Taubes's ``Gr=SW'' result, and applies to several $4$-manifolds beyond $S^2\times S^2$. 

And at this point is where Prasad's main original contributions are put into use. Rather than appealing to an SFT-compactness style theorem, or a Fish--Hofer procedure of extracting a feral curve, he uses a completely different argument. Roughly speaking, he introduces the notion of a \emph{stretched limit set}, a topological object in nature, which arises when keeping track of slices of a fixed size of the curve as it is being stretched, and whose elements are roughly speaking limits of these windows but where their height, i.e.\ the $\mathbb R$-coordinate at which they are centered, is recorded. The stretched limit set is not necessarily a holomorphic curve (it could be a fractal even, or a set with non-empty interior), but contains information from all parts of the stretching curve, as opposed to simply well-chosen parts of the necks (which is what happens in the Fish--Hofer approach). 

For completeness, we outline a definition of stretched limit set is what follows. Assume that $M=H^{-1}(0)\subset (W,\Omega)$ is a regular level set, fix an $\Omega$-compatible almost complex structure $J$, and assume that $(W_k,\Omega_k,J_k)$ is a sequence of symplectic manifolds obtained by the standard stretching the neck procedure along $M$ (see e.g.\ \cite{EGH}). Denote by $\mathcal{K}(X)$ the set of closed subsets of a space $X$, with the Hausdorff topology.

\begin{definition}[Stretched limit set, \cite{Pr24}]
    For $a_0 \in \mathbb R$, denote the shift map $$\tau_{a_0}:\mathbb R\times M \rightarrow \mathbb R\times M,$$ $$(a,x)\mapsto (a-a_0,x).$$ For a subsequence $\{J_{k_j}\}$ as above, fix a closed, connected Riemann surface $\Sigma_j$ and a $J_{k_j}$-holomorphic curve $u_j:\Sigma_j\rightarrow W_k$. The \emph{stretched limit set} $\chi \subset \mathcal{K}((-1,1)\times M)\times (-1,1)$ is the collection of pairs $(\zeta,s)$ for which there exists a sequence $a_j \in (-k_j,k_j)$ satisfying:
\begin{enumerate}
    \item $a_j/k_j\rightarrow s;$
    \item A subsequence of the slices
    $$
    \tau_{a_j}\cdot (u_j(\Sigma_j)\cap (a_j-1,a_j+1)\times Y)\subset (-1,1)\times M
    $$
    converge in $\mathcal{K}((-1,1)\times M)$ to $\zeta$.
\end{enumerate}

\end{definition}

For each degree $d\geq 1$, the strategy is then to stretch a degree $d$ curve, of which there are plenty in $S^2\times S^2$ by work of Gromov. While its stretched limit set could be large and complicated, it can be shown to contain a connected subset which consists of nearly invariant sets. Moreover, the elements $(\zeta,s)$ of this subset satisfy that $\zeta$ is $\epsilon$ close to a cylinder $(-1,1)\times \Lambda$, where $\Lambda$ is a closed invariant set, and where $\epsilon$ goes to zero as the degree goes to infinity. With this basic idea in mind, after another limiting procedure taking $d\rightarrow \infty$ and further work, the desired infinite family of proper, distinct, compact invariants subsets is obtained. The fact that they are proper appeals to positivity of intersections, with respect to an additional moduli space of suitably chose holomorphic curves, whereas the fact that there are infinitely many uses moduli spaces with point constraints (and the fact that one can add sufficiently many of the latter, combined with careful quantitative estimates). For more details (of which we have omitted many), the reader is of course referred to \cite{FH23,Pr24}. 

In a nutshell, we have sketched a proof of the following result.

\begin{thm}[M.-Prasad, 2024]\label{thm:goal}
    Consider the Hamiltonian $H$ for the \emph{planar} CR3BP. Choose any mass ratio $\mu$ and any energy $c \in (-\infty, H(L_1))\cup (H(L_1),H(L_2))$. Then, for the near-primary dynamics, there exists an infinite collection of pairwise distinct, proper, either properly embedded or closed, invariant subsets $\Lambda \subset H^{-1}(c)$ whose union is dense in $H^{-1}(c)$. 
\end{thm}

Note here that the conclusion holds for the \emph{unregularized} dynamics, as one obtains it first for the regularized case, and the conclusion holds directly for the unregularized system, with the observation that any closed subset intersecting the collision locus now simply becomes a properly embedded invariant subset (here, do \emph{not} confuse proper --i.e.\ nonempty nor the whole space-- with properly embedded --i.e.\ preimage of compact is compact under the embedding).

\part{Practical aspects}

\chapter{Symplectic data analysis}\label{ch:symp_data_analysis}

This chapter is motivated by the study of periodic orbits of Hamiltonian systems, in families, for practical purposes. A given Hamiltonian system usually depends on parameters (e.g.\ energy or mass parameters), which one may vary. Under such deformation, periodic orbits may undergo \emph{bifurcation}, a mechanism by which new families of periodic trajectories arise, occurring when one orbit in the family becomes degenerate. One may be interested in following those new families, as they might entail practical interest. For instance, in the preliminary stages of space mission design, one is interested in mapping out the largest possible data base of periodic orbits of a concrete system, modelling the region of space to which one will send a spacecraft (e.g.\ the Earth--Moon system). The way different families can connect to each other is encoded in the topology of a \emph{bifurcation graph}. These may be thought of as ``highways'' representing families of orbits, and the ``exits'' from one highway to another are represented by bifurcating orbits.  

One then wishes to optimize over all orbits, by taking into account a large number of practical considerations (e.g.\ minimize fuel usage, risk of collisions, station-keeping and time flight, and maximize safety, as well as design trajectories which fulfill the requirements of a given mission). While stable orbits are used for parking spacecraft, unstable ones are usually used for transfer trajectories by intersecting their stable and unstable manifolds. The vast amount of data imposes the need to introduce tools and methods in order to keep track of all the information, which leads naturally into the real of data analysis, where computationally cheap methods are central. The aim of this chapter is then to introduce a ``symplectic toolkit'', extracted from the modern methods of symplectic geometry, and designed to study periodic orbits of arbitrary Hamiltonian systems, together with their bifurcations in familes, and their stability. The emphasis is on visual, easy and resource--efficient implementation.

After introducing the basic mechanism of bifurcations, we will review our toolkit, consisting of the following elements:

\begin{enumerate}
    
    \item[(1)] \textbf{Floer numerical invariants:} Numbers which stay invariant before and after a bifurcation, and so can help predict the existence of orbits, as well as being easy to implement. There is one invariant for arbitrary periodic orbits, and another for \emph{symmetric} periodic orbits \cite{FKM}. They are defined as the Euler characteristic of local Floer homology, and of local Lagrangian Floer homology of the fixed-point locus, respectively.

\medskip
    
    \item[(2)] \textbf{The B-signs} \cite{FM}: a $\pm$ sign associated to each elliptic or hyperbolic Floquet multiplier of an orbit, which helps predict bifurcations. This is generalization of the classical Moser--Krein signature \cite{Kre2,Kre3,M78}, which originally applies only to elliptic Floquet multipliers, to also include the case of hyperbolic multipliers, whenever the corresponding orbit is \emph{symmetric}. 

\medskip
    
    \item[(3)] \textbf{Global topological methods:} the \emph{GIT-sequence} \cite{FM}, a sequence of spaces whose global topology encodes (and sometimes forces) bifurcations, and refines Broucke's stability diagram \cite{Br69} by adding the $B$-signs.

\medskip
    
    \item[(4)]\textbf{Conley-Zehnder index \cite{CZ84,RS93}:} a winding number associated to each non-degenerate orbit, extracted from the topology of the symplectic group, which does not change unless a bifurcation occurs. Therefore it can be used to determine which families connect to which.
\end{enumerate}

\begin{figure}
    \centering
    \includegraphics[width=1\linewidth]{toolkit.png}
    \caption{The diagram, due to Bhanu Kumar, summarizes the symplectic toolkit and its uses.}
    \label{fig:placeholder}
\end{figure}

There are currently two ``CZ-index calculators'' available in the literature, one implemented by Otto van Koert (in Python\footnote{Available at https://github.com/ovkoert/cz-index}) and one by Bhanu Kumar (in MATLAB\footnote{Available at https://github.com/bhanukumar314}). They provide a means of directly computing the CZ-index of an orbits given its initial conditions.

\section{Bifurcations}\label{sec:bifurcations} We begin with a standard definition. Let $(M,\omega)$ be a symplectic manifold, and $H: M\rightarrow \mathbb R$ an autonomous Hamiltonian.

\begin{definition}
    A \emph{regular orbit cylinder} is a smooth map $\Gamma: [a,b] \times S^1 \rightarrow M$ such that $\gamma_c=\Gamma(c,\cdot): S^1 \rightarrow M$ is a non-degenerate periodic orbit for $H$, and $H(\Gamma(c,t))=c$ is a regular value of $H$.
\end{definition}

In other words, a regular orbit cylinder is a family of non-degenerate periodic orbits parametrized by the energy, moving along regular level sets of the Hamiltonian. One could also use other parameters to parametrized a family, but most of our applications will involve the energy as the underlying parameter being varied. An easy consequence of the implicit function theorem is that every non-degenerate periodic orbit can be extended locally to a regular orbits cylinder, i.e.\ it may be continued to a family. But there is a limit to continuing the family, as the orbits in the cylinder might converge to a \emph{degenerate} orbit, or to a critical point of the Hamiltonian. A \emph{bifurcation} occurs precisely when this happens. The first case corresponds to the situation when the eigenvalue $1$ appears in the spectrum of the matrix of the bifurcating orbit, and is crossed transversely in the $1$-parameter family. In what follows, rather than give a formal definition, we give a number of examples for Hamiltonian systems with two degrees of freedom.

\begin{figure}
    \centering
    \includegraphics[width=0.5\linewidth]{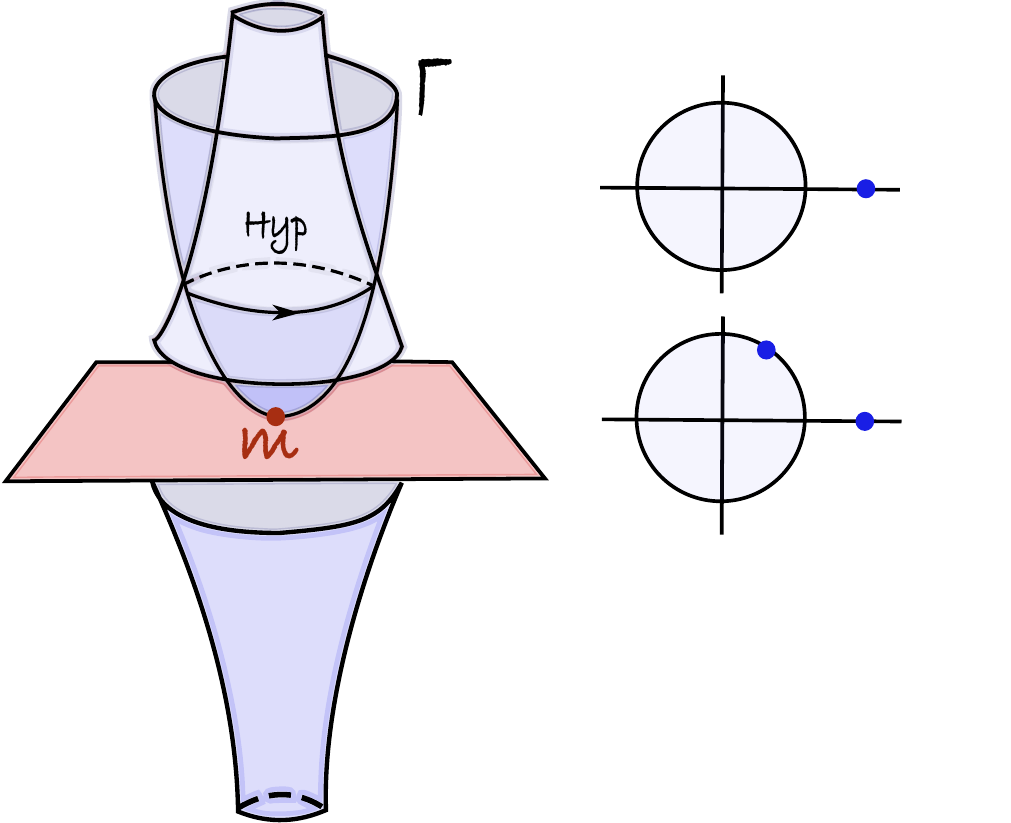}
    \caption{The phantom burst. The principal Floquet multipliers of the critical point $m$ are shown on the bottom right (elliptic-positive hyperbolic), and those of the hyperbolic cylinder forming the center manifold, on the top right.}
    \label{fig:phantom-burst}
\end{figure}

\begin{figure}[h]
    \centering
    \includegraphics[width=0.6\linewidth]{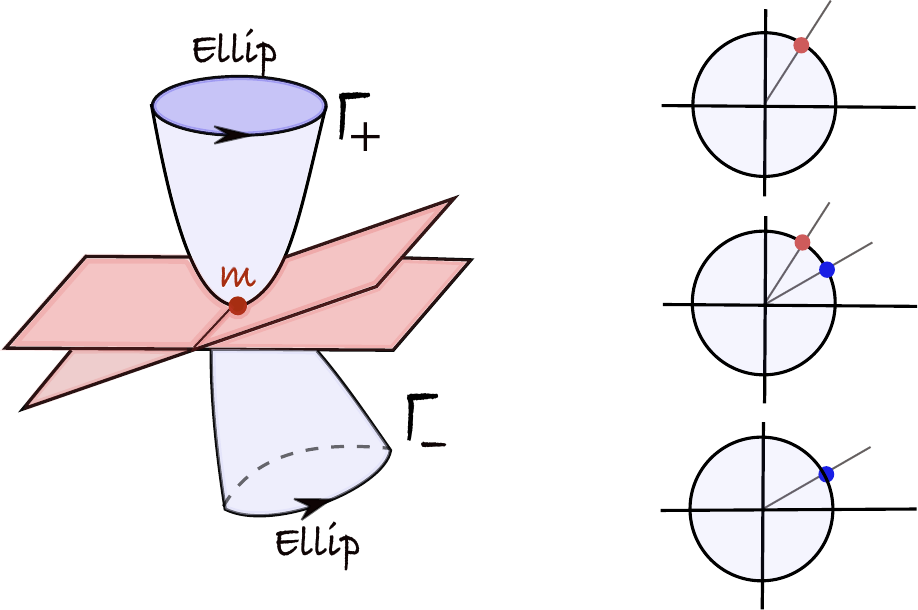}
    \caption{The stable burst. The critical point is doubly elliptic.}
    \label{fig:stable-burst}
\end{figure}

\subsection{A zoo of bifurcations} We will follow the list given in page 595 of the classical book of Abraham--Marsden \cite{AM78}. The following are the ``generic'' cases of bifurcations for a Hamiltonian in dimension $4$ (we refer to \cite{AM78} for a precise meaning of genericity, for each case discussed below). In practice, when considering a concrete system like the CR3BP, genericity is not something that can be counted on, and more complicated bifurcations can occur. This is also the case in the presence of symmetry.

\medskip

\textbf{The burst.} This bifurcation occurs when an orbit family originates or terminates at a critical point $m$ of the Hamiltonian. There are two possibilities: either $m$ is hyperbolic-elliptic (a \emph{phantom} burst), or elliptic-elliptic (a \emph{stable} burst). In both cases, the family of periodic orbits approaches the critical point along a center manifold, which exists by Lyapunov's center theorem. The center manifold consists of a family of periodic orbits which close up at the critical point, where it is tangent to the eigenspace of an elliptic Floquet multiplier, and so forms a disk near $m$ with $m$ at the origin. In the first case, the family consists of hyperbolic orbits, with Floquet multipliers near the hyperbolic pair of multipliers of $m$. In the second case, there are two (sub-)center manifolds, one for each elliptic Floquet multiplier $e^{2\pi i \theta_j}$ (here we need to assume that there are no resonances between the two \emph{frequencies} $\theta_j$ --i.e.\ the Floquet exponents--, which holds generically), each consisting of elliptic orbits, with Floquet multipliers close to those of $m$. In the latter case, there are three further sub-cases. if $m$ is a local minimum, the energy level set is locally a three-sphere containing two orbits (one for each sub-center manifold) which come together at $m$ as the energy decreases to the critical value. If $m$ is a local maximum, the parameter is reversed. If $m$ is hyperbolic, then the level set is locally a hyperboloid, containing one orbit which collapses to $m$, and then is reborn into a new orbit. One speaks of \emph{reincarnation} in the last case.

\medskip

\textbf{Birth-death.} In this case, as the parameter varies, a positive hyperbolic orbit comes together with an elliptic periodic orbit along a degenerate periodic orbit, which is tangent to the energy level set. After the critical value is passed, the orbits disappear (\emph{death}, or \emph{annihilation}). If the parameter is reversed, then the two orbits get created out of nothing (\emph{birth} or \emph{creation}). This can be viewed either as two distinct families (one elliptic, the other hyperbolic) which cancel each other, or as a single family which is tangent to the level set. In the latter case, the orbit transitions transversely from elliptic to positive hyperbolic, along the eigenvalue $1$.  

\begin{figure}
    \centering
    \includegraphics[width=0.7\linewidth]{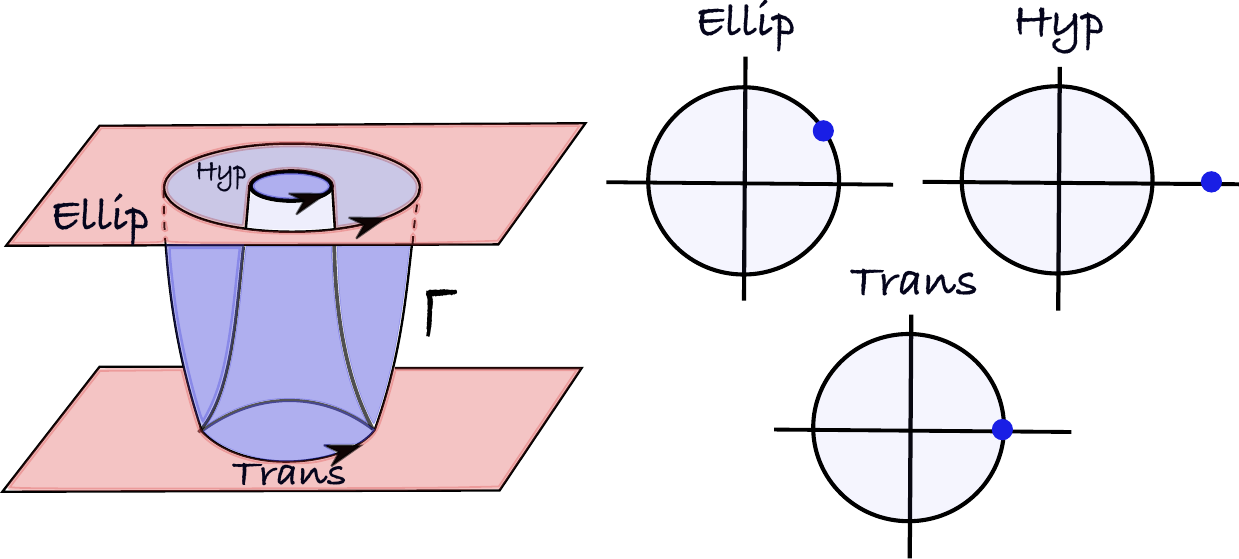}
    \caption{Birth-death.}
    \label{fig:birth-death}
\end{figure}

\medskip

\textbf{Period-doubling bifurcation.} This case occurs when an elliptic orbit $\gamma$ transversely transitions to negative hyperbolic, along the eigenvalue $-1$. The orbit itself does not bifurcate, but its double cover $\gamma^2$ does. There are two possibilities, depending on whether the new orbit $\beta$ that appears out of $\gamma^2$ arises before or after the critical energy value. In the first case, $\beta$ is positive hyperbolic orbit and doubles over itself, converging to $\gamma^2$ at the critical value. This case is sometimes called \emph{murder}, as $\gamma$ becomes unstable and hence qualitatively ``dies'', at the hand of the ``murderer'' $\beta$. When the parameter is reversed, this is called \emph{materialization}, as the unstable $\gamma$ ``materializes'' into existence by becoming stable, and hence qualitatively visible. In the second case, $\beta$ is elliptic, and its period is close to twice the period of the original orbit, converging to $\gamma^2$ as the parameter approaches the critical value. This is sometimes called \emph{subtle division}, or \emph{subtle doubling} in one direction, and \emph{subtle halving} in the other.

\begin{figure}
    \centering
    \includegraphics[width=0.6\linewidth]{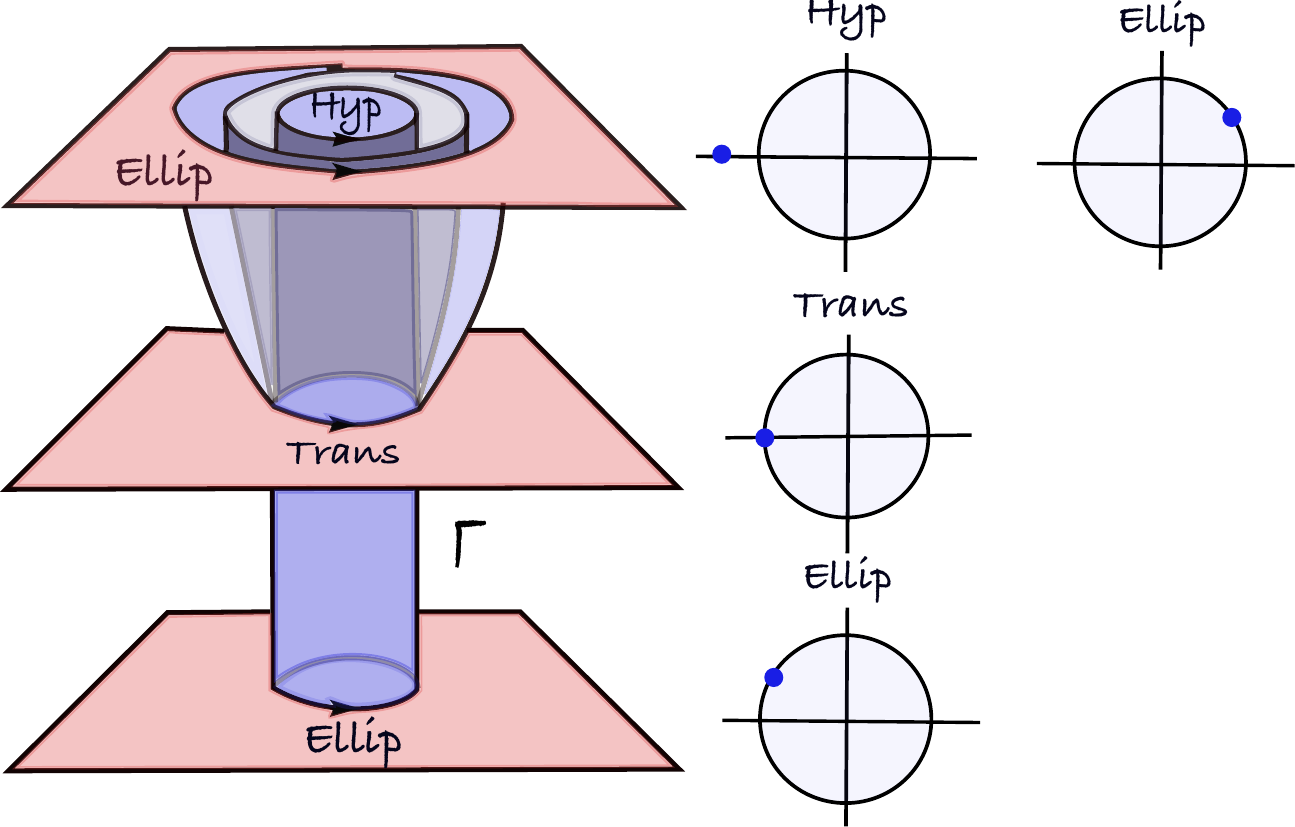}
    \caption{Period-doubling bifurcation.}
    \label{fig:enter-label}
\end{figure}

\begin{figure}
    \centering
    \includegraphics[width=0.6\linewidth]{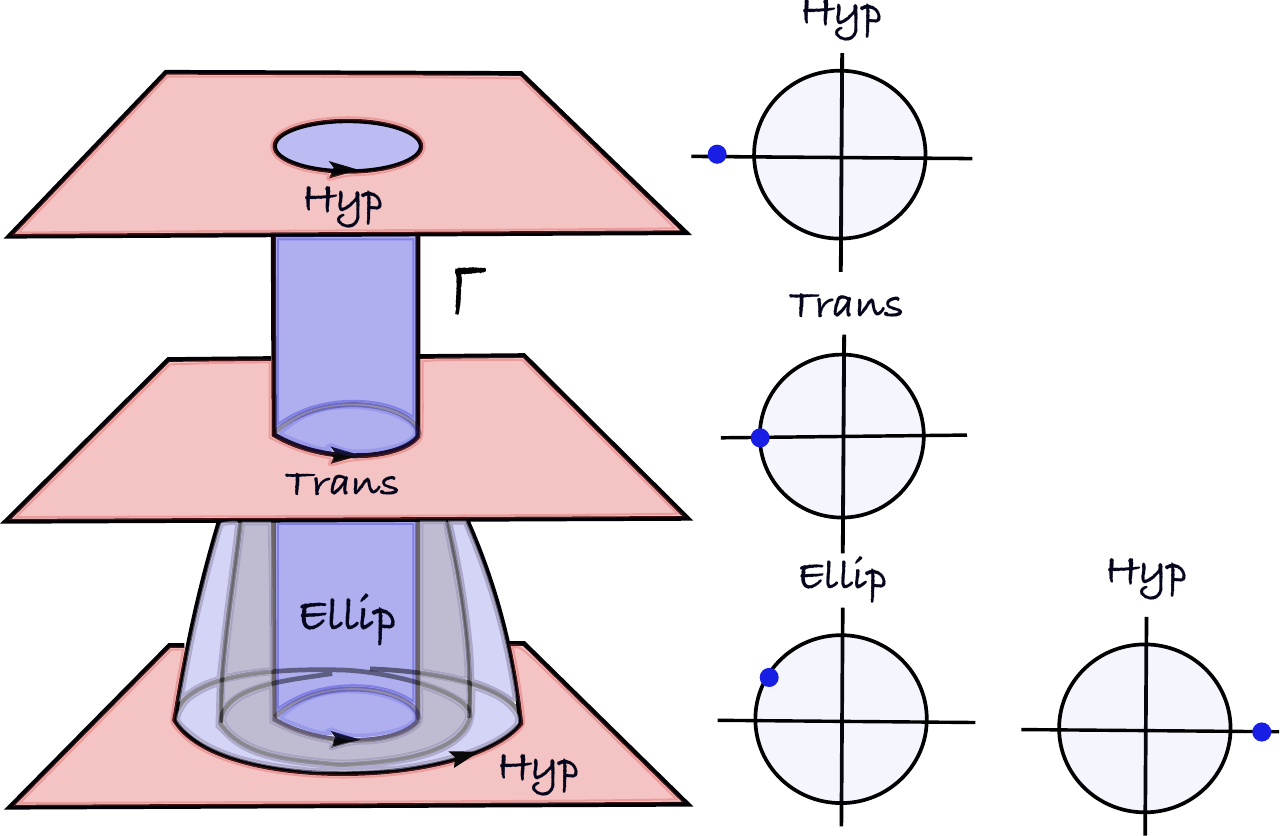}
    \caption{Murder.}
    \label{fig:murder}
\end{figure}

\medskip

\textbf{Phantom kisses.} This bifurcation takes place when an orbit cylinder $\Gamma$ transversely crosses an eigenvalue which is either a third or fourth root of unity. Both cases are analogous. While the family $\Gamma$ does not bifurcate (as it stays elliptic), there are two hyperbolic orbit cylinders, one in either direction with respect to the parameter. These cylinders ``kiss'' $\Gamma$ at the critical energy value, where they end, converging to the third (respectively the fourth) cover of the underlying simple orbit of $\Gamma$. 

\begin{figure}
    \centering
    \includegraphics[width=0.6\linewidth]{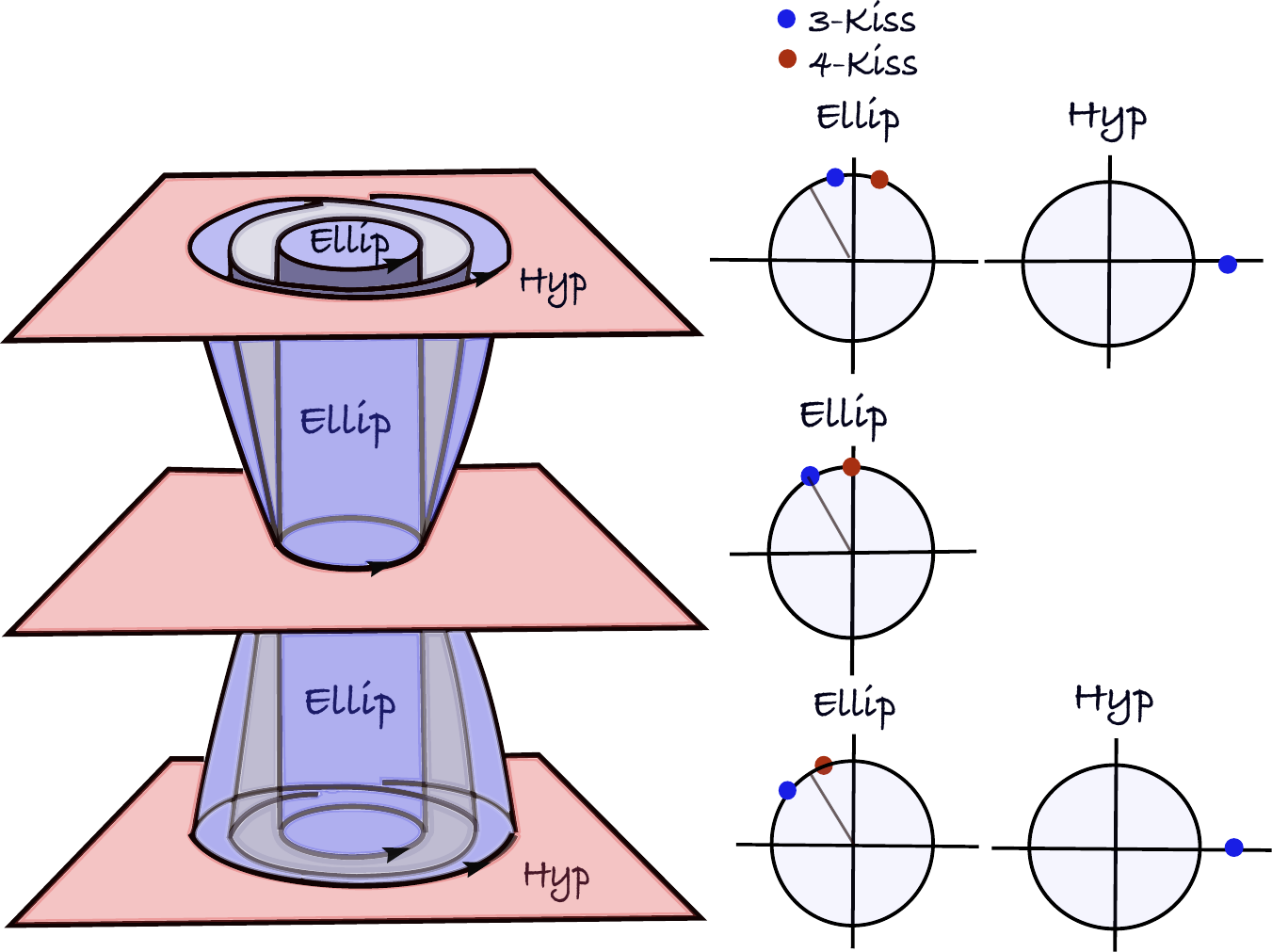}
    \caption{Phantom kisses, for the case of a $3$-kiss.}
    \label{fig:phantom_kisses}
\end{figure}

\medskip

\textbf{$k$-fold bifurcation.} When $k\geq 4$, this case occurs when a family of elliptic orbits crosses a $k$-th root of unity. The cylinder continues without bifurcating, but the $k$-fold cover of the orbit bifurcates. At the critical energy value, there appear two periodic orbits out of the $k$-fold cover, one which is elliptic, and the other, hyperbolic, only defined in the positive direction of the parameter. One sometimes talks about \emph{emission} in this direction, and \emph{absorption} in the opposite direction (i.e.\ two families get absorbed and disappear).

\begin{figure}
    \centering
    \includegraphics[width=0.6\linewidth]{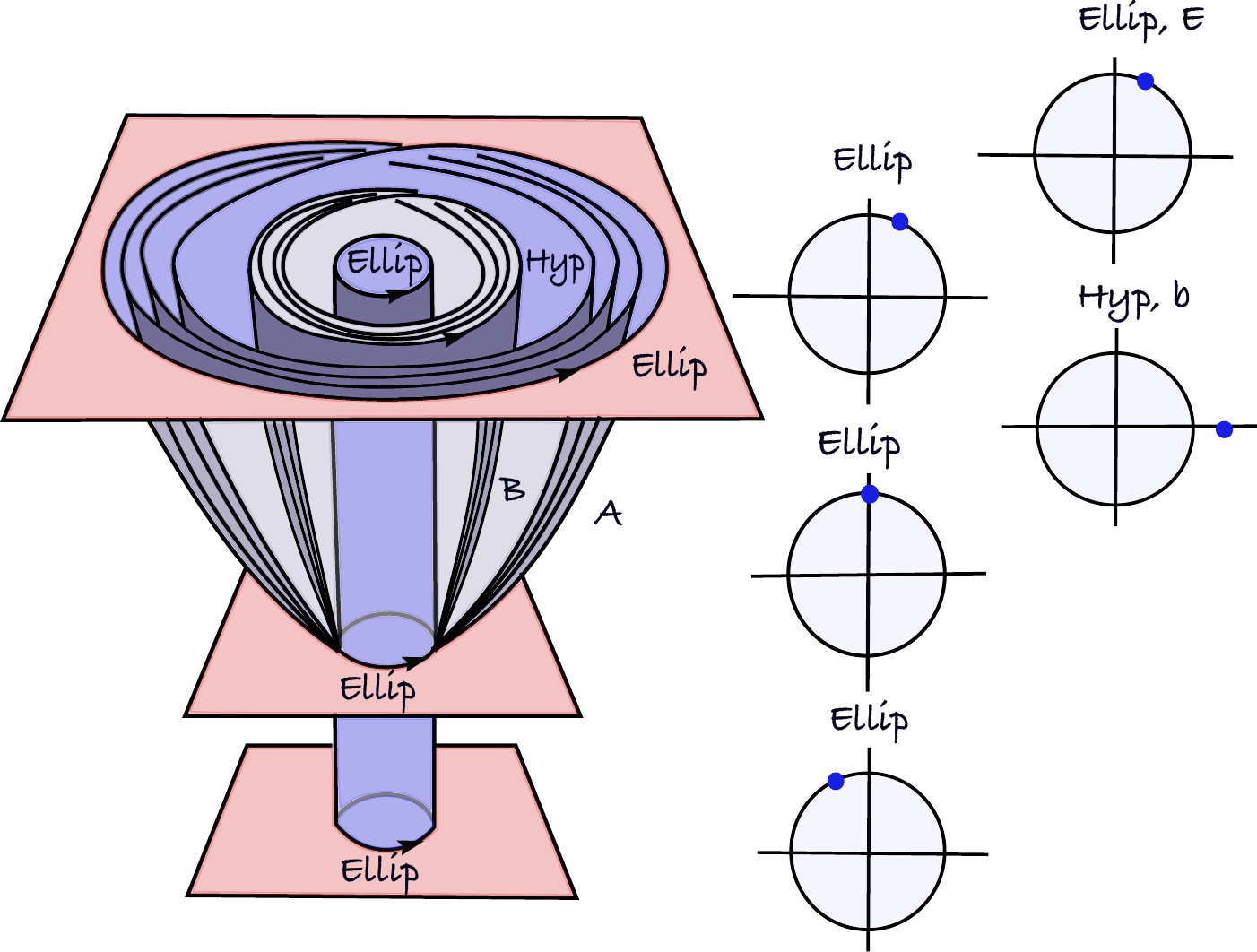}
    \caption{A $k$-fold bifurcation, with $k=4$.}
    \label{fig:emission}
\end{figure}

\section{The B-signs} In this section, we discuss the notion of \emph{$B$-signature}, or \emph{$B$-signs}, introduced by the author and Urs Frauenfelder in \cite{FM}, associated to a symmetric periodic orbit (see also \cite{Ay22}).

Recall from Chapter \ref{ch:basic_notions} that the monodromy matrix of a symmetric periodic orbit, at a symmetric point, and after choosing a basis for the Lagrangian fixed-point locus of the symmetry, is a Wonenburger matrix, i.e.\ of the form 
\begin{equation*}\label{symsymp}
M=M_{A,B,C}=\left(\begin{array}{cc}
A & B\\
C & A^t
\end{array}\right)\in Sp^\mathcal{I}(2n)\subset Sp(2n),
\end{equation*}
where
\begin{equation}\label{eq:Wonenburger2}
B=B^t,\quad C=C^t,\quad AB=BA^t,\quad
A^tC=CA,\quad A^2-BC=\mathds 1.
\end{equation}

Recall that the spectrum of $M$ is completely determined by the spectrum of $A$. Let $\lambda,1/\lambda$ be an elliptic or hyperbolic pair of eigenvalues of $M$, which is simple (i.e.\ each eigenvalue is of multiplicity $1$). Then the stability index $a=a(\lambda)=\frac{1}{2}(\lambda+\frac{1}{\lambda})$ is an eigenvalue of $A$ (which has the same eigenvalues as $A^t)$, and is also simple. Let $v$ an eigenvector of $A$ and $w$ an eigenvector of $A^t$ with eigenvalue $a$. Assume that the pair $\lambda,1/\lambda$ is ordered so that $\lambda$ has positive imaginary part in the elliptic case, or has absolute value which is greater than $1$ in the hyperbolic case (i.e.\ $\lambda$ is the \emph{principal} eigenvalue).

\begin{definition}[\textbf{B/C-signs}] The \emph{B-sign} of $\lambda$ is the sign of $w^tBw$, and the \emph{C-sign} of $\lambda$ is the sign of $v^tCv$, denoted
$$
\mathrm{sign}_B(\lambda)=\mathrm{sign}(w^tBw)= \pm, 
$$
$$
\mathrm{sign}_C(\lambda)=\mathrm{sign}(v^tCv)= \pm. 
$$  
By definition, we let $\mathrm{sign}_{B/C}(1/\lambda)=-\mathrm{sign}_{B/C}(\lambda)$.
\end{definition}

We can think that we have attached a sign to the \emph{pair} $\lambda,1/\lambda$, ordered as above. Alternatively, we can attach these signs to the eigenvalues of $A$, as they correspond to the principal eigenvalues of $M$, after choosing the complex square root with positive imaginary part in the expression $\lambda=a+\sqrt{1-a^2}$. That the numbers $w^tBw$ and $v^tCv$ are non-vanishing follows from Equations (\ref{eq:Wonenburger2}), and the fact that we are assuming the eigenvalues are elliptic or hyperbolic. It is easily checked that these definitions are independent of the eigenvectors $v,w$, and of the basis chosen. Therefore, to each symmetric point of a symmetric periodic orbit we have associated, for each simple elliptic or hyperbolic pairs of eigenvalues, a $B/C$-sign. For a given Wonenburger matrix, after ordering the simple and real eigenvalues of its $A$ block in strictly increasing order (which gives an order of the corresponding pairs of elliptic or hyperbolic pairs of eigenvalues of $M$), we obtain an ordered tuple of $B/C$ signs of the form $(\pm, \dots, \pm)$, respectively called the \emph{$B/C$-signature} of the matrix. Therefore a symmetric periodic orbit has \emph{two} such signatures, one for each of the two symmetric points. Moreover, the $C$-sign is completely determined by the $B$-sign, and viceversa, so they provide the same information. Namely, the $B$-sign agrees with the $C$-sign if $\lambda$ is hyperbolic, and they disagree if $\lambda$ is elliptic (this can be seen by inspecting the normal forms provided in \cite{FM}; see Example \ref{ex:simple_case} for the case $n=1$).

\begin{example}\label{ex:simple_case}
For the simplest case $n=1$, each Wonenburger matrix
\begin{equation*}
M=\left(\begin{array}{cc}
a& b\\
c & a
\end{array}\right)\in Sp^\mathcal{I}(2)\subset SL(2,\mathbb R),
\end{equation*}
with $a^2-bc=1$, has one $B/C$-sign, corresponding respectively to the sign of the entries $b$ and $c$, which are non-vanishing if $\mathrm{tr}(M)\neq \pm 2$, i.e.\ for the non-parabolic case. 

In particular, if $M$ is elliptic, it is symplectically conjugated to
a rotation matrix
$$
M\sim \left(\begin{array}{cc}
  \cos(\theta)   &  -\sin(\theta)\\\
   \sin(\theta)  &  \cos(\theta)
\end{array}\right),
$$
with Floquet multipliers $e^{\pm 2\pi i \theta}$ and stability index $\cos(\theta)$, so the $B$-sign of this pair is $$\mathrm{sign}_B(\lambda)=-\mathrm{sign}(\sin(\theta)),$$ whereas the $C$-sign is
$$\mathrm{sign}_C(\lambda)=\mathrm{sign}(\sin(\theta))=-\mathrm{sign}_B(\lambda).$$

If $M$ is hyperbolic, it is symplectically conjugated to a matrix of the form

$$
M\sim \left( \begin{array}{cc}
\pm \cosh(u)    & \sinh(u) \\
\sinh(u)     & \pm \cosh(u)
\end{array}\right),
$$
with Floquet multipliers $\lambda=\pm e^u, 1/\lambda=\pm e^{-u}$. Then
$$
\mathrm{sign}_B(\lambda)=\mathrm{sign}(\sinh(u))=\mathrm{sign}_C(\lambda).
$$
\end{example}

\begin{example} As a simple example with $n=2$, consider the Wonenburger matrix
$$
M=\left(\begin{array}{cccc}
\cos \theta_1 &      0       & -\sin \theta_1 & 0\\
0             & \cos\theta_2 &            0   & -\sin \theta_2\\
\sin \theta_1 &     0        & \cos \theta_1  & 0\\
0             & \sin\theta_2 &            0   & \cos \theta_2
\end{array}\right)\in Sp^\mathcal{I}(4),
$$
which has two pairs of elliptic eigenvalues $\lambda_j^\pm=e^{\pm 2\pi i \theta_j}$, for $j=1,2$. Assume wlog that $\theta_j\in (0,\pi)$, and if $a_j=\cos\theta_j$, then $a_1<a_2$. Then the $B$-sign of $a_j$ as an eigenvalue of the $A$-block, i.e.\ the $B$-sign of $\lambda_j^+$, is 
$$
\mathrm{sign}_B(a_j)=\mathrm{sign}_B(\lambda^+_j)=-\mathrm{sign}(\sin \theta_j)=-,
$$
for $j=1,2$, so that the $B$-signature of $M$ is $\mathrm{sign}_B(M)=(-,-)$.

Similarly, if we consider the Wonenburger matrix
$$
N=\left(\begin{array}{cccc}
\cos \theta_1 &      0       & \sin \theta_1 & 0\\
0             & \cos\theta_2 &            0   & -\sin \theta_2\\
-\sin \theta_1 &     0        & \cos \theta_1  & 0\\
0             & \sin\theta_2 &            0   & \cos \theta_2
\end{array}\right)\in Sp^\mathcal{I}(4),
$$
one checks as above that its $B$-signature is $\mathrm{sign}_B(N)=(+,-)$. Analogously, the $B$-signature of 
$$
P=\left(\begin{array}{cccc}
\cos \theta_1 &      0       & -\sin \theta_1 & 0\\
0             & \cos\theta_2 &            0   & \sin \theta_2\\
\sin \theta_1 &     0        & \cos \theta_1  & 0\\
0             & -\sin\theta_2 &            0   & \cos \theta_2
\end{array}\right)
$$
is $\mathrm{sign}_B(P)=(-,+)$, and that of
$$
R=\left(\begin{array}{cccc}
\cos \theta_1 &      0       & \sin \theta_1 & 0\\
0             & \cos\theta_2 &            0   & \sin \theta_2\\
-\sin \theta_1 &     0        & \cos \theta_1  & 0\\
0             & -\sin\theta_2 &            0   & \cos \theta_2
\end{array}\right).
$$
is $\mathrm{sign}_B(R)=(+,+)$. The matrices $M,N,P,R$ are the four normal forms for Wonenburger $4\times 4$ matrices with a pair of distinct elliptic eigenvalues, up to the action of $GL_n(\mathbb R)$; see \cite{FM}.
\end{example}

The following result \cite{FMb} illustrates the uses of these signs, as they give information on the type of periodic orbit.

\begin{thm}[\cite{FMb}]
A symmetric periodic orbit of a Hamiltonian system with two degrees of freedom is negative hyperbolic if and only if its two $B$-signs are different.
\end{thm}

The $B/C$-signs can be also defined in a straightforward way to the case where the eigenvalues are not necessarily simple. Indeed, let $\lambda$ be an elliptic or hyperbolic eigenvalue of $M\in Sp^{\mathcal{I}}(2n)$. Let $E_{\lambda}$ and $E^t_{\lambda}$ denote the $a(\lambda)$-eigenspaces of $A$ and $A^t$ respectively. We can then view $B$ as a bilinear form on $E^t_\lambda$, via
$$
B(v,w)=v^t B w,
$$
for $v,w \in E^t_\lambda$. Similarly, we can view $C$ as a bilinear form on $E_\lambda$.

\begin{definition}[\textbf{B/C-signature}] The \emph{B-signature} of $\lambda$ is the signature of $B\vert_{E^t_\lambda}$, and the \emph{C-signature} of $\lambda$ is the signature of $C\vert_{E_\lambda}$, denoted
$$
\mathrm{sign}_B(\lambda)=\mathrm{sign}(B\vert_{E^t_\lambda}), 
$$
$$
\mathrm{sign}_C(\lambda)=\mathrm{sign}(C\vert_{E_\lambda}). 
$$
We define the $B/C$-signature of $1/\lambda$ as the $B/C$-signature of $\lambda$.
\end{definition}

Recall that the signature of a non-degenerate bilinear form $G$ is the pair $(p,q)$, where $p$ is the dimension of a maximal subspace where $G$ is positive definite, and $q$ is the dimension of a maximal subspace where $G$ is negative definite. The fact that bilinear forms above are non-degenerate follows from Equations (\ref{eq:Wonenburger2}) and ellipticity/hyperbolicity of the eigenvalues (see Section \ref{sec:GIT_any_dim} for more details). Given a Wonenburger matrix, we order the real eigenvalues of $A$ in (non-strictly) increasing order, and this gives an ordered tuple $((p_1,q_1),\dots,(p_m,q_m))$ of $B/C$-signatures, which we call the $B/C-signature$ of the matrix. In the case where the eigenvalues are simple as above, one replaces $(1,0)$ with a $+$, and $(0,1)$ with a $-$ (as these are the only possibilities).

\begin{example}
    As a simple example with $n=2$, consider the Wonenburger matrices

$$M=\left(\begin{array}{cccc}
\cosh \theta &      0       & -\sinh \theta & 0\\
0             & \cosh\theta &            0   & -\sinh \theta\\
-\sinh \theta &     0        & \cosh \theta & 0\\
0             & -\sinh\theta &            0   & \cosh \theta
\end{array}\right), 
N=\left(\begin{array}{cccc}
\cosh \theta&      0       & \sinh \theta & 0\\
0             & \cosh\theta &            0   & -\sinh \theta\\
\sinh \theta &     0        & \cosh \theta  & 0\\
0             & -\sinh\theta &            0   & \cosh \theta
\end{array}\right),$$ 
$$
P=\left(\begin{array}{cccc}
\cosh \theta &      0       & \sinh \theta & 0\\
0             & \cosh\theta &            0   & \sinh \theta\\
\sinh \theta &     0        & \cosh \theta  & 0\\
0             & \sinh\theta &            0   & \cosh \theta
\end{array}\right),$$
with $\theta\in (0,+\infty)$. As the $A$-block has eigenvalue $a=\cosh \theta$ with multiplicity two, $M,N,P$ have a pair of positive hyperbolic eigenvalues $\lambda=e^\theta,1/\lambda=e^{-\theta}$ with multiplicity two. Their respective $B$-signatures are $\mathrm{sign}_B(M)=((0,2)),\mathrm{sign}_B(N)=((1,1)),\mathrm{sign}_B(P)=((2,0)).$ These are the three normal forms for Wonenburger matrices which are doubly positive hyperbolic and have eigenvalue $\lambda$ with multiplicity two, see \cite{FM}.
\end{example}

\section{Global topological methods: the GIT sequence}\label{sec:GIT} We now discuss global topological methods in the study of periodic orbits, following the exposition in \cite{AFvKKM}. These methods encode: bifurcations; stability; eigenvalue configurations; obstructions to existence of regular families; and $B$-signs, in a visual and resource-efficient way. The main tool is the \emph{GIT sequence}, introduced in \cite{FM}, as a refinement of the \emph{Broucke stability diagram} \cite{Br69} via implementing the $B$-signs. This is a sequence of three branched spaces (or \emph{layers}), arranged into top, middle, and bottom (or base), together with two maps between them, which collapse certain branches together. Each branch is labelled by the $B$-signs. A symmetric orbit gives a point in the top layer, and an arbitrary orbit, in the middle layer. The base layer is $\mathbb R^n$ (the space of coefficients of the characteristic polynomial of the first block $A$ of $M_{A,B,C}$). Then a family of orbits gives a path in these spaces, so that their topology encodes valuable information, as it may sometimes enforce bifurcations, i.e.\ provide obstructions to the existence of regular families. The details are as follows. 

\subsection{GIT sequence: 2D} We first consider the simplest case, i.e.\ the case of an autonomous Hamiltonian of two degrees of freedom (like the \emph{planar} CR3BP), so that the reduced monodromy matrix is an element in $SL(2,\mathbb R)$.

Let $\lambda$ eigenvalue of $M\in Sp(2)$, with associated stability index $a(\lambda)=\frac{1}{2}(\lambda+1/\lambda)$. Then we have
\begin{itemize}
    \item $\lambda=\pm 1$ if and only if $a(\lambda)=\pm 1$;
    \item $\lambda$ positive hyperbolic if and only if $a(\lambda)>1$; 
    \item $\lambda$ negative hyperbolic if and only if $a(\lambda)<-1$; and 
    \item $\lambda$ elliptic (i.e.\ stable) if and only if $-1<a(\lambda)<1$.
\end{itemize}
The Broucke stability diagram is then the real line, split into three components; see Figure \ref{fig:GIT_sequence}. If two orbits lie in different components of this diagram, then there are always bifurcations in any family joining them, as the topology of the diagram implies that any path between them has to cross the $\pm 1$ eigenvalues (which correspond respectively to bifurcation or period-doubling bifurcation).

Topologically, the stability index ``collapses'' the two elliptic branches in the middle layer of Figure \ref{fig:GIT_sequence} together, to form an interval. These two branches are distinguished by the $B$-signs, which coincide with the Krein signs \cite{Kre2,Kre3}. There is an extra top layer corresponding to symmetric orbits, where now each hyperbolic branch separates into two, and there is a collapsing map from the top to middle layer similarly as before. To go from one branch to the other (say from the positive hyperbolic branch I to the positive hyperbolic branch II), the topology of the top layer implies that the eigenvalue 1 is necessarily crossed. This then means that there are always bifurcations in any (symmetric) family joining them, even if they project to the same component of the Broucke diagram. Therefore the information given by the diagram is much more refined for the case of symmetric orbits. If we declare two orbits to be \emph{qualitatively equivalent} if they can be joined by a regular orbit cylinder, then the topology of the spaces in the GIT sequence give (topological) criteria in order to determine whenever two orbits are \emph{not} qualitatively equivalent. In conclusion:

\begin{figure}[h]
    \centering
    \includegraphics[width=0.8\linewidth]{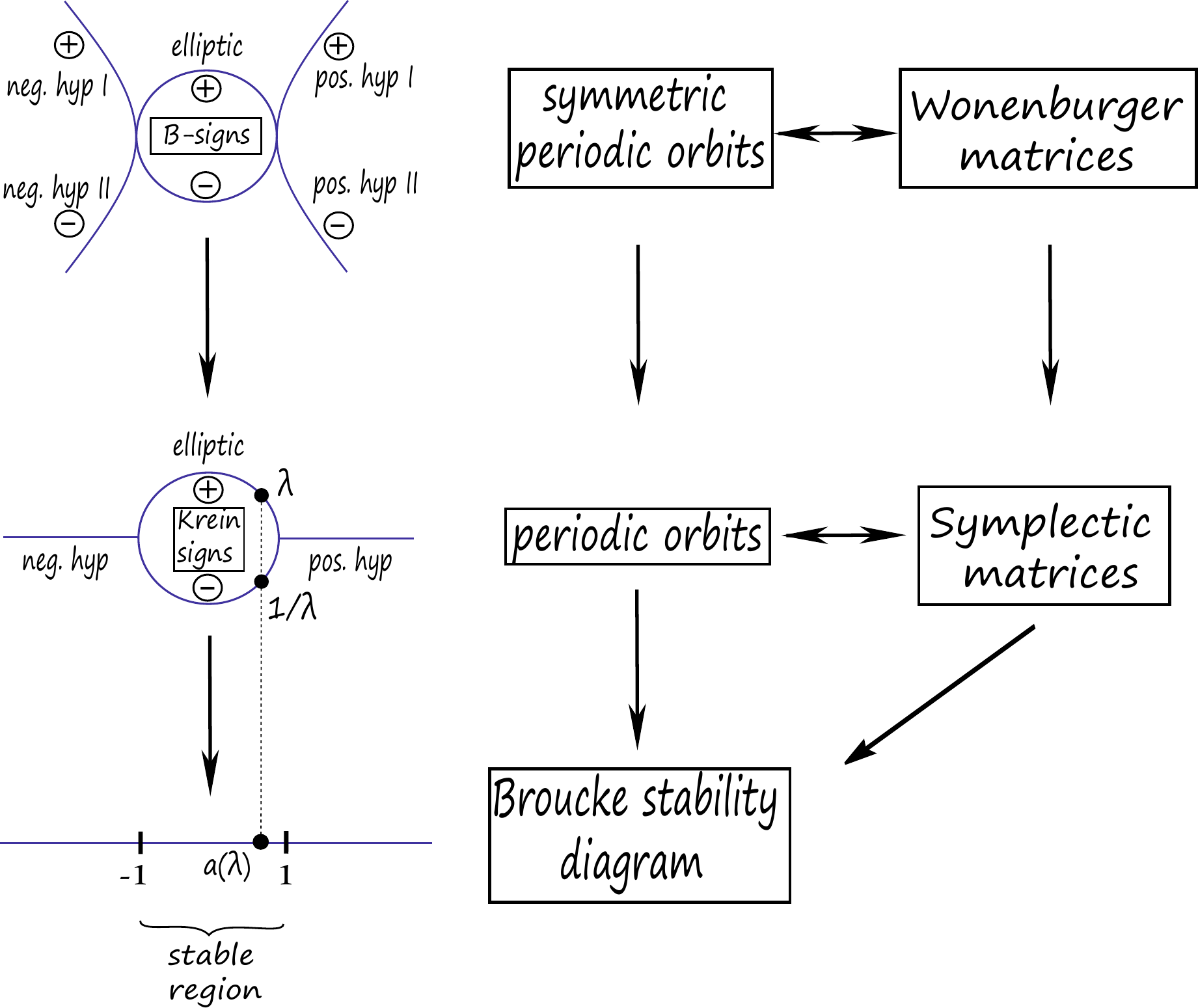}
    \caption{The 2D GIT sequence. The information is more refined for symmetric orbits.}
    \label{fig:GIT_sequence}
\end{figure}

\begin{itemize}
    \item $B$-signs ``separate'' hyperbolic branches, for the case of symmetric orbits.

\medskip
    
    \item If two orbits lie in different components of the Broucke stability diagram, there are always bifurcations in any path joining them.

 \medskip
 
    \item If two symmetric orbits lie in the same component of the Broucke diagram, but if their $B$-signs differ, there are always bifurcations in any (symmetric) path joining them. 
    \end{itemize}

\subsection{GIT sequence: 3D} Now we use the same idea, but for the case of autonomous Hamiltonian systems with three degrees of freedom, for which reduced monodromy matrices are elements in $Sp(4)$, e.g.\ the \emph{spatial} CR3BP.

\begin{figure}
    \centering
    \includegraphics[width=0.85\linewidth]{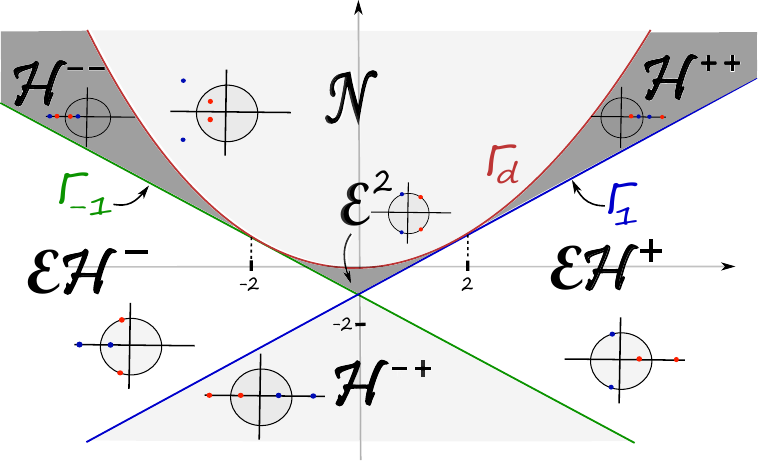}
    \caption{The 3D Broucke stability diagram, where $\Gamma_{\pm 1}$ corresponds to eigenvalue $\pm 1$, $\Gamma_d$ to double eigenvalue, $\mathcal{E}^2$ to doubly elliptic (the \emph{stable region}), and so on; see \cite{FM}.}
    \label{fig:Broucke_3D}
\end{figure}

Given a Wonenburger matrix $M=M_{A,B,C}\in Sp^\mathcal{I}(4)$, its \emph{stability point} is defined to be $$p=(\mbox{tr}(A),\det(A)) \in \mathbb R^2.$$ This point lies in the plane, which splits into regions corresponding to the eigenvalue configuration of $M$ (see Figure~\ref{fig:Broucke_3D}, which depicts the Broucke stability diagram for $n=2$). Each component of the plane is labelled according to the eigenvalue configuration. As an example, $\mathcal{E}^2$ (the \emph{doubly elliptic} component) corresponds to two pairs of elliptic eigenvalues; $\mathcal{EH}^{+}$ (the \emph{elliptic-positive hyperbolic} component), corresponds to one elliptic pair and a positive hyperbolic pair; $\mathcal{N}$, corresponds to complex quadruples, and so on. The parabola $\Gamma_d=\{y=x^2/4\}$ represents double eigenvalues, that is, when two eigenvalues come together. The lines $\Gamma_{\pm 1}$ tangent to $\Gamma_d$ and with slope $\pm 1$ corresponds to symplectic matrices with eigenvalue $\pm 1$ in their spectrum. 

\begin{figure}
    \centering
    \includegraphics[width=1.13\linewidth]{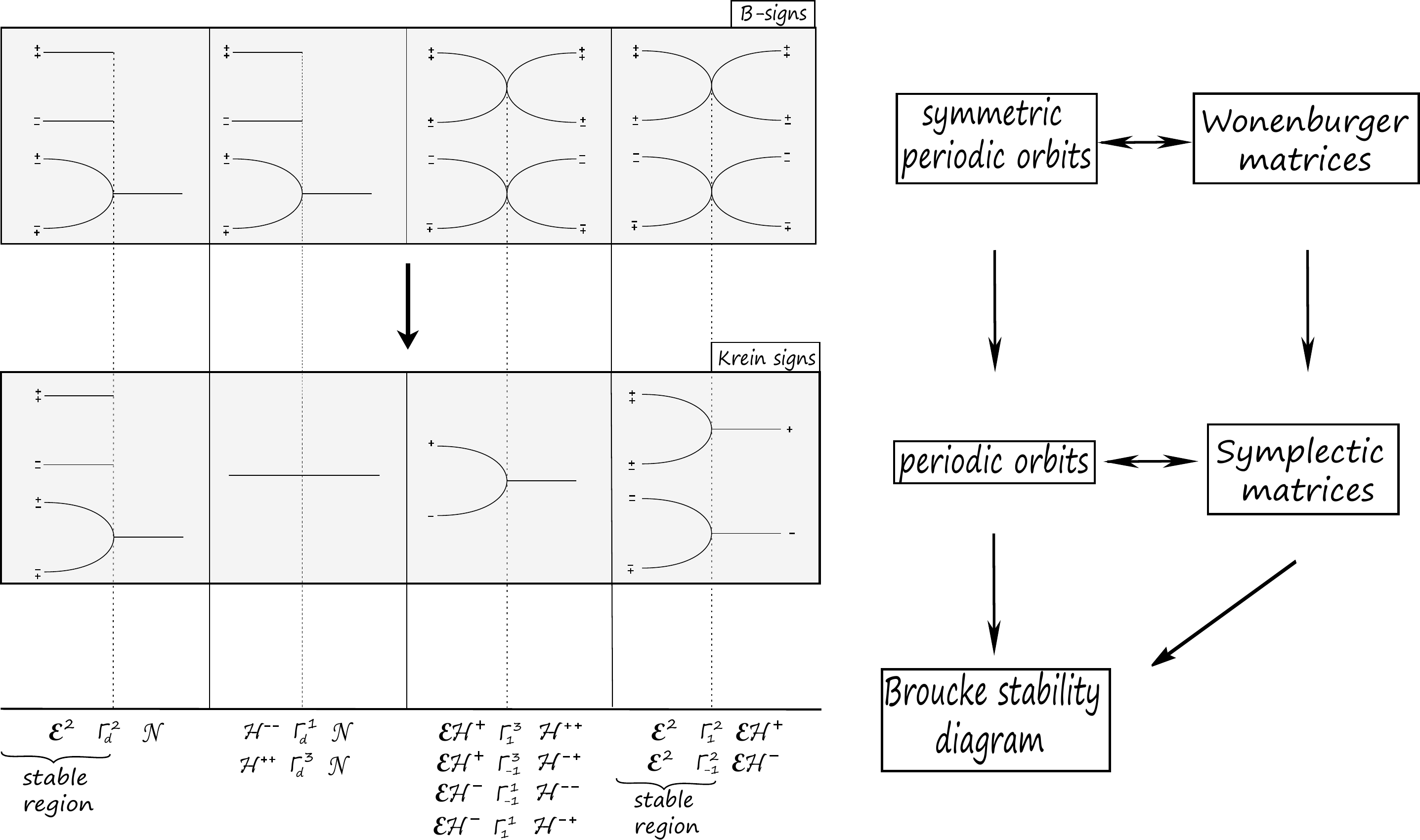}
    \caption{The branches (represented as lines) are two-dimensional, and come together at the 1-dimensional ``branching locus'' (represented as points), where we cross from one region to another of the Broucke diagram. The $1$-dimensional loci collapse to points over each of the three singular points $(2,1),(0,-1), (-2,1) \in \mathbb R^2$.}
    \label{fig:GIT_sequence_3D}
\end{figure}

The GIT sequence \cite{FM} adds two more layers to this diagram, as shown in Figure~\ref{fig:GIT_sequence_3D}. The top layer comes with two extra branches than the middle one, one for each hyperbolic eigenvalue. While the combinatorics and the global topology of the spaces involved is more complicated than the 2D case, the intuitive idea is still the same: namely, that the amount of information for symmetric orbits is much richer, and that one can distinguish more orbits up to qualitative equivalence. In this dimension we have two \emph{pairs} of eigenvalues, and so the $B$-signature is a pair $(\pm,\pm)$ of signs, so that the top layer has four branches over each component of the Broucke diagram (except the nonreal component, which has only one).

\subsection{Bifurcations in the Broucke diagram.} 
\begin{figure}
    \centering
\includegraphics[width=0.65\linewidth]{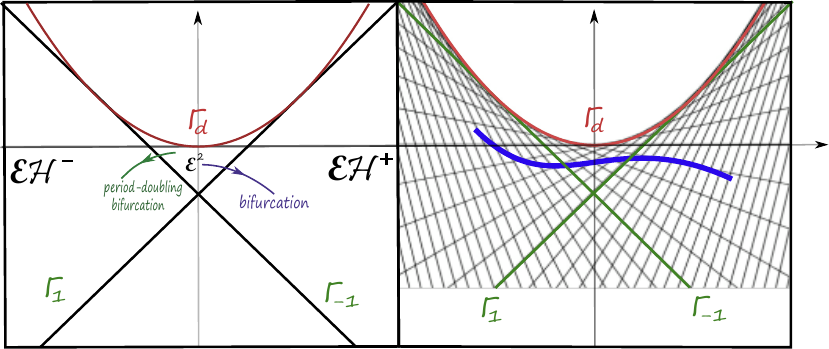}
\caption{Bifurcations are encoded by a pencil of lines.}
   \label{fig:complete_bifurcation}
\end{figure}

Any family of symmetric periodic orbits $c\mapsto \gamma_c$ yields a path $c\mapsto p_c\in \mathbb R^2$ of stability points in the plane. The family bifurcates if and only if $p_c$ crosses $\Gamma_1$. More generally, denote by $\Gamma^e_{\varphi}$ the line with slope $\cos(2\pi \varphi)\in [-1,1]$ tangent to $\Gamma_d=\{y=x^2/4\}$, which corresponds to matrices having eigenvalue $e^{2\pi i\varphi}$; and let $\Gamma^h_\lambda$ be the tangent line with slope $a(\lambda)\in \mathbb R\backslash [-1,1]$, corresponding to matrices with eigenvalue $\lambda$. Then a $k$-fold bifurcation occurs precisely when crossing $\Gamma^e_{l/k}$ for some $l$, i.e.\ the eigenvalue $e^{2\pi i l/k}$ is crossed in a family. In other words, higher order bifurcations are encoded by a pencil of lines tangent to a parabola, as shown in Figure~\ref{fig:complete_bifurcation}. Two such lines intersect at a point, which lies in a component determined by the lines (for instance, $\Gamma^e_\varphi \cap \Gamma^h_\lambda$ lies in $\mathcal{EH}^+$ if $\lambda>1$, and similarly for the other cases).
    
\subsection{Mathematical treatment} Now that we have given an idea of the salient features of the GIT sequence, we present some of the mathematical details. In order to deal with topological technicalities, we need to recall the definition of a GIT quotient, a notion borrowed from geometric invariant theory.

\begin{definition}[\textbf{GIT quotient}] Let $G$ be a group acting on a topological space $X$ by homeomorphisms. The \emph{GIT quotient} is the quotient space $X//G$ defined by $x\sim y $ if the closures of the $G$-orbits of $x$ and $y$ intersect, endowed with the quotient topology.     
\end{definition}
The condition $x\sim y$ is reflexive and symmetric, but not always transitive, but when it is, it defines an equivalence relation (all the examples considered below will satisfy this). By definition, it means that there exists sequences $g_n,h_n \in G$ such that $\lim_n g_n \cdot x =\lim_n h_n \cdot y \in X$. The point is that the naive quotient space $X/G$ might not be Hausdorff in general, and $X//G$ always is, in a universal way. Intuitively, to obtain the GIT quotient from the naive quotient, one gets rid of those bad points where the Hausdorff condition does not hold. 

Now we are ready to define the GIT sequence. This is the sequence of maps and spaces given by
$$
Sp^\mathcal{I}(2n)// GL_n(\mathbb R) \rightarrow Sp(2n)//Sp(2n) \rightarrow M_{n\times n}(\mathbb R)//GL_n(\mathbb R),
$$
$$
[M_{A,B,C}]\mapsto [M_{A,B,C}]\mapsto [A].
$$

Recall that the action of $GL_n(\mathbb R)$ on $Sp^{\mathcal{I}}(2n)$ is given by

\begin{equation*}
R_*\big(A,B,C\big)=\Big(RAR^{-1},RBR^t,(R^t)^{-1}CR^{-1}\Big),
\end{equation*}
whereas the remaining actions are via conjugation. We then see that the above maps are well-defined, by checking that $M_{A,B,C}$ and $M_{R_*(A,B,C)}$ are symplectically conjugated. We are also implicitly using Theorem \ref{thm:Wonenburger} to define the second map, which is independent of the Wonenburger representative. Moreover, for the base of the sequence, we have the following nice fact.

\begin{proposition}
    The base of the GIT sequence $M_{n\times n}(\mathbb R)//GL_n(\mathbb R)$ is homeomorphic to $\mathbb R^n$, where the homeomorphism maps a matrix to the coefficients of its characteristic polynomial, i.e.\
$$
M_{n\times n}(\mathbb R)//GL_n(\mathbb R)\rightarrow \mathbb R^n,
$$
    $$
    [A] \mapsto (c_{n-1},c_n,\dots,c_0),
    $$
    whenever
    $$
    p_A(t)=\det(A-t\cdot\mathds 1)=(-1)^nt^n+c_{n-1}t^{n-1}+\dots+c_0.
    $$
\end{proposition}

The map is well-defined as conjugated matrices have the same characteristic polynomial. Indeed, if $A\sim B$, by definition then there exists $D$, and a sequence of matrices $R_m, S_m \in GL_n(\mathbb R)$ such that 
$$
\lim_{m} R_m A R_m^{-1}= \lim_m S_m B S_m^{-1}=D.
$$
But then $p_A=\lim_m p_{R_m A R_m^{-1}}=p_D = \lim_m p_{S_m B S_m^{-1}}=p_B$. The basic observation for the proof is then that every matrix is conjugated to its Jordan form. And moreover, by a simple calculation, every Jordan block is equivalent, in the GIT quotient, to a diagonal matrix with the same eigenvalues. In other words, the passage from the naive quotient $M_{n\times n}(\mathbb R)/GL_n(\mathbb R)$ to the GIT quotient $M_{n\times n}(\mathbb R)//GL_n(\mathbb R)$ basically means ignoring Jordan blocks, and pretending that everything is diagonalizable.
Therefore each equivalence class has a canonical representative consisting of a real Jordan matrix with no super-diagonal entries. One then observes that the Jordan decomposition is determined by the characteristic polynomial, and that every polynomial with $(-1)^n$ as the leading term is the characteristic polynomial of some matrix. See Appendix A in \cite{FM} for further details. 

\begin{example}
    If $n=2$, the coefficients of the characteristic polynomial of $A\in M_{2\times 2}(\mathbb R)$ are $p=(\mathrm{tr}(A),\det(A))$, i.e.\ our definition of the stability point of $M=M_{A,B,C}\in Sp^{\mathcal{I}}(4).$
\end{example}

\begin{definition} In arbitrary dimension, we define the \emph{stability point} of the Wonenburger matrix $M_{A,B,C}\in Sp^{\mathcal I}(2n)$ as
$$
p=(s_1(\mu_1,\dots,\mu_n),\dots,s_n(\mu_1,\dots,\mu_n))\in \mathbb R^n,
$$
where $\mu_1,\dots,\mu_n$ are the eigenvalues of $A$, and $s_j$ is the $j$-th elementary symmetric polynomial, given by
$$
s_j(\mu_1,\dots,\mu_n)=\sum_{1\leq i_1 < \dots < i_j \leq n} \mu_{i_1}\dots\mu_{i_j}.
$$
\end{definition}

The stability point is then the result of applying the GIT sequence of maps to the given matrix.

\subsection{Formulas for GIT sequence: 2D} The low dimensional case of the GIT sequence can be explicitly described as follows. This case has also been studied in \cite{BZ}, where it plays an important role when trying to define a real version of Embedded Contact Homology (ECH).

The GIT sequence for $n=1$ is 
$$
Sp^{\mathcal{I}}(2)//GL_1(\mathbb R)\rightarrow Sp(2)//Sp(2) \rightarrow \mathrm{M}_{1\times 1}(\mathbb{R})//\mathrm{GL}_1(\mathbb{R})\cong \mathbb R,
$$
$$
[M_{A,B,C}]\mapsto [M_{A,B,C}]\mapsto [A]=A=\mathrm{tr}(M_{A,B,C})/2.
$$
The action of $GL_1(\mathbb{R})=\mathbb{R}^+$ on $\mathrm{Sp}^\mathcal{I}(2)$ is simply 
$$
\epsilon \cdot \left(\begin{array}{cc}
    a & b \\
    c & a
\end{array} \right)
=\left(\begin{array}{cc}
    a & \epsilon^2 b \\
    \frac{1}{\epsilon^2}c & a
\end{array} \right),$$ where $a^2-bc=1$, $\epsilon>0$. We have $\mathrm{Sp}(2)=SL(2,\mathbb{R})$, and a matrix $A\in \mathrm{Sp}(2)$ is either hyperbolic (i.e.\ $\vert \mathrm{tr}(A) \vert >2$, in which case it has two real eigenvalues $r,1/r$ with $\vert r \vert >1$), elliptic (i.e.\ $\vert \mathrm{tr}(A) \vert <2$, in which case it has two conjugate complex eigenvalues in the unit circle), or parabollic (i.e.\ $\vert \mathrm{tr}(A) \vert =2$, in which case it has eigenvalue $\pm 1$ with algebraic multiplicity two). From the discussion in \cite[Section 10.5]{FvK18}, we gather that $\mathrm{Sp}(2)//\mathrm{Sp}(2)$ admits a homeomorphism 
$$
\mathrm{Sp}(2)//\mathrm{Sp}(2)=\{z\in \mathbb{C}: \vert z\vert =1\}\cup \{r\in \mathbb{R}: \vert r\vert\geq 1\}\subset \mathbb{C},
$$
via the identification
$$
s(e^{i\theta})=\left( \begin{array}{cc}
\cos(\theta)     & -\sin(\theta) \\
\sin(\theta)     & \cos(\theta)
\end{array}\right), \; s(r)= \left( \begin{array}{cc}
r    & 0 \\
0     & \frac{1}{r}
\end{array}\right).
$$
The hyperbolic locus consists of closed orbits and corresponds to $\{\vert r \vert >1\}$; the elliptic locus also consists of closed orbits, and corresponds to $\{\vert z\vert=1\}\backslash\{\pm 1\}$; and the parabollic locus is $\{\pm 1\}$, where $\{+1\}$ corresponds to the three different Jordan forms with eigenvalue $1$ of algebraic multiplicity two, and, similarly $\{-1\}$ corresponds to the three Jordan forms with eigenvalue $-1$ of algebraic multiplicity two.

Similarly, the GIT quotient $\mathrm{Sp}^\mathcal{I}(2)//\mathrm{GL}_1(\mathbb{R})$ admits an identification \cite[Appendix B]{BZ}
$$
\mathrm{Sp}^\mathcal{I}(2)//\mathrm{GL}_1(\mathbb{R})=\{z\in \mathbb{C}:\vert z\vert=1\}\cup \{(\pm \cosh(u),\sinh(u)):u\in \mathbb{R}\}\subset \mathbb{C},
$$
via
$$
t(e^{i\theta})=s(e^{i\theta})=\left( \begin{array}{cc}
\cos(\theta)     & -\sin(\theta) \\
\sin(\theta)     & \cos(\theta)
\end{array}\right), \; t(u)= \left( \begin{array}{cc}
\pm \cosh(u)    & \sinh(u) \\
\sinh(u)     & \pm \cosh(u)
\end{array}\right).
$$
The matrix $t(u)$ has eigenvalues $\pm e^u, \pm e^{-u}$. Moreover, the matrices $t(u)$ and $t(-u)$ are both symplectically conjugate to diag$(\pm e^u,\pm e^{-u})$, hence to each other, and therefore define the same element in $\mathrm{Sp}(2)//\mathrm{Sp}(2)$.
After these identifications, the GIT sequence becomes
$$
e^{i\theta}\mapsto e^{i\theta} \mapsto \cos(\theta),
$$
$$
(\pm \cosh(u),\sinh(u))\mapsto r=\pm e^{\vert u \vert} \mapsto r=\pm e^{\vert u \vert}.
$$
    
See Figure \ref{fig:GIT_sequence}.

\subsection{Stability and the Krein--Moser theorem} We now discuss how the GIT sequence \emph{topologically} encodes the (linear) stability of periodic orbits, and compare it to the notions of Krein theory, and the classical Krein--Moser stability theorem. For the latter, we follow the exposition in Ekeland's book \cite{Eke90} (cf.\ \cite{Ab01}).

Consider a linear symplectic ODE, given by
$$
\dot x=M(t)x,
$$
where $M(t)=JA(t)$, with $A(t)$ a symmetric and $T$-periodic matrix, i.e.\ $A(t+T)=A(t)$ for all $t$, and where $J=\left(\begin{array}{cc}
    0 &  \mathds 1\\
    -\mathds 1 & 0
\end{array}\right)$ is the standard complex multiplication. The solutions are given by $x(t)=R(t)x(0)$, where
$R(t)\in \mathrm{Sp}(2n)$ is symplectic and solves the ODE $\dot R(t)=M(t)R(t)$, $R(0)=\mathds 1$. This type of linear ODE naturally arises when linearizing the Hamiltonian flow along a periodic orbit.

\smallskip

\begin{definition}\textbf{(stability)}
The ODE $\dot x=JA(t)x$ is \emph{stable} if all solutions remain bounded for all $t\in \mathbb{R}$. It is called \emph{strongly} stable if there is $\epsilon>0$ such that, if $B(t)$ is symmetric and satisfies $\Vert A(t)- B(t)\Vert <\epsilon$, then the corresponding ODE $\dot x=JB(t)x$ is stable. 

Similarly, a symplectic matrix $R$ is said to be \emph{stable} if all its iterates $R^k$ remain bounded for $k\in \mathbb{Z}$, and it is called \emph{strongly} stable if there is $\epsilon>0$ such that all symplectic matrices $S$ satisfying $\Vert R-S \Vert<\epsilon$ are also stable.    
\end{definition}

From Floquet theory, one obtains that the ODE $\dot x=JA(t)x$ is (strongly) stable if and only if $R(T)$ is (strongly) stable; see e.g.\ \cite[Section 2, Proposition 3]{Eke90}. Moreover, stability is equivalent to $R(T)$ being diagonalizable (i.e.\ all eigenvalues are semi-simple), with its spectrum lying in the unit circle \cite[Section 1, Proposition 1]{Eke90}.

Let $\{\lambda,\overline{\lambda}\}$ be an elliptic pair of eigenvalues of a symplectic matrix $R$. Then any other symplectic matrix close to $R$ will also have simple eigenvalues in the unit circle different from $\pm 1$. Therefore $R$ is strongly stable. The case of eigenvalues with higher multiplicity is the subject of Krein theory: whenever two elliptic eigenvalues come together, this gives a criterion for when they cannot possibly escape the circle and transition into a complex quadruple. This works as follows.

Consider the nondegenerate bilinear form $G(x,y)=x^t \cdot (-iJ) \cdot \overline{y}$ on $\mathbb{C}^{2n}$, associated to the Hermitian matrix $-iJ$; then every real symplectic matrix $R$ preserves $G$. Moreover, if $\lambda,\mu$ are eigenvalues of $R$ which satisfy $\lambda \overline{\mu}\neq 1$, then the corresponding eigenspaces are $G$-orthogonal, since we have
$$
G(x,y)=G(Rx,Ry)=\lambda \overline{\mu}G(x,y), 
$$
if $x,y$ are the corresponding eigenvectors. Moreover, if we consider the generalized eigenspaces 
$$
E_\lambda=\bigcup_{m\geq 1}\ker(R-\lambda I)^m,
$$
then $E_\lambda,E_\mu$ are also $G$-orthogonal if $\lambda \overline{\mu}\neq 1$ \cite[Section 2, Proposition 5]{Eke90}. This implies that if $\vert\lambda\vert \neq 1$, then $E_\lambda$ is $G$-isotropic, i.e.\ $G\vert_{E_\lambda}=0$. If $\sigma(R)$ denotes the spectrum of $R$, we then have a $G$-orthogonal decomposition
$$
\mathbb{C}^{2n}=\bigoplus_{\substack{\lambda \in \sigma(R)\\\vert \lambda \vert\geq 1}}F_\lambda,
$$
where $F_\lambda=E_\lambda$ if $\vert \lambda\vert=1$, and $F_\lambda=E_\lambda \oplus E_{\overline{\lambda}^{-1}}$ if $\vert\lambda\vert >1$. Since $G$ is non-degenerate, and this splitting is $G$-orthogonal, the restriction $G_\lambda=G\vert_{F_\lambda}$ is also non-degenerate. Moreover, if $\vert \lambda \vert \neq 1$, with algebraic multiplicity $d$, then the $2d$-dimensional space $F_\lambda$ has $E_\lambda$ as a $d$-dimensional isotropic subspace, and therefore the signature of $G_\lambda$ is $(d,d)$. But if $\vert\lambda\vert=1$, then $G_\lambda$ can have any signature, which justifies the following notion.

\begin{definition}\textbf{(Krein-positivity/negativity)}
If $\lambda$ is an eigenvalue of the symplectic matrix $R$ with $\vert \lambda \vert=1$, then the signature $(p,q)$ of $G_\lambda$ is called the \emph{Krein-type} or \emph{Krein signature} of $\lambda$. If $q=0$, i.e.\ $G_\lambda$ is positive definite, $\lambda$ is called \emph{Krein-positive.} If $p=0$, i.e.\ $G_\lambda$ is negative definite, $\lambda$ is called \emph{Krein-negative.} If $\lambda$ is either Krein-negative or Krein-positive, it is said to be \emph{Krein-definite}. Otherwise, it is \emph{Krein-indefinite}.
\end{definition}

If $\lambda$ is of Krein-type $(p,q)$, then $\overline{\lambda}$ is of Krein-type $(q,p)$ \cite[Section 2, Lemma 9]{Eke90}. If $\lambda$ satisfies $\vert \lambda \vert=1$ and it is not semi-simple, then it is Krein-indefinite \cite[Section 2, Proposition 7]{Eke90}. Moreover, $\pm 1$ are always Krein-indefinite if they are eigenvalues. The following, originally proved by Krein in \cite{Kre1,Kre2,Kre3,Kre4} and independently rediscovered by Moser in \cite{M78}, gives a characterization of strong stability in terms of Krein signatures:
\begin{thm}[\textbf{Krein--Moser}]\label{Kreinthm}
 $R$ is strongly stable if and only if it is stable and all its eigenvalues are Krein-definite. 
\end{thm}
See \cite[Section 2, Theorem 3]{Eke90} for a full proof. This generalizes the case where all eigenvalues are simple, different from $\pm 1$ and in the unit circle, discussed above. 

The way that the GIT sequence relates to Krein theory is the following.

\begin{proposition}\label{prop:Krein_vs_B}
    For a Wonenburger matrix, the Krein signature coincides with the $B$-signature, for the case of elliptic eigenvalues.
\end{proposition}

\begin{example}
    As a simple example, to illustrate Proposition \ref{prop:Krein_vs_B}, let us consider the Wonenburger matrices
$$M=\left(\begin{array}{cccc}
\cos \theta &      0       & -\sin \theta& 0\\
0             & \cos\theta &            0   & -\sin \theta\\
\sin \theta&     0        & \cos \theta  & 0\\
0             & \sin\theta&            0   & \cos \theta
\end{array}\right), 
N=\left(\begin{array}{cccc}
\cos \theta &      0       & \sin \theta & 0\\
0             & \cos\theta &            0   & -\sin \theta\\
-\sin \theta &     0        & \cos \theta & 0\\
0             & \sin\theta &            0   & \cos \theta
\end{array}\right),$$ 
$$
P=\left(\begin{array}{cccc}
\cos \theta&      0       & \sin \theta & 0\\
0             & \cos\theta &            0   & \sin \theta\\
-\sin \theta &     0        & \cos \theta  & 0\\
0             & -\sin\theta &            0   & \cos \theta
\end{array}\right),$$

with $\theta \in (0,\pi)$. These are the three normal forms for $4\times 4$ Wonenburger matrices which are doubly elliptic and have eigenvalues $e^{\pm i \theta}$ of multiplicity two, see \cite{FM}. For the matrix $M$, the eigenspace for the eigenvalue $e^{i \theta}$ is spanned by the two eigenvectors
$$
v_1=(1,0,i,0),\; v_2=(0,1,0,i),
$$
and the eigenspace of $e^{-i \theta}$, by their corresponding conjugates 
$$
w_1=(1,0,-i,0),\; w_2=(0,1,0,-i).
$$
With $G=\left(\begin{array}{cc}
   0  &  -i\mathds 1\\
   i\mathds 1  & 0
\end{array} \right)$, we can compute that
$$
G(w_1,w_1)=G(w_2,w_2)=-G(v_1,v_1)=-G(v_2,v_2)=2,
$$
$$
G(w_1,w_2)=G(v_1,v_2)=0.
$$
Then the Krein signature of $e^{i \theta}$ is $(0,2)$, and that of $e^{-i\theta}$ is $(2,0)$. By inspection, we see that these coincide with the $B$-signatures of these eigenvalues. The remaining matrices are dealt with completely analogously.   
\end{example}

As a corollary of the Krein--Moser theorem and of Proposition \ref{prop:Krein_vs_B}, we obtain the following result.

\begin{thm}
 Let $R$ be a Wonenburger matrix. Then $R$ is strongly stable if and only if it is stable and all its eigenvalues are $B$-definite. 
\end{thm}

The Krein--Moser theorem is therefore detected topologically by the GIT sequence in a very simple, visual manner. Indeed, e.g.\ for the case $n=2$, if the stability point lies in the interior of the doubly elliptic region $\mathcal{E}^2$, then the corresponding (equivalence class of) matrices are strongly stable. This follows e.g.\ because they cannot be perturbed away from $\mathcal{E}^2$, but also because of the Krein--Moser theorem. At the boundary of $\mathcal{E}^2$ is where things becomes interesting, and in particular, along the boundary component of $\mathcal{E}^2$ lying in the parabola $\Gamma_d$. As shown in Figure \ref{fig:GIT_sequence_3D}, the $++$ and the $--$ branches over $\mathcal{E}^2$, in either the middle layer $Sp(4)//Sp(4)$ or the top layer $Sp^{\mathcal{I}}(2n)/GL_n(\mathbb R)$, do \emph{not} cross from $\mathcal{E}^2$ to $\mathcal{N}$. This means that the boundary of these branches over the correspnding portion of $\Gamma_d$ corresponds to two elliptic eigenvalues coming together, but such that the corresponding matrices are \emph{also} strongly stable. This is compatible with the Krein--Moser theorem, as the $B$-signature or equivalentely the Krein signature is positive (respectively negative) definite along the $++$ (respectively the $--$) branch. Note here that we need to use the full definition of $B$-signature rather than of $B$-sign, as along $\mathcal{E}^2\cap \Gamma_d$ the eigenvalues are no longer simple.

Moreover, the same phenomenon happens when two positive/negative-hyperbolic eigenvalues of a Wonenburger matrix come together. Indeed, from Figure \ref{fig:GIT_sequence_3D}, we detect that we cannot cross from $\mathcal{H}^{\pm}$ to $\mathcal{N}$ along the $++$ or $--$ branch of the top layer. In other words, we obtain the following as a corollary.

\begin{proposition}
Consider a Wonenburger matrix $M\in Sp^\mathcal{I}(4)$ with a hyperbolic eigenvalue of multiplicity $2$. Then $M$ cannot be perturbed to a Wonenburger matrix with a complex quadruple if and only if its $B$-signature is definite. 
\end{proposition}

\begin{remark} Note that the GIT sequence is designed to study \emph{linear} stability, and not the stronger \emph{non-linear} notion of stability for periodic orbits (i.e.\ trajectories that start near the given orbit stay near the orbit for all times), sometimes called \emph{Lyapunov} stability. Moreover, as Jordan blocks are trivial in the GIT quotient, the GIT sequence does not distinguish between linear stability, and the weaker notion of \emph{spectral} stability (i.e.\ that eigenvalues lie in the unit circle, but are not necessarily semi-simple).
\end{remark}

\section{The Floer numerical invariants}

\begin{figure}
    \centering
    \includegraphics[width=0.42\linewidth]{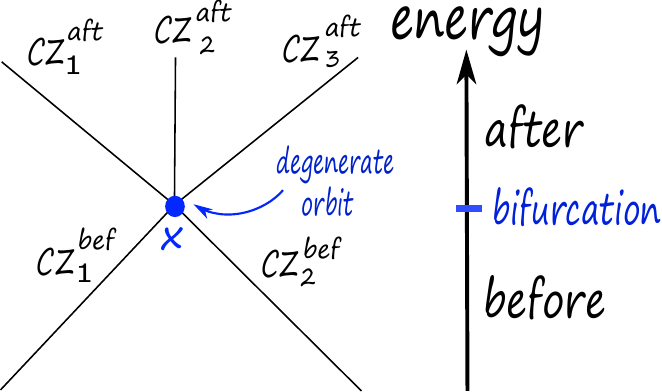}
    \caption{A sketch of a bifurcation at a degenerate orbit, with the before/after orbits determined by the deformation parameter (the energy), each branch with its own CZ-index. The Floer number is a signed count of orbits which stays invariant.}
    \label{fig:bifurcation_sketch}
\end{figure}

Recall that bifurcations occurs when studying families $t\mapsto x_t$ of periodic orbits, as a mechanism by which at some parameter time $t=t_0$ the orbit $x_{t_0}$ becomes degenerate, and several new families may bifurcate out of it; see Figure~\ref{fig:bifurcation_sketch}. The Floer numbers are meant to give a simple test to keep track of all new families. 

\subsection{Floer numbers} We will first need the following technical definition: a periodic orbit $x=y^k$, where $y$ is its underlying simple orbit, is \emph{bad} if $k$ is even and $$\mu_{CZ}(x)\neq \mu_{CZ}(y) \mbox{ mod }2.$$ Otherwise, it is \emph{good}. In fact, a planar orbit is bad iff it is an even cover of a negative hyperbolic orbit. And a spatial orbit is bad iff it is an even cover of either an elliptic-negative hyperbolic or a positive-negative hyperbolic orbit. Note that a good planar orbit can be bad if viewed in the spatial problem.

\begin{definition}[\textbf{Floer number}]
Given a bifurcation at $x$, the \emph{SFT-Euler characteristic} (or the \emph{Floer number}) of $x$ is

$$
\chi(x)=\sum_i (-1)^{CZ^{bef}_i}=\sum_j (-1)^{CZ^{aft}_j}. 
$$    

The sum on the LHS is over \textbf{good} orbits \emph{before} bifurcation, and RHS is over \textbf{good} orbits \emph{after} bifurcation.     
\end{definition}

The Floer number is simply the Euler characteristic of local Floer homology. As these numbers only involve the parity of the CZ-index, one has simple formulas which bypass the computation of this index, as they only involve the Floquet multipliers:

\medskip

\begin{itemize}
    \item \textbf{Planar case.} $\chi(x)=\#\big\{\textrm{ good}\,\,\mathcal H^{+}\big\}-\#\big\{\mathcal E,\,\mathcal{H}^-\big\}.$

\medskip
    
\item \textbf{Spatial case.} $\chi(x)=\#\big\{\mathcal H^{--}\,,\mathcal{EH}^-\,,\mathcal E^2\,,\textrm{ good}\,\,\mathcal H^{++}\,,\mathcal N\big\}-\#\big\{\mathcal H^{-+},\,\textrm{good}\,\,\mathcal{EH}^+\big\}.$
\end{itemize}

\medskip

Here, $\mathcal{E}$ denotes \emph{elliptic}, $\mathcal{H}^\pm$ denotes \emph{positive/negative hyperbolic}, and $\mathcal{N}$ denotes \emph{complex quadruples} $\lambda,1/\lambda,\overline{\lambda},1/\overline{\lambda}$. The above simply tells us which type of orbit comes with a plus or a minus sign (the formula should be interpreted as either before or after).

The fact that the sums agree before and after --\emph{invariance}-- follows from invariance of local Floer homology. In practice, when carrying out numerical work, the Floer number can be used as a \textbf{test}: 

\medskip

\textbf{Motto.} \emph{If the sums do \textbf{not} agree, we know the algorithm missed an orbit}.

\begin{example} Let $t\mapsto \gamma_t$ be a family of planar periodic orbits in the spatial CR3BP. We claim that it cannot go from doubly elliptic to complex quadruple. Indeed, planar orbits cannot have non-real quadruples as Floquet multipliers, as the CR3BP admits the reflection along the ecliptic as symmetry. Then if the family becomes complex quadruple, it has to become spatial. But then applying the symmetry we obtain a new family, also complex quadruple. As there is no bifurcation in the family (the eigenvalue $1$ is never crossed), the (spatial) Floer number of the planar orbit before becoming spatial is $1$, but it jumps to $2$ after. This is a contradiction to invariance. 
\end{example}

For the above example, we remark that the Floer number is also defined for orbits which are non-degenerate, although it always coincides with $\pm 1$, depending on the orbit type.

\subsection{Real Floer number} The invariant above works for arbitrary periodic orbits. There is a similar Floer invariant for \emph{symmetric} orbits \cite{FKM}, given by the Euler characteristic of the local Lagrangian Floer homology of the fixed-point locus of the corresponding involution. This works as follows.

Recall that a symmetric periodic $x$ can be seen both as a periodic orbit, as well as a chord between the Lagrangian fixed-point locus of the involution. Therefore, it has a CZ-index $\mu_{CZ}(x)$ and a Lagrangian Maslov index $\mu_L(x)$, a \emph{half}-integer, i.e.\ taking values in $\frac{1}{2}\mathbb{Z}$. The difference of these two, as introduced in \cite{FvK14}, is the \emph{Hörmander index} $$s(x) = \mu_{CZ}(x) - \mu_L(x)\in \frac{1}{2}\mathbb{Z},$$ also a half-integer. One can use this index to detect when $x$ bifurcates as a chord, even when it does not bifurcate as an orbit.

Before or after a bifurcation of a symmetric orbit, one obtains a collection of non-degenerate symmetric orbits for which one may compute the parity of the Maslov index $\mu_L(x)$. By this, if $\mu_L(x)=\frac{1}{2}m_L(x)$, we mean the parity of $m_L(x)\in \mathbb Z$. 

\begin{definition}[\textbf{Real Floer number}]

The real Euler characteristic $\chi_L(x)$ is then defined as
$$
\chi_L(x)=\sum_j(-1)^{\mu_L(x_j)}=\sum_ii^{m_L(x_j)} \in \mathbb{C},
$$
where the sum runs over the collection $x_j$ of non-degenerate chords arising after perturbation of $x$ (e.g.\ before or after a bifurcation). 

\end{definition}

Note that, $\chi_L(x)$ is complex-valued by definition. Its invariance under bifurcation follows from invariance of the local Lagrangian Floer homology of $x$.

\subsection{Formulas for computation} In order to compute the real Floer numbers, we need only to know the parity of the Maslov index. We now give useful formulas for this.

The Chebyshev polynomials of the first kind are recursively defined by
$$
T_0(x) = 1
$$
$$
T_1(x) = x
$$
$$
T_{k+1}(x) = 2xT_k(x) - T_{k-1}(x).
$$
The Chebyshev polynomials of the second kind are similarly defined by
$$
U_0(x) = 1
$$
$$
U_1(x) = 2x
$$
$$
U_{k+1}(x) = 2xU_k(x) - U_{k-1}(x).
$$
The following gives a formula for computing the Hörmander index of the iterates of a symmetric orbit in terms of the monodromy matrix.

\begin{thm}[\cite{FvK14}]\label{thm:sindex} Let $x$ be a nondegenerate, symmetric periodic orbit with monodromy
$$
M=M_{A,B,C}=\left(\begin{array}{cc}
A & B\\
C & A^T
\end{array}\right),
$$
a Wonenburger matrix. Then the Hörmander indices of its iterates are given by
\begin{equation}\label{eq:sindex}
s(x^k)= \frac{1}{2}\mathrm{sign}\left((Id -T_k(A))U_{k-1}(A)^{-1}C^{-1}\right),
\end{equation}
$k \in \mathbb N$. For $k=1$, we have in particular that
\begin{equation}\label{eq:sindex1}
s(x) = \frac{1}{2}\mathrm{sign}\left((Id-A)C
^{-1}\right).
\end{equation}
\end{thm}
Here, $\mathrm{sign}$ denotes the signature. Note that the covers of a symmetric orbit are also symmetric. We have also used the fact that $C$ is invertible if $x$ is non-degenerate \cite[Lemma 3.2]{FvK2}. 

\medskip

Note that, the parity of the Maslov index can be determined from the following: 
\begin{itemize}

\medskip

    \item The monodromy matrix;

\medskip

    \item the formula $\mu_L(x)=\mu_{CZ}(x)-s(x)$;

\medskip

    \item Formula (\ref{eq:sindex}) (or (\ref{eq:sindex1})) of Theorem \ref{thm:sindex}, which in particular gives the parity of $s$;

\medskip
    
    \item The parity of the CZ-index.
\end{itemize}

\begin{figure}
    \centering
    \includegraphics[width=1\linewidth]{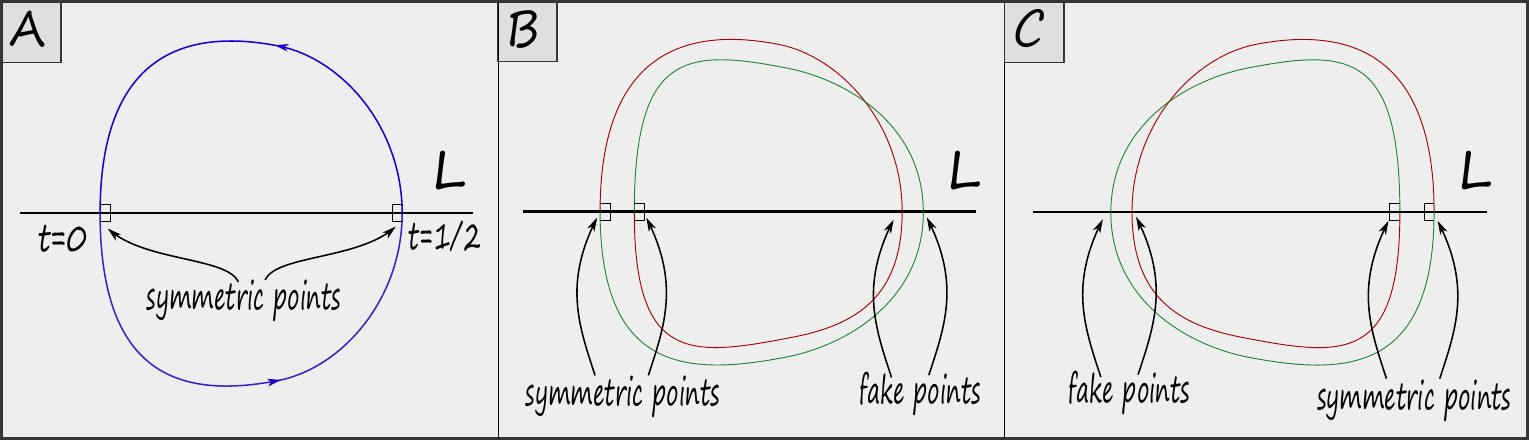}
    \caption{Symmetric period doubling bifurcation. The \emph{fake} symmetric points, while close to intersection points, do \emph{not} intersect the fixed-point loci.}
    \label{fig:period_doubling}
\end{figure}

\begin{example}[\textbf{Symmetric period doubling bifurcation}] \label{ex:sym_double}We finish this section with an example where our invariants give new information. Consider a symmetric orbit $x$ going from elliptic to negative hyperbolic. A priori there could be two bifurcations, one near each symmetric point (B or C in Figure~\ref{fig:period_doubling}). However, invariance of $\chi(x^2)$ implies only \emph{one} can happen (note $x^2$ is \emph{bad}). And where the bifurcation happens is determined by the $B$-sign, occurring at the symmetric point in which the $B$-sign does \emph{not} jump; or alternatively, where the $C$-sign jumps (recall they differ for the elliptic case, and agree for the hyperbolic case). Near the symmetric point where bifurcation does not occur, it splits into two points which, although close to the fixed-point locus, do not lie in it, i.e.\ they are \emph{fake symmetric}.

As a simple orbit, there is no bifurcation since there is no eigenvalue $1$ in the reduced monodromy matrix. However, we can interpret this orbit as a chord from the Lagrangian to itself, where  $x(0)$ happens to agree with $x(1)$; now, as a \emph{chord}, it might bifurcate.\footnote{Such a chord bifurcation takes place if and only if the Lagragian Maslov index jumps, which happens if and only if the Hörmander index jumps.} If this happens, we can apply the symmetry again to the red chord in Figure \ref{fig:period_doubling} (B)/(C) to obtain the green chord in the same figure. So, \emph{two} chords bifurcate. This is compatible with the real Euler characteristic. Indeed, the Lagrangian Maslov index of $x$ before and after bifurcation (thought of as a chord) differ by one. The green and red chords have the same Maslov index, say $k$, as they are symmetric to each other, and this coincides with that of $x$ before bifurcation. This makes sure that the real Euler characteristic stays invariant. Indeed, before bifurcation we have $\mu_{L}(x)=k$ and so $\chi_L(x)=(-1)^k$; and after bifurcation, $\mu_{L}(x)=k+1$, so $\chi_L(x)=2(-1)^k+(-1)^{k+1}=2(-1)^k-(-1)^k=(-1)^k$, which explicitly shows invariance.
\end{example}

\section{Digression: GIT sequence and the associahedron}\label{sec:GIT_any_dim} We will now further delve into the details of the GIT sequence, in arbitrary dimension. This section is based on \cite{MR23}.

We will study the topology of the GIT quotients in the sequence, that is, we will determine their branching structure. In particular, we will find normal form representatives in $Sp^{\mathcal{I}}(2n)//GL_n(\mathbb R)$, i.e.\ of Wonenburger matrices up to the natural action of $GL_n(\mathbb R)$, and therefore compute the degrees of the maps in the GIT sequence, which vary over the the components of a suitable decomposition of the base $\mathbb R^n$, each of them labelled by eigenvalue configurations. Among these components, there is a special one, the \emph{stable} component, which corresponds to stable periodic orbits. We will show that its combinatorics is governed by a quotient of the \emph{associahedron}. This will be a generalization to arbitrary dimension of the results in \cite{FM}, which have appeared in \cite{MR23} (see also Howard--Mackay \cite{HM87}, Howard--Dullin \cite{HD98}), and will be rather technical, so can be skipped on a first read.

\medskip

\subsection{Some real algebraic geometry} We shall consider the space of monic polynomials with real coefficients, and of degree $n$, that is, of the form 
$$
p(t)=(-1)^nt^n +c_{n-1}t^{n-1}+\dots + c_0,
$$
with $c_i\in \mathbb R$. We identify this space with $\mathbb R^n$ via
$$
p \longleftrightarrow (c_0,\dots,c_{n-1}).
$$
Recall that the \emph{discriminant} of a polynomial is defined as 
$$
\Delta(p)=\prod_{i<j} (\lambda_i-\lambda_j)^2,
$$
where $\lambda_1\dots,\lambda_n$ are the complex roots of $p$. Then by definition $p$ has a multiple root if and only if $\Delta(p)=0$. If these roots are all real and simple, then $\Delta(p)>0$. Moreover, as $\Delta(p)$ is a symmetric polynomial in the roots of $p$, it follows from the fundamental theorem of symmetric polynomials that it is a polynomial in the basic symmetric polynomials of the roots, i.e.\ a polynomial in the coefficients $c_i$ of $p$. Then
$$
V(\Delta)=\{\Delta=0\}\subset \mathbb R^n
$$
is a (in general, singular) real algebraic variety, which we call the \emph{discriminant variety}. For instance, the discriminant of a degree $2$ polynomial $p(x)=x^2+bx+c$ is $\Delta(p)=b^2-4c$, and so $V(\Delta)=\{b^2=4c\}$ is a parabola in $\mathbb R^2$, and that of a degree $2$ polynomial $p(x)= -x^3+bx^2+cx+d$ is $\Delta(p)=b^2c^2+4c^3-4b^3d-27d^2-18bcd,$ so that $V(\Delta)$ is a quadric hypersurface in $\mathbb R^3$. 

For $a\in \mathbb R$, we denote $\Gamma_a=\{p: p(a)=0\}\subset \mathbb R^n$. Note that the equation $p(a)=0$ is linear in the coefficients of $p$, and so $\Gamma_a$ is a linear hyperplane in $\mathbb R^n$, given by $\Gamma_a=\{P_a=0\}$ for some linear equation $P_a$. Similarly, for $\alpha\in \mathbb C\backslash \mathbb R$, we let $V_\alpha=\{p: p(\alpha)=0\}\subset \mathbb R^n$. Note that $V_\alpha=V_{\overline{\alpha}}$, as complex roots come in conjugate pairs. We then define the regions of $\mathbb R^n$ given by
$$
\mathcal{E}=\bigcup_{a\in[-1,1]}\Gamma_a,\; \mathcal{H}^+=\bigcup_{a\geq 1}\Gamma_a,\; \mathcal{H}^-=\bigcup_{a\leq -1}\Gamma_a, \;\mathcal{N}=\bigcup_{\alpha \in \mathbb C\backslash \mathbb R} V_\alpha,
$$
which we respectively call the \emph{elliptic, positive/negative hyperbolic, nonreal} regions. The intersection $\Gamma_{aa'}=\Gamma_a\cap \Gamma_{a'}$ of two hyperplanes $\Gamma_a,\Gamma_{a'}$ for $a\neq a'$ is a codimension-$2$ affine linear subspace of $\mathbb R^n$ corresponding to polynomials $$p(x)=(x-a)(x-a')q(x)$$ for some $q\in \mathbb R^{n-2}$. Similarly, $\Gamma_{a_1\dots a_n}=\Gamma_{a_1}\cap \dots \cap \Gamma_{a_n}$ is a point, corresponding to the polynomial
$$
p(x)=(-1)^n(x-a_1)\dots (x-a_n).
$$
Moreover, $V_{\alpha_1\dots \alpha_n}=V_{\alpha_1}\cap\dots\cap V_{\alpha_n}$ is a point. If we let $a'\rightarrow a$, then $\Gamma_{aa'}$ converges to $\Gamma_a$, but also becomes tangent to the discriminant variety along the subspace $$V_a(\Delta)=\{(x-a)^2q(x): q \in \mathbb R^{n-2}\}=V(\Delta)\cap \Gamma_a\subset V(\Delta).$$ 
The \emph{regular} points of $V(\Delta)$ are those points $(x-a)^2q(x)$ which satisfy $q(a)\neq 0$, and so $V(\Delta)$ has a well-defined tangent space at these points. By comparing dimensions, we obtain that $\Gamma_a$ coincides with this tangent space at each of the regular points. Moreover, $\Gamma_a$ is tangent to $V(\Delta)$ along $V_a(\Delta)$. Varying $a$, we see that the discriminant variety is the \emph{envelope} of the family of hyperplanes $\{\Gamma_a\}_{a\in \mathbb R}$, i.e.\ it is the locus of points $V(\Delta)=\bigcup_a V_a(\Delta)$ tangent to each $\Gamma_a$, and coincides with the union of these tangency loci (the \emph{characteristic subspaces} $V_a(\Delta)$). The envelope $V(\Delta)$ can therefore be described as the set of coefficients $c\in \mathbb R^n$ satisfying
$$
P_a(c)=\partial_aP_a(c)=0
$$
for some $a$.

We also have a decomposition of $\mathbb R^n$ into regions labelled by the root configurations of the corresponding polynomials. Namely, for $k,l,m,r$ with $k+l+m+2r=n$, we define
$$
\mathcal{M}^{klmr}_0:=(\mathcal{H}^-)^k\mathcal{E}^l(\mathcal{H}^+)^m\mathcal{N}^r:=\bigcup_{\substack{a_1,\dots,a_l\in [-1,1]\\ b_1,\dots, b_m\geq 1 \\ c_1,\dots, c_k\leq -1\\ \alpha_1,\dots, \alpha_r\in \mathbb C\backslash\mathbb R}}\Gamma_{a_1\dots a_l}\cap \Gamma_{b_1\dots b_m} \cap \Gamma_{c_1 \dots c_k} \cap V_{\alpha_1\dots\alpha_r}.
$$
Then we have
$$
\mathbb R^n=\bigcup_{\substack{k,l,m,r\geq 0\\ k+l+m+2r=n}}\mathcal{M}^{klmr}_0.
$$
See Figures \ref{fig:Broucke_3D} and \ref{fig:bifurcation_sketch} for the case $n=2$ ($V(\Delta)$ is denoted $\Gamma_d$ there). 

Among the regions in this decomposition, there is a unique compact one, given by
$$
\mathcal{E}^n:=(\mathcal{H}^-)^0\mathcal{E}^n(\mathcal{H}^+)^0\mathcal{N}^0.
$$
This is the \emph{stable} region. Below, we will give a combinatorial description of this component.

\begin{figure}
    \centering
    \includegraphics[width=0.8\linewidth]{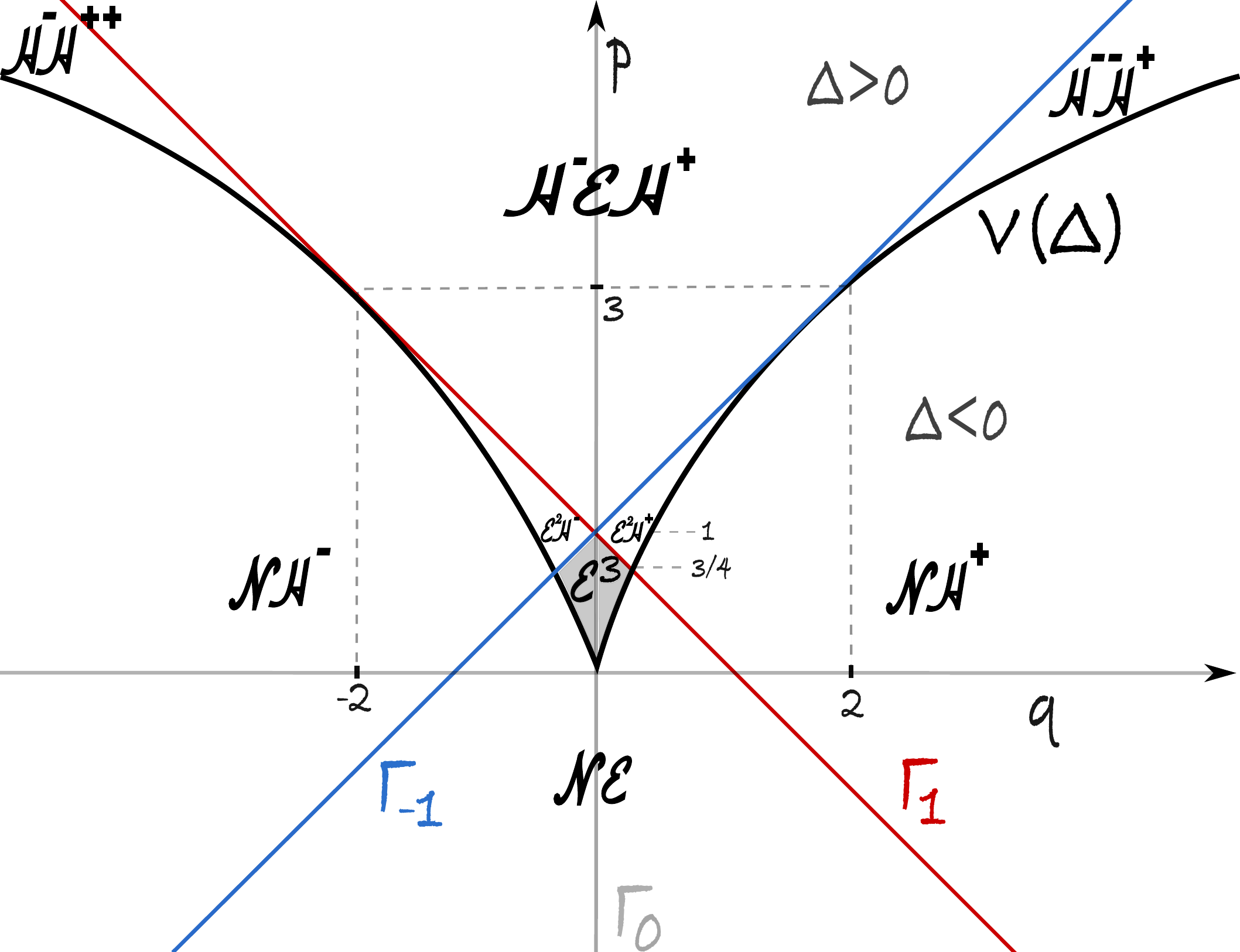}
    \caption{The stability diagram for depressed cubics.}
    \label{fig:depressed_diagram}
\end{figure}

\begin{example} For the case $n=3$, every polynomial
$$
A(x)=-x^3+bx^2+cx+d
$$  
can be transformed via the change of variables $y=x-b/3$ to a polynomial with no degree $2$ term (which we call a \emph{depressed} cubic polynomial), that is, of the form
$$
B(y)=-y^3+py+q
$$
with $p=-c-\frac{b^2}{3},\; q=-\frac{2b^3+9bc}{27}-d$. Its discriminant is 
$$
\Delta(B)=4p^3-27q^2.
$$
We identify the space of depressed cubics with $\mathbb R^2$ via $B\leftrightarrow (q,p)$, and so the vanishing locus of $\Delta$ is a cubic in the plane singular at the origin (corresponding to the polynomial $y^3$, and so having triple root). The decomposition of the plane according to the root configuration is depicted in Figure \ref{fig:depressed_diagram}. There are no $(\mathcal{H}^{+})^3$ or $(\mathcal{H}^{-})^3$ components, as no depressed cubic can have roots which are all greater than $1$ in absolute value. The pencil of lines $\{\Gamma_a\}_{a\in \mathbb R}$ has $V(\Delta)$ as envelope, with $\Gamma_a=\{P_a=-a^3+pa+q=0\}$ of slope $-1/a$, corresponding to polynomials with root $a$. The discriminant variety is the loci of points satisfying 
$$
P_a=-a^3+pa+q=\partial_aP_a=-3a^2+p=0
$$
for some $a$.

\end{example}

\subsection{The GIT sequence}

Recall that the GIT sequence is by definition the sequence of maps and spaces given by
$$
Sp^\mathcal{I}(2n)//GL_n(\mathbb R)\rightarrow Sp(2n)//Sp(2n)\rightarrow M_{n\times n}// GL_n(\mathbb R)=\mathbb R^n
$$
$$
[M_{A,B,C}]\mapsto [M_{A,B,C}]\mapsto [A].
$$
We will now study the topology and combinatorics of these spaces and maps.

\medskip

\textbf{Regular cases.} We start with the regular cases, corresponding to the case when the matrix has only simple eigenvalues different from $\pm 1$. Let $\mathcal{M}^{klmr}$ be the stratum of matrices $[M]\in Sp^\mathcal{I}(2n)//GL_n(\mathbb R)$ with $M=M_{A,B,C}$, and with the $A$-block having an eigenvalue configuration with
$$
\mu^-_1< \dots < \mu^-_k< -1< \nu_{1}< \dots < \nu_l < 1 < \mu^+_{1}<\dots<\mu_m^+ 
$$
the real eigenvalues of $A$, and 
$$
\alpha_1,\overline{\alpha}_1,\dots, \alpha_r,\overline{\alpha}_r
$$
the complex conjugate eigenvalues of $A$, which satisfy
$$
k+l+m+2r=n.
$$

Then there are unique $\beta^\pm_j, r_j \in (0,+\infty)$ and $\theta_j,\gamma_j\in (0,\pi)$ such that
$$
\mu_j^-=-\cosh \beta_j^-,\;\mu_j^+=\cosh \beta_j^+,\; \nu_j=\cos \theta_j,\;\alpha_j=r_je^{i\gamma_j},\;\overline{\alpha}_j=r_je^{-i\gamma_j}. 
$$ 
Moreover, up to conjugation $A$ can be assumed to split as direct sum
$$
A=\mathrm{diag}\left(\mu^-_1,\dots,\mu^-_k,\nu_{1},\dots,\nu_l,\mu^+_{1},\dots,\mu_m^+\right)\oplus\bigoplus_{j=1}^rR(r_j,\gamma_j), 
$$
where 
$$
R(r,\gamma)=\left(\begin{array}{cc}
   r\cos \gamma & -r\sin \gamma \\
   r\sin \gamma  & r\cos \gamma
\end{array} \right)
$$
is the composition of a dilation with a rotation. We then see that each summand can be treated separately, and so we can appeal to \cite{FM}. The result is the normal form
\begin{equation}
    \begin{split}
B= &\; \mathrm{ diag}\left(\epsilon^-_1\sinh(\beta_1^-),\dots,\epsilon^-_k\sinh(\beta_k^-),\delta_1\sin \theta_j,\dots,\delta_l\sin \theta_l,\epsilon^+_1\sinh(\beta_1^+),\dots,\epsilon^+_m\sinh(\beta_m^+)\right) \\
&\oplus\bigoplus_{j=1}^r\mathrm{diag}(1,-1),\\
C= &\; \mathrm{ diag}\left(\epsilon^-_1\sinh(\beta_1^-),\dots,\epsilon^-_k\sinh(\beta_k^-),-\delta_1\sin \theta_j,\dots,-\delta_l\sin \theta_l,\epsilon^+_1\sinh(\beta_1^+),\dots,\epsilon^+_m\sinh(\beta_m^+)\right)\\ 
&\oplus\bigoplus_{j=1}^rS(r_j,\gamma_j),\\ 
    \end{split}
\end{equation}
where $\epsilon_j^\pm,\delta_j$ are $\pm$ signs, and
$$
S(r,\gamma)=\left(\begin{array}{cc}
r^2 \cos 2\gamma -1 & -r^2\sin 2\gamma\\
-r^2 \sin 2\gamma & -r^2 \cos 2\gamma +1
\end{array}\right).
$$
The $B$-signature of $M$ is therefore
$$
\mathrm{sign}_B(M)=\left(\epsilon^-_1,\dots,\epsilon^-_k,\delta_1,\dots,\delta_l,\epsilon^+_1,\dots,\epsilon^+_m\right),
$$
and therefore they are distinct as elements in the GIT quotient $Sp^\mathcal{I}(2n)//GL_n(\mathbb R)$. Moreover, it is easy to check that if
we flip some of the $\epsilon_j^\pm$ signs, we obtain matrices that are symplectically conjugated to the original one.

By varying the $B$-signature we conclude that there are $2^{n-2l}=2^{k+l+m}$ connected components of the loci $\mathcal{M}^{klmr}$, i.e.\
$$
\mathcal{M}^{klmr}=\bigsqcup_{\epsilon}\mathcal{M}^{klmr; \epsilon},
$$
where $\epsilon$ is a $B$-signature with $n$ entries. If we denote by $\mathcal{M}^{klmr}_0$ the region of $\mathbb R^n$ given by the image of $\mathcal{M}^{klmr}$ under the GIT map $Sp^\mathcal{I}(2n)//GL_n(\mathbb R)\rightarrow \mathbb R^n$, the covering degree of this map over $\mathcal{M}^{klmr}_0$ is exactly the number of connected components (or \emph{branches}) of $\mathcal{M}^{klmr}$, i.e.\ given by
$$
\mathrm{deg}\;\mathcal{M}^{klmr}=2^{n-2l}=2^{k+l+m}.
$$
The GIT map then collapses the different $\mathcal{M}^{klmr; \epsilon}$ to the region $\mathcal{M}^{klmr}_0$.

We also consider the stratum $\mathcal{M}_1^{klmr}$ of matrices $[M]\in Sp(2n)//Sp(2n)$ with the above eigenvalue configuration. By varying the Krein-signature $\kappa$, we obtain a splitting
$$
\mathcal{M}_1^{klmr}=\bigsqcup_{\kappa}\mathcal{M}_1^{klmr; \kappa},
$$
into connected components. The GIT map $Sp^\mathcal{I}(2n)//GL_n(\mathbb R)\rightarrow Sp(2n)//Sp(2n)$ collapses certain branches $\mathcal{M}_1^{klmr; \kappa}$ together, i.e.\ whenever two such signatures $\kappa,\kappa'$ coincide up to removing all those signs corresponding to hyperbolic eigenvalues. More precisely, we define a \emph{collapsing} operation on signatures, determined by
$$
\mathrm{Col}\left(\epsilon^-_1,\dots,\epsilon^-_k,\delta_1,\dots,\delta_l,\epsilon^+_1,\dots,\epsilon^+_m\right)=\left(\delta_1,\dots,\delta_l\right),
$$
that only keeps the signs of the elliptic eigenvalues. Then the GIT sequence maps
$$
\mathcal{M}^{klmr;\epsilon} \mapsto \mathcal{M}_1^{klmr;\mathrm{Col}(\epsilon)} \mapsto \mathcal{M}_0^{klmr}
$$
homeomorphically. The degree of $\mathcal{M}_1^{klmr}$, defined as the covering degree of the second map in the GIT sequence over the component $\mathcal{M}_0^{klmr}$ (that is, the number of branches $\mathcal{M}_1^{klmr;\kappa}$ that project to it), is then given by
$$
\mathrm{deg}\;\mathcal{M}_1^{klmr}=2^{l}.
$$
The degree of the first map in the GIT sequence over a branch $\mathcal{M}_1^{klmr;\kappa}$, which we call the \emph{relative} degree (the number of branches $\mathcal{M}^{klmr;\epsilon}$ that collapse to it), is 
$$
\mathrm{deg}\;\left(\mathcal{M}^{klmr}\Big\vert\mathcal{M}_1^{klmr}\right)=2^{k+m},
$$
which is in particular independent of $\kappa$. Moreover, we have the obvious multiplicative formula
$$ \mathrm{deg}\;\left(\mathcal{M}^{klmr}\Big\vert\mathcal{M}_1^{klmr}\right)\cdot \mathrm{deg}\;\mathcal{M}_1^{klmr}=\mathrm{deg}\;\mathcal{M}^{klmr}.
$$
\medskip

\textbf{Nonregular cases.} The nonregular cases can be dealt with analogously; the combinatorics just gets more involved. Assume that $A$ has real eigenvalues 
$$
\mu^-_1\leq\dots \leq \mu^-_k\leq -1\leq \nu_{1}\leq \dots \leq \nu_l \leq 1 \leq \mu^+_{1}\leq \dots\leq \mu_m^+, 
$$
where we also allow $\pm 1$ as an eigenvalue, and complex eigenvalues
$$
\alpha_1,\overline{\alpha}_1,\dots,\alpha_r,\overline{\alpha}_r.
$$
We denote the multiplicities by
$$
m_j^\pm=\mathrm{mult}(\mu_j^\pm),\; o_j=\mathrm{mult}(\nu_j), \;M_\pm= \mathrm{mult}(\pm 1),\; p_j= \mathrm{mult}(\alpha_j)=\mathrm{mult}(\overline{\alpha_j}).
$$
If we let $m^\pm=(m_1^\pm,\dots,m_k^\pm), o=(o_1,\dots,o_l), p=(p_1,\dots,p_r)$, we have
$$
\vert m^+ \vert + M_-+\vert o \vert + M_+ +\vert m^- \vert +2\vert p\vert=n, 
$$
where we define
$$
\vert(a_1,\dots,a_n)\vert=\sum_{j=1}^na_i.
$$
We denote by 
$$\mathcal{M}^{klmr}_{m^-M_-oM_+m^+p}=(\mathcal{H}^-)_{m^-}^k(-\mathcal I)_{M_-}\mathcal{E}_o^l\;\mathcal{I}_{M_+}(\mathcal{H}^+)_{m^+}^m\mathcal{N}_p^r$$ the stratum of matrices $[M]\in Sp^{\mathcal{I}}(2n)//GL_n(\mathbb R)$ with $A$-block with an eigenvalue configuration as above.

\begin{figure}
    \centering
    \includegraphics[width=0.65\linewidth]{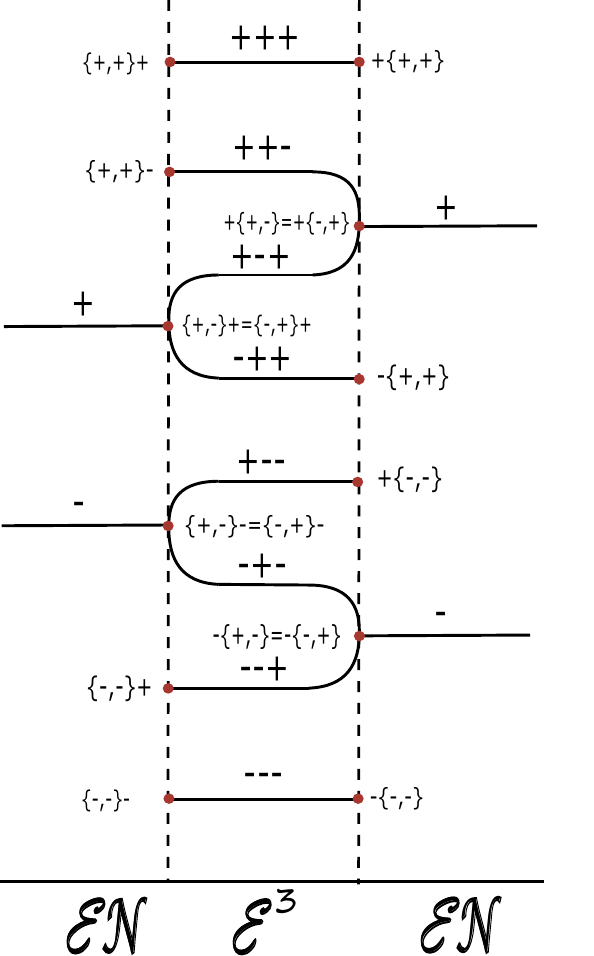}
    \caption{The branching structure of $Sp^{\mathcal{ I}}(6)//GL_3(\mathbb R)$, for the case $n=3$, corresponding to a transition from $\mathcal{E}^3$ to $\mathcal{EN}$. The right and left hand sides should be identified. The branching structure for $Sp(6)//Sp(6)$ coincides for this transition, as no hyperbolic eigenvalues are involved. Note that the branches do \emph{not} cross along a contraction which has definite signature inside the brackets, as predicted by the Krein--Moser theorem.}
    \label{fig:EN}
\end{figure}

\begin{figure}
    \centering
    \includegraphics[width=0.8\linewidth]{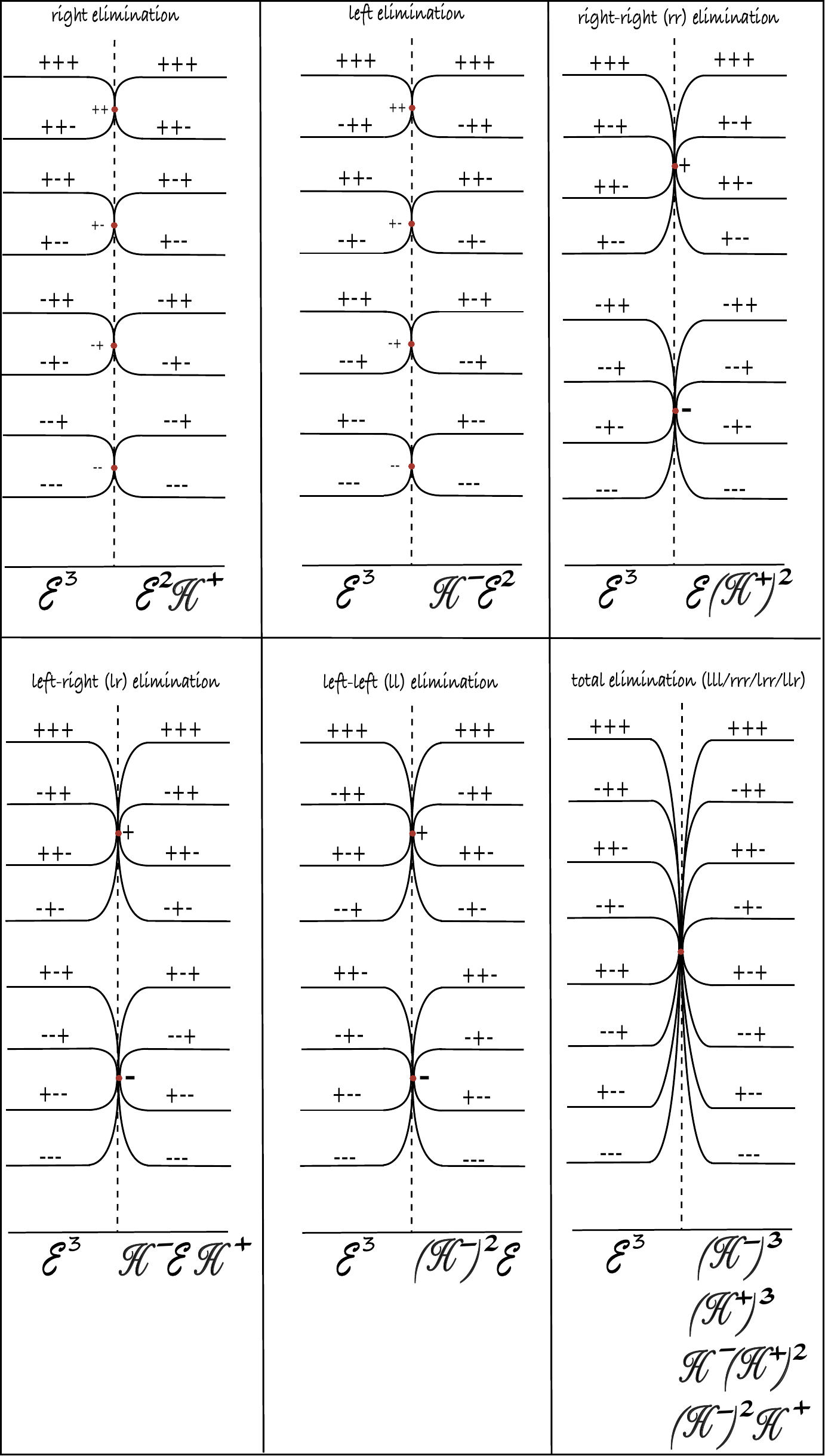}
    \caption{Case $n=3$, for $Sp^{\mathcal{I}}(6)//GL_3(\mathbb R)$. The transitions out of the stable region which do not involve $\mathcal{N}$ correspond to all possible elimination contractions.}
    \label{fig:eliminations}
\end{figure}

Then there exist unique $\beta^\pm_j, r_j \in [0,+\infty)$ and $\theta_j,\gamma_j\in [0,\pi]$ satisfying
$$
\mu_j^-=-\cosh \beta_j^-,\;\mu_j^+=\cosh \beta_j^+,\; \nu_j=\cos \theta_j,\;\alpha_j=r_je^{i\gamma_j},\;\overline{\alpha}_j=r_je^{-i\gamma_j}. 
$$ 
As we can ignore Jordan blocks in the GIT quotient, we may assume that 
$$
A=\bigoplus_{j=1}^k \mu_j^- \mathds 1_{m_j^-} \oplus -\mathds 1_{M_-} \oplus \bigoplus_{j=1}^l \nu_j \mathds 1_{o_j} \oplus \mathds 1_{M_+} \oplus \bigoplus_{j=1}^m \mu_j^+ \mathds 1_{m_j^+} \oplus \bigoplus_{j=1}^r \bigoplus_{i=1}^{p_j} R(r_j,\gamma_j), 
$$
where $\mathds 1_n$ denotes the identity matrix of size $n$. Applying \cite{FM} to each summand as before, we deduce the normal form
\begin{equation}
\begin{split}
B=&\bigoplus_{j=1}^k \left(\sinh(\beta_j^-) \mathds 1_{a_j^-}  \oplus -\sinh(\beta_j^-) \mathds 1_{b_j^-}\right)\oplus 0_{M_-} \oplus \bigoplus_{j=1}^l \left(\sin(\theta_j) \mathds 1_{c_j}  \oplus -\sin(\theta_j) \mathds 1_{d_j}\right) \oplus 0_{M_+}\oplus\\
& \bigoplus_{j=1}^m  \left(\sinh(\beta_j^+) \mathds 1_{a_j^+}  \oplus -\sinh(\beta_j^+)\mathds 1_{b_j^+}\right) \oplus \bigoplus_{j=1}^r \bigoplus_{i=1}^{p_j} \mathrm{diag}(1,-1), 
\end{split}
\end{equation}
\begin{equation}
\begin{split}
C=&\bigoplus_{j=1}^k \left(\sinh(\beta_j^-) \mathds 1_{a_j^-}  \oplus -\sinh(\beta_j^-) \mathds 1_{b_j^-}\right)\oplus 0_{M_-} \oplus \bigoplus_{j=1}^l \left(-\sin(\theta_j) \mathds 1_{c_j}  \oplus \sin(\theta_j) \mathds 1_{d_j}\right) \oplus 0_{M_+}\oplus\\
& \bigoplus_{j=1}^m  \left(\sinh(\beta_j^+) \mathds 1_{a_j^+}  \oplus -\sinh(\beta_j^+)\mathds 1_{b_j^+}\right) \oplus \bigoplus_{j=1}^r \bigoplus_{i=1}^{p_j} S(r_j,\gamma_j), 
\end{split}
\end{equation}
with $0_m$ the zero matrix of size $m$, and with the constraints
$$
a_j^\pm+b_j^\pm=m_j^\pm,\; c_j+d_j=o_j.
$$
The $B$-signature of $M$ is therefore given by
$$
\mathrm{sign}_B(M)=((a_1^-,b_1^-),\dots,(a_k^-,b_k^-), (c_1,d_1),\dots,(c_l,d_l),(a_1^+,b_1^+),\dots,(a_m^+,b_m^+)).
$$
If we let $\mathcal{M}^{klmr}_{m^-M_-oM_+m^+p;0}\subset \mathbb R^n$ the image of $\mathcal{M}^{klmr}_{m^-M_-oM_+m^+p}$ under the maps defining the GIT sequence, then by counting the normal forms for varying $B$-signature, we obtain that the degree of the GIT map over this component (that is, the number of connected components $\mathcal{M}^{klmr;\epsilon}_{m^-M_-oM_+m^+p}$, or the \emph{branches} of $\mathcal{M}^{klmr}_{m^-M_-oM_+m^+p}$) is given by
$$
\mathrm{deg}\;\mathcal{M}^{klmr}_{m^-M_-oM_+m^+p}=\prod_{j=1}^k(m_j^-+1)\prod_{j=1}^l (o_j+1)\prod_{j=1}^m(m_j^++1).
$$
Note that if $m_j^\pm=o_j=1$ then we recover the formula for the regular case, given above. And moreover if $m_j^\pm=o_j=0$ (i.e.\ only complex quadruples or $\pm 1$ appear as eigenvalues), then the above degree is $1$. The reader should compare the above formula to Figures \ref{fig:GIT_sequence_3D}, \ref{fig:EN}, \ref{fig:eliminations}, as it precisely gives the number of branches above each component of the base.

In a similar fashion, we denote by $\mathcal{M}^{klmr}_{m^-M_-oM_+m^+p;1}$ the stratum of matrices $[M]\in Sp(2n)//Sp(2n)$ that have eigenvalue configuration as above, and also by $\mathcal{M}^{klmr;\kappa}_{m^-M_-oM_+m^+p;1}$ its connected components, which are labelled by the Krein signature $\kappa$. By counting these connected components, we obtain that its degree is
$$
\mathrm{deg}\;\mathcal{M}^{klmr}_{m^-M_-oM_+m^+p;1}=\prod_{j=1}^l (o_j+1).
$$
We also have an analogous collapsing map
$$
\mathrm{Col}((a_1^-,b_1^-),\dots,(a_k^-,b_k^-), (c_1,d_1),\dots,(c_l,d_l),(a_1^+,b_1^+),\dots,(a_m^+,b_m^+))=((c_1,d_1),\dots,(c_l,d_l)),
$$
and the GIT sequence sends
$$
\mathcal{M}^{klmr;\epsilon}_{m^-M_-oM_+m^+p}\mapsto\mathcal{M}^{klmr; \mathrm{Col}(\epsilon)}_{m^-M_-oM_+m^+p;1}\mapsto \mathcal{M}^{klmr}_{m^-M_-oM_+m^+p;0}.
$$
The relative degree is 
$$
\mathrm{deg}\;\left(\mathcal{M}^{klmr}_{m^-M_-oM_+m^+p}\Big\vert \mathcal{M}^{klmr}_{m^-M_-oM_+m^+p;1}\right)=\prod_{j=1}^k(m_j^-+1)\prod_{j=1}^m(m_j^++1),
$$
so that the multiplicative formula holds:
$$
\mathrm{deg}\;\left(\mathcal{M}^{klmr}_{m^-M_-oM_+m^+p}\Big\vert \mathcal{M}^{klmr}_{m^-M_-oM_+m^+p;1}\right)\mathrm{deg}\;\mathcal{M}^{klmr}_{m^-M_-oM_+m^+p;1}=\mathrm{deg}\;\mathcal{M}^{klmr}_{m^-M_-oM_+m^+p}.
$$

\medskip

\textbf{Stable region.} We now relate the above computations to stability of periodic orbits. The stable region is by definition
$$
\mathcal{E}^n=\bigcup_{l=0}^n\bigcup_{\substack{M_+,M_-,o\\M_++M_-+\vert o \vert=n}}\mathcal{E}^l_{M_-oM_+}
$$
where we have simplified the notation by
$$
\mathcal{E}^l_{M_-oM_+}:=\mathcal{M}^{0l00}_{\mathbf{0}M_-oM_+\mathbf{0}\mathbf{0}},
$$
with $\mathbf{0}=(0,\dots,0)$. This is naturally a stratified space: its top open stratum is $\mathcal{E}^n_{0\mathbf{1}0}$ with $\mathbf{1}=(1,\dots,1)$, which is $n$-dimensional, and represented by matrices with only simple elliptic eigenvalues. Its boundary is given by
$$
\partial \mathcal{E}^n_{0\mathbf{1}0}=\bigcup_{l=0}^{n-1}\bigcup_{\substack{M_+,M_-,o\\M_++M_-+\vert o \vert=n}}\mathcal{E}^l_{M_-oM_+},
$$
and is represented by matrices having only elliptic eigenvalues, at least one with multiplicity higher than $1$, and/or with $\pm 1$ in the spectrum. The dimension of the boundary stratum $\mathcal{E}^l_{M_-oM_+}$ is given by
$$
\dim \;\mathcal{E}^l_{M_-oM_+}=n - M_--M_+-\sum_{j=1}^l(o_j-1)=l.
$$
In particular, $\mathcal{E}^0_{M_-\mathbf{0}M_+}$ is simply a point (represented by the matrix $\mathds 1_{M_-}\oplus \mathds 1_{M_+}$). In order to describe the boundary of $\mathcal{E}^l_{M_-oM_+}$, we need to introduce a \emph{basic} contraction operator $C_j: \mathbb R^k\rightarrow \mathbb R^{k-1}$ for each $j=1,\dots, l-1$. Given an ordered tuple of real numbers $(a_1,\dots,a_k) \in \mathbb R^{k}$, its contraction on consecutive entries is defined by
$$
C_j(a_1,\dots,a_k):=(a_1,\dots,\wick{\c a_j,\c a_{j+1}}, \dots, a_n):=(a_1,\dots, a_j+a_{j+1},\dots, a_k)\in \mathbb R^{k-1}.
$$
A \emph{contraction} of $a=(a_1,\dots,a_k)$ is by definition the result of applying a sequence of basic contractions to $a$. 

The boundary of $\mathcal{E}^l_{M_-oM_+}$ is then given by
$$
\partial \mathcal{E}^l_{M_-oM_+}=\bigcup_{l'=0}^{l-1}\bigcup_{\substack{(M'_-,o',M'_+)\in\\ C(M_-,o,M_+)}}\mathcal{E}^{l'}_{M'_-o'M'_+},
$$
where $C(M_-,o,M_+)$ denotes the (finite) set of all possible contractions of $(M_-,o,M_+)\in \mathbb R^{l+2}$.

\medskip

\textbf{The associahedron.} The boundary combinatorics of the stable region can be, alternatively, encoded in the following manner. We identify the simple eigenvalues
$$
-1<\nu_1<\dots<\nu_l<1
$$
with the word $-1\nu_1\dots \nu_l 1$. When moving from one stratum to the next, we add a bracket, only allowing to group two consecutive elements (i.e.\ letters) at a time. For instance
$$
-1\nu_1\dots \nu_l 1\mapsto -1\nu_1\dots \{\nu_j,\nu_{j+1}\}\dots\nu_l 1
$$
indicates that the eigenvalues $\nu_j$ and $\nu_{j+1}$ come together into a multiplicity two eigenvalue, and therefore corresponds to the contraction of multiplicities given by
\begin{equation*}
(1,\dots, 1) \mapsto (1,\dots,2,\dots,1).   \end{equation*}
Analogously, a further parenthesis
$$
-1\nu_1\dots \nu_{j-1}\{\nu_j,\nu_{j+1}\}\dots\nu_l 1\mapsto -1\nu_1\dots \{\nu_{j-1},\nu_j,\nu_{j+1}\}\dots\nu_l 1
$$
indicates that the eigenvalue $\nu_{j-1}$ came together with the previous multiplicity two eigenvalue, forming a multiplicity three eigenvalue, and therefore corresponding to the contraction
$$
(1,\dots, 1,2,\dots,1)\mapsto (1,\dots,3,\dots,1).
$$
This construction iterates in the obvious way. Note that here we should also allow eigenvalues to come together with $\pm 1$, i.e.\ $\{-1,\nu_1\}\nu_2\{\nu_3,\nu_4,1\}$ is a valid expression. We use the bracket notation to indicate that the order of the elements in the bracket is irrelevant.
Iterating this construction results in a poset of strings (in which all open brackets are accompanied by a corresponding closed one, and there are no nested brackets), and where two strings $a,b$ satisfy $a\leq b$ if and only if $b$ is obtained by a sequence of brackets operations from $a$. By construction, this poset encodes the boundary combinatorics of the stable region.

\begin{figure}
    \centering
    \includegraphics[width=1\linewidth]{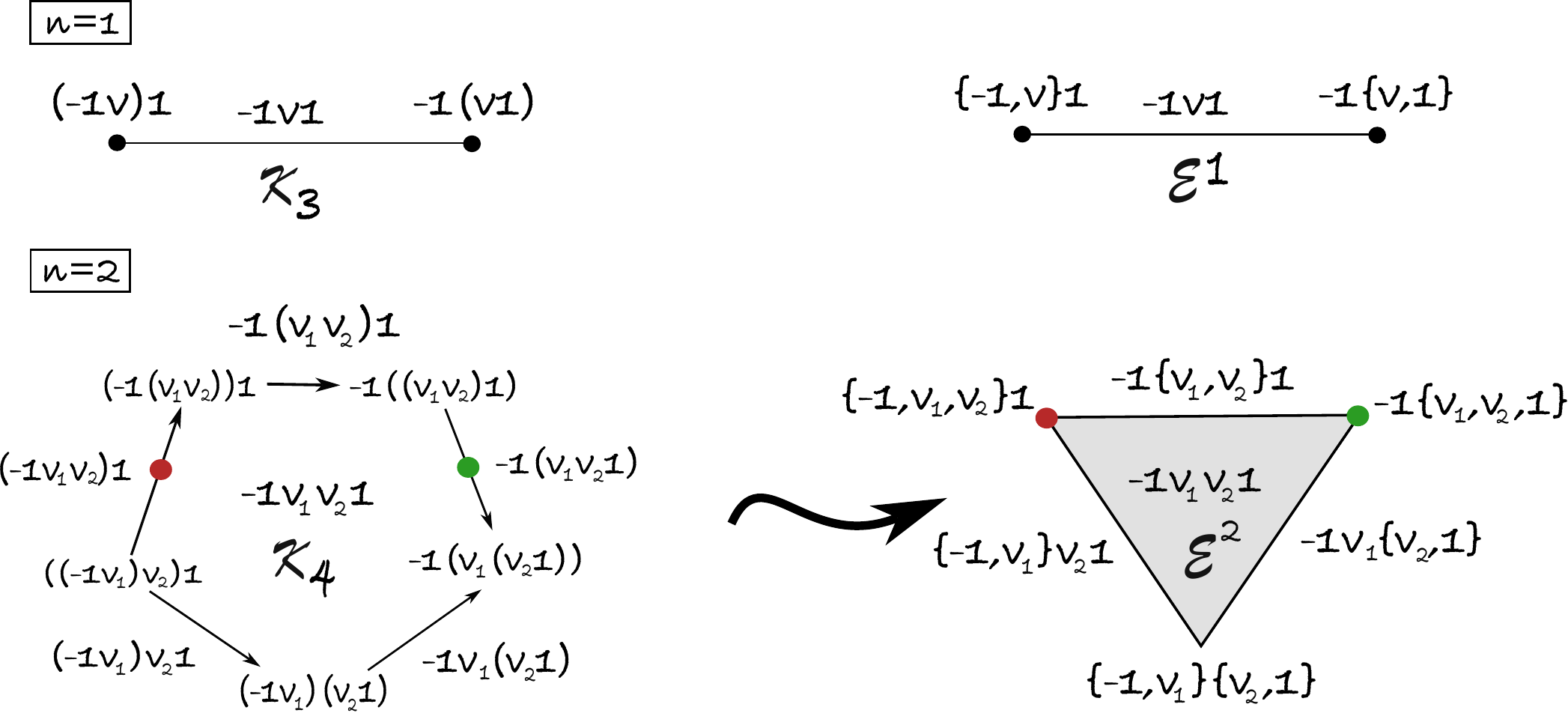}
    \caption{The associahedron $K_3$ agrees with the $1$-dimensional stable region $\mathcal{E}^1$ (an interval). The $2$-dimensional stable region, a triangle, is obtained up to homeomorphism by collapsing two edges in $K_4$, a pentagon, to a point, indicated by points in the figure.}
    \label{fig:stable_regionlowd}
\end{figure}

We should note that the above is related to the operation of taking parenthesis
$$
-1\nu_1\dots \nu_l 1\mapsto -1\nu_1\dots (\nu_j\nu_{j+1})\dots\nu_l 1
$$
and iterating them, similarly as above, e.g.\ as
$$
-1\nu_1\dots \nu_{j-1}(\nu_j\nu_{j+1})\dots\nu_l 1\rightarrow -1\nu_1\dots (\nu_{j-1}(\nu_j\nu_{j+1}))\dots\nu_l 1,
$$
and so on, where now a valid expression is e.g.\ $((-1\nu_1)\nu_2)\nu_3(\nu_41)$. The bracket is then obtained as the result of removing all interior parenthesis in an expression, symbolically via $(\dots (\dots )\dots)\mapsto (\dots)$, and modding out by the action of the corresponding permutation group (which acts on the number of elements inside the bracket), symbolically via
$$
(\underbrace{\dots}_m)\mapsto \{\underbrace{\dots}_m\}=(\underbrace{\dots}_m)/S_m.
$$ As an example, the above expression becomes $\{-1,\nu_1,\nu_2\}\nu_3\{\nu_4,1\}$, where now the order of the elements inside the bracket is irrelevant.

But the combinatorics of expressions with parenthesis is goverened by a very-well known polytope: the \emph{associahedron}. This is by definition the $(m -2)$-dimensional convex polytope $K_m$ in which each vertex corresponds to a way of \emph{correctly} inserting opening and closing parentheses in a string of $m$ letters (this means that it uniquely determines the order of the product operations), and the edges correspond to single application of the associativity rule. This can also be viewed as a poset, when the arrow indicates that the parentheses have been moved to the right (this is called the \emph{Tamari lattice}). Moreover, one can label the edges with ``\emph{incorrect}'' expressions, containing the common parentheses to each of its boundary vertices. And one can further label the faces by also ``incorrect'' expressions containing the parentheses common to all its boundary edges. This process iterates in the obvious way, ending in the top strata, which is labelled by the string with no brackets $-1\nu_1\dots\nu_l1$.  

\begin{figure}
    \centering
    \includegraphics[width=1.02\linewidth]{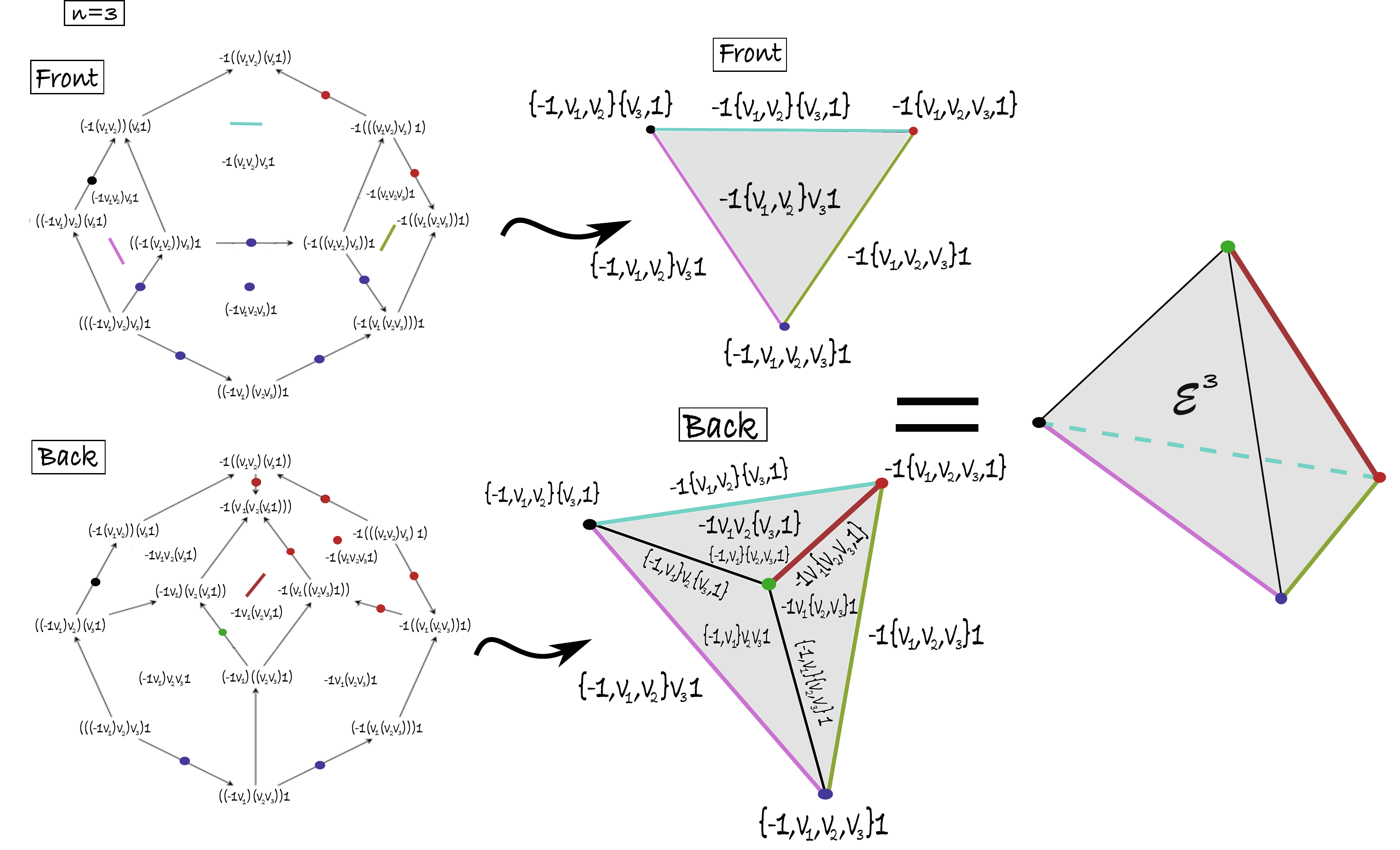}
    \caption{The $3$-dimensional stable region, a tetrahedron, is obtained up to homeomorphism from $K_5$, an enneahedron, by collapsing 2 faces and 12 edges of $K_5$ to a point (labelled with a point), and 4 faces of $K_5$ to a line (labelled with a line).}
    \label{fig:associahedron}
\end{figure}

In order to obtain the stable region from the associahedron, we note that many labels in this polytope are actually equivalent when written with the bracket notation. We then conclude that \emph{the boundary combinatorics of the stable region is determined by the associahedron $K_{n+2}$ in $n+2$ letters, up to collapsing a suitable collection of strata to lower dimensional strata.} In other words, the stable region is homeomorphic to a quotient of the associahedron, where we identify those strata whose label become equivalent when written in the bracket notation. The low dimensional cases ($n=1,2,3$) are depicted in Figure \ref{fig:stable_regionlowd} and Figure \ref{fig:associahedron}.

\medskip

\textbf{Branching structure.} We will denote the branches of $\mathcal{E}^l_{M_-oM_+}$, labelled by the $B$-signature $$\epsilon=((c_1,d_1),\dots,(c_l,d_l)),$$ by $\mathcal{E}^{l;\epsilon}_{M_-oM_+},$ so that
$$
\mathcal{E}^l_{M_-oM_+}=\bigsqcup_\epsilon
\mathcal{E}^{l;\epsilon}_{M_-oM_+}.
$$
By the computation of degrees above, the interior of the top stratum $\mathcal{E}^n_{0\mathbf{1} 0}$ has $2^n$ branches. Some of these come together at the boundary strata, and which branches can come together over each boundary stratum is determined by the corresponding $B$-signatures, as these behave continuously whenever defined. To explain this, we first define a few operations. A \emph{middle contraction} operation on consecutive entries of a $B$-signature, symbolically, is defined via the formula
\begin{equation*}
((c_1,d_1),\dots,\wick{(\c c_j,d_j), (c_{j+1},\c d_{j+1})},\dots,(c_l,d_l)):=((c_1,d_1),\dots,(c_j+c_{j+1},d_j+d_{j+1}),\dots,(c_l,d_l)).
\end{equation*}
We also allow the \emph{elimination} contraction of the first and last entries, defined by
$$
(\wick{(\c c_1,\c d_1)},\dots,(c_l,d_l))=((c_2,d_2),\dots, (c_{l},d_{l})),
$$
$$
(( c_1, d_1),\dots,\wick{\c( c_l,\c d_l)})=((c_1,d_1),\dots, (c_{l-1},d_{l-1})).
$$
The first operation is the \emph{left} elimination contraction, and the second, the \emph{right} elimination contraction. A \emph{contraction} of a $B$-signature is then by definition a sequence of middle contractions of consecutive entries, composed with elimination contractions. The \emph{right/left elimination order} of the contraction is defined as the number of right/left elimination contractions in the sequence, the \emph{middle order}, as the number of middle contractions in the sequence, and the \emph{total order}, as the sum of right/left/middle contractions. A middle contraction represents elliptic eigenvalues of possibly high multiplicity coming together but staying elliptic; a left elimination contraction represents an elliptic eigenvalue of possibly high multiplicity becoming $-1$; and a right elimination contraction represents an elliptic eigenvalue of possibly high multiplicity becoming $1$. In the first case, signatures simply add up, while in the second cases, since the $B$-signature is not defined (as the $B$-block becomes degenerate), it is simply suppressed. We note that a contraction of $B$-signatures determines an underlying contraction of multiplicities, by replacing the pair $(c_j,d_j)$ in the above tuple by $o_j=c_j+d_j$. 

In keeping with the bracket notation, we will sometimes also use the notation $(c,d)=\{\underbrace{+,\dots,+}_{c},\underbrace{-,\dots,-}_{d}\}$, and write contractions as
\begin{equation*}
\begin{split}
(\{\underbrace{+,\dots,+}_{c_1},\underbrace{-,\dots,-}_{d_1}\},\dots, \wick{\{\underbrace{+,\dots,+}_{c_j}\c,\underbrace{-,\dots,-}_{d_j}\}, \{\underbrace{+,\dots, +}_{c_{j+1}}\c,\underbrace{-,\dots,-}_{d_{j+1}}\}}, \dots, \{\underbrace{+,\dots,+}_{c_l},\underbrace{-,\dots,-}_{d_l}\})
\end{split}
\end{equation*}
$$
:=(\{\underbrace{+,\dots,+}_{c_1},\underbrace{-,\dots,-}_{d_1}\},\dots, \{\underbrace{+,\dots,+}_{c_j+c_{j+1}},\underbrace{-,\dots,-}_{d_j+d_{j+1}}\}, \dots, \{\underbrace{+,\dots,+}_{c_l},\underbrace{-,\dots,-}_{d_l}\}).
$$

\begin{figure}
    \centering
    \includegraphics[width=0.8\linewidth]{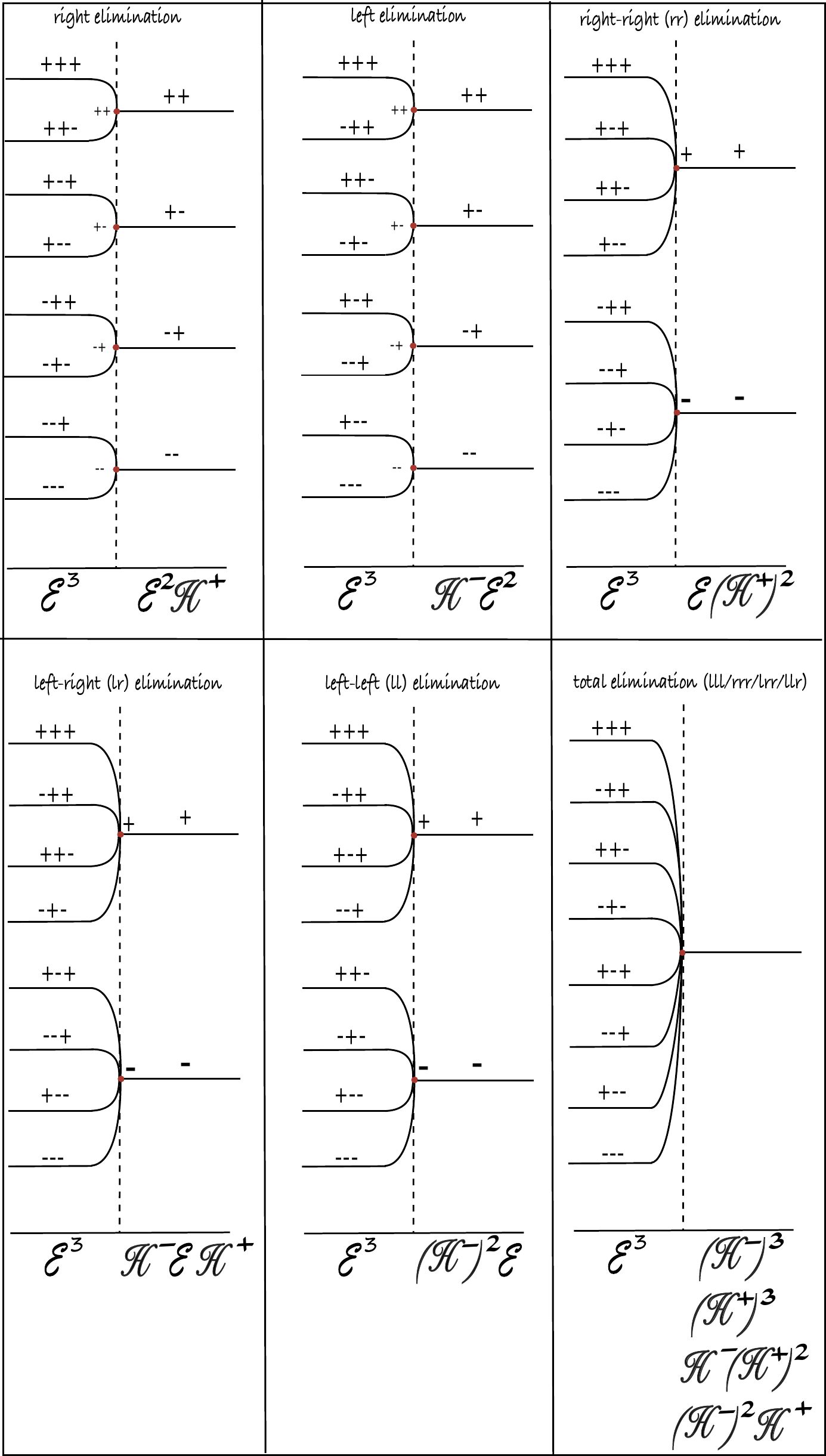}
    \caption{The branching structure for $Sp(6)//Sp(6)$ is obtained from the ones in Figure \ref{fig:eliminations} by collapsing all branches corresponding to hyperbolic eigenvalues together.}
    \label{fig:eliminations_middle}
\end{figure}

The branch $\mathcal{E}^{l;\epsilon}_{M_-oM_+}$ can then meet the branch $\mathcal{E}^{l';\epsilon'}_{M_-'o'M_+'}$ at the boundary, where $M'_-o'M'_+=c(M_-,o,M_+)$ for some contraction $c$, if and only if 
$$
\epsilon'=((c'_1,d'_1),\dots, (c'_{l'},d'_{l'}))
$$
is obtained by the underlying contraction of $\epsilon$ determined by $c$, with $$
l'=l-\mathrm{tord}(c),
$$ 
$$
M_-'=M_- + \mathrm{lord}(c),
$$ 
$$
M_+'=M_+ + \mathrm{rord}(c),
$$ 
$$
\vert o'\vert =\vert o \vert-\mathrm{lord}(c)-\mathrm{mord}(c),
$$ 
where $\mathrm{mord}(c),\mathrm{lord}(c),\mathrm{rord}(c), \mathrm{tord}(c)$ are respectively the middle/left/right/total order of $c$. It follows that, in particular, each $c_j'$ is a sum of consecutive $c_j$'s, and each $d_j'$ is a sum of consecutive $d_j$'s. We then conclude that two branches can come together at the boundary if and only if their corresponding $B$-signatures have a common contraction, and the branch at which they meet (and therefore the boundary strata at which they meet) is completely determined by this contraction. The co-dimension of this strata is the total order of the contraction, as follows from the first of the above formulae. The low-dimensional cases are depicted in Figure \ref{fig:GIT_sequence} ($n=1$) and Figure \ref{fig:GIT_sequence_3D} ($n=2$). For $n=3$, the branching structure corresponding to transitions from the stable region is shown in Figure \ref{fig:EN} and Figure \ref{fig:eliminations}. The combinatorics for the branching structure for all other transitions between components is obtained analogously; several of the diagrams are the same. 

\chapter{Astrodynamics: numerical work}\label{ch:num_work}

This chapter is devoted to numerical work carried out for the CR3BP, by the author's collaborators Otto van Koert, Cengiz Aydin, Bhanu Kumar and Dayung Koh. The aim is to illustrate the use of all the theory we discussed so far, in the context of numerical studies of periodic orbits for various systems, in the context of astroynamics, and with a view towards space mission design. The systems of interest are:

\medskip

\begin{itemize}
    \item Hill's lunar problem;

\medskip
    
    \item Saturn--Enceladus (CR3BP with $\mu=1.9002485658670E$-$7$);

\medskip
    
    \item Jupiter--Europa (CR3BP with $\mu=2.5266448850435E$-$5$);

\medskip
    
    \item Earth--Moon (CR3BP with $\mu=1.215058560962404E$-$2$).
\end{itemize}

\medskip

\begin{figure}[h]
    \centering
    \includegraphics[width=0.66\linewidth]{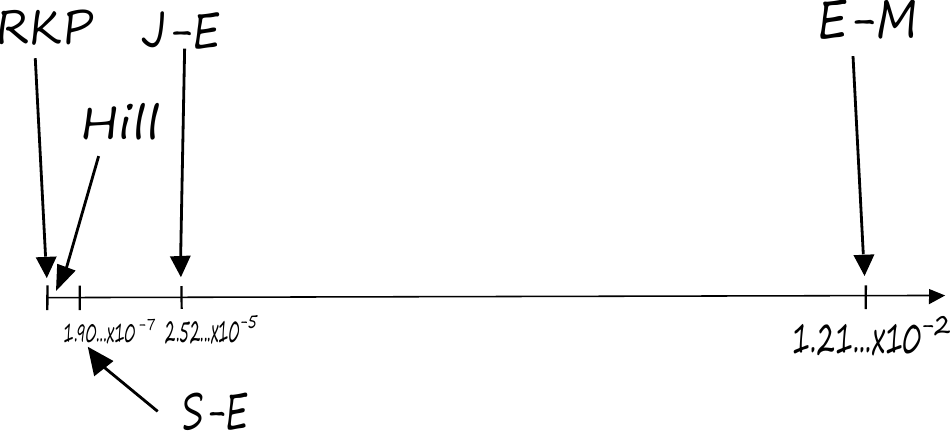}
    \caption{The values of $\mu$.}
    \label{fig:mus}
\end{figure}

These systems are currently very popular due to the ambitious Artemis program (centered around the Earth--Moon system), and the fact that Europa and Enceladus are icy moons which might harbor conditions for the existence of life. Note that the Earth--Moon system is much further away from the integrable case than the rest of the above systems, as sketched in Figure \ref{fig:mus}. Other systems of current interest include Jupiter---Ganymede, which e.g.\ is the target system of the JUICE mission of ESA, although we will not be carrying out studies for this particular system in what follows.

\section{Bifurcation graphs} This section is based on the papers \cite{Ayd, AFvKKM}. We will study the Jupiter--Europa system as a deformation of Hill's lunar problem, i.e.\ by deforming the $\mu$ parameter (althoug this is fomally obtained by performing a suitable Taylor expansion and truncation with respect to $\mu$, and a change of coordinates). In order to do this, we first need to study periodic orbits in the Hill problem. This is done in \cite{Ay22}, which is our starting point. The main tool are the \emph{bifurcation graphs}, which represent families of orbits with varying Jacobi constant, and bifurcations between them. We will be interested in studying \emph{spatial} orbits bifurcating out of planar orbits.

\section{Hill's lunar problem: numerical work} We now summarize the main results of Aydin's thesis \cite{Ay22} (where we refer the reader for more details), following his exposition.

\medskip

Based on foundational work of Hill \cite{H77}, Hénon \cite{He69, He70, He74, He03} numerically studied families of planar direct and retrograde periodic orbits, which are referred to as family $g$ and $f$, respectively. These two families are the fundamental families of symmetric planar periodic orbits in the Hill lunar problem. They are both doubly symmetric with respect to reflection at the $q_1$- and $q_2$-axis (see Section \ref{sec:symmetries_Hill} for the symmetries in Hill's problem). We will be interested in studying families that bifurcate out of these basic ones.

Moreover, we also have the northern and southern polar collision orbits $n$ and $s$, both spatial, bouncing on the primary along the upper half-space and the lower-half space respectively, which arise by continuation from the RKP, see \cite{BFvK}.

The following results are obtained \emph{analytically} for very low energies, by deforming to the Kepler problem. 

\begin{thm}[\cite{Ay22}] The CZ-indices of the orbits $f,g,n,s$ are given by
$$
\mu_{CZ}=\left\{\begin{array}{cc}
    6 & \mbox{ for } g \\
    4 & \mbox{ for } n\\
    4 & \mbox{ for } s\\
    2 & \mbox{ for } f
\end{array}\right.
$$
For the planar orbits $f,g$, their planar CZ-indices are
$$
\mu_{CZ}^p=\left\{\begin{array}{cc}
    3 & \mbox{ for } g \\
    1 & \mbox{ for } f\\
\end{array}\right.
$$
and their spatial CZ-indices are
$$
\mu_{CZ}^s=\left\{\begin{array}{cc}
    3 & \mbox{ for } g \\
    1 & \mbox{ for } f\\
\end{array}\right.
$$
so that indeed $\mu_{CZ}=\mu_{CZ}^p+\mu_{CZ}^s$.
    
\end{thm}

These families were then followed numerically, as well as other families which bifurcate from them, which is detected via jumps in the indices. If the planar index jumps, this is a planar-to-planar bifurcation, and if the planar one jumps, it is planar-to-spatial. The list of orbits relevant to Aydin's work are listed in Table \ref{orbits}.

\begin{table}[h]\fontsize{12}{12} \centering
\caption{Families of orbits studied in \cite{Ay22}.}
	\begin{tabular}{|c|c|c|}
 \hline
		Reference & family & type \\
		\hline  
H\'enon (1969) \cite{He69}	& $g,f$ & planar \\
 & $g'$ & planar ($\mu_{CZ}^p$ jump from $g$)\\
H\'enon (1970) \cite{He70}, (2003) \cite{He03} & $g_3$ & planar\\
Batkhin–Batkhina \cite{BB} (2009) & $g_{2\nu}$ & spatial ($\mu_{CZ}^s$ jump from $g$) \\
 & $g^{YOZ}_{1\nu}$ & spatial (from the 2nd cover of $g$) \\
 Michalodimitrakis (1980) & $g_{1\nu}$ & spatial (from the 2nd cover of $g$ and $g'$)\\
Kalantonis \cite{Ka20} (2020) & $f_g^{(2,3)}$ , $f_g^{(2cut,3)}$ & spatial (from the 3rd cover of $g$)\\
 & $f_{g'}^{(2,3)}$ , $f_{g'}^{(2cut,3)}$ & spatial (from the 3rd cover of $g'$) \\
 & $f_g^{(1,4)}$ , $f_g^{(1cut,4)}$ & spatial (from the 4th cover of $g$)\\
 & $f_{g'}^{(1,4)}$ , $f_{g'}^{(1cut,4)}$& spatial (from the 4th cover of $g'$)\\
\hline
 \end{tabular}
	
\label{orbits}
\end{table}

We now describe the numerical results. In what follows, the deformation parameter will be $\Gamma = -2c$, where $c$ is the Jacobi constant (note that $\Gamma$ is also traditionally called the Jacobi constant, but not here).

\begin{figure}
    \centering
    \includegraphics[width=1\linewidth]{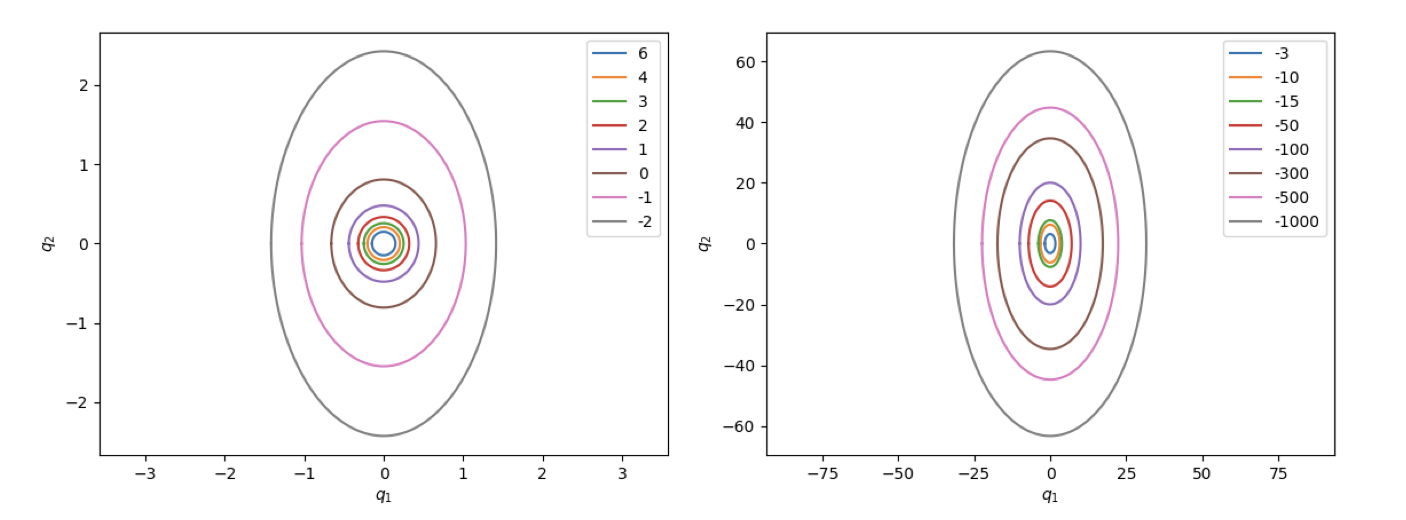}
    \caption{The $f$ family (plot from \cite{Ay22}).}
    \label{fig:f_family}
\end{figure}

\begin{figure}[h]
    \centering
    \includegraphics[width=0.6\linewidth]{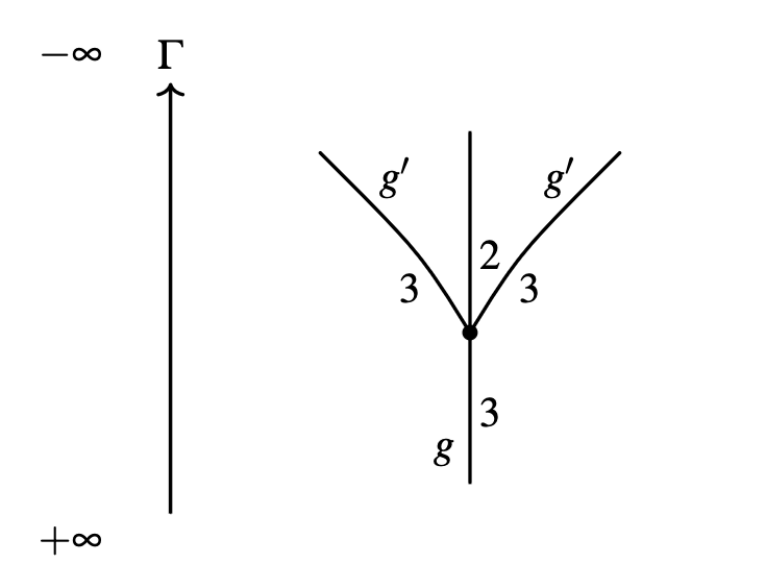}
    \caption{Bifurcation graph for the $g,g'$ orbits with \emph{planar} CZ-index.}
    \label{fig:pitchfork}
\end{figure}

\subsection{Planar direct/prograde orbits} The family $f$ stays spatial and planar elliptic for all times, thus their indices do not change, and they do not undergo bifurcation. This family is plotted in Figure \ref{fig:f_family}. However, the family $g$ is more interesting, since, as described by H\'enon, it undergoes a non-generic pitchfork bifurcation, going from elliptic to positive hyperbolic. Two new families of elliptic orbits, called $g'$, appear. These new families are still invariant under reflection at the $x$-axis, but not under reflection at the $y$-axis. Reflection at the $y$-axis maps one branch of the $g'$-family to the other branch. Figure \ref{fig:pitchfork} shows the bifurcation graph which is constructed as follows. Each vertex denotes a degenerate orbit at which bifurcation happens and each edge represents a family of orbits with varying energy, labeled by the corresponding CZ-index. From this data, it is easy to determine the associated Floer number. For instance in Figure \ref{fig:pitchfork} on the left, the (planar) Floer number is $(-1)^3 = -1$ before bifurcation, and $(-1)^2 + 2(-1)^3 = -1$ after bifurcation; they coincide, as they should.

\begin{figure}
    \centering
    \includegraphics[width=1 \linewidth]{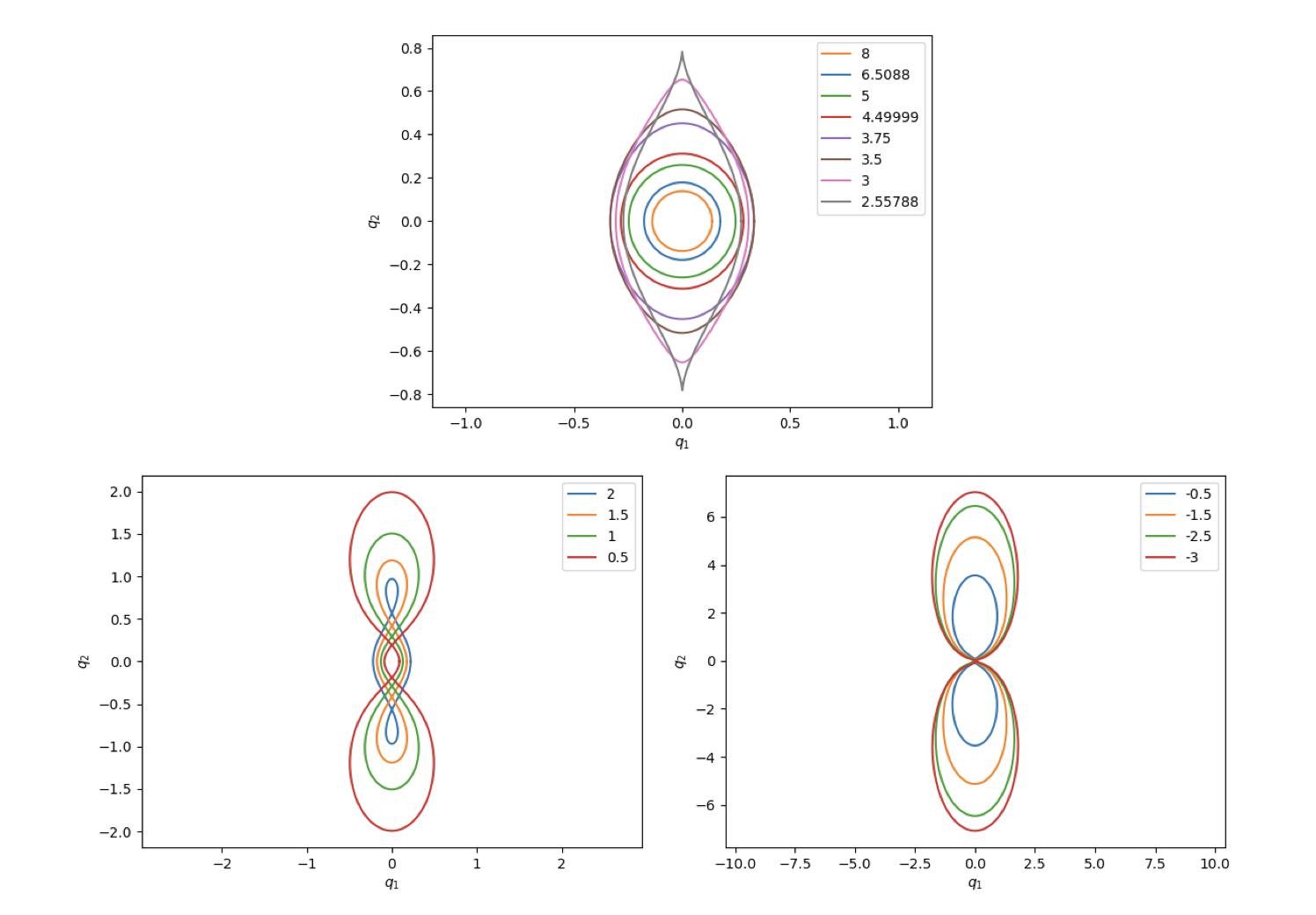}
    \caption{The $g$ family (plot from \cite{Ay22}).}
    \label{fig:g_family}
\end{figure}

The data for the family $g$ is presented in Table \ref{table:g_family}, and the family is plotted in Figure \ref{fig:g_family}.

\begin{table}[h]\fontsize{12}{12} \centering
\caption{Data for the $g$ family.}
	\begin{tabular}{|c|c|c|c|c|c|}
 \hline
		Values of $\Gamma$ & planar & spatial & $\mu_{CZ}^p$ & $\mu_{CZ}^s$ & $\mu_{CZ}$  \\
		\hline  
$(+\infty, 4.49999)$ & elliptic & elliptic & 3 & 3 & 6 \\
$(4.49999, 1.3829)$  & pos.\ hyperbolic & elliptic & 2 & 3 & 5 \\
$(1.3829, -\infty)$ & pos.\ hyperbolic & pos.\ hyperbolic & 2 & 4 & 6 \\
\hline
 \end{tabular}
\label{table:g_family}
\end{table}

The jump in the planar CZ-index at $\Gamma=4.49999$ corresponds to the appearance of $g'$. The data for this family is collected in Table \ref{table:gprime_family}. As a sanity check, we can also check invariance of the spatial Floer numbers at the $g'$ bifurcation, which are $(-1)^6=1$ before bifurcation, and $(-1)^5+2(-1)^6=1$ after bifurcation.

\begin{table}[h]\fontsize{12}{12} \centering
\caption{Data for the $g'$ family.}
	\begin{tabular}{|c|c|c|c|c|c|}
 \hline
		Values of $\Gamma$ & planar & spatial & $\mu_{CZ}^p$ & $\mu_{CZ}^s$ & $\mu_{CZ}$  \\
		\hline  
$(4.49999, 4.2851)$ & elliptic & elliptic & 3 & 3 & 6 \\
$(4.2851, 4.2806)$  & elliptic & neg.\ hyperbolic & 3 & 3 & 6 \\
$(4.2806, 4.2714)$ & elliptic & elliptic & 3 & 3 & 6 \\
$(4.2714, 3.3901)$ & neg.\ hyperbolic & elliptic & 3 & 3 & 6 \\
$(3.3901, 0.4771)$ & neg.\ hyperbolic & pos.\ hyperbolic & 3 & 4 & 7 \\
$(0.4771, -0.2195)$ & neg.\ hyperbolic & elliptic & 3 & 5 & 8 \\
$(-0.2195, -4.6921)$ & neg.\ hyperbolic & neg.\ hyperbolic & 3 & 5 & 8 \\
$(-4.6921, -4.7047)$ & elliptic & neg.\ hyperbolic & 3 & 5 & 8 \\
$(-4.7047, -\infty)$ & pos.\ hyperbolic & neg.\ hyperbolic & 4 & 5 & 9 \\

\hline
 \end{tabular}
\label{table:gprime_family}
\end{table}

\begin{figure}
    \centering
    \includegraphics[width=1\linewidth]{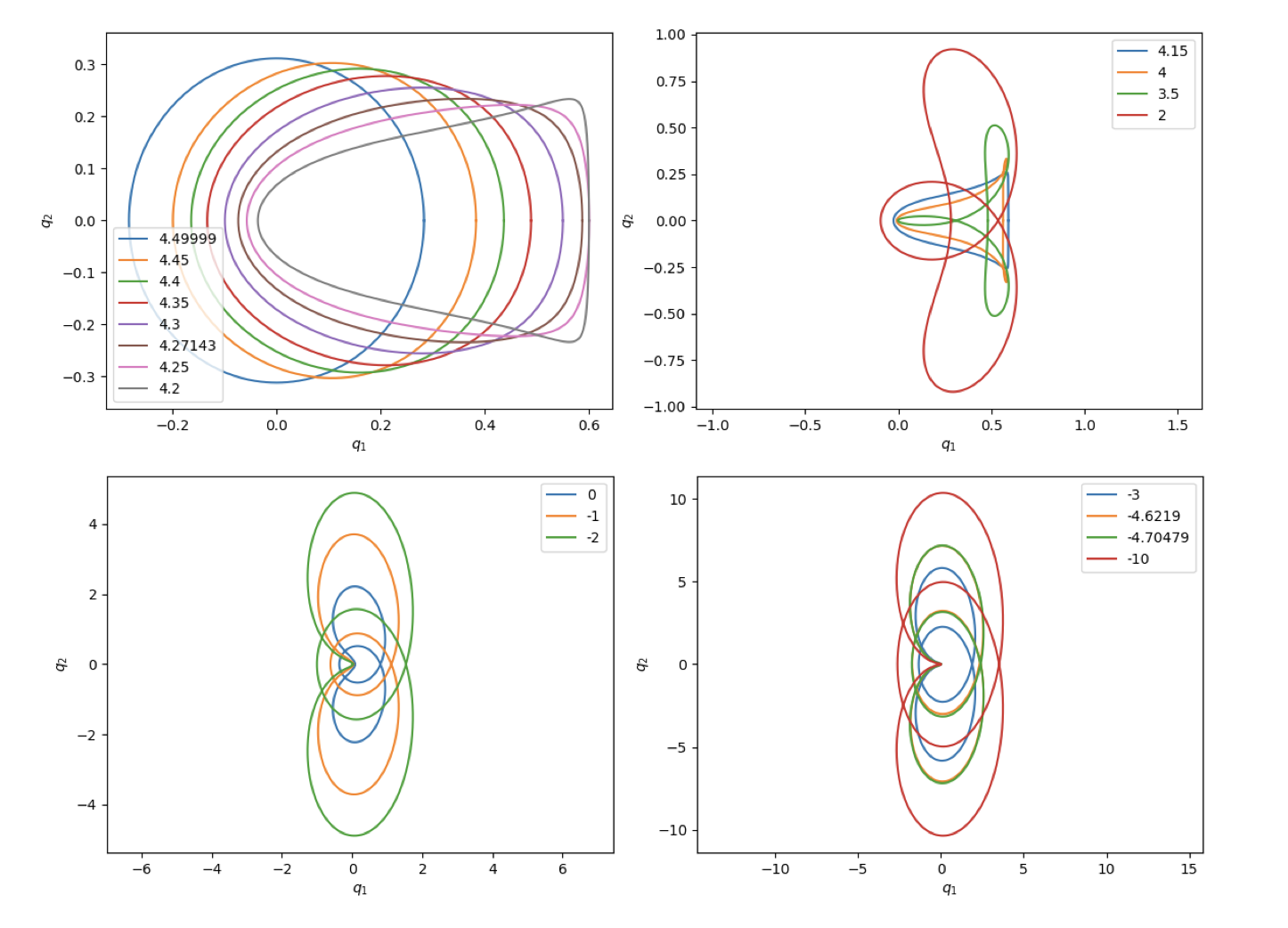}
    \caption{The $g'$ family (plot from \cite{Ay22}).}
    \label{fig:enter-label}
\end{figure}

\begin{figure}
    \centering
    \includegraphics[width=1\linewidth]{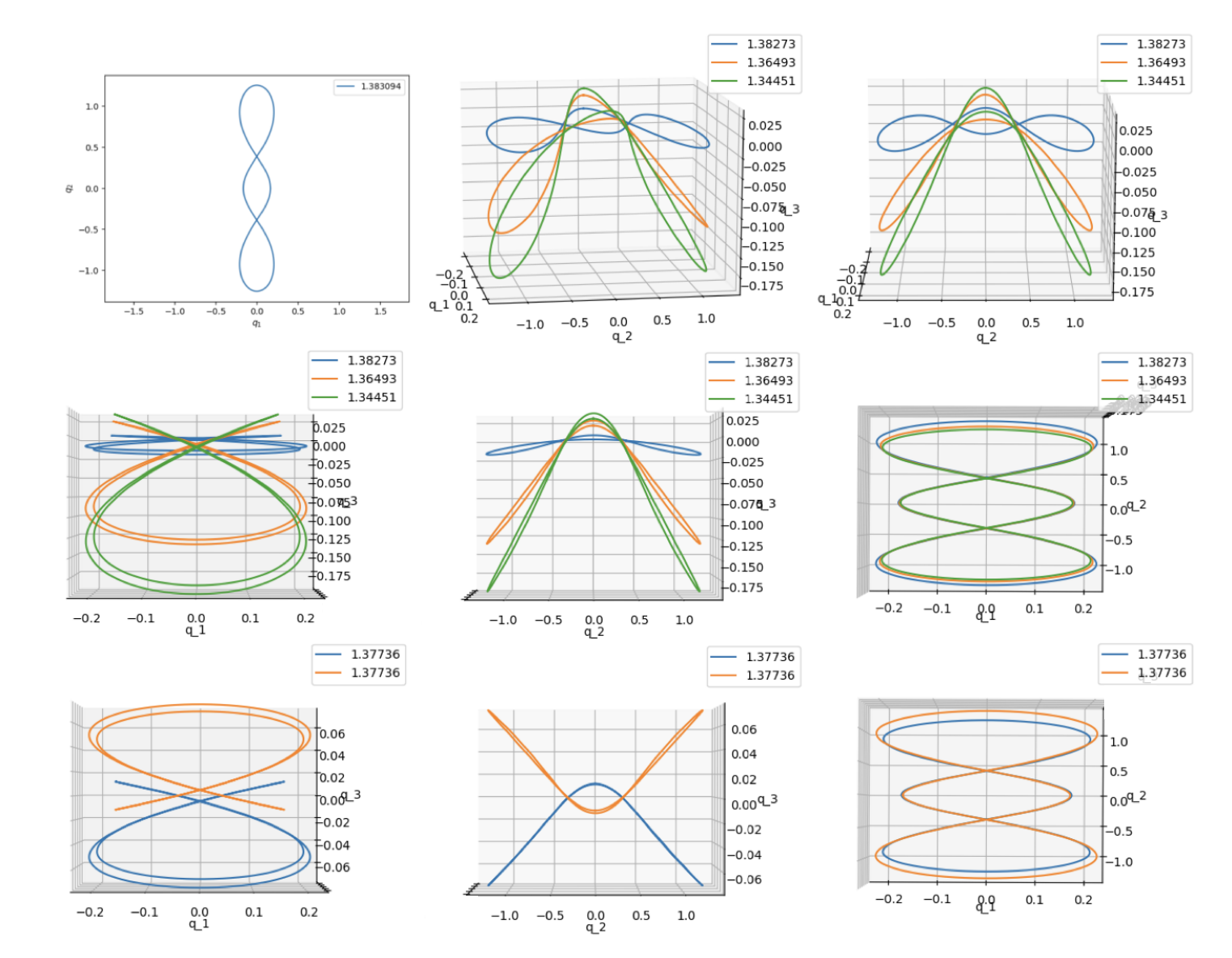}
    \caption{The $g_{2\nu}$ family (plot from \cite{Ay22}).}
    \label{fig:g2nu_family}
\end{figure}

From Table \ref{table:g_family}, we see that the family $g$ undergoes a spatial transition from elliptic to positive hyperbolic at $\Gamma=1.3289$, where $\mu_{CZ}^s$ jumps from $3$ to $4$. At this point, the \emph{spatial} family $g_{2\nu}$ appears, which is plotted in Figure \ref{fig:g2nu_family}. Its orbits are doubly-symmetric with respect to $\rho_1$ and $\rho_2$ (see Section \ref{sec:symmetries_Hill}), and by reflection along the ecliptic we obtain its symmetric family, which is also doubly-symmetric. The CZ-index for this family is $\mu_{CZ}=5$.

\medskip

\begin{figure}
    \centering
    \includegraphics[width=0.9\linewidth]{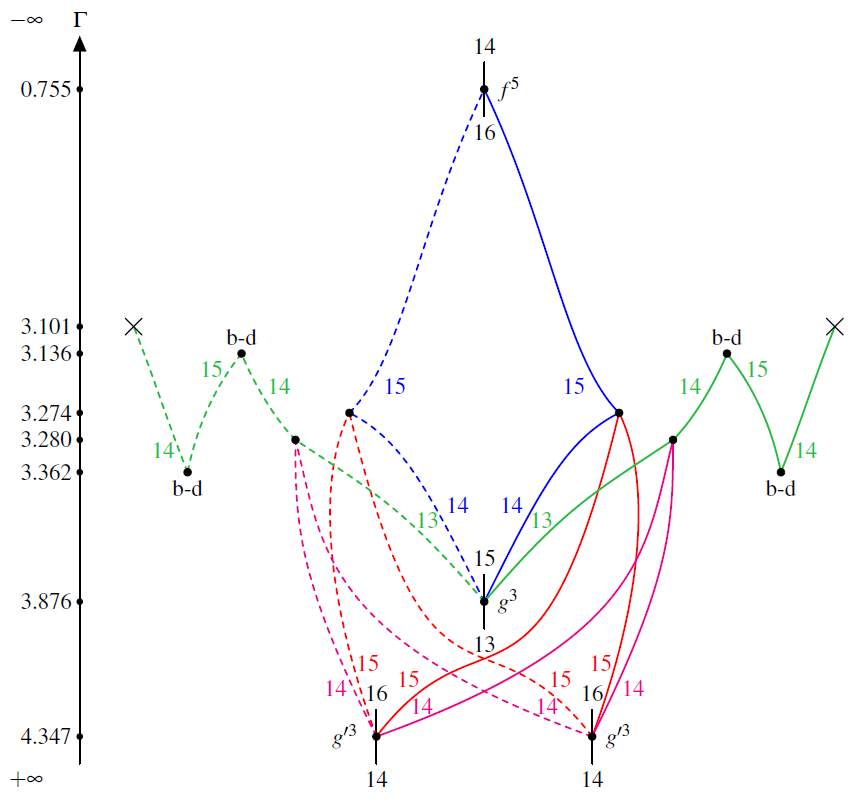}
    \caption{Bifurcation graph between the 3rd cover of $g$, the 3rd cover of $g'$, and the 5th cover of $f$ with the families $f_g^{(2,3)}$ (blue) , $f_g^{(2cut,3)}$ (green) , $f_{g'}^{(2,3)}$ (red), $f_{g'}^{(2cut,3)}$ (pink) \cite{Ay22}.}
    \label{fig:bif_diagram_Hill}
\end{figure}

The remaining families of Table \ref{orbits} bifurcate from the corresponding cover of the underlying planar orbit, as explained there. The easiest way to visualize this information is via bifurcation graphs. For instance, see Figure \ref{fig:bif_diagram_Hill}. There, the direction of increasing energy is from bottom to top. Each family has an attached CZ-index. The edges in black represent the underlying simple planar families where bifurcation occurs. The dashed families are obtained from the non-dashed ones by reflection along the ecliptic, which explains the horizontal symmetry. A cross means collision, and b-d means birth-death (note that the Floer number is always zero at birth-death).

At $\Gamma=3.274$, the blue non-dashed family is doubly symmetric with respect to $\overline{\rho}_1,\overline{\rho}_2$, and the two red non-dashed families are related to each other by $\overline{\rho}_2$, each being symmetric with respect to $\overline{\rho}_1$. In other words, the situation is completely analogous as to the pitchfork bifurcation relating $g$ with $g'$. The same happens at the bifurcation relating the pink and red families. This point will be relevant for Jupiter--Europa later.

We emphasize that the graphs organize the local bifurcations, and help to check that the Floer numbers match up. For instance, at $\Gamma=3.876$, and at $g^3$, the Floer numbers before and after are respectively $(-1)^{13}=-1$, and $2(-1)^{13}+2(-1)^{14}+(-1)^{15}=-1$. But note that at $\Gamma=0.755$, and at $f^5$, the Floer number number before is $2(-1)^{15}+(-1)^{16}=-1$, but after it is $(-1)^{14}=1$. Therefore we conclude that there are orbits still missing, that bifurcate out of $f^5$. 

We refer to \cite{Ay22} for many more plots and bifurcation graphs between the families in Table \ref{orbits}.

\section{Jupiter--Europa system: numerical work}

\begin{figure}
    \centering
    \includegraphics[width=0.8\linewidth]{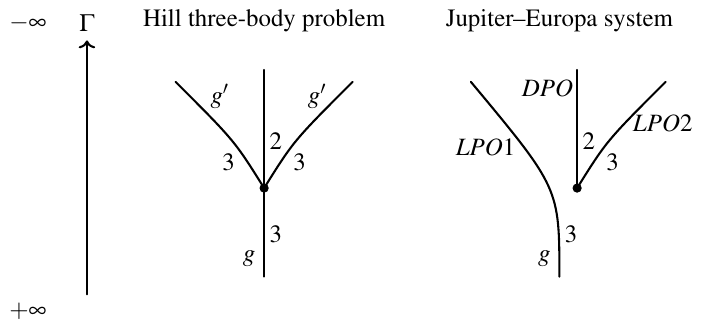}
    \caption{Deforming the pitchfork bifurcation in Hill's lunar problem to a ``broken'' bifurcation in Jupiter-Europa.}
    \label{fig:pitchfork2}
\end{figure}

By deforming the mass parameter $\mu$, we may go from Hill's lunar problem to the Jupiter-Europa system. We wish to understand how the periodic orbits from the next section deform with varying $\mu$. The reference for this section is \cite{AFvKKM}.

\begin{figure}[h]
    \centering
\includegraphics[width=0.8\linewidth]{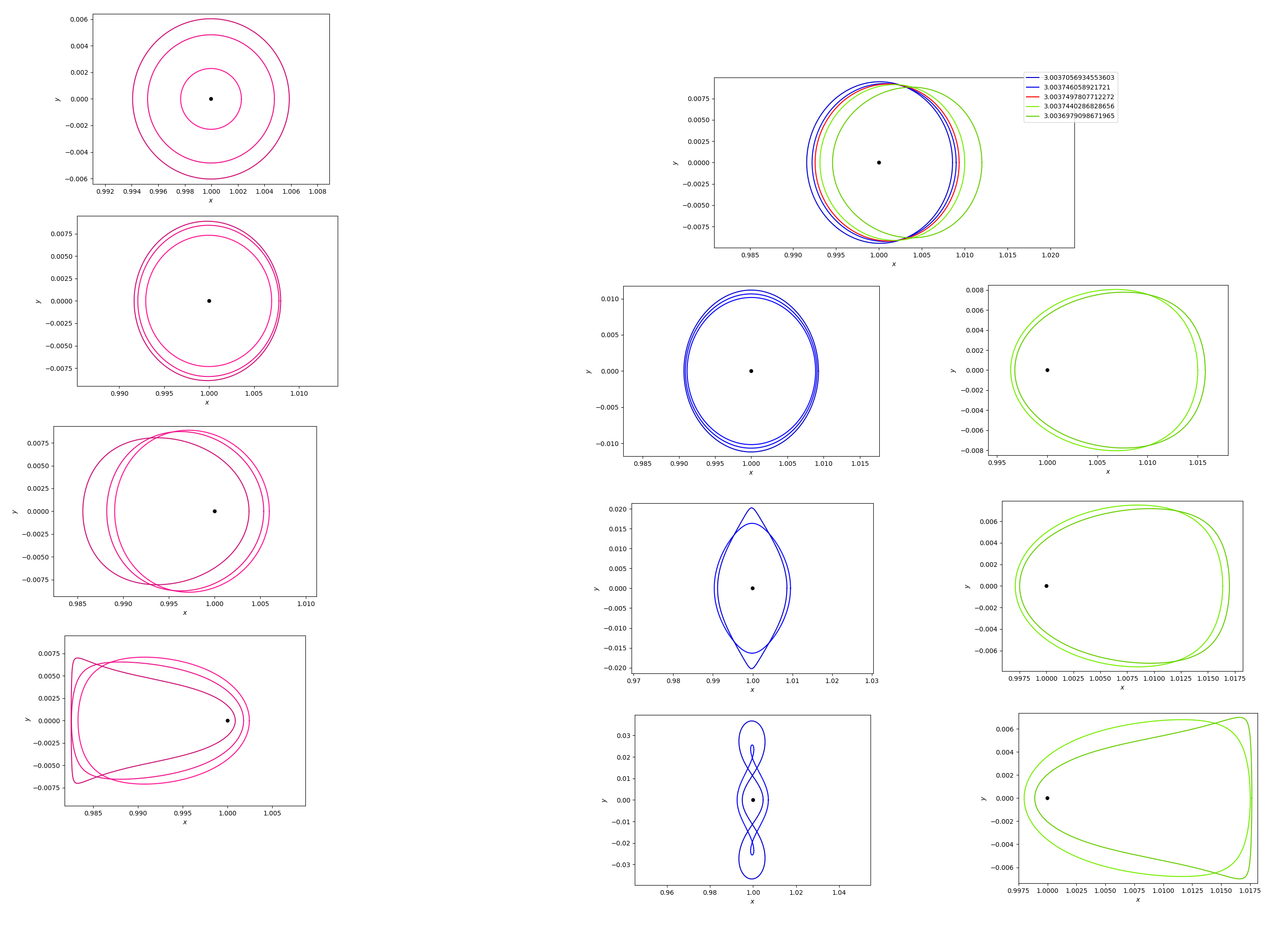}
    \caption{The ``broken'' bifurcation in Jupiter-Europa. Left: the g-LPO1 branch. Right: the DPO-LPO2 branch, split into the DPO sub-branch (left) and of the LPO2 sub-branch (right).}
    \label{fig:broken_bifurcation}
\end{figure}

\subsection{A ``broken'' bifurcation} The pitchfork bifurcation in Figure \ref{fig:pitchfork} deforms to a \emph{generic} situation, where one of the $g'$ branches glues to the before-bifurcation part of the $g$ branch, the result of which we call the \emph{g-LPO1} branch, and where the other $g'$ branch glues to the after-bifurcation part of the $g$ branch, which we call the \emph{DPO-LPO2} branch (undergoing birth-death bifurcation); see Figure \ref{fig:pitchfork2}. The $DPO$-orbits are planar positive hyperbolic and the $LPO2$-orbits are planar elliptic. As the symmetry with respect to the $y$-axis is lost, the new orbits will be \emph{approximately} symmetric with respect to the $y$-axis, but not exactly symmetric; similarly, the $y$-symmetric relation between the $g'$ branches persists only approximately for the corresponding deformed orbits. See Figure \ref{fig:broken_bifurcation}, where this phenomenon is manifest. Via this bifurcation analysis, one may predict the existence of the DPO-LPO2 branch, which a priori is not straightforward to find. While these families are already known and appear e.g.\ in page 12 of \cite{RR17}, this suggests a general mechanism which we will exploit.

\subsection{Spatial pitchfork bifurcation} 
\begin{figure}[h]
	\centering
	\includegraphics[]{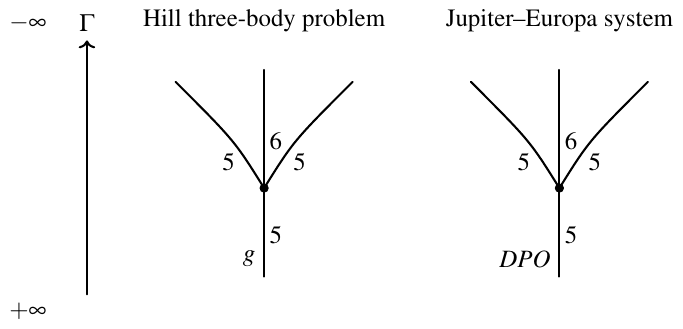}
	\caption{Left:\ The bifurcation graph between simple closed $g$-orbit and the new families of spatial orbits generated by the spatial index jump in Hill's system. Right: In the Jupiter--Europa system. The horizontal symmetry corresponds to the reflection at the ecliptic.}
	\label{bifurcation_graph_g_spatial_index_jump}
\end{figure}

\begin{figure}[h]
    \centering
    \includegraphics[width=0.9\linewidth]{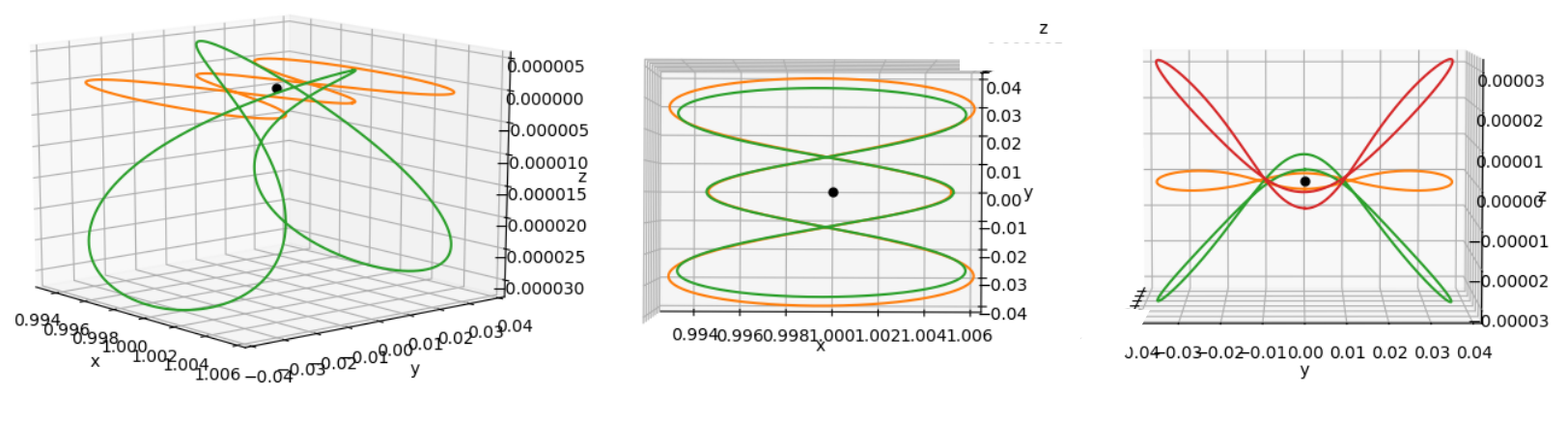}
    \caption{Jupiter-Europa: A planar-to-spatial bifurcation of a simple closed planar DPO orbit, from the side, from above and with its symmetric family by using the reflection at the ecliptic.}
    \label{fig:DPO}
\end{figure}

The spatial CZ-index of the simple closed DPO-orbit jumps by $+1$ at around $\Gamma=3.0011$, and therefore generates a planar-to-spatial bifurcation; see Figure \ref{fig:DPO}. As in Hill's problem, this new family of spatial orbits appears twice by using the reflection at the ecliptic. Surprisingly, compared to Figure \ref{fig:pitchfork2}, because the symmetry is preserved, the bifurcation graph has the same topology after deformation and is still non-generic; see Figure \ref{bifurcation_graph_g_spatial_index_jump}.

\subsection{Spatial bifurcation graph in Jupiter-Europa} 
\begin{figure}[h]
    \centering
    \includegraphics[width=0.65\linewidth]{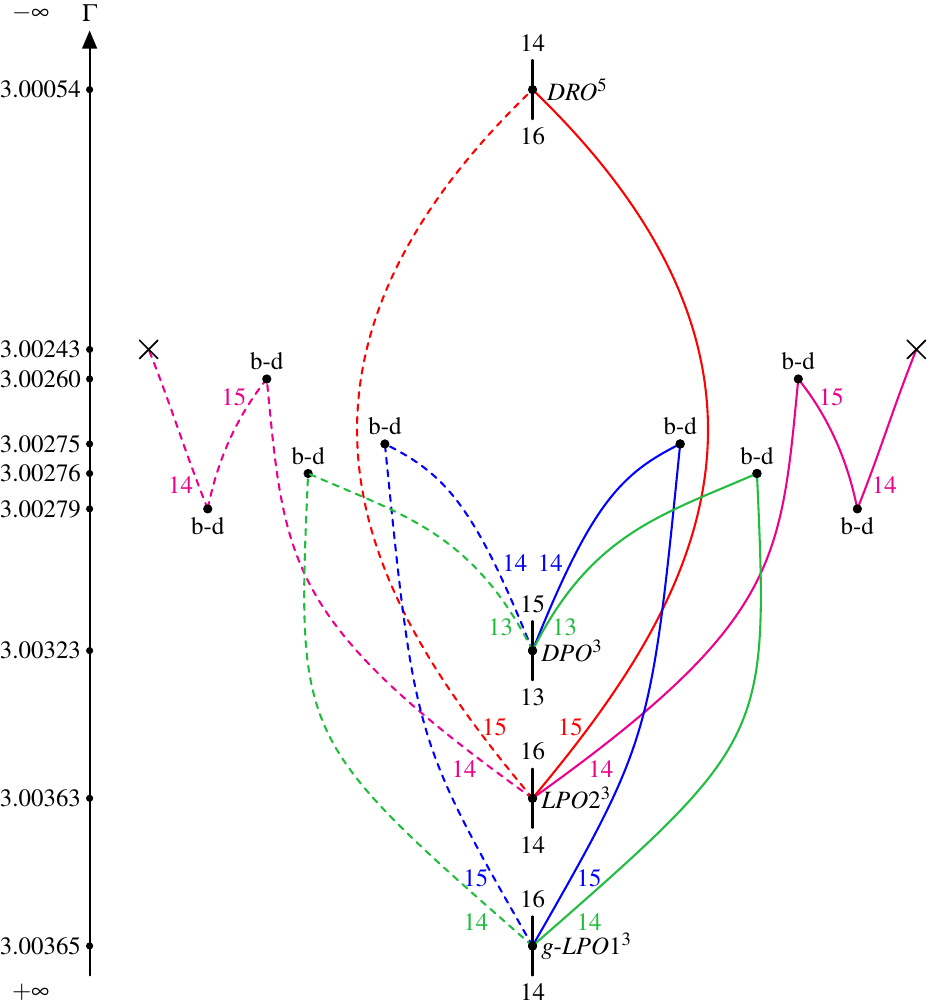}
    \caption{Bifurcation graph for the Jupiter--Europa system, between the 3rd cover of $g$-LPO1, the third cover of DPO,  the third cover of LPO2, and the 5th cover of DRO. The horizontal symmetry is reflection along the ecliptic.}
    \label{fig:bifurc_diagram_new}
\end{figure}

Taking Figure \ref{fig:bif_diagram_Hill} as a starting point, we deform it to the Jupiter--Europa system. The result is plotted in Figure \ref{fig:bifurc_diagram_new}. 

Let us focus on the two vertices on the center-right of Figure \ref{fig:bif_diagram_Hill} which are not of birth-death type. After deformation, the (red) family starting at $g'^3$ on the right of CZ-index 15 glues to the (blue) family of the same index ending in $f^5$, resolving the vertex at which they meet, as in Figure \ref{fig:bifurc_diagram_new}; note that $f$ is replaced with DRO, and $g'$, with LPO2. The two other families meeting at the same vertex coming from $g'^3$ and $g^3$ now glue to a family undergoing birth-death, where now $g'$ is replaced by $g$-LPO1, and $g$, with DPO. A similar phenomenon happens at the other vertex, where the (pink) family starting at $g'^3$ with CZ-index 14 on the right glues to the (green) family of the same index, and the other two families now undergo birth-death. This is exactly the ``broken'' bifurcation phenomenon which we observed in Figure \ref{fig:broken_bifurcation} (the picture is simply ``upside down''). This is the way in which we predicted the existence of the blue and green birth-death branches, which might have been hard to detect via different methods. 


\begin{figure}[h]
    \centering
    \includegraphics[width=0.86\linewidth]{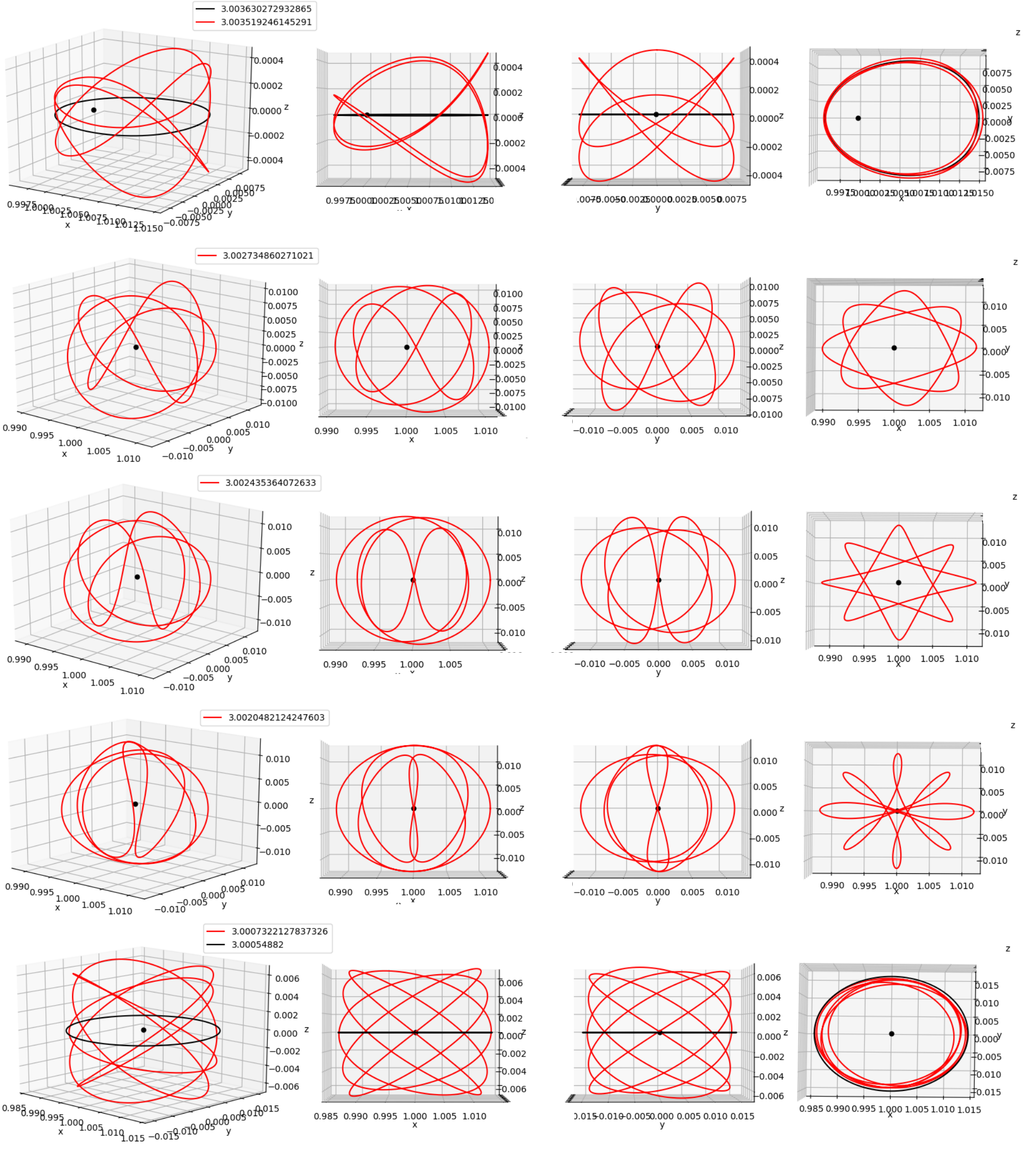}
    \caption{The red prograde to retrograde spatial connection, with CZ-index 15, for varying energy. The black planar orbit is an LPO2, and the red family bifurcates from its third cover. A row corresponds to the same orbit from different angles.}
    \label{fig:pro-ret}
\end{figure}

Another notable feature is the (red) family between LPO2$^3$ and DRO$^5$ of CZ-index 15, a \emph{spatial} family connecting two \emph{planar} orbits, one retrograde (DRO$^5$), and the other, prograde (LPO2$^3$). This family is plotted in Figure \ref{fig:pro-ret}.

\begin{figure}[h]
    \centering
    \includegraphics[width=1\linewidth]{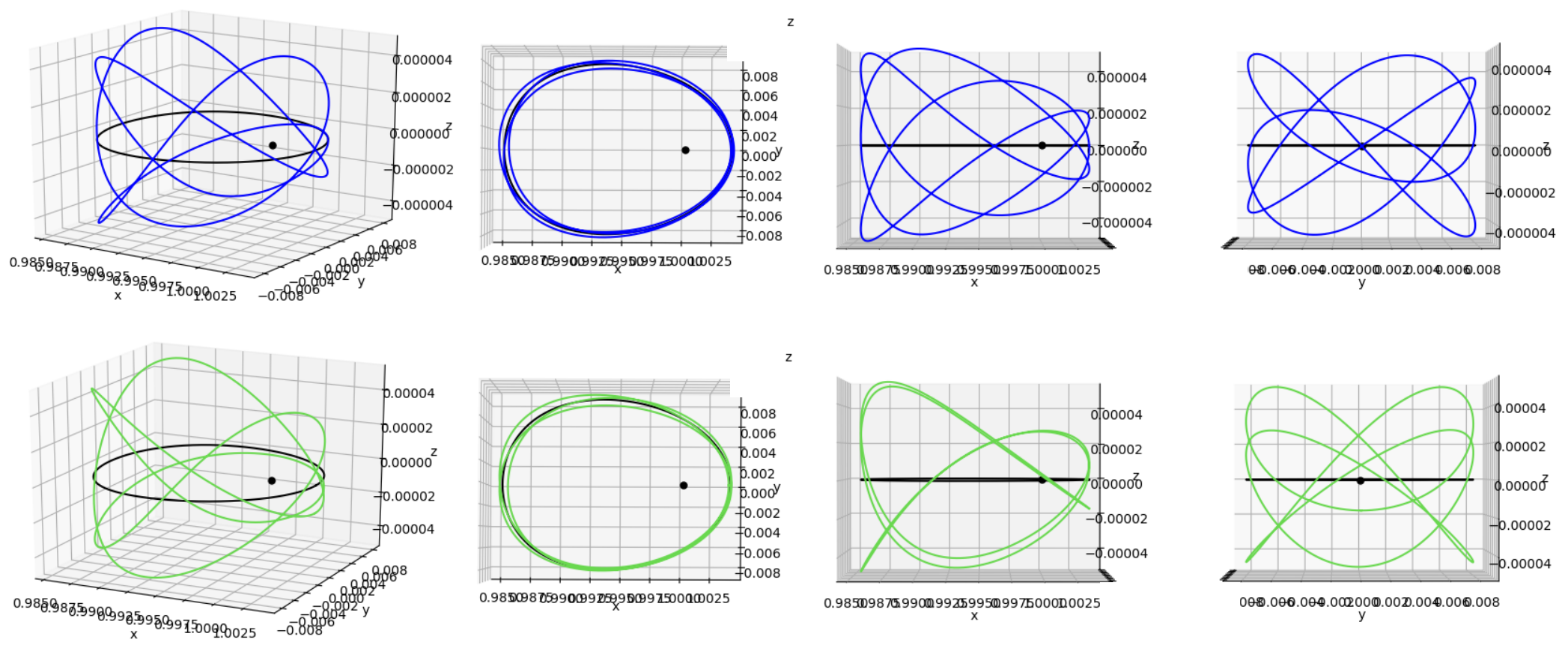}
    \caption{Jupiter-Europa: Two families of spatial orbits branching out from the $g$-LPO1$^3$ orbit; above: these orbits are symmetric wrt the $x$-axis. This is the (blue) family of CZ-index 15 in Figure \ref{fig:bifurc_diagram_new}; below: these orbits are symmetric wrt the $xz$-plane. This is the (green) family of CZ-index 14 in Figure \ref{fig:bifurc_diagram_new}. Each family has a symmetric family by using the reflection at the ecliptic.}
    \label{fig:LPO1_3rd}
\end{figure}

\newpage

\subsection{Spatial bifurcation graph for Saturn-Enceladus} 
\begin{figure}
    \centering
    \includegraphics[width=0.6\linewidth]{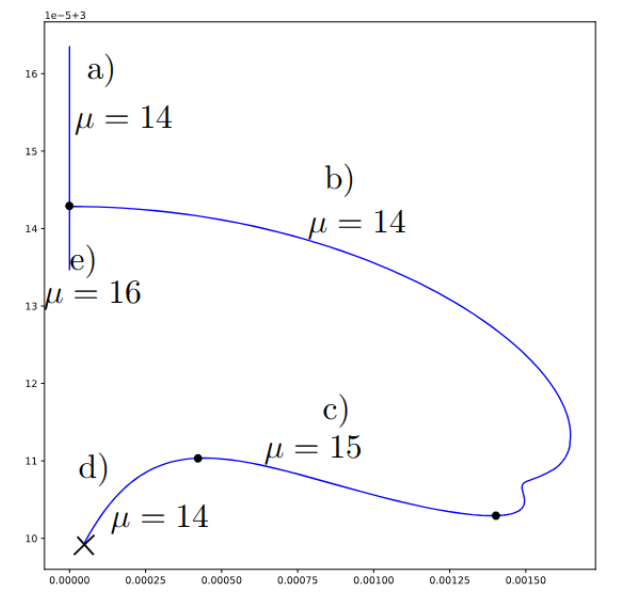}
    \caption{A bifurcation graph for Saturn-Enceladus of $xz$-plane symmetric orbits, which deforms to the pink families in Jupiter-Europa. Horizontal axis is $z$ starting value. Vertical axis is energy.}
    \label{fig:SE}
\end{figure}

\begin{figure}
        \centering
        \includegraphics[width=0.9\linewidth]{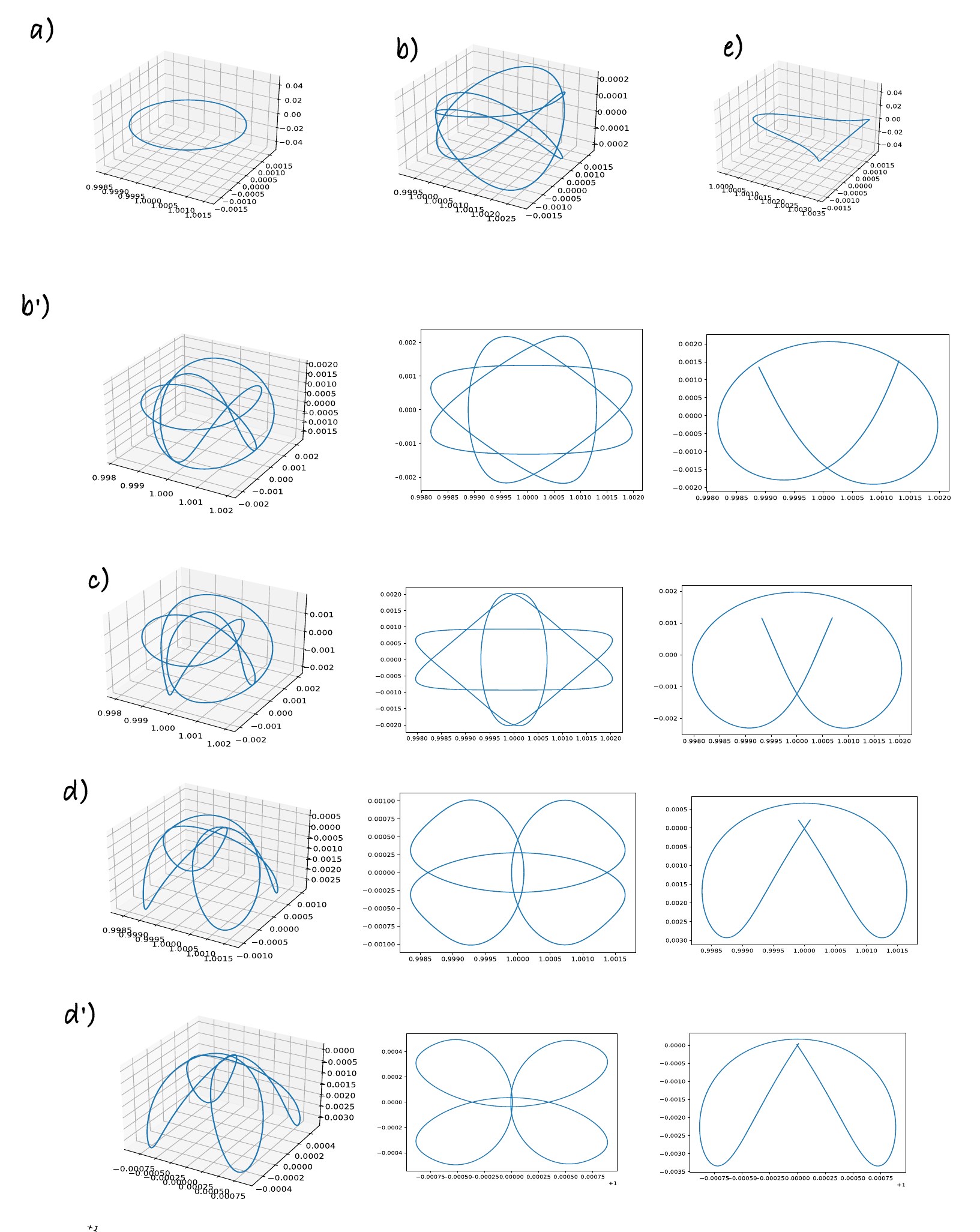}
        \caption{Plots of the orbits represented by the bifurcation graph of Figure \ref{fig:SE}. The b') orbits are part of the b) family, only for different energy: similarly for d) and d'). The second to last rows correspond to the same orbit from different angles.}
        \label{fig:SEplots}
\end{figure}

The bifurcation graph for Saturn-Enceladus corresponding to the one shown in Figure \ref{fig:bifurc_diagram_new} has exactly the same topology (but different energy values). For instance, Figure \ref{fig:SE} shows a bifurcation graph corresponding to the pink families of Figure \ref{fig:bifurc_diagram_new} (but drawn upside down). Note that it is not only topological, as we also record the starting value along the $z$ axis. The corresponding families of orbits are plotted in Figure \ref{fig:SEplots}.

\newpage

\section{GIT plots}

In this section, we illustrate the practical uses of the GIT sequence, with a number of numerical plots, produced by Dayung Koh. This section is based on \cite{FKM}.

\subsection{Planar to spatial period-doubling bifurcation in Jupiter-Europa.} 
\begin{figure}[h]
    \centering
    \includegraphics[width=0.6\linewidth]{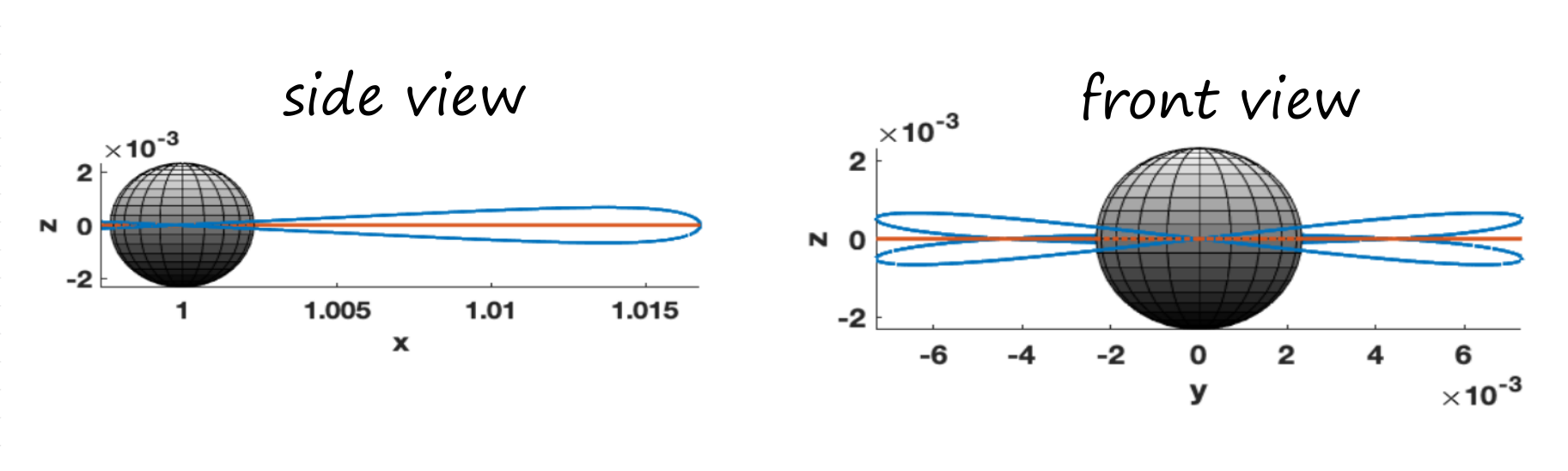}
    \includegraphics[width=0.3\linewidth]{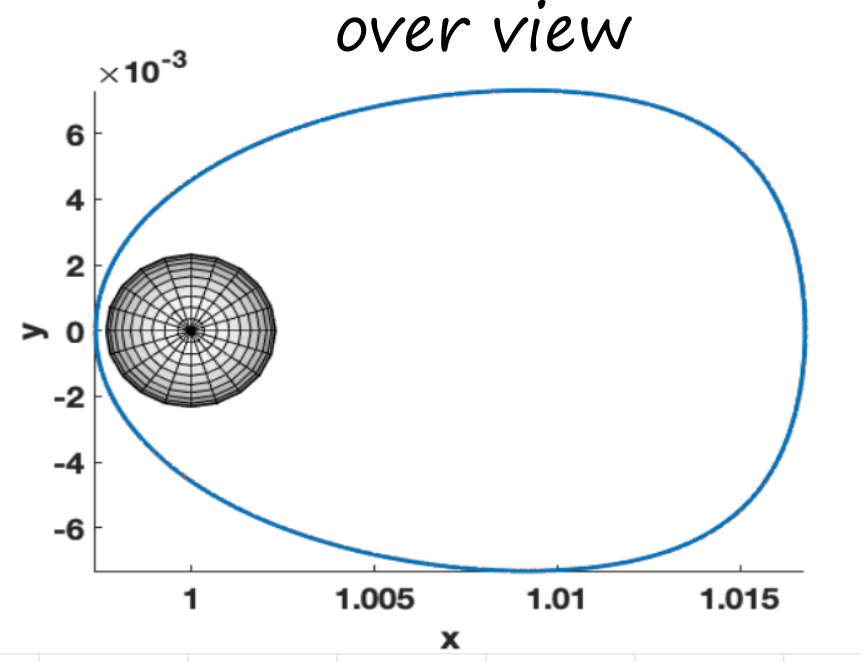}
    \includegraphics[width=0.4\linewidth]{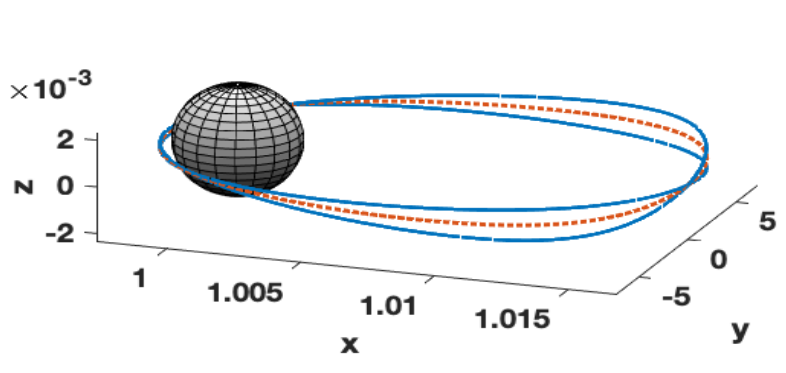}
    \caption{The prograde planar orbit $\gamma_{bef}$ is the dotted orange line; the Jacobi constant is $c=3.00357414$, and its period is $T_0=2.1215$. The spatial orbit $\beta$ of double period (after bifurcation) is the blue one; the Jacobi constant is now $c=3.003571774$, with period $T_1=4.245\approx 2T_0$. We call this the ``snitch'' configuration.}
    \label{fig:prograde}
\end{figure}

The $H_2$ family \cite{KABM} orbit depicted in Figure \ref{fig:prograde} undergoes planar to spatial period-doubling bifurcation. We denote by $\gamma_{bef}$ and $\gamma_{aft}$ the simple orbit before and after the bifurcation and by $\beta$ the orbit with double period appearing after bifurcation. This is a \emph{doubly} symmetric period-doubling bifurcation, as all orbits are doubly symmetric, i.e.\ symmetric with respect to two involutions $\rho_1,\overline{\rho}_1$. Here, $\gamma_{bef}$ is of type $\mathcal{E}^2$, and $\gamma_{aft}$ is of type $\mathcal{EH^-}$. The symmetric points with respect to one of the symmetries are the \emph{fake} points with respect to the other symmetry. Which one is which is determined by the $B$-sign, as explained in Example \ref{ex:sym_double}.

\begin{figure}[h]
    \centering
    \includegraphics[width=1\linewidth]{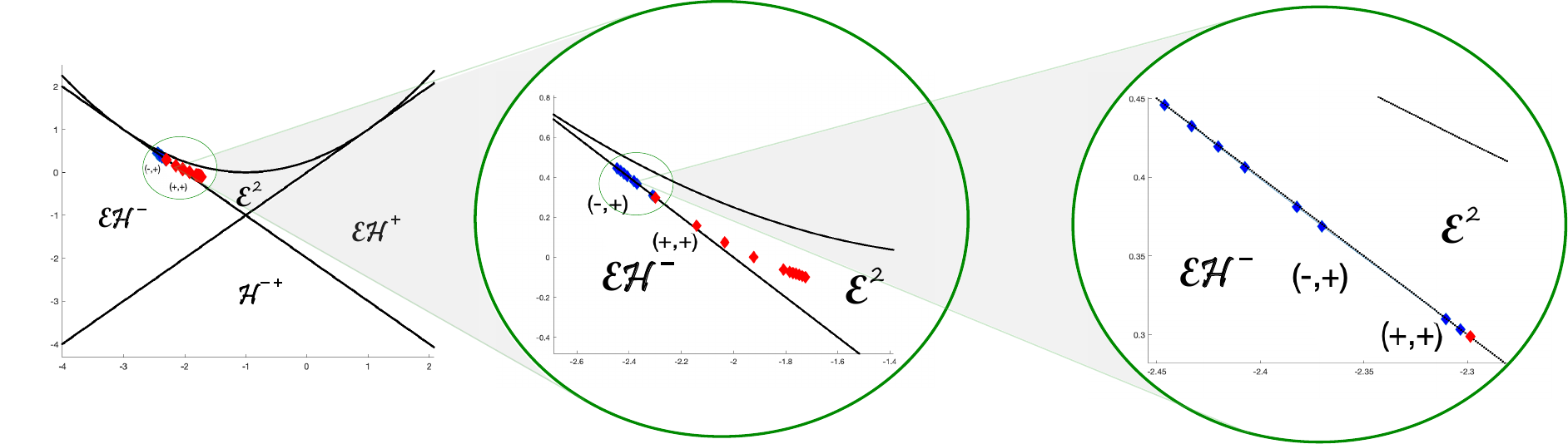}
    \caption{GIT plot of the period-doubling bifurcation of the snitch configuration.}
    \label{fig:GIT}
\end{figure}

Figure \ref{fig:GIT} shows a numerical plot of this period-doubling bifurcation, as seen in the base $\mathbb{R}^2$ of the GIT sequence, in three different scales. The time parameter is the Jacobi constant. Red dots correspond to $\gamma_{bef}$, and blue dots, to $\gamma_{aft}$. The bifurcation takes place when the period-doubling line of slope $-1$ separating the doubly-elliptic region $\mathcal{E}^2$ and the elliptic-negative hyperbolic region $\mathcal{EH}^-$ is crossed. The plot also contains the $B$-signature of the simple orbit, before and after the bifurcation.

\begin{figure}[h]
    \centering
    \includegraphics[width=0.6\linewidth]{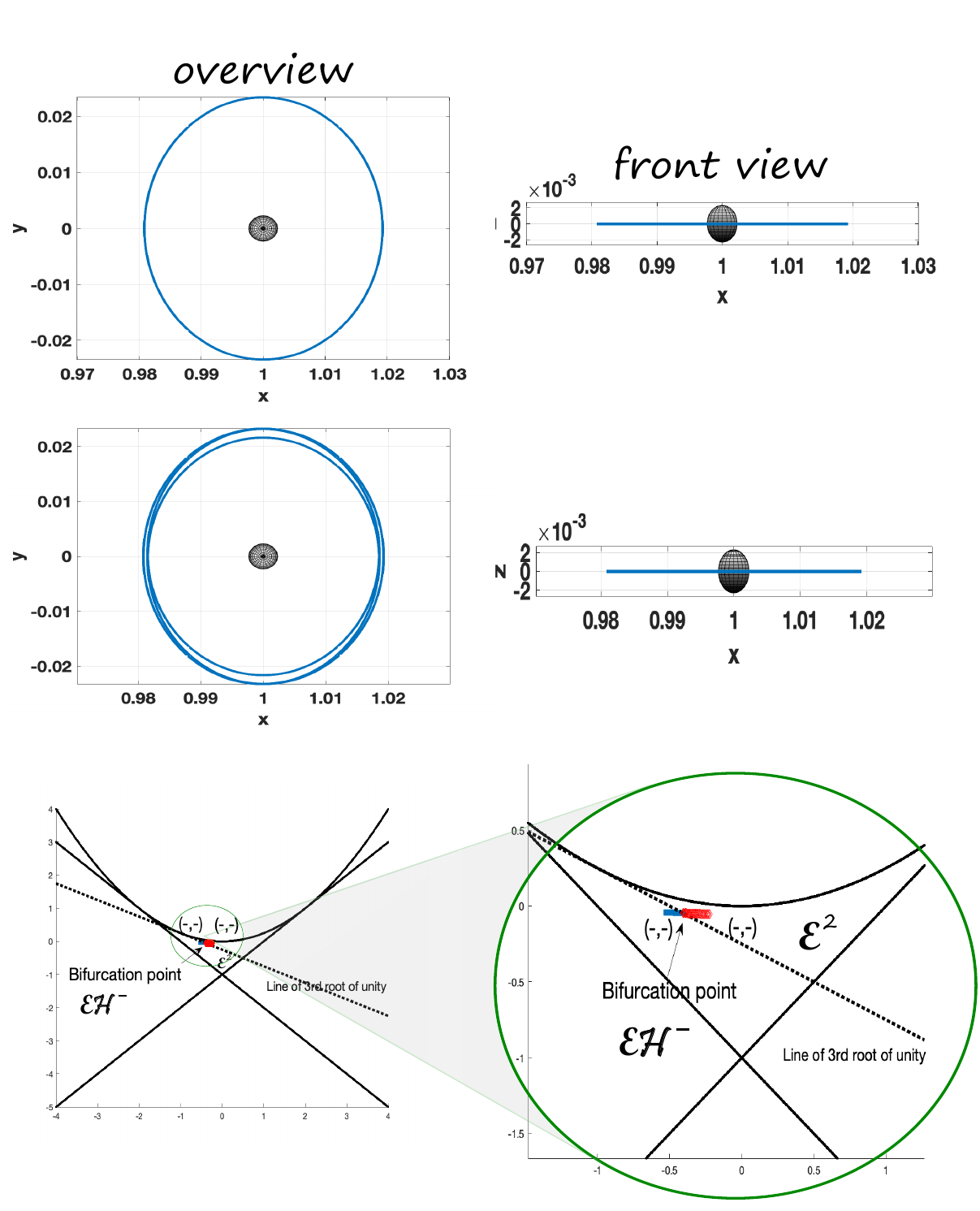}
    \caption{Jupiter-Europa system: a symmetric planar to planar period-tripling bifurcation of DRO family. Above: the planar simple orbit at bifurcation, the \emph{distant retrograde orbit}, $c=2.9999$, $T_0=2.504$. Middle: the planar triple period orbit after bifurcation $c\gtrsim 2.9999$, $T= 7.3\approx 3T_0$. Below: GIT plot, including $B$-signs.}
    \label{fig:JEDRO3}
\end{figure}

\begin{figure}[h]
    \centering
    \includegraphics[width=0.6\linewidth]{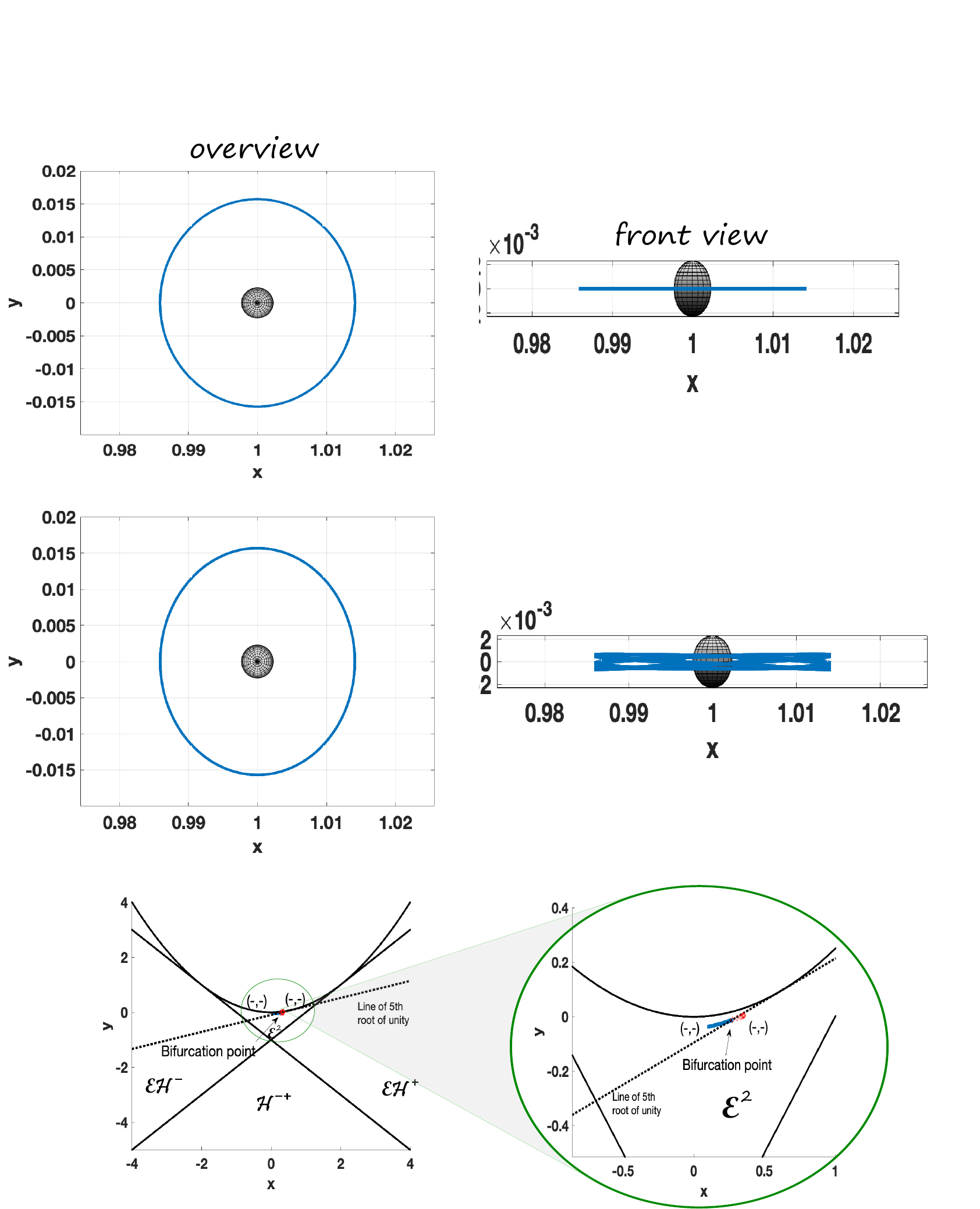}
    \caption{Jupiter-Europa system: a symmetric planar to spatial $5$-fold bifurcation of DRO family. Above: the planar simple orbit at bifurcation, a distant retrograde orbit in the same family of Figure \ref{fig:JEDRO3}, but with  $c=3.0005,$ $T_0=1.705$. Middle: the spatial $5$-fold period orbit after bifurcation, $c\gtrsim 3.0005,$ $T=8.52\approx 5T_0$. Below: GIT plot, including $B$-signs.}
    \label{fig:JE_GIT}
\end{figure}

\begin{figure}[h]
    \centering
    \includegraphics[width=0.5\linewidth]{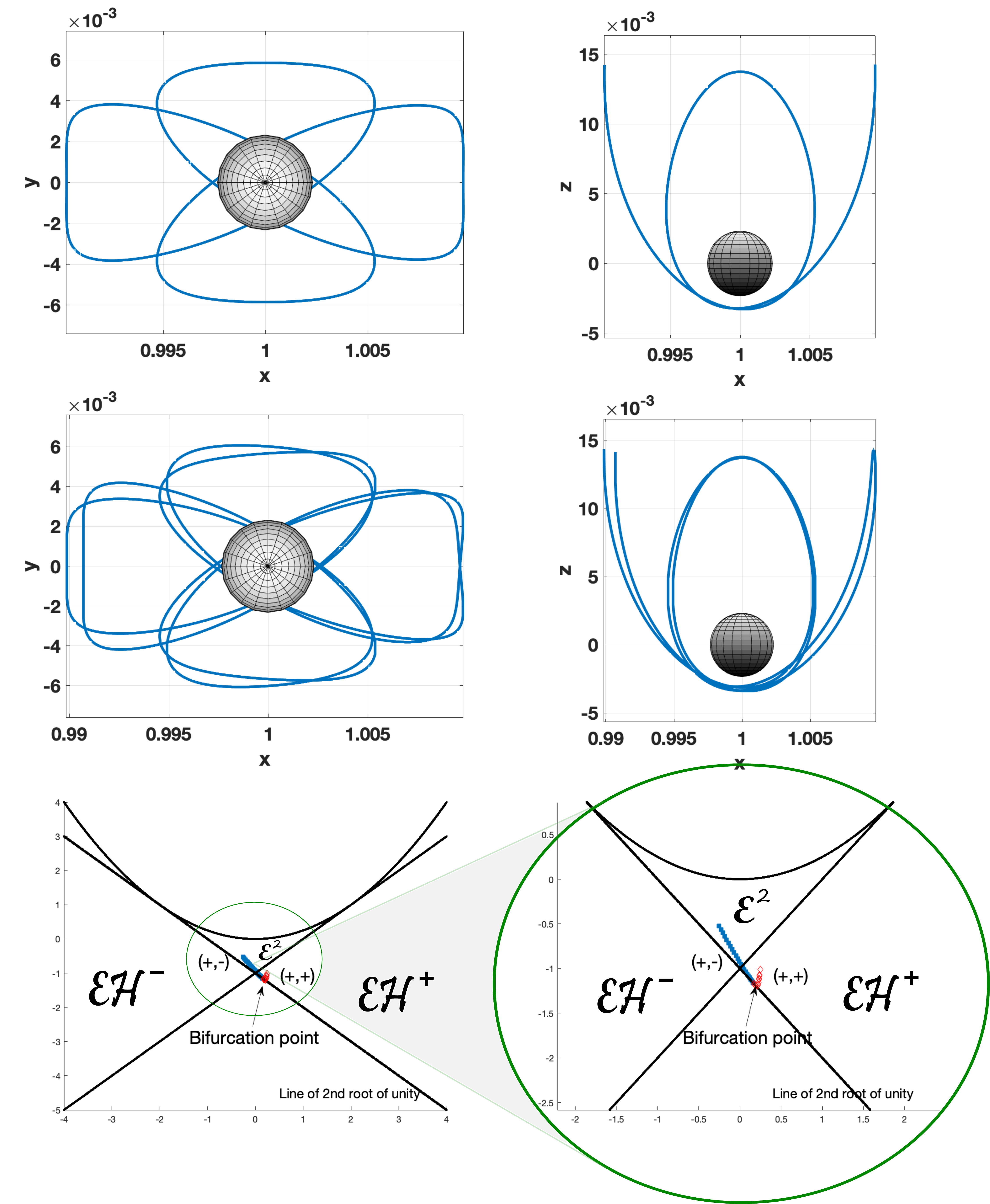}
    \caption{Jupiter-Europa system: a symmetric spatial to spatial period-doubling bifurcation. Above: the spatial simple orbit at bifurcation, $c=3.0028$, $T_0=4.62$. Middle: the spatial period-doubling orbit after bifurcation $c\gtrsim 3.0028,$ $T=9.23\approx 2T_0$. Below: GIT plot, including $B$-signs.}
    \label{fig:JE_GIT2}
\end{figure}

Figures \ref{fig:JEDRO3}, \ref{fig:JE_GIT}, \ref{fig:JE_GIT2} show various examples of bifurcations for the Jupiter-Europa and Saturn-Enceladus systems, with their corresponding GIT plots.

\begin{figure}[h]
    \centering
    \includegraphics[width=0.4\linewidth]{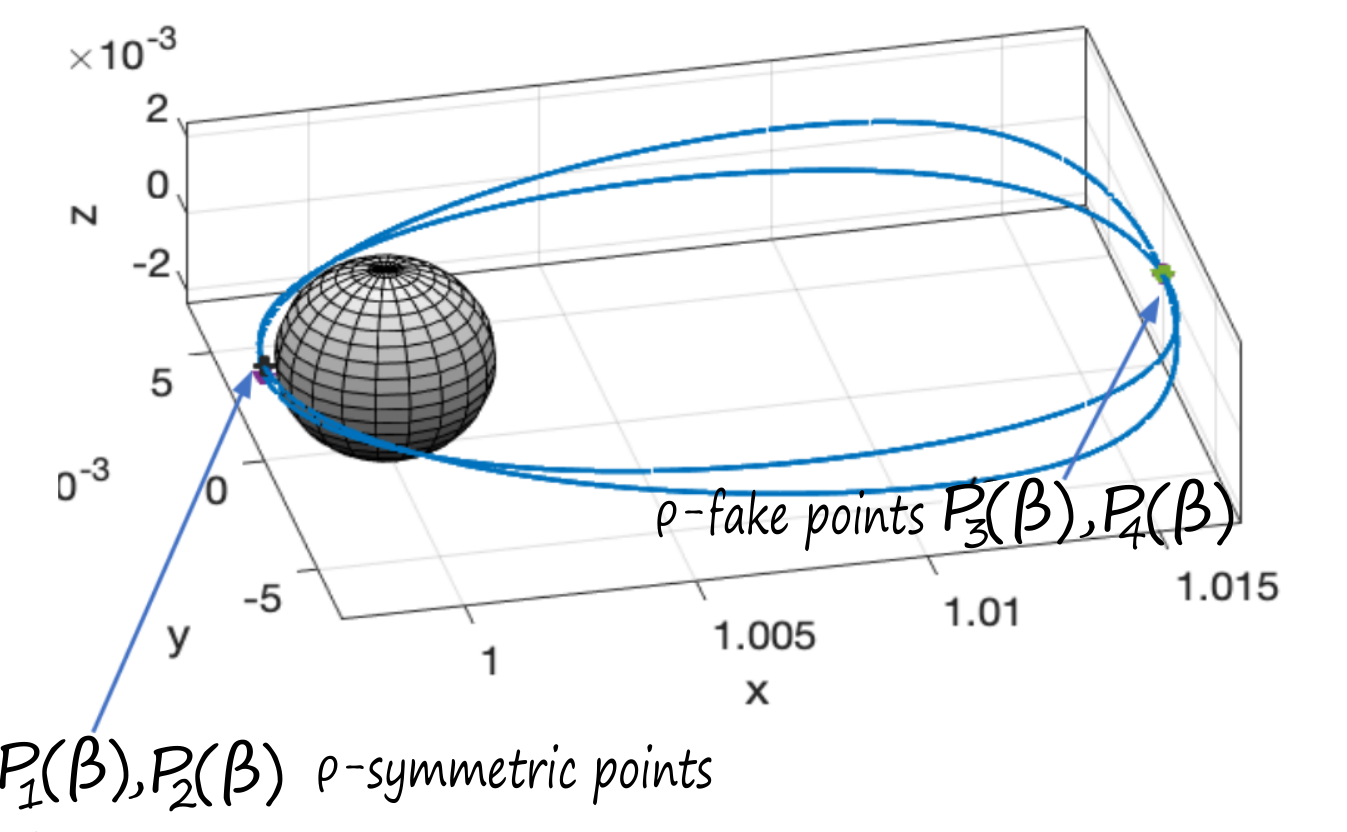}
    \caption{The symmetric points $P_1(\beta),P_2(\beta)$, and the fake points $P_3(\beta),P_4(\beta)$, of the spatial orbit $\beta$, with respect to $\rho=\overline{\rho}_1$. The roles are reversed when $\overline{\rho}_1$ is replaced by $\rho_1$.}
    \label{fig:sym_vs_fake}
\end{figure}

\newpage\newpage\newpage\newpage

\section{The Earth--Moon system: a study of the Halo orbits} The Earth--Moon system (or the study of \emph{cislunar dynamics}) is currently a very popular system of study, among other reasons, because of NASA's ambitious ARTEMIS program. In particular, the Lunar Gateway, a space station meant to serve as communication hub, science laboratory and habitation module, falls withing the umbrella of this agenda. Its implementation, well underway at the time of writing, is also intended as a first step towards deeper space exploration. Moreover, the plan is to put the Gateway in orbit along a specific orbit, which forms part of a family called the \emph{Halo} orbits. 

In general, there is one such family associated to the critical points $L_1,L_2$ and $L_3$ in every CR3BP system. Moreover, because of the symmetry of the CR3BP with respect to reflection along the ecliptic, one always has a Northern family of Halo orbits, passing close to the Northern pole of the small primary, but also a Southern one, passing close to the South pole. As these orbits tend to be unstable, from a practical point of view, they require station-keeping and thruster correction maneuvers, in order to stay in orbit. In practice, these orbits are also amenable for putting scientific instruments like telescopes in orbit, as they have a particular geometry which does not obstruct communication from Earth (the spacecraft is never blocked by the Moon, when observed from Earth, as it ``halos'' the Moon from its opposite side with respect to Earth).

\medskip

\textbf{A brief history of the halo orbits.} The history of these orbits is quite remarkable in it own right. Robert W.\ Farquhar first coined the name ``halo'' in 1966 for orbits around L2 which were forced to be periodic via burners. He also advocated using spacecraft in such an orbit for an Earth–Moon L2, as a communications relay station for an Apollo mission to the far side of the Moon. This is exactly because of its geometry, as explained above. In the end, no relay satellite was launched for Apollo, as landings were on the near side of the Moon.

In 1973 Farquhar and Ahmed Kamel found that a Lissajous orbit ceased to be a Lissajous orbit and became approximately an ellipse, i.e.\ a halo orbit, for which they used analytical expressions. In 1984, Kathleen Howell numerically computed more precisely such trajectories, and found that for most values of the mass ratio there was a range of stable orbits.

The first spacecraft to utilize a halo orbit was ISEE-3, a collaborative mission between ESA and NASA, launched in 1978. It traveled to the Sun–Earth L1 point, where it remained for several years. The next mission to employ a halo orbit was the Solar and Heliospheric Observatory (SOHO), another joint ESA/NASA project dedicated to studying the Sun. SOHO arrived at the Sun–Earth L1 point in 1996, adopting an orbit similar to that of ISEE-3. Since then, while several missions—such as the Gaia astrometric observatory—have ventured to Lagrange points, they typically use related, non-periodic paths called Lissajous orbits rather than true halo orbits.

Although halo orbits were well understood within the framework of the CR3BP, finding such orbits for the actual Earth-Moon system posed significant challenges. In 1998, M.A.\ Andreu successfully computed translunar halo orbits by introducing a novel model for spacecraft motion in the Earth-Moon-Sun system, known as the Quasi-Bicircular Problem (QBCP).

In May 2018, Farquhar's original concept came to fruition when China launched the first communications relay satellite, Queqiao, into a halo orbit around the Earth-Moon L2 point. On January 3, 2019, the Chang'e 4 spacecraft successfully landed in the Von Kármán crater on the far side of the Moon, relying on the Queqiao satellite to maintain communication with Earth. Later, on January 24, 2022, the James Webb Space Telescope entered a halo orbit around the Sun-Earth L2 point. Euclid followed suit, entering a similar orbit around the same point in August 2022.

India's space agency, ISRO, launched Aditya-L1 to study the Sun from a halo orbit around the L1 point. On January 6, 2024, the Aditya-L1 spacecraft, India's first solar mission, successfully reached its final orbit. This orbit has a period of about 180 days and is located around the first Sun-Earth Lagrangian point (L1), approximately 1.5 million kilometers from Earth.

\medskip

\textbf{Results: new connections.} With this motivation in mind, the goal of this section is to carry out a self-contained study of periodic orbits families which are near or around the Halo orbits. To the best of the author's knowledge, several of the connections we present here are novel, and can potentially be of use when navigating in and out of the Halo families. All plots in this section were produced by Cengiz Aydin and Bhanu Kumar.

\begin{figure}
    \centering
    \includegraphics[width=1\linewidth]{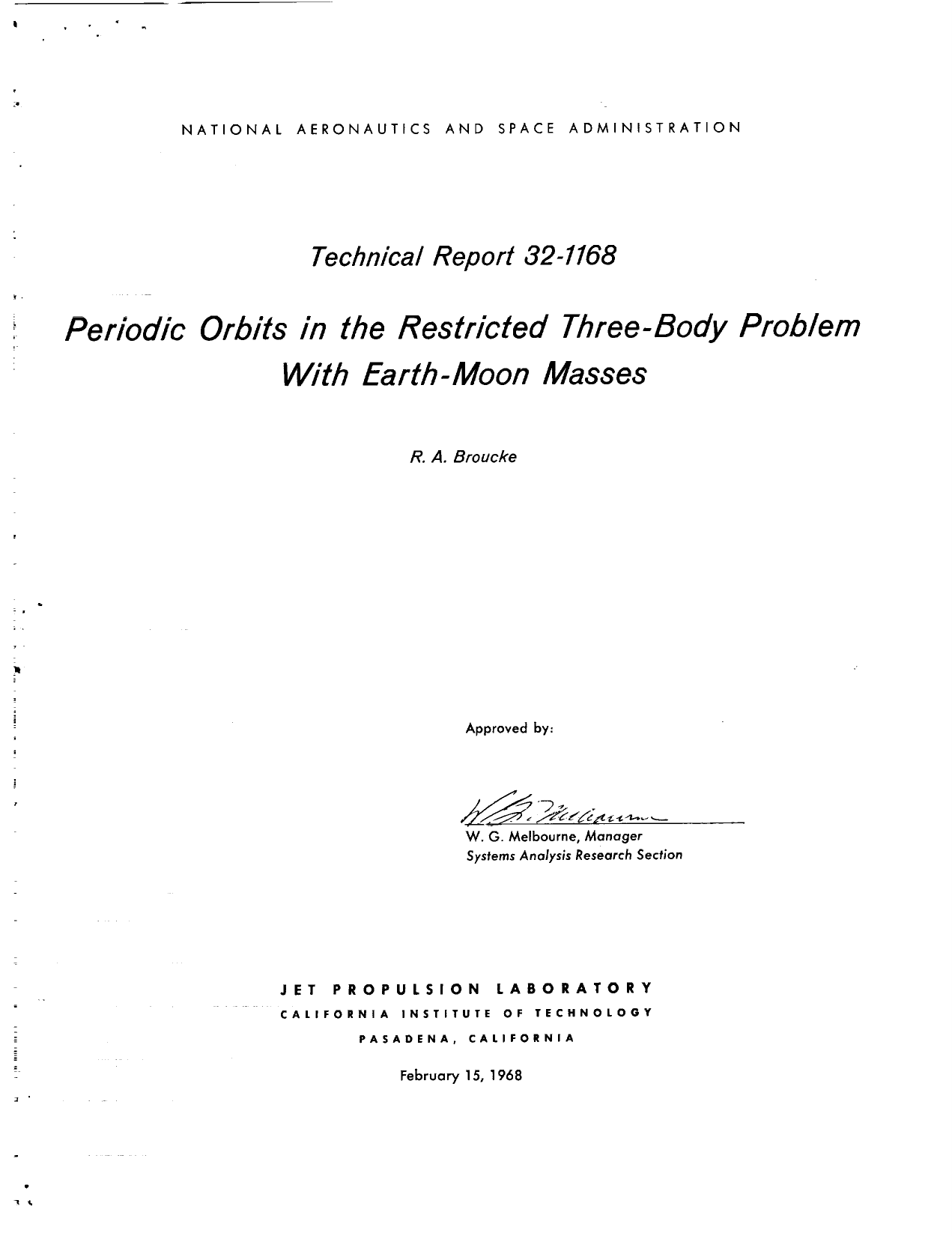}
    \caption{Broucke's original JPL Technical Report, 1968.}
\label{Broucke}
\end{figure}

The starting point for the study is the family of planar orbits for the Earth-Moon system originally studied by Broucke, all of which can be found in a Technical Report of JPL NASA dating back to 1968 (see Figure \ref{Broucke}). We will then be interested in finding out-of-plane bifurcations from these families, from which new spatial families appear, some of which connect to the well-known Halo orbits. The pattern of families which arises is quite rich, and brings in a completely novel web of possibilities for moving in and out regions of phase-space containing the Halo orbits, using natural dynamics in combination with small $\Delta V$ (if maneuvers are carried out at bifurcation points).

\begin{figure}[t]
    \centering
    \includegraphics[width=0.5\linewidth]{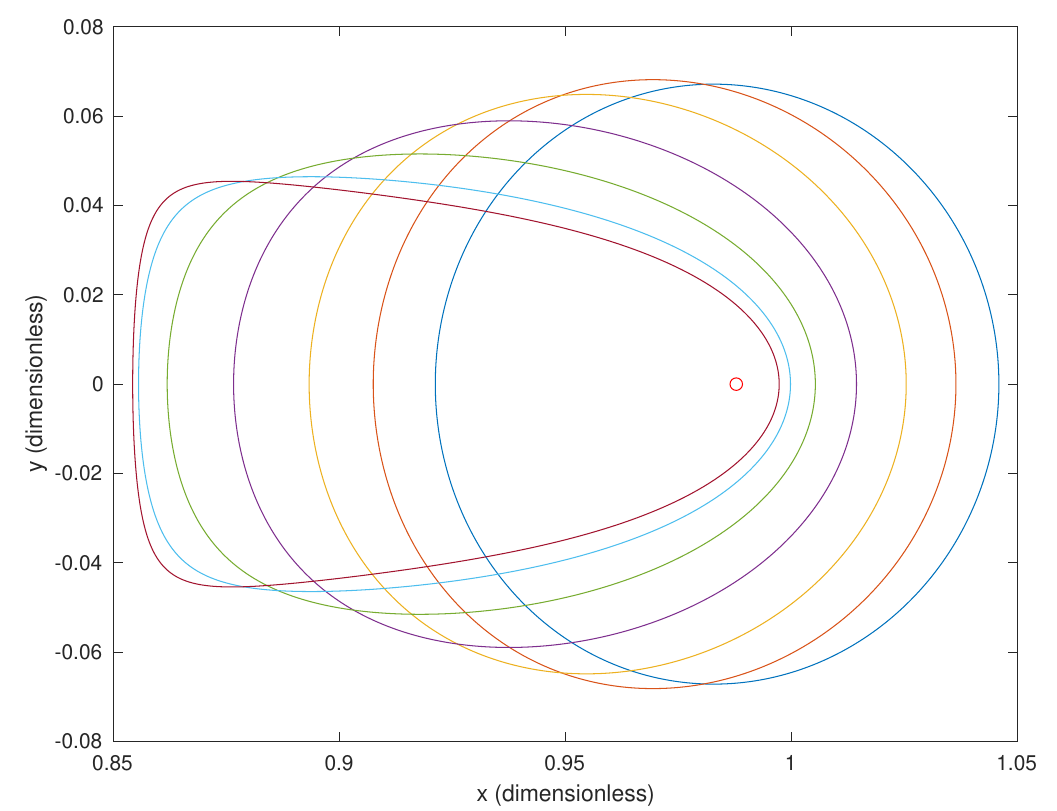}
    \caption{The H1 family of planar lunar orbits found by Broucke (Part I).}
    \label{H1I}
\end{figure}

\begin{figure}[b]
    \centering
    \includegraphics[width=0.5\linewidth]{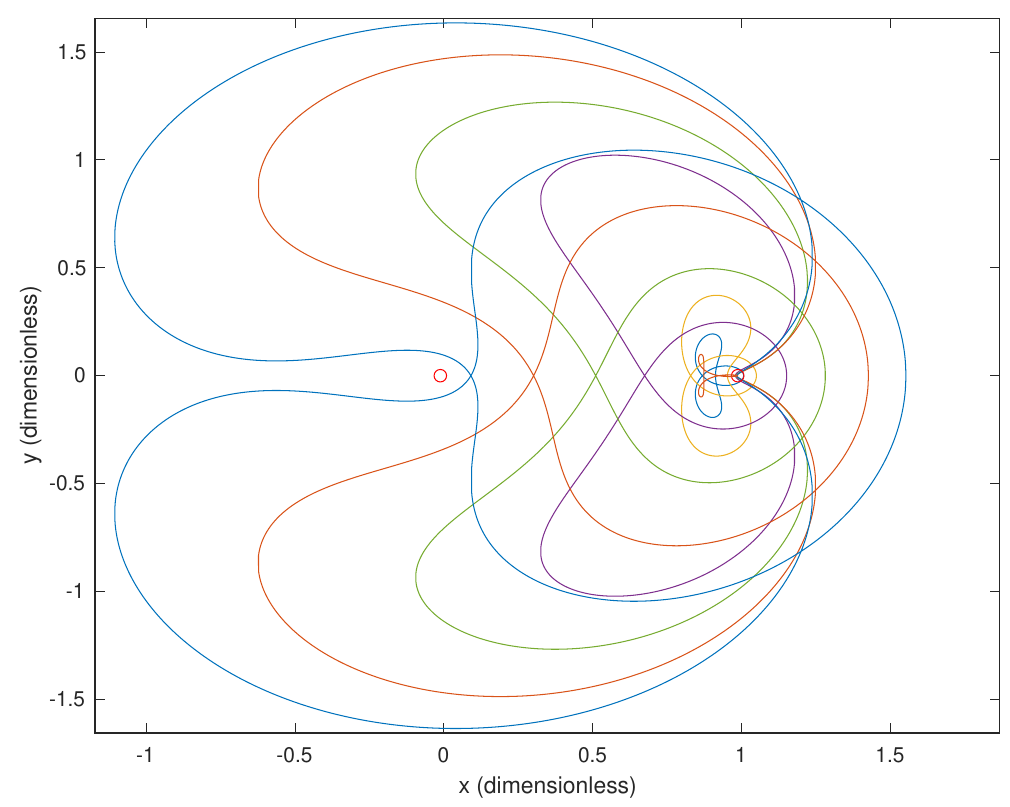}
    \caption{The H1 family of planar lunar orbits found by Broucke (Part II).}
    \label{H1II}
\end{figure}

\begin{figure}
    \centering
    \includegraphics[width=0.5\linewidth]{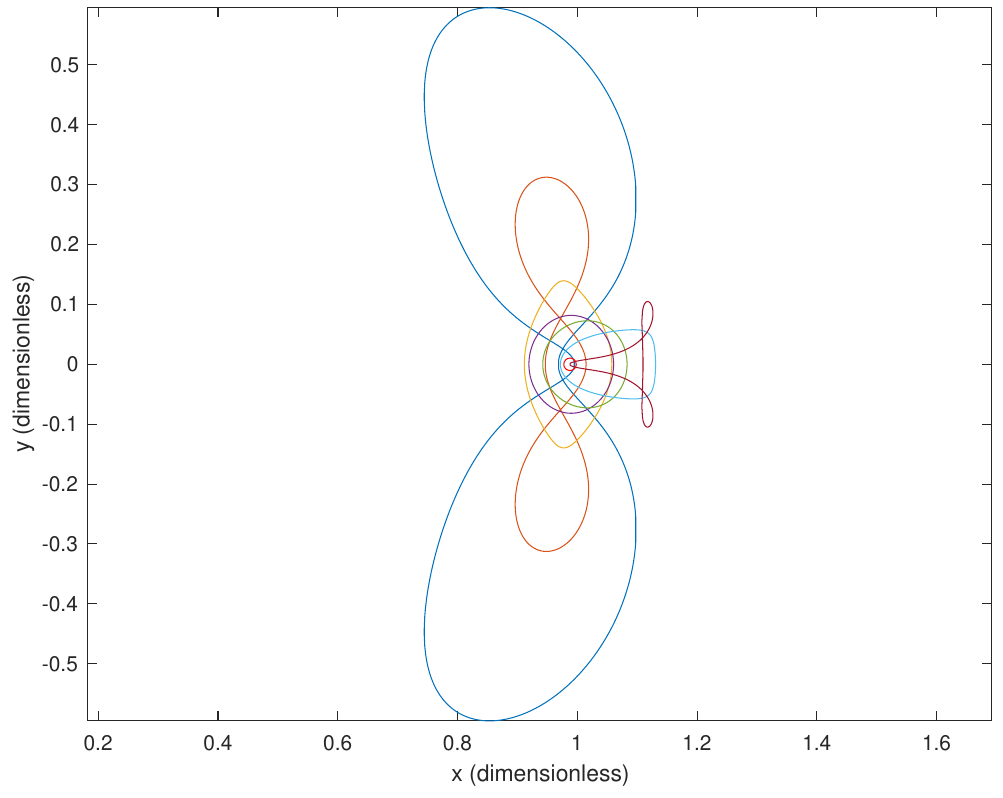}
    \caption{The H2 family of planar lunar orbits found by Broucke (Part I).}
    \label{H2I}
\end{figure}

\begin{figure}
    \centering
    \includegraphics[width=0.5\linewidth]{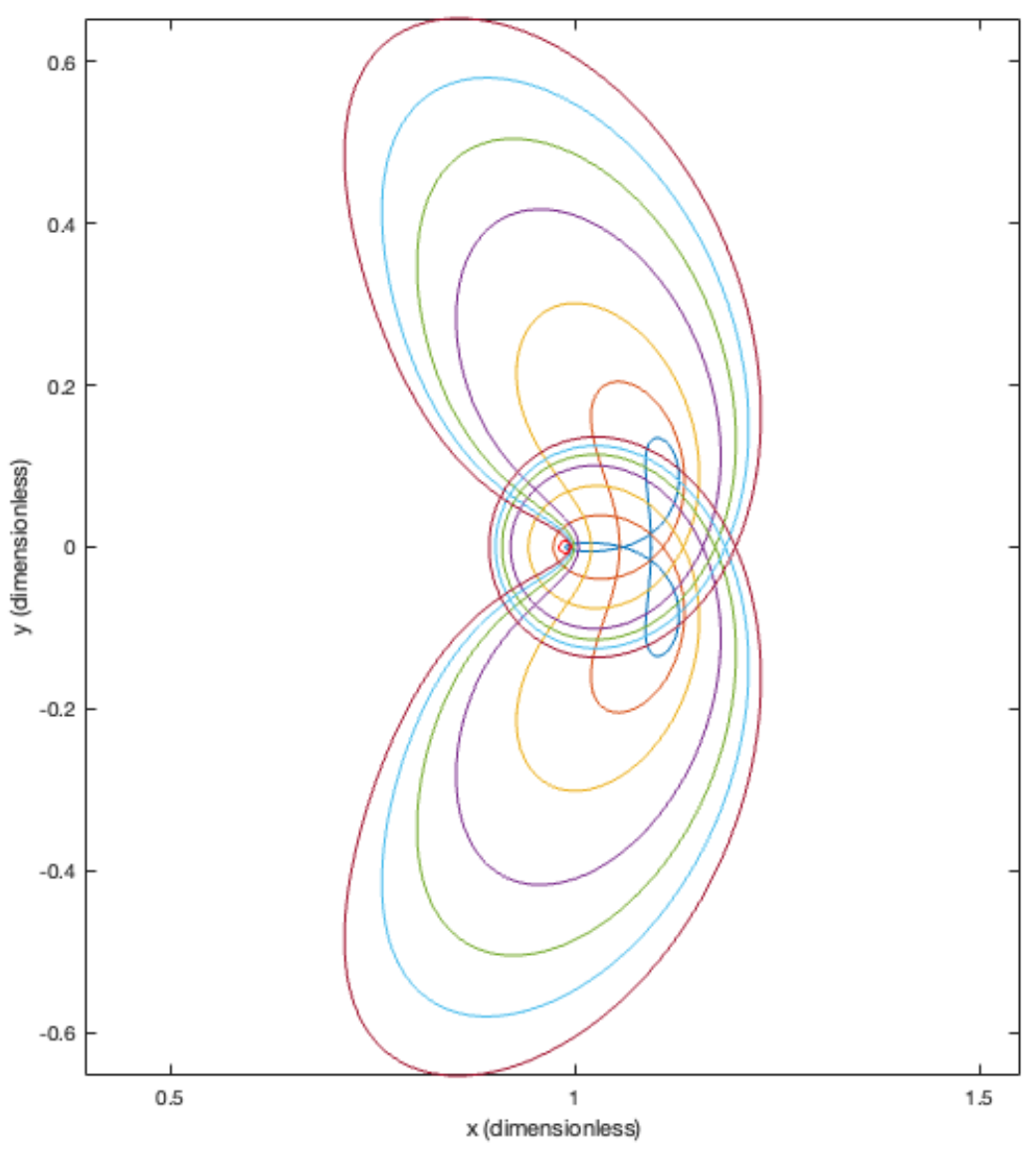}
    \caption{The H2 family of planar lunar orbits found by Broucke (Part II).}
    \label{H2II}
\end{figure}

For starters, Figures \ref{H1I}, \ref{H1II} form one family, called H1, and Figures \ref{H2I}, \ref{H2II}, form another one, called H2, as found by Broucke. We should emphasize that H1 is exactly the g-LPO1 branch of Figure \ref{fig:pitchfork2} (but for Earth-Moon rather than Jupiter-Europa), while H2 is the Birth-Death DPO-LPO2 branch in the same figure.

\begin{figure}
    \centering
    \includegraphics[width=0.8\linewidth]{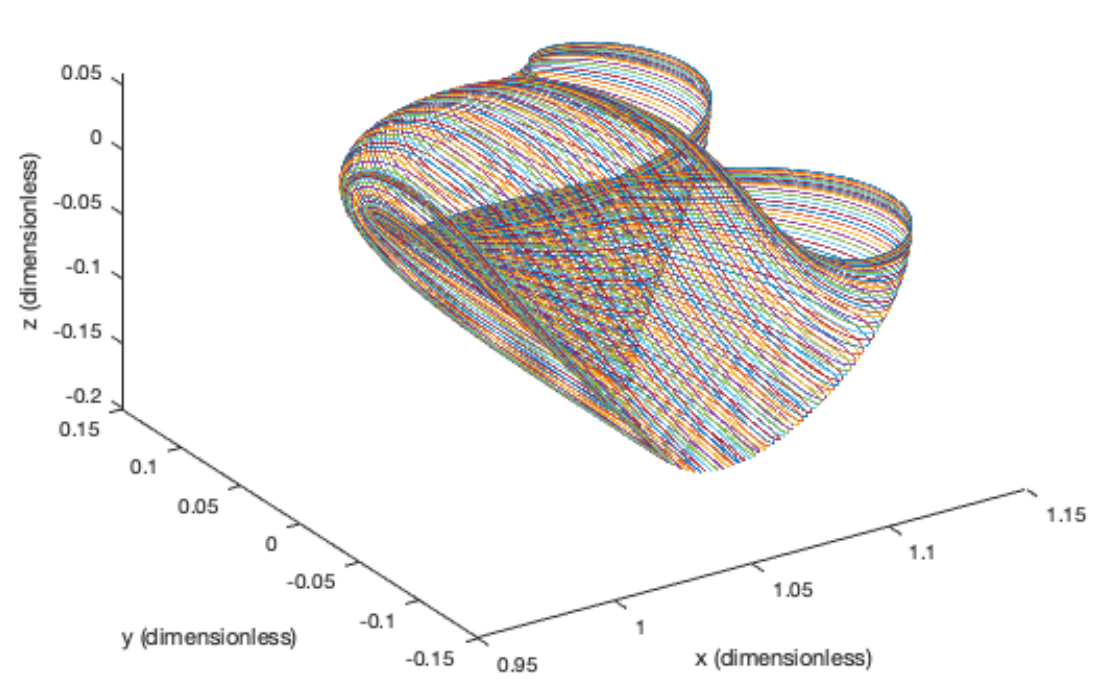}
    \caption{A family connecting H2 to the second cover of a Southern L2 Halo.}
    \label{fig:H2L2}
\end{figure}

\begin{figure}
    \centering
    \includegraphics[width=0.8\linewidth]{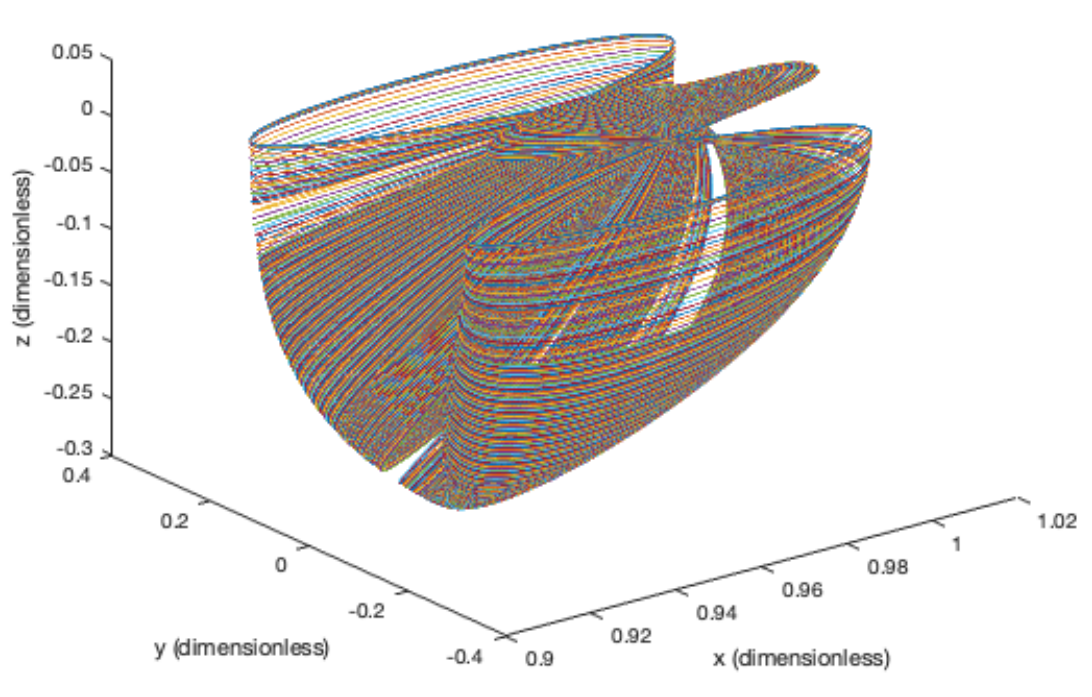}
    \caption{A family connecting H2 to a Southern L1 Halo.}
    \label{fig:non-sparse}
\end{figure}

\begin{figure}
    \centering
    \includegraphics[width=0.8\linewidth]{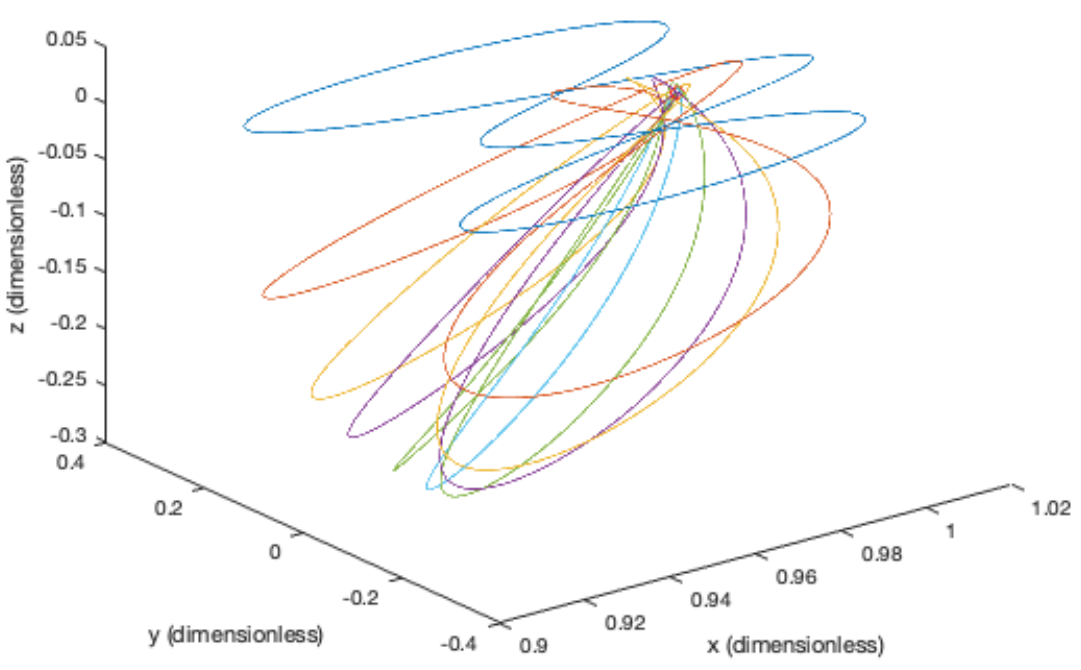}
    \caption{The same family from Figure \ref{fig:non-sparse} but sparser, to aid visualization. The two ``wings'' of the H2 orbit come together to doubly cover the Halo orbit.}
    \label{fig:sparse}
\end{figure}

The remaining Figures show out of plane connections. For instance, Figures \ref{fig:non-sparse} and \ref{fig:sparse} show (the same, not sparse and sparse, respectively) family that starts at one of Broucke’s H2 orbits, goes out of plane, passes through the singularity at the Moon (a single collision), and then ends at a double cover of the L1 Halo orbit at Jacobi constant value of $c = 2.94344$. Computing the orbit family as it passed through the singularity at the Moon required KS regularization. The same Broucke H2 family connects (at a different energy value) to a doubly covered L2 Halo orbit (see Figures \ref{fig:H2L2}), which in turn joins H1.
In summary:
\begin{itemize}
\item H2 joins L2 halo double cover (no collision);
\item H1 joins L1 halo double cover (no collision);
\item H2 joins L1 halo double cover after passing through moon collision.
\end{itemize}

\begin{figure}
    \centering
    \includegraphics[width=0.7\linewidth]{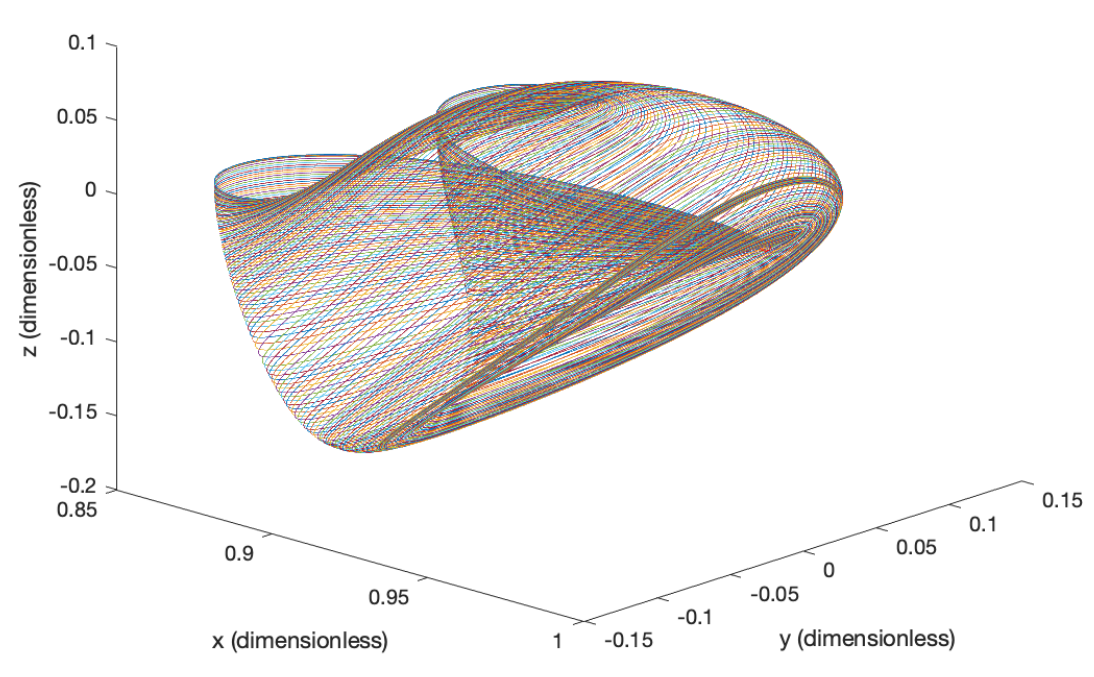}
    \caption{A family connecting H1 to Southern L1 Halo.}
    \label{fig:H1SL1HALO}
\end{figure}

\begin{figure}
    \centering
    \includegraphics[width=0.7\linewidth]{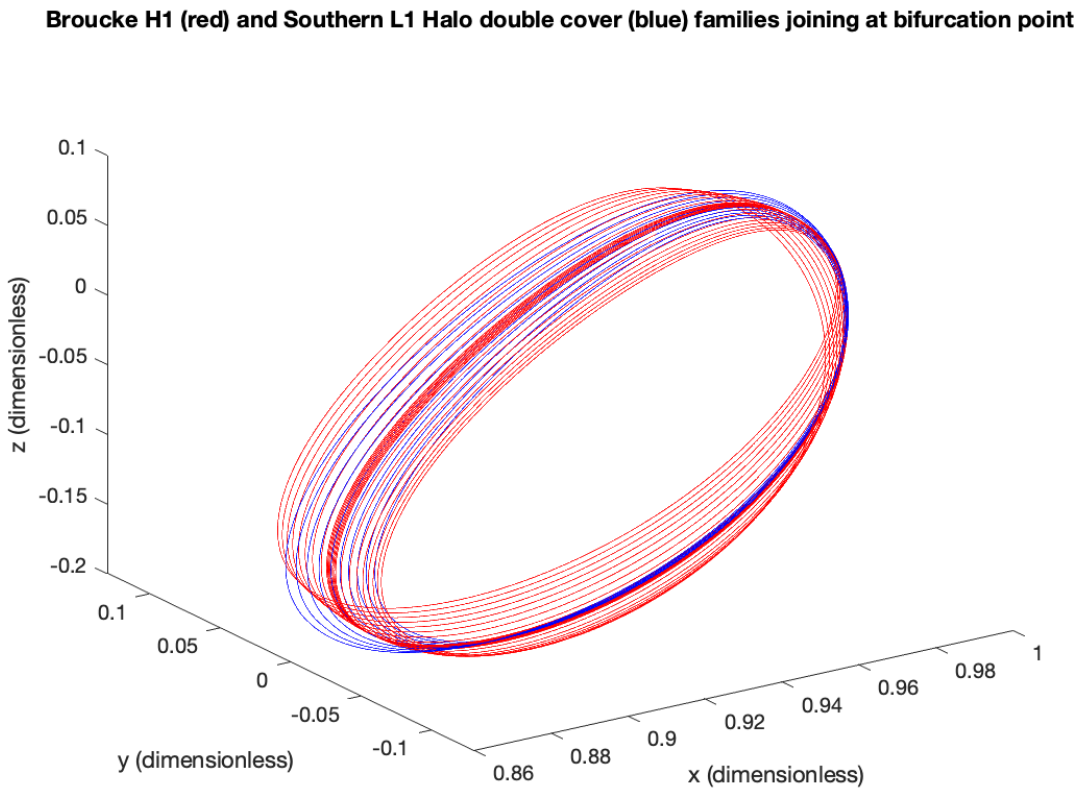}
    \caption{A portion of the same family from Figure \ref{fig:H1SL1HALO} and a portion (in blue) of the halo orbits, showing more closely how the two families merge together.}
    \label{fig:H12NDHALO}
\end{figure}


\newpage

\begin{figure}
    \centering
    \includegraphics[width=0.77\linewidth]{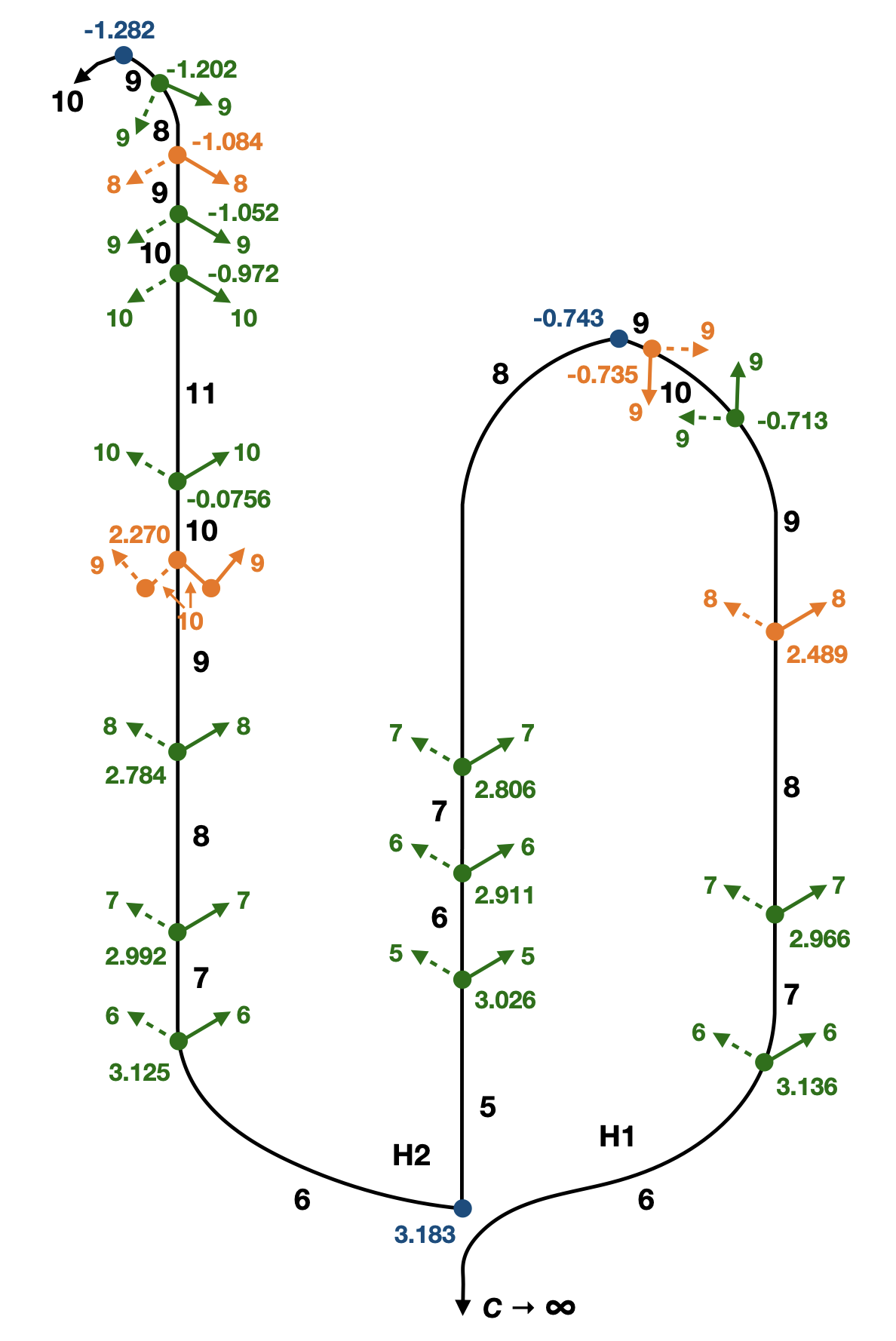}
    \caption{A bifurcation diagram for the H1-H2 prograde orbit family. Green dots indicate
out-of-plane bifurcations, orange represent planar symmetry-breaking bifurcations, and dark blue
are fold bifurcations. Decimal numbers are Jacobi constants $c$ of critical bifurcation orbits; $c$ decreases (energy increases) as one moves up the page. The integers are CZ indices. Dotted lines
indicate that the corresponding families have been obtained by reflection along the xy-plane (for
spatial orbits) or xz-plane (for planar asymmetric orbits) of the families represented by solid lines. Picture courtesy of Bhanu Kumar.}
    \label{fig:Bhanu1}
\end{figure}

\begin{figure}
    \centering
    \includegraphics[width=1\linewidth]{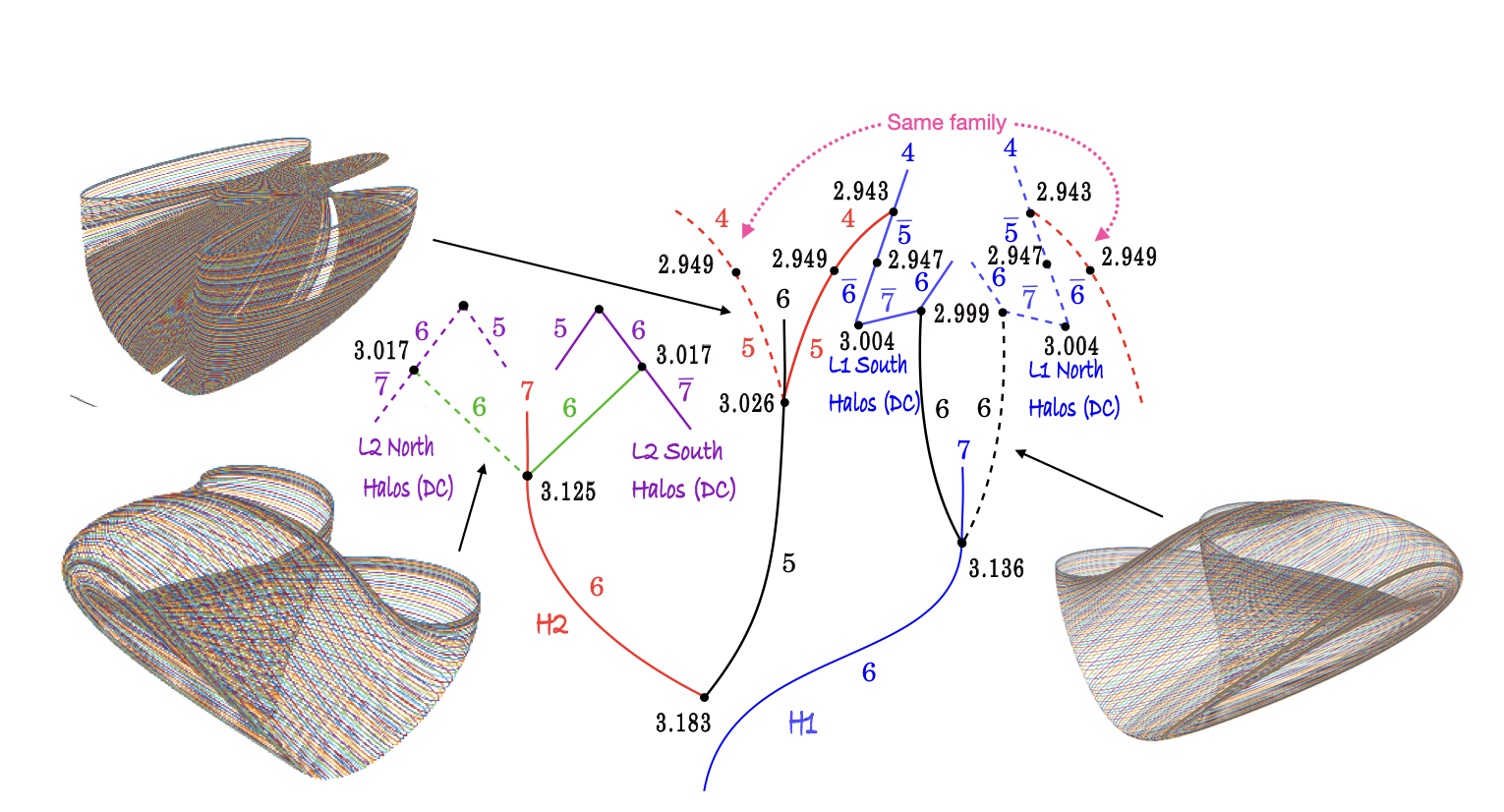}
    \caption{Bifurcation graph relating H1 \& H2 orbits with L1 \& L2 Halo double cover orbits. Orbits from 3 new families are shown. Dotted lines indicate orbits obtained by reflection.
Bars over CZ indices indicate that the corresponding orbits are homologically-bad (do not affect the
Floer invariant). The two dotted red curves represent the same family and should be identified. Floer numbers all check out.}
    \label{fig:Bhanu2}
\end{figure}

\newpage

\textbf{Further new connections in the Earth--Moon system.} The following is based on \cite{KM25}, a paper which carries a thorough investigation of the bifurcation network structure of several Earth-Moon CR3BP periodic orbit families. Among others, orbits near the Moon are considered, such as low and distant prograde, retrograde, and Halo orbits, as well as Earth-centered prograde, retrograde, and resonant periodic orbits. In that paper, several previously unnoticed families and connections between them by combining numerical continuation with another regularization scheme due to Kustaanheimo--Steifel, together with the symplectic toolkit. Regularization reveals that certain families previously thought of as disjoint form the same family (going through collision), confirming a conjecture by R. Broucke from 1968. The symplectic toolkit is useful to characterize and catalogue the families that arise, and to detect bifurcations. We will content ourselves with showing plots of some of the bifurcation diagrams in that paper; see Figures \ref{fig:Bhanu1} and \ref{fig:Bhanu2}.

\end{document}